\documentclass[12pt,twoside]{article}
\usepackage{multirow}
\usepackage{amsmath,amssymb}
\usepackage{graphicx}
\usepackage{subcaption}

\pdfoutput=1

\allowdisplaybreaks[4]

\newtheorem{lm}{Lemma}[section]
\newtheorem{thm}{Theorem}[section]
\newtheorem{rmk}{Remark}[section]
\newtheorem{prop}{{\bf Proposition}}[section]
\newtheorem{cor}{Corollary}[section]
\newtheorem{pro}{Proposition}[section]
\newtheorem{eg}{Example}[section]

\newcommand{\bm}[1]{\mbox{\boldmath{$#1$}}}

\newcounter{saveeqn}%

\makeatletter
\evensidemargin
\oddsidemargin
\makeatother

\title{\large\bf
Orbits and tsectors
in irregular exceptional directions of full-null degenerate singular point
}
\author{
{\sc Jun Zhang}\,$^{a}$,~~
{\sc Xingwu Chen}\,$^{b}$,~~
{\sc Weinian Zhang}\,$^{b}$
\footnote{Corresponding author. E-mail addresses: wnzhang@scu.edu.cn, matzwn@126.com (W. Zhang)}
\\
$^{a}$
{\small College of Mathematics and Physics \& Sichuan Geomath Key Lab}
\\
{\small Chengdu University of Technology, Sichuan 610059, P. R. China}
\\
$^{b}$
{\small School of Mathematics, Sichuan University}
\\
{\small Chengdu, Sichuan 610064, P. R. China}
}
\date{}

\begin{document}
\maketitle

\begin{abstract}
Near full-null degenerate singular points of analytic vector fields,
asymptotic behaviors of orbits are not given by eigenvectors but totally decided by nonlinearities.
Especially, in the case of high full-null degeneracy, i.e.,
the lowest degree of nonlinearities is high,
such a singular point may have irregular exceptional directions and
the blow-up technique can be hardly applied,
which leaves a problem how to determine
numbers of orbits and (elliptic, hyperbolic and parabolic) tangential sectors in this case.
In this paper we work on this problem.
Using Newton polygons to decompose nonlinearities into principal parts and remainder parts,
we convert the problem to the numbers of nonzero real roots
of edge-polynomials of principal parts.
Computing Newton polygons for multiplication and differentiation of analytic functions
and
giving Newton polygons for addition, which was not found in literatures,
we determine semi-definiteness of the Lie-bracket of principal parts and therefore
obtain criteria for those numbers.

\vskip 0.2cm


\vskip 0.2cm
{\bf AMS (2020) subject classification:}
34A26; 
34C05. 

\end{abstract}


\baselineskip 15pt    
\parskip 10pt         

\thispagestyle{empty}
\setcounter{page}{1}

\tableofcontents

\newpage

\section{Introduction}
\setcounter{equation}{0}
\setcounter{lm}{0}
\setcounter{thm}{0}
\setcounter{rmk}{0}
\setcounter{df}{0}
\setcounter{cor}{0}
\setcounter{pro}{0}

Let $\mathfrak{X}$ be an analytic vector field on a 2-dimensional real analytic manifold.
An isolated singular point $P$ of $\mathfrak{X}$ is said to be {\it degenerate}
if the Jacobian matrix at $P$ has a vanished eigenvalue.
The case of exact one vanished eigenvalue can be discussed on a center manifold,
where the most elementary bifurcations of codimension 1
such as saddle-node bifurcation, transcritical bifurcation and pitchfork bifurcation
are investigated in \cite[p.177]{CLW} or \cite[p.145]{GH}.
When eigenvalues are both vanished,
there are two possibilities: nilpotent case and (linear) full-null case.
In the nilpotent case the singular point may be cusp,
focus, center, saddle, node, saddle-node,
or a singular point with exactly one hyperbolic sector and one elliptic sector
as indicated in \cite[p.116]{DLA} or \cite[p.132]{ZZF}
and bifurcations such as the Bogodanov-Takens bifurcation (\cite{Bog, Takens}) may occur.

Along with the deep research on vector fields,
especially for the full-null case (where the Jacobian matrix is the zero matrix),
delicate investigation to nonlinearities becomes important
because the linear part is not decisive
and nonlinear terms are all resonant by the Poincar\'e normal form theory.
In this case, efforts have been made to versal unfoldings within specific classes, for example,
Krauskopf and Rousseau (\cite{K-R}) discussed
the versal unfolding for a quartic polynomial differential system of full-null degeneracy
within the class of reflectionally symmetric systems,
Ruan, Tang and Zhang (\cite{RTZ-JDE}) gave a versal unfolding for
a predator-prey system of full-null degeneracy
within the class of generalized Lotka-Volterra systems, and
Tang and Zhang (\cite{T-Z2}) investigated the versal unfolding
for an analytic Hamiltonian system of full-null degeneracy
within the class of Hamiltonian systems.

Before investigating bifurcations,
a fundamental work is to give qualitative properties of orbits near degenerate singular points.
In a small neighborhood of the singular point $P$,
the vector field $\mathfrak{X}$ can be expressed equivalently as
\begin{equation}
{\cal X}(x,y)\frac{\partial }{\partial x}+{\cal Y}(x,y)\frac{\partial }{\partial y}
\label{equ:initial}
\end{equation}
for $(x,y)\in\mathbb{R}^2$ near the origin $O: (0,0)$ with
$$
{\cal X}(x,y):=\sum_{k\ge m} X_k(x,y),~~~
{\cal Y}(x,y):=\sum_{k\ge m} Y_k(x,y),
$$
where $X_k$ and $Y_k$ are both homogeneous polynomials of degree
$k \ge m\ge 2$ and $X_m^2+Y_m^2\not\equiv0$.
For convenience we refer to the case of full-null linear part as the {\it full-null degeneracy}
and the case of $m>2$ as {\it high full-null degeneracy}.
If $O$ is neither a center nor a focus,
then there is at least one orbit connecting with $O$ in a definite direction,
called an {\it exceptional orbit} or {\it characteristic orbit} in \cite[p.17]{DLA}.
Further,
a neighborhood of $O$ can be decomposed by some exceptional orbits
into finitely many sectors,
which are elliptic, hyperbolic or parabolic sectors as indicated in
\cite[pp.17-19]{DLA} and \cite[pp.103-113]{ZZF}.
The qualitative properties near $O$ are decided by the number of those sectors
and their arrangement.
Bendixson (\cite{Ben}) concluded that
the numbers of elliptic sectors and hyperbolic sectors
are bounded by $2m$ and $2m+2$ respectively
and that the bound for hyperbolic sectors is sharp.
Berlinskii (\cite{Ber}) improved the bound for elliptic sectors to be $2m-1$ and
proved it to be sharp.
Sagalovich (\cite{Saga}) and Schecter and Singer (\cite{SS})
proved that
the number of separatrices at $O$,
exceptional orbits being boundaries of heperbolic sectors,
is bounded by $4m-2$ sharply.
However,
{\it exact number} of exceptional orbits
and {\it exact numbers} of various sectors remain to be answered.

\begin{figure}[h]
    \centering
    \subcaptionbox{%
     }{\includegraphics[height=1.6in,width=1.6in]{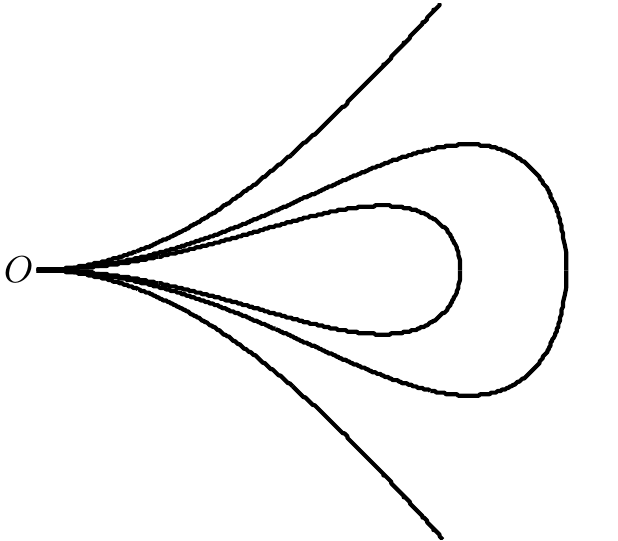}}~
     \subcaptionbox{%
     }{\includegraphics[height=1.6in,width=1.6in]{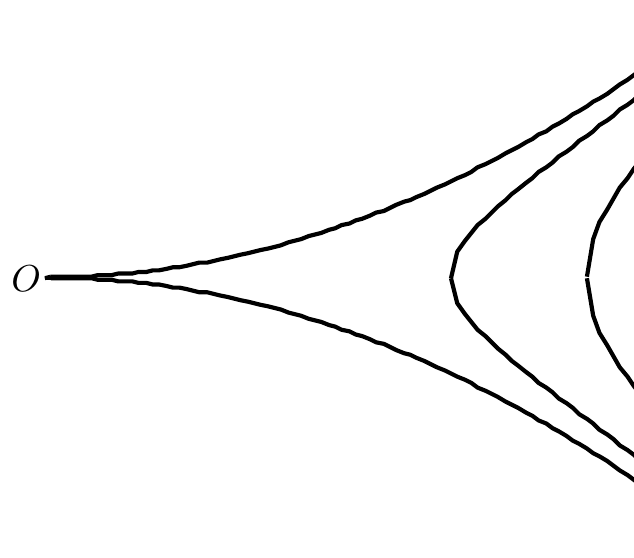}}~~~
     \subcaptionbox{%
     }{\includegraphics[height=1.6in,width=1.6in]{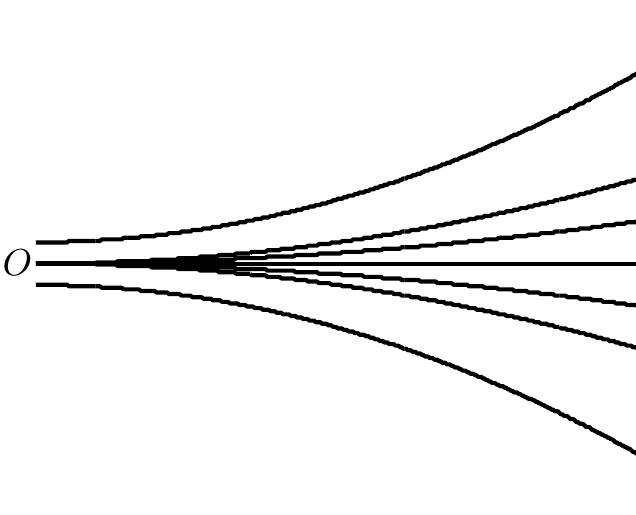}}
    \caption{Three types of tsector: (a) e-tsector, (b) h-tsector and (c) p-tsector.}
    \label{figPHE}
\end{figure}

Rewriting vector field \eqref{equ:initial} in the polar coordinates
$x=\rho\cos\theta$ and $y=\rho\sin\theta$ and
eliminating the common factor $\rho^{m-1}$,
we obtain
\begin{equation}
\rho{\cal H}(\rho,\theta)\frac{\partial }{\partial \rho}
+{\cal G}(\rho,\theta)\frac{\partial }{\partial \theta}
\label{equ:polar system}
\end{equation}
with
\begin{equation*}
{\cal H}(\rho,\theta)
:=H_0(\theta)+\sum_{k\ge 1} H_k(\theta)\rho^k,~~~
{\cal G}(\rho,\theta)
:=G_0(\theta)+\sum_{k\ge 1} G_k(\theta)\rho^k,
\end{equation*}
where for all $k\ge 0$,
\begin{eqnarray*}
&&G_{k}(\theta):= Y_{k+m}(\cos\theta,\sin\theta)\cos\theta - X_{k+m}(\cos\theta,\sin\theta)\sin\theta,
\\
&&H_{k}(\theta):= X_{k+m}(\cos\theta,\sin\theta)\cos\theta + Y_{k+m}(\cos\theta,\sin\theta)\sin\theta.
\end{eqnarray*}
In what follows, we assume that
neither ${\cal G}$ nor ${\cal H}$ is identical to 0;
otherwise, the singular point $O$ is either a center or a node.
As indicated in \cite{Frommer28, Hartman, NS60, SC, ZZF},
each zero $\theta=\theta_*$ of $G_0$ is called an {\it exceptional direction} (or {\it characteristic direction}, or {\it critical direction})
and exceptional orbits approach or leave the singular point $O$
only in exceptional directions.
The exceptional direction $\theta=\theta_*$ is referred to as
a {\it regular} exceptional direction by Frommer~(\cite{Frommer28})
if $H_0(\theta_*)\ne 0$.
Otherwise,
the direction is called an {\it irregular} (or {\it singular}) exceptional direction.
We call a (elliptic, hyperbolic or parabolic) sector $S$ a {\it tangential sector} or
simply {\it tsector} in an exceptional direction $\theta=\theta_*$
if for arbitrarily small $\delta>0$ the intersection $S\cap {\cal C}_\delta$,
where
$$
{\cal C}_\delta:=\{(\rho,\theta)\in\mathbb{R}^2:\rho\ge 0, |\theta-\theta_*|< \delta\}
$$
is a neighborhood of the direction $\theta=\theta_*$,
always contains a sector $S_\delta$ of the same type as $S$.
As shown in Figure~\ref{figPHE}, there are also three types of tsector,
called {\it e-tsector}, {\it h-tsector} and {\it p-tsector} respectively for short.
Similarly to \cite[pp.17-19]{DLA},
we keep the number of p-tsectors as small as possible
both by joining two adjacent p-tsectors,
and by adding a p-tsector to an e-tsector if they are adjacent.


Note that the complexity of a degenerate singular point lies in exceptional directions
because exceptional orbits connecting with the singular point in exceptional directions
and moreover orbits topologically rotate in a small circular sector with the vertex at $O$
and containing no exceptional directions,
as indicated in \cite[Theorem~3.1, p.60]{ZZF}.
Moreover,
singular points attempt to arise in exceptional directions as parameters vary,
correspondingly generating rich bifurcation phenomena,
because branches of isoclines gather in exceptional directions and therefore
might generate intersection points as parameters vary.
Thus,
we focus on an exceptional direction of the degenerate singular point $O$ and consider
\begin{description}
  \item[Question 1:] {\it How many orbits connect with $O$ in an exceptional direction?}
  \item[Question 2:] {\it How many tsectors are there in an exceptional direction?}
\end{description}

In order to answer the above two questions,
one method is the blowing-up technique (see \cite[Chapter~3]{DLA}),
performing the Briot-Bouquet transformation repeatedly
to decompose a degenerate singular point into several simpler ones.
Another method is applying Frommer's idea (\cite{Frommer28}), i.e.,
using isoclines and some constructed lines (curves) to
partition the angular neighborhood of the exceptional direction so that
orbits going forth or back to the singular point $O$ can be determined
in those partitioned sub-regions.
Moreover, Z-sector (\cite{Lef57,SC}),
normal sector (\cite{Frommer28,Hartman,Lef57,NS60,SC,ZZF}) and
generalized normal sectors (\cite{T-Z1}) were proposed to improve Frommer's idea.
As indicated in \cite[Theorem~3.3]{DLA},
a degenerate singular point of a given analytic vector field can be
blown-up into elementary singular points (with at most one zero eigenvalue)
in finitely many steps.
However, the blowing-up process gets rather complicated
when the nonlinearities possess higher degree or even unspecified degree (see e.g. \cite{T-Z2})
because one needs too many times of blowing-up
or even does not know how many times of blowing-up.
Moreover,
if we consider a family of vector fields having a full-null degenerate singular point $O$ with parameters,
for each fixed tuple of parameters
 one can use the blowing-up technique to
obtain the phase portrait in a small neighborhood of $O$
but the neighborhood depends on the parameters,
i.e.,
the blowing-up does not provide complete information of the degenerate singular point $O$
in a fixed neighborhood as parameters vary,
for example, does not answer whether another singular point arises from $O$.
So we will use Frommer's idea to investigate {\bf Questions 1} and {\bf 2}
because the idea neither requires to blow-up many times
nor demands to inverse so many times.

There were made efforts to the above questions (see e.g. \cite{Frommer28,Hartman,ZZF}).
Attentions are still paid to the lowest order term $G_0$, which has three cases:
$G_0(\theta)\ne 0$, $G_0(\theta)\equiv0$, or
$G_0$ has finitely many zeros in the interval $[0,2\pi)$.
In the first case,
there are no exceptional directions.
By Theorem 3.2 of \cite[p.61]{ZZF}
system \eqref{equ:initial} is {\it monodromy},
that is, no orbits connect with $O$ in a definite direction,
which leads to works of identifying focus from center
(e.g. \cite{GLMM, Medvedeva}).
In the second case,
every $\theta_*\in[0,2\pi)$ is an exceptional direction
but only finitely many of them are irregular and,
by Theorem 3.3 of \cite[p.63]{ZZF},
there is exactly one orbit connecting with $O$,
one p-tsector and no other tsectors in each regular direction.
In the third case, without loss of generality,
let $\theta=0$ be an exceptional direction.
Then $G_0(\theta)=g_0{\theta}^{\ell_0}+O(\theta^{\ell_0+1})$,
where $\ell_0\in\mathbb{N}:=\{1,2,3,...\}$ and $g_0\in\mathbb{R}\backslash\{0\}$.
If $H_0(0)\ne 0$,
Theorems 3.4, 3.7 and 3.8 of \cite[pp.68-75]{ZZF} ensure that
in the regular exceptional direction $\theta=0$,
\begin{description}

\item[(K1)]
there are infinitely many orbits connecting with $O$,
one p-tsector and no other tsectors
as either $\ell_0$ is even, or $\ell_0$ is odd and $g_0H_0(0)>0$, or

\item[(K2)]
there is one orbit connecting with $O$ and no tsectors
as $\ell_0$ is odd and $g_0H_0(0)<0$.

\end{description}
Oppositely,
if $H_0(0)=0$
then $H_0(\theta)=h_0\theta^{\tilde{\ell}_0}+O(\theta^{\tilde{\ell}_0+1})$,
where $\tilde{\ell}_0\in\mathbb{N}\cup\{+\infty\}$, $h_0\in\mathbb{R}\backslash\{0\}$,
and $\tilde{\ell}_0=+\infty$ implies that $H_0(\theta)\equiv0$.
Theorems 4.2 and 4.4 of \cite[pp.221-222]{Hartman} show that
in the irregular exceptional direction $\theta=0$,
\begin{description}
\item[(K3)]
at least one orbit connects with $O$
as $\ell_0$ is odd,
or

\item[(K4)]
infinitely many orbits connect with $O$ in the region $R:=\{(\rho,\theta)\in\mathbb{R}^2:0<\theta(\rho)<\theta<\delta, 0<\rho<\epsilon\}$,
where $\epsilon$ and $\delta$ are sufficiently small constants and
$\theta(\rho)$ is a continuously differentiable function,
as $\ell_0>\tilde{\ell}_0+1$ and $g_0h_0>0$.
\end{description}
However,
in the second case
neither the number of orbits nor the number of tsectors
in the irregular exceptional direction $\theta=\theta_*$, i.e., $H_0(\theta_*)=0$, are given yet.
In {\bf(K3)} of the third case,
no results show the exact number of orbits connecting with $O$,
not mentioning tsectors.
In {\bf(K4)},
it is unknown if there are orbits connecting with $O$ and tsectors
in the complement of the region $R$ near the direction $\theta=0$.
If we simply assume that $\theta=0$ is an exceptional direction,
then {\it the above three situations with unknown things} can be referred to the circumstance
\begin{eqnarray}
G_0(0)=H_0(0)=0,
\label{G0H00}
\end{eqnarray}
i.e., the circumstance that the exceptional direction $\theta=0$ is irregular.
This is actually the unsolved situation given in \cite[section~7]{Frommer28}.
Although there are examples showing the existence of several tsectors
in an irregular exceptional direction
(see e.g. \cite[Section~10]{Frommer28} and \cite[Example~4]{SS}),
we hardly find a criterion for those numbers in {\bf Questions 1} and {\bf 2}
for irregular exceptional directions.

In this paper
we answer to {\bf Questions 1} and {\bf 2}
in the above circumstance \eqref{G0H00}.
We use Newton polygons to decompose those nonlinearities of ${\cal G}$ and ${\cal H}$
into principal parts and remainder parts separately and then,
determine those numbers mentioned in {\bf Questions 1} and {\bf 2}
by numbers of nonzero roots of edge-polynomials of principal parts.
The main results Theorems \ref{th:finite} and \ref{th:finiteJ2},
devoted to the ``one-above case'' and the ``two-below case''
(defined at the beginning of section 2.2) respectively,
are stated in section 2.
Furthermore,
in the case of no parallel edges (section~2.3),
we simplify conditions given in Theorems \ref{th:finite} and \ref{th:finiteJ2}
to determine numbers of orbits and tsectors in an exceptional diection
with sign changes in the coefficient sequences corresponding to
vertices of Newton polygons of ${\cal G}$ and ${\cal H}$
(Corollaries \ref{cor:nonpara}-\ref{cor-n2}).
Our main results are proved in section 5,
where we divide a small neighborhood of the exceptional direction into Z-sectors
by real branches of the equation ${\cal H}(\rho,\theta)=0$ and
determine the number of orbits in each Z-sector by
the semi-definiteness of the function
$$
[{\cal G},{\cal H}]_\theta:={\cal G}'_\theta{\cal H}-{\cal G}{\cal H}'_\theta,
$$
the Lie-bracket of ${\cal G}$ and ${\cal H}$ in the variable $\theta$.
For preparations,
in section 3,
we first use Newton polygons to desingularize ${\cal H}$
so as to find all real branches of the equation ${\cal H}(\rho,\theta)=0$,
which are employed to divide Z-sectors.
Then,
we determine the classes for those Z-sectors
by desingularizing ${\cal G}$ to determine the signs of ${\cal G}$ on those real branches.
Since
a Z-sector of Class II (resp. III)
does not answer whether 1 (resp. 0) or infinitely many orbits connect with $O$ definitely
by known result~{\bf(c)} (resp. {\bf(d)}) of \cite[pp.81-83]{SC},
finally
we reduce the problem on the exact number of orbits in the Z-sector
to determining the semi-definiteness of
the Lie-bracket $[{\cal G},{\cal H}]_\theta$ in the Z-sector,
for which we need to apply desingularization on each edge of its Newton polygon.
In order to present the Newton polygon of $[{\cal G},{\cal H}]_\theta$
in terms of the Newton polygons of the given ${\cal G}$ and ${\cal H}$,
we study {\it addition, multiplication} and {\it differentiation} of Newton polygons
in section 4.
We hardly found a theory of addition in literature but
only theories of the latter two operations from references \cite{BK86, CA00, JP00, Tei}.
In section 5 we complete the proofs of the main theorems.
At last,
in section 6,
we first apply Theorem~\ref{th:finiteJ2} to a quartic polynomial system of full-null degeneracy,
showing exact 1 h-tsector or e-tsector in the irregular exceptional direction,
which corrects a result given in \cite{K-R}.
Then,
we use Theorems~\ref{th:finite} and \ref{th:finiteJ2} and Corollary~\ref{cor-n2} to
display five polynomial differential systems
which have 1 h-tsector, 2 h-tsectors, 3 h-tsectors, 2 e-tsectors and 3 e-tsectors
in an irregular exceptional direction separately.

\section{Results in exceptional directions}
\setcounter{equation}{0}
\setcounter{lm}{0}
\setcounter{thm}{0}
\setcounter{rmk}{0}
\setcounter{df}{0}
\setcounter{cor}{0}
\setcounter{pro}{0}

For convenience,
we let $\eta_O$ denote the number of orbits of system~\eqref{equ:initial}
connecting with $O$ in the exceptional direction $\theta=0$,
and let ${\cal S}_O^e$, ${\cal S}_O^h$ and ${\cal S}_O^p$
denote the numbers of e-tsectors, h-tsectors and p-tsectors
in the direction $\theta=0$ respectively.
In order to determine the numbers
$\eta_O$, ${\cal S}_O^e$, ${\cal S}_O^h$ and ${\cal S}_O^p$
in the direction $\theta=0$ which satisfies \eqref{G0H00},
we need the tool ``Newton polygon'' (\cite{BK86,CA00,JP00,Fischer}) and
some related concepts.

\subsection{Concepts on Newton polygon}

It is assumed that
neither ${\cal G}$ nor ${\cal H}$ is identical to 0
just below \eqref{equ:polar system}.
Then we expand nonzero analytic functions
${\cal G}$ and ${\cal H}$ at the point $(\rho,\theta)=(0,0)$ as
\begin{eqnarray}
{\cal G}(\rho,\theta)=\sum_{i+j\ge 0} a_{i,j}({\cal G})\rho^i\theta^j
~~~\mbox{and}~~~
{\cal H}(\rho,\theta)=\sum_{i+j\ge 0} a_{i,j}({\cal H})\rho^i\theta^j
\label{expGH}
\end{eqnarray}
respectively.
The set $\Delta({\cal G}):=\{(i,j)\in\mathbb{Z}_+^2:a_{i,j}({\cal G})\ne 0\}$,
where $\mathbb{Z}_+^2:=\{(i,j)\in\mathbb{Z}^2:i\ge0,j\ge0\}$,
is called the {\it valid index set}
(or {\it carrier} in \cite[p.380]{BK86} and \cite[p.197]{Fischer},
or {\it Newton diagram} in \cite[p.15]{CA00},
or {\it support} in \cite[p.176]{JP00}) of ${\cal G}$.
Further,
we embed the set $\Delta({\cal G})$ into the $(u,v)$-plane and
consider its {\it lower convex semi-hull}
$$
\wp(\Delta({\cal G})):={\rm conv}\{\Delta({\cal G})+\mathbb{R}_+^2\},
$$
where $\mathbb{R}_+^2:=\{(u,v)\in\mathbb{R}^2:u\ge0,v\ge0\}$ and
the Minkowski-sum
$\Delta({\cal G})+\mathbb{R}_+^2:=
\{(u,v)\in\mathbb{R}^2:u=u_1+u_2,v=v_1+v_2,
(u_1,v_1)\in\Delta({\cal G}),(u_2,v_2)\in\mathbb{R}_+^2\}$
as defined in \cite[p.3]{Soltan}.
Clearly, the boundary $\partial\wp(\Delta({\cal G}))$ consists of
one closed vertical ray, one closed horizontal ray and one compact polygon,
denoted by $L^v_{\cal G}$, $L^h_{\cal G}$ and ${\cal N}_{\cal G}$ respectively.
The compact polygon ${\cal N}_{\cal G}$, i.e.,
the principal boundary of $\wp(\Delta({\cal G}))$,
is called the {\it Newton polygon} of ${\cal G}$.
Then
the function ${\cal G}$ can be expressed as
${\cal G}(\rho,\theta)={\cal G}_P(\rho,\theta)+{\cal G}_R(\rho,\theta)$,
where
$$
{\cal G}_P(\rho,\theta):=
\sum_{(i,j)\in\Delta({\cal G})\cap{\cal N}_{\cal G}}a_{i,j}({\cal G})\rho^i\theta^j
~~~\mbox{and}~~~
{\cal G}_R(\rho,\theta):=
\sum_{(i,j)\in\Delta({\cal G})\setminus{\cal N}_{\cal G}}a_{i,j}({\cal G})\rho^i\theta^j,
$$
called
the {\it principal part} ${\cal H}_P$ and the {\it remainder part} ${\cal H}_R$ respectively.
We similarly define $\Delta({\cal H})$, $\wp(\Delta({\cal H}))$,
$L^v_{\cal H}$, $L^h_{\cal H}$, ${\cal N}_{\cal H}$, ${\cal H}_P$ and ${\cal H}_R$.

\begin{figure}[!h]
  \centering
  \includegraphics[height=1.6in,width=1.6in]{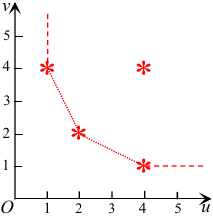}\\
  \caption{The Newton polygon of polynomial ${\cal G}$ defined in \eqref{egNP}.}
  \label{NGf}
\end{figure}

For example,
the polynomial
\begin{align}
{\cal G}(\rho,\theta)=-\rho\theta^4+\rho^2\theta^2-\rho^4\theta+\rho^4\theta^4
\label{egNP}
\end{align}
has valid indices $(1,4)$, $(2,2)$, $(4,1)$ and $(4,4)$,
and the Newton polygon ${\cal N}_{\cal G}$ is the compact polygon linking
points $(1,4)$, $(2,2)$ and $(4,1)$ successively,
as shown in Figure~\ref{NGf}.

Suppose that Newton polygon ${\cal N}_{\cal G}$ has $s({\cal G})+1$ vertices
$V_0({\cal G}):(i_0,j_0), ..., V_{s({\cal G})}({\cal G}):(i_{s({\cal G})},j_{s({\cal G})})$
with $i_0<\cdots<i_{s({\cal G})}$.
Then we call
$$
\Delta^V({\cal G}):=\{V_0({\cal G}),...,V_{s({\cal G})}({\cal G})\}
~~~\mbox{and}~~~
\vec\Delta^V({\cal G}):=(V_0({\cal G}),...,V_{s({\cal G})}({\cal G}))
$$
the {\it vertex set} and {\it vertex sequence} of ${\cal N}_{\cal G}$ respectively.
Let
\begin{equation}
\begin{aligned}
{\cal A}({\cal G})&:=(a_{i_0,j_0}({\cal G}),..., a_{i_{s({\cal G})},j_{s({\cal G})}}({\cal G})),
\\
{\cal A}_j({\cal G})&:=((-1)^{j_{s({\cal G})}}a_{i_{s({\cal G})},j_{s({\cal G})}}({\cal G}),...,
(-1)^{j_0}a_{i_0,j_0}({\cal G})),
\end{aligned}
\label{defAA}
\end{equation}
and call them the {\it vertex coefficient sequence} and
{\it $j$-algebraic vertex coefficient sequence} of ${\cal G }$ respectively.
Similarly, suppose that Newton polygon ${\cal N}_{\cal H}$ has $s({\cal H})+1$ vertices
$
V_0({\cal H}):(\tilde{i}_0,\tilde{j}_0),
...,
V_{s({\cal H})}({\cal H}):(\tilde{i}_{s({\cal H})},\tilde{j}_{s({\cal H})})
$
with $\tilde{i}_0<\cdots<\tilde{i}_{s({\cal H})}$.
We similarly define $\Delta^V({\cal H})$, $\vec\Delta^V({\cal H})$,
${\cal A}({\cal H})$ and ${\cal A}_j({\cal H})$.

For each $k=1,...,s({\cal G})$,
let $E_k({\cal G}):=\underline{V_{k-1}({\cal G})V_{k}({\cal G})}$,
the edge formed by the closed linear segment
linking $V_{k-1}({\cal G})$ with $V_k({\cal G})$.
Let $\zeta(E_k({\cal G}))$ denote the slope of the edge $E_k({\cal G})$.
Then we call
$$
{\mathfrak{S}}({\cal G}):=\{\zeta(E_1({\cal G})),...,\zeta(E_{s({\cal G})}({\cal G}))\}
~~\mbox{and}~~
\vec{\mathfrak{S}}({\cal G}):=(\zeta(E_1({\cal G})),...,\zeta(E_{s({\cal G})}({\cal G})))
$$
the {\it slope set} and {\it slope sequence} of ${\cal N}_{\cal G}$ respectively.
For each $k=0,...,s({\cal G})$,
let
\begin{equation}
\zeta^-(V_k({\cal G})):=\zeta(E_k({\cal G}))~~~\mbox{and}~~~
\zeta^+(V_k({\cal G})):=\zeta(E_{k+1}({\cal G}))
\label{SP+-}
\end{equation}
be the {\it left-sided slope} and
the {\it right-sided slope} of the vertex $V_k({\cal G})$ respectively,
where we complementarily define
$\zeta(E_0({\cal G})):=-\infty$ and $\zeta(E_{s({\cal G})+1}({\cal G})):=0$.
For any given $\xi<0$ we refer to the edge or vertex
\begin{equation}
\sigma({\cal G},\xi):=\left\{
\begin{array}{llll}
E_k({\cal G})~\mbox{if}~\xi=\zeta(E_k({\cal G}))
~\mbox{for a}~k\in\{1,2,...,s({\cal G})\},
\\
V_k({\cal G})~\mbox{if}~\xi\in(\zeta^-(V_k({\cal G})),\zeta^+(V_k({\cal G})))
~\mbox{for a}~k\in\{0,1,...,s({\cal G})\}
\end{array}
\right.
\label{defscomp}
\end{equation}
as the {\it $\xi$-component} of ${\cal N}_{\cal G}$.
Further,
we call
$$
u({\cal G},\xi):=i-j/\xi~~~\mbox{for}~(i,j)\in \sigma({\cal G},\xi)
$$
the {\it $u$-intercept} of the $\xi$-component $\sigma({\cal G},\xi)$
and call
\begin{eqnarray}
{\cal K}_\xi({\cal G})(\theta):=\sum_{(i,j)\in\sigma({\cal G},\xi)\cap\Delta({\cal G})} a_{i,j}({\cal G})\theta^j
\label{defkxip}
\end{eqnarray}
the {\it $\xi$-componential polynomial} of ${\cal G}$.
As above,
we similarly define $E_k({\cal H})$ for all $k=1,...,s({\cal H})$,
$\zeta^\pm(V_k({\cal H}))$ for all $k=0,...,s({\cal H})$,
${\mathfrak{S}}({\cal H})$, $\vec{\mathfrak{S}}({\cal H})$,
$\sigma({\cal H},\xi)$, $u({\cal H},\xi)$ and ${\cal K}_\xi({\cal H})$.

Take polynomial \eqref{egNP} as an example.
It has the vertex sequence $\vec\Delta^V({\cal G})=((1,4),(2,2),(4,1))$,
coefficient sequence ${\cal A}({\cal G})=(-1,1,-1)$,
the $j$-algebraic coefficient sequence ${\cal A}_j({\cal G})=(1,1,-1)$
and the slope sequence $\vec{\mathfrak{S}}({\cal G})=(-2,-1/2)$.
For $\xi=-2$,
the $\xi$-component $\sigma({\cal G},\xi)$ is the edge
linking vertex $(1,4)$ with vertex $(2,2)$,
and the $\xi$-componential polynomial
${\cal K}_\xi({\cal G})(\theta)=-\theta^4+\theta^2$.
For $\xi=-1/2$,
$\sigma({\cal G},\xi)$ is the edge linking vertex $(2,2)$ with vertex $(4,1)$
and ${\cal K}_\xi({\cal G})(\theta)=\theta^2-\theta$.
For $\xi\in(-\infty,-2)$,
$\sigma({\cal G},\xi)$ is the vertex $(1,4)$
and ${\cal K}_\xi({\cal G})(\theta)=-\theta^4$.
For $\xi\in(-2,-1/2)$,
$\sigma({\cal G},\xi)$ is the vertex $(2,2)$
and ${\cal K}_\xi({\cal G})(\theta)=\theta^2$.
For $\xi\in(-1/2,0)$,
$\sigma({\cal G},\xi)$ is the vertex $(4,1)$
and ${\cal K}_\xi({\cal G})(\theta)=-\theta$.

Define the {\it sequence union}
$\vec{\mathfrak{S}}({\cal G})\cup\vec{\mathfrak{S}}({\cal H})$
to be the sequence obtained by arranging those elements in the union
$\mathfrak{S}({\cal G})\cup\mathfrak{S}({\cal H})$
in ascending order.
Let $s({\cal G},{\cal H})$ denote the cardinality of the set
$\mathfrak{S}({\cal G})\cup\mathfrak{S}({\cal H})$.
Then, assume that
\begin{eqnarray}
\vec{\mathfrak{S}}({\cal G})\cup\vec{\mathfrak{S}}({\cal H})
=(\xi_1,...,\xi_{s({\cal G},{\cal H})}).
\label{unionseq}
\end{eqnarray}
Similarly,
we define the {\it sequence intersection}
$\vec{\mathfrak{S}}({\cal G})\cap\vec{\mathfrak{S}}({\cal H})$
to be the sequence obtained by arranging those elements in the intersection
$\mathfrak{S}({\cal G})\cap\mathfrak{S}({\cal H})$
in ascending order.
Obviously,
the sequence intersection
$\vec{\mathfrak{S}}({\cal G})\cap\vec{\mathfrak{S}}({\cal H})$
is a subsequence of the sequence union
$\vec{\mathfrak{S}}({\cal G})\cup\vec{\mathfrak{S}}({\cal H})$.
By \eqref{unionseq}, we assume that
\begin{eqnarray}
\vec{\mathfrak{S}}({\cal G})\cap\vec{\mathfrak{S}}({\cal H})
=(\xi_{s_1},...,\xi_{s_{r({\cal G},{\cal H})}}),
\label{interseq}
\end{eqnarray}
where $r({\cal G},{\cal H})\ge0$ and
$1\le s_1<\cdots<s_{r({\cal G},{\cal H})}\le s({\cal G},{\cal H})$.
By \eqref{interseq},
we divide the vertex sequence $\vec{\Delta}^V({\cal G})$
into $r({\cal G},{\cal H})+1$ sub-sequences
\begin{eqnarray}
\Lambda_\varrho({\cal G})
:=(V_{k_{\varrho-1}+1}({\cal G}),...,V_{k_{\varrho}}({\cal G})),
~~~\varrho=0,...,r({\cal G},{\cal H}),
\label{def-Lambda}
\end{eqnarray}
such that $\zeta^+(V_{k_{\varrho}}({\cal G}))=\xi_{s_{\varrho+1}}$
for all $\varrho=0,...,r({\cal G},{\cal H})-1$,
where $-1=k_{-1}<k_0<\cdots<k_{r({\cal G},{\cal H})}=s({\cal G})$.
Especially, $\Lambda_0({\cal G})=\vec{\Delta}^V({\cal G})$ if $r({\cal G},{\cal H})=0$.
Similarly, we can divide the vertex sequence $\vec{\Delta}^V({\cal H})$ into
$r({\cal G},{\cal H})+1$ subsequences
$\Lambda_0({\cal H}),...,\Lambda_{r({\cal G},{\cal H})}({\cal H})$.

\subsection{Main theorems}

\begin{figure}[h]
    \centering
     \subcaptionbox{%
     }{\includegraphics[height=1.6in]{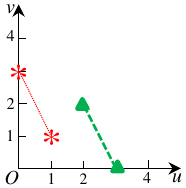}}~~~~~~~~~~~~
     \subcaptionbox{%
     }{\includegraphics[height=1.6in]{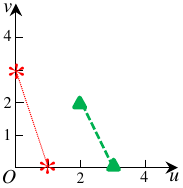}}
    \caption{Illustrations for (a) one-above case {\bf(J1)} and (b) two-below case {\bf(J2)},
    with dots for ${\cal N}_{\cal G}$ and dashes for ${\cal N}_{\cal H}$.}
    \label{fig:J1J2}
\end{figure}

Having known the above concepts and notations,
we investigate the number $\eta_O$ of orbits and
the numbers ${\cal S}_O^e$, ${\cal S}_O^h$ and ${\cal S}_O^p$
of tsectors in the direction $\theta=0$
satisfying \eqref{G0H00} in the following two cases:
\begin{description}
\item[(J1)]
$j_{s({\cal G})}+\tilde{j}_{s({\cal H})}>0$, i.e.,
either ${\cal N}_{\cal G}$ or ${\cal N}_{\cal H}$ ends above the $u$-axis;

\item[(J2)]
$j_{s({\cal G})}=\tilde{j}_{s({\cal H})}=0$, i.e.,
neither ${\cal N}_{\cal G}$ nor ${\cal N}_{\cal H}$ ends above the $u$-axis.
\end{description}
We simply refer to {\bf (J1)} and {\bf (J2)} as the {\it one-above case} and
the {\it two-below case} respectively, as illustrated in Figure~\ref{fig:J1J2}.


In the one-above case {\bf (J1)},
we make the following hypotheses:
\begin{description}

\item[(P1)]
Ordinates of lattice points in $\Lambda_\varrho({\cal G})$
(and $\Lambda_\varrho({\cal H})$),
defined in \eqref{def-Lambda}, have the same parity, and
the parity of $\Lambda_\varrho({\cal G})$ is different from
the one of $\Lambda_\varrho({\cal H})$
for each $\varrho=0,...,r({\cal G},{\cal H})$.
This is referred to as {\bf uniformly opposite parities} in {\bf(J1)}.

\item[(Q)]
Polynomial
$
[{\cal K}_{\xi_k}({\cal G}),{\cal K}_{\xi_k}({\cal H})]
:=({\cal K}_{\xi_k}({\cal G}))'_\theta{\cal K}_{\xi_k}({\cal H})
-{\cal K}_{\xi_k}({\cal G})({\cal K}_{\xi_k}({\cal H}))'_\theta,
$
the Lie-bracket of ${\cal K}_{\xi_k}({\cal G})$ and ${\cal K}_{\xi_k}({\cal H})$,
has no nonzero real roots for each $k=1,...,$ $s({\cal G},{\cal H})$,
where $\xi_k$ and $s({\cal G},{\cal H})$ are given in \eqref{unionseq},
${\cal K}_{\xi_k}({\cal G})$ and ${\cal K}_{\xi_k}({\cal H})$ are
$\xi_k$-componential polynomials,
i.e.,
all Lie-brackets of $\xi$-componential polynomials have {\bf definite signs}.
\end{description}
Moreover, let
\begin{eqnarray}
\chi({\cal G}):=\sum_{\xi\in\mathfrak{S}({\cal G})}
\sharp\{\theta\in\mathbb{R}\backslash\{0\}:{\cal K}_\xi({\cal G})(\theta)=0\},
\label{defChi}
\end{eqnarray}
where $\mathfrak{S}({\cal G})$ is the slope set,
${\cal K}_\xi({\cal G})$ is the $\xi$-componential polynomial,
and $\sharp$ is the cardinality of a set.
Define
\begin{eqnarray*}
N({\cal G}):=\left\{\begin{array}{ll}
\chi({\cal G})     & \mbox{ as }j_{s({\cal G})}=0,
\\
\chi({\cal G})+1   & \mbox{ as }j_{s({\cal G})}>0
\end{array}\right.
\end{eqnarray*}
and similarly define $\chi({\cal G})$ and $N({\cal H})$.
Let
$$
C_0:=(j_0-\tilde{j}_0)a_{i_0,j_0}({\cal G})a_{\tilde{i}_0,\tilde{j}_0}({\cal H}),
$$
where $a_{i_0,j_0}({\cal G})$ and $a_{\tilde{i}_0,\tilde{j}_0}({\cal H})$
are coefficients of ${\cal G}$ and ${\cal H}$ corresponding to their left-most vertices
$V_0({\cal G}):(i_0,j_0)$ and $V_0({\cal H}):(\tilde{i}_0,\tilde{j_0})$ respectively.

\begin{thm}
In the one-above case {\bf (J1)} of the circumstance \eqref{G0H00},
suppose that conditions {\bf(P1)} and {\bf(Q)} hold.
Then
in the subcase $G_0(\theta)\not\equiv0$,
\begin{description}
\item[(ia)]
if $C_0<0$,
$\eta_O=N({\cal G})$,
${\cal S}_O^e=0$, ${\cal S}_O^h=\max\{0,N({\cal G})-1\}$ and ${\cal S}_O^p=0$;

\item[(ib)]
if $C_0>0$,
either $\eta_O={\cal S}_O^e={\cal S}_O^h={\cal S}_O^p=0$ when $N({\cal G})=0$,
or $\eta_O=+\infty$, ${\cal S}_O^e=N({\cal G})-1$ and ${\cal S}_O^h={\cal S}_O^p=0$
when $N({\cal G})>0$;
\end{description}
and
in the subcase $G_0(\theta)\equiv0$,
\begin{description}
\item[(iia)]
if $C_0<0$,
$\eta_O=N({\cal H})+1$,
${\cal S}_O^e=0$, ${\cal S}_O^h=N({\cal H})$ and ${\cal S}_O^p=2$;

\item[(iib)]
if $C_0>0$,
either
$\eta_O=1$, ${\cal S}_O^e={\cal S}_O^h=0$ and ${\cal S}_O^p=1$ when $N({\cal H})=0$,
or $\eta_O=0$ or $+\infty$, ${\cal S}_O^e=N({\cal H})$ and ${\cal S}_O^h={\cal S}_O^p=0$
when $N({\cal H})\ge 1$.
In particular,
$\eta_O=0$ if and only if $\tilde{j}_{s({\cal H})}>0$ and $N({\cal H})=1$.
\end{description}
\label{th:finite}
\end{thm}

This theorem determines the numbers $\eta_O, {\cal S}_O^e,{\cal S}_O^h$ and ${\cal S}_O^p$, defined at the beginning of this section, in the case {\bf (J1)}.
For example, result {\bf (ia)} indicates that
system~\eqref{equ:initial} has $N({\cal G})$ orbits connecting with $O$
in the direction $\theta=0$,
neither e-tsectors nor p-tsectors,
and $\max\{0,N({\cal G})-1\}$ h-tsectors
in this direction if $G_0(\theta)\not\equiv 0$ and $C_0<0$.
The proof of Theorem~\ref{th:finite} is long and will be given in section 5.
Note that we state results {\bf (ia)} and {\bf (ib)} in terms of $N({\cal G})$
but {\bf (iia)} and {\bf (iib)} in terms of $N({\cal H})$,
which will be explained after the proof of this theorem.

\begin{rmk}
{\rm
In Theorem~\ref{th:finite},
we need to determine the numbers of nonzero real roots of
the polynomial $[{\cal K}_{\xi_k}({\cal G}),{\cal K}_{\xi_k}({\cal H})]$,
given in {\bf(Q)},
and polynomials ${\cal K}_\xi({\cal G})$ and ${\cal K}_\xi({\cal H})$,
which can be settled down by applying the theory of Hankel form
given in \cite[Chapter~XV]{Gant}
or the Complete Discrimination System given in \cite{YL}
(see also in section~3.3).
}
\label{Rk:poly}
\end{rmk}

\begin{figure}[h]
    \centering
     \subcaptionbox{%
     }{\includegraphics[height=1.6in]{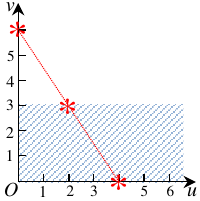}}~~~~~~~~~~~~
     \subcaptionbox{%
     }{\includegraphics[height=1.6in]{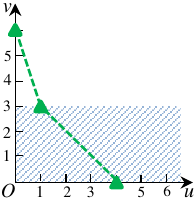}}
    \caption{Illustrations for (a) ${\cal N}_{\cal G}$ satisfying {\bf(S)} and
    (b) ${\cal N}_{\cal H}$ satisfying {\bf(S)}.}
    \label{fig:GHstripe}
\end{figure}

In the two-below case {\bf(J2)},
we need the following hypotheses:
\begin{description}
\item[(P2)]
Ordinates of lattice points in $\Lambda_\varrho({\cal G})$ (and $\Lambda_\varrho({\cal H})$),
defined in \eqref{def-Lambda}, have the same parity, and
the parity of $\Lambda_\varrho({\cal G})$ is different from
the one of $\Lambda_\varrho({\cal H})$
for each $\varrho=0,...,r({\cal G},{\cal H})-1$.
Moreover,
ordinates of lattice points in
$
\Lambda_{r({\cal G},{\cal H})}({\cal G})\backslash$ $(V_{s({\cal G})}({\cal G}))
:=(V_{k_{r({\cal G},{\cal H})-1}+1}({\cal G}),...,V_{s({\cal G})-1}({\cal G}))
$
are all odd, and those in $\Lambda_{r({\cal G},{\cal H})}({\cal H})$ are all even.
Like {\bf(P1)}, these are referred to as {\bf uniformly opposite parities} in {\bf(J2)}.

\item[(S)]
The Newton polygon $\mathcal{N}_{\cal G}$ (and $\mathcal{N}_{\cal H}$) has at least two valid indices,
the last two of which are
$V_{s({\cal G})}({\cal G}):(i_{s({\cal G})}, 0)$
(and $V_{s({\cal H})}({\cal H}):(\tilde{i}_{s({\cal H})}, 0)$)
and $V_*({\cal G}):(i_*, j_*)$
(and $V_*({\cal H}):(\tilde{i}_*, \tilde{j}_*)$)
such that
$\{(i,j)\in\Delta({\cal G}):0< j< j_*\}=\emptyset$
(and $\{(\tilde{i},\tilde{j})\in\Delta({\cal H}):0< \tilde{j}< \tilde{j}_*\}=\emptyset$),
as shown in Figure~\ref{fig:GHstripe}.
This means that there are
{\bf no valid indices} in the horizontal stripes between the last two valid indices.
\end{description}


\begin{thm}
In the two-below case {\bf (J2)} of the circumstance \eqref{G0H00},
suppose that
{\bf(P2)}, {\bf(Q)} and {\bf(S)} are satisfied and
one of the following conditions holds:
\begin{description}
\item[(H1)]
$\zeta(E_{s({\cal G})}({\cal G}))>\zeta(E_{s({\cal H})}({\cal H}))$ and
$j_*\le \tilde{j}_*$,

\item[(H2)]
$\zeta(E_{s({\cal G})}({\cal G}))>\zeta(E_{s({\cal H})}({\cal H}))$,
$j_*>\tilde{j}_*$, $\tilde{j}_*$ is odd and
$a_{i_*,j_*}({\cal G})a_{\tilde{i}_{s({\cal H})},\tilde{j}_{s({\cal H})}}({\cal H})$
$a_{i_{s({\cal G})},j_{s({\cal G})}}({\cal G})a_{\tilde{i}_*,\tilde{j}_*}({\cal H})<0$,

\item[(H3)]
$\zeta(E_{s({\cal G})}({\cal G}))=\zeta(E_{s({\cal H})}({\cal H}))$ and $j_*< \tilde{j}_*$,

\item[(H4)]
$\zeta(E_{s({\cal G})}({\cal G}))=\zeta(E_{s({\cal H})}({\cal H}))$,
$j_*=\tilde{j}_*$ and
$
a_{i_*,j_*}({\cal G})a_{\tilde{i}_{s({\cal H})},\tilde{j}_{s({\cal H})}}({\cal H})\ne$
$a_{i_{s({\cal G})},j_{s({\cal G})}}({\cal G})$
$a_{\tilde{i}_*,\tilde{j}_*}({\cal H})$.
\end{description}
Then
conclusions of subcases {\bf (ia)}, {\bf(ib)} and {\bf (iia)} in Theorem~\ref{th:finite} are all true
and moreover,
in the subcase {\bf(iib)},
i.e., $G_0(\theta)\equiv 0$ and $C_0>0$,
either
$\eta_O=1$, ${\cal S}_O^e={\cal S}_O^h=0$ and ${\cal S}_O^p=1$
when $N({\cal H})=0$,
or
$\eta_O=0$ or $+\infty$, ${\cal S}_O^e=N({\cal H})$ and ${\cal S}_O^h={\cal S}_O^p=0$
when $N({\cal H})\ge 1$.
In particular,
$\eta_O=+\infty$ if either $N({\cal H})>1$,
or
$N({\cal H})=1$ and $a_{i_0,j_0}({\cal G})a_{i_{s({\cal G})},j_{s({\cal G})}}({\cal G})>0$,
or
$N({\cal H})=1$ and $a_{i_0,j_0}({\cal G})a_{i_{s({\cal G})},j_{s({\cal G})}}({\cal G})<0$ and $u({\cal G},\xi)>u({\cal H},\xi)-1/\xi$,
but
$\eta_O=0$ if
$N({\cal H})=1$ and $a_{i_0,j_0}({\cal G})a_{i_{s({\cal G})},j_{s({\cal G})}}({\cal G})<0$ and $u({\cal G},\xi)<u({\cal H},\xi)-1/\xi$,
where $u({\cal G},\xi)$ is defined just below \eqref{defscomp} and
$\xi$ is the slope of the edge $E_q({\cal H})$,
the only edge of ${\cal N}_{\cal H}$ such that
ordinates of its two end-points have different parities.
\label{th:finiteJ2}
\end{thm}

Similar to Theorem~\ref{th:finite},
the proof of Theorem~\ref{th:finiteJ2} also needs the techniques of desingularization and
computation of Newton polygons, which are given in sections 3 and 4 respectively.
So we still give the proof in section~5.

\begin{rmk}
{\rm
The reason why we cannot obtain a necessary and sufficient condition
for `in particular' as in Theorem~\ref{th:finite} is that
the number $\eta_O$ is not determined yet
in the critical situation
$N({\cal H})=1$,
$a_{i_0,j_0}({\cal G})a_{i_{s({\cal G})},j_{s({\cal G})}}({\cal G})<0$
and
$u({\cal G},\xi)=u({\cal H},\xi)-1/\xi$.
In this situation the result cannot be determined only by
principal parts ${\cal G}_P$ and ${\cal H}_P$,
but remainder parts ${\cal G}_R$ and ${\cal H}_R$ will be involved.
\label{Rk:th22u}
}
\end{rmk}

\begin{rmk}
{\rm
Unlike Theorem~\ref{th:finite},
we further need {\bf (H1)}-{\bf (H4)} to guarantee that
the Newton polygon of the Lie-bracket $[{\cal G},{\cal H}]_\theta$ can be determined by
Newton polygons of ${\cal G}$ and ${\cal H}$
since both of them end on the $u$-axis.
Note that
{\bf (H1)}-{\bf (H4)} are only concerning the case that
$\zeta(E_{s({\cal G})}({\cal G}))\ge \zeta(E_{s({\cal H})}({\cal H}))$.
In the opposite case,
this theorem still holds if we exchange ${\cal G}$ with ${\cal H}$
in all conditions of the theorem.
}
\label{Rk:Sym}
\end{rmk}

\begin{rmk}
{\rm
In the two-below case {\bf(J2)},
condition {\bf(S)} ensures that
the shapes of ${\cal N}_{{\cal G}'_\theta}$ and ${\cal N}_{{\cal H}'_\theta}$
are completely determined by ${\cal N}_{\cal G}$ and ${\cal N}_{\cal H}$ respectively.
If condition {\bf(S)} is not true,
that is, there are valid indices in the horizontal stripes
(see the shadow regions given in Figure~\ref{fig:GHstripe} for example)
between the last two valid indices of ${\cal N}_{\cal G}$ and ${\cal N}_{\cal H}$,
then the shapes of
${\cal N}_{{\cal G}'_\theta}$ and ${\cal N}_{{\cal H}'_\theta}$
are determined by not only
${\cal N}_{\cal G}$ and ${\cal N}_{\cal H}$
but also valid indices in the horizontal stripes,
which implies that
the qualitative properties of system~\eqref{equ:initial}
in the direction $\theta=0$ cannot be determined
by those principal parts ${\cal G}_P$ and ${\cal H}_P$ only
and the discussion in this case will be more complicated.
}
\end{rmk}

\subsection{Corollaries without parallel edges}

Theorems~\ref{th:finite} and \ref{th:finiteJ2} seem complicated with
computing the $\xi$-componential polynomials ${\cal K}_\xi({\cal G})$ and
${\cal K}_\xi({\cal H})$ for $N({\cal G})$ and $N({\cal H})$
in the process given in Remark~\ref{Rk:poly}.
Actually, this computation can be simplified in the case
\begin{description}
\item[(NP)]
none of edges of ${\cal N}_{\cal G}$ is parallel to an edge of ${\cal N}_{\cal H}$,
i.e., $\mathfrak{S}({\cal G})\cap\mathfrak{S}({\cal H})=\emptyset$,
where $\mathfrak{S}({\cal G})$ and $\mathfrak{S}({\cal H})$ are the slope sets of
${\cal N}_{\cal G}$ and ${\cal N}_{\cal H}$ respectively.
\end{description}
In this case, we simply give $N({\cal G})$ and $N({\cal H})$
in terms of the sign changes in the coefficient sequences
${\cal A}({\cal G})$, ${\cal A}_j({\cal G})$,
${\cal A}({\cal H})$ and ${\cal A}_j({\cal H})$
given in \eqref{defAA}.
For a short statement,
let $(A,B)$ present the {\it conjunctive sequence} $(a_1,..., a_n, b_1,..., b_m)$ of the
two sequences $A:=(a_1,...,a_n)$ and $B:=(b_1,..., b_m)$, and
let $\Xi(A)$ denote the number of sign changes in a sequence $A$ not containing zero.
For example, $\Xi(A)=2$ if $A=(-1,-2, 1,1, 5,-3)$.

\begin{cor}
In the case {\bf (NP)}, under the same conditions of Theorem~\ref{th:finite} or Theorem~\ref{th:finiteJ2},
$N({\cal G})=\Xi({\cal A}({\cal G}),{\cal A}_j({\cal G}))$ and
$N({\cal H})=\Xi({\cal A}({\cal H}),{\cal A}_j({\cal H}))$.
\label{cor:nonpara}
\end{cor}


{\bf Proof.}
In the case {\bf(NP)},
we claim that for each $\xi\in\mathfrak{S}({\cal G})$
the polynomial ${\cal K}_\xi({\cal G})$
has at most one positive root and at most one negative root.
For a reduction to absurdity,
let $\theta_1<\theta_2$ be
two adjacent positive roots of ${\cal K}_\xi({\cal G})$.
By condition {\bf(Q)}, the Lie-bracket
$[{\cal K}_\xi({\cal G}),{\cal K}_\xi({\cal H})]$
has no nonzero real roots,
which implies that both $\theta_1$ and $\theta_2$ are simple roots of
${\cal K}_\xi({\cal G})$;
otherwise, one of them is a common real root of ${\cal K}_{\xi}({\cal G})$ and
its derivative.
Additionally,
\begin{eqnarray}
({\cal K}_\xi({\cal G}))'_\theta(\theta_1)({\cal K}_\xi({\cal G}))'_\theta(\theta_2)<0
\label{kgt1kht2}
\end{eqnarray}
because $\theta_1$ and $\theta_2$ are two adjacent positive roots (simple).
On the other hand,
condition {\bf(NP)} ensures that
$\xi\in \mathfrak{S}({\cal G})\backslash\mathfrak{S}({\cal H})$.
For ${\cal H}$, it implies that
the $\xi$-component of ${\cal N}_{\cal H}$ takes the second option in \eqref{defscomp}.
Hence, the summation in \eqref{defkxip} has exactly one term and
${\cal K}_\xi({\cal H})$ is a monomial of the form $\alpha \theta^\ell$,
where $\alpha\ne 0$ and $\ell$ is a positive integer.
Then
$$
[{\cal K}_\xi({\cal G}),{\cal K}_\xi({\cal H})]
=\alpha\theta^{\ell-1}F(\theta),
$$
where
$F(\theta):=({\cal K}_\xi({\cal G}))'_\theta(\theta)\theta
-\ell{\cal K}_\xi({\cal G})(\theta)$.
It follows from \eqref{kgt1kht2} that $F(\theta_1)F(\theta_2)<0$,
which implies by the continuity that
the Lie-bracket
$[{\cal K}_\xi({\cal G}),{\cal K}_\xi({\cal H})]$
has a positive root between $\theta_1$ and $\theta_2$,
a contradiction to {\bf(Q)}.
This proves the claim with positive roots.
We can similarly prove the claim with negative roots.
In addition,
for the same reason as $\theta_1$ and $\theta_2$ being simple roots
as mentioned before \eqref{kgt1kht2},
the unique positive (and negative) root of ${\cal K}_\xi({\cal G})$ is simple if exists.

Let $\xi\in\mathfrak{S}({\cal G})$ be the slope of the $k$-th edge $E_k({\cal G})$,
linking vertex $V_{k-1}({\cal G}):(i_{k-1},j_{k-1})$ with vertex $V_k({\cal G}):(i_k,j_k)$,
where $k=1,...,s({\cal G})$.
Then, by definitions \eqref{defscomp} and \eqref{defkxip},
$$
{\cal K}_\xi({\cal G})(\theta)=
\{a_{i_{k-1},j_{k-1}}({\cal G})\theta^{j_{k-1}-j_k}+\cdots
+a_{i_k,j_k}({\cal G})\}\theta^{j_k}.
$$
By continuity, ${\cal K}_\xi({\cal G})$ has a positive root if
the coefficient $a_{i_{k-1},j_{k-1}}({\cal G})$ of the highest degree and
the coefficient $a_{i_k,j_k}({\cal G})$ of the lowest degree have different signs.
The above claim further ensures that
${\cal K}_\xi({\cal G})$ has one positive root and no positive roots
if
$$
\Xi(a_{i_{k-1},j_{k-1}}({\cal G}),a_{i_k,j_k}({\cal G}))=1~\mbox{and}~0
$$
respectively,
where $\Xi$ is defined just before Corollary~\ref{cor:nonpara}.
Similarly,
${\cal K}_\xi({\cal G})$
has one negative root and no negative roots if
$$
\Xi((-1)^{j_k}a_{i_k,j_k}({\cal G}),(-1)^{j_{k-1}}a_{i_{k-1},j_{k-1}}({\cal G}))
=1~\mbox{and}~0
$$
respectively.
Then
$\chi({\cal G})=\Xi({\cal A}({\cal G}))+\Xi({\cal A}_j({\cal G}))$
by definition~\eqref{defChi}.
It follows that
\begin{align*}
N({\cal G})
&=\Xi(a_{i_0,j_0}({\cal G}),...,a_{i_{s({\cal G})},j_{s({\cal G})}}({\cal G}),
(-1)^{j_{s({\cal G})}}a_{i_{s({\cal G})},j_{s({\cal G})}}({\cal G}),...,
(-1)^{j_0}a_{i_0,j_0}({\cal G}))
\\
&=\Xi({\cal A}({\cal G}),{\cal A}_j({\cal G})).
\end{align*}
For the same reason,
$N({\cal H})=\Xi({\cal A}({\cal H}),{\cal A}_j({\cal H}))$.
Thus, the proof of this corollary is completed.
\qquad$\Box$


Although the computation of $N({\cal G})$ and $N({\cal H})$
can be simplified in Corollary~\ref{cor:nonpara},
there are still difficulties in checking condition {\bf (Q)},
where we need to determine the number of nonzero real roots of the Lie-bracket
$[{\cal K}_{\xi_k}({\cal G}),{\cal K}_{\xi_k}({\cal H})]$.
In what follows, we additionally consider the restriction
\begin{description}
\item[(T)]
${\cal N}_{\cal G}\cap\Delta({\cal G})=\Delta^V({\cal G})$ and
${\cal N}_{{\cal H}}\cap\Delta({\cal H})=\Delta^V({\cal H})$,
\end{description}
i.e.,
there are no valid indices on Newton polygons
${\cal N}_{\cal G}$ and ${\cal N}_{\cal H}$ except for vertices.
Then, condition {\bf (Q)} can be realized by signs of some constants.
For this purpose,
we need a sequence composed of those vertices
$V_0({\cal G}),..., V_{s({\cal G})}({\cal G})$,
which is defined by $D_0({\cal G}):=V_0({\cal G})$ and
\begin{equation}
D_k({\cal G}):=V_{n(k)}({\cal G})~~~\mbox{if}~~~
\xi_k \in [\zeta^-(V_{n(k)}({\cal G})),\zeta^-(V_{n(k)}({\cal G})))
\label{Dkf}
\end{equation}
for all $k=1,...,s({\cal G},{\cal H})$,
where $n(k)\in\{1,...,s({\cal G})\}$ and
$\xi_k$ is given in \eqref{unionseq}.
Similarly, we define $D_0({\cal H}),...,D_{s({\cal G},{\cal H})}({\cal H})$.
It follows that
$C_0$, defined just before Theorem \ref{th:finite},
is equal to
$$
(j(D_0({\cal G}))-j(D_0({\cal H})))a_{D_0({\cal G})}a_{D_0({\cal H})},
$$
where $j(D_0({\cal G}))$ denotes the ordinate of the vertex $D_0({\cal G})$
and $a_{D_0({\cal G})}$ simply denotes the coefficient of the function ${\cal G}$
corresponding to the vertex $D_0({\cal G})$.
We further define
\begin{equation}
C_k:=(j(D_k({\cal G}))-j(D_k({\cal H})))a_{D_k({\cal G})}a_{D_k({\cal H})}
\label{defCk}
\end{equation}
for all $k=1,...,s({\cal G},{\cal H})$.
Clearly,
these $C_k$\,s can be determined easily by
ordinates of vertices of ${\cal N}_{\cal G}$ and ${\cal N}_{\cal H}$
and coefficients of those vertices.
These $C_k$\,s can be used to simplify Corollary~\ref{cor:nonpara}.

\begin{cor}
Under the assumptions {\bf(NP)} and {\bf(T)},
results of Corollary~\ref{cor:nonpara} hold if either
\\
{\bf (i)}
condition {\bf(P1)}  holds
and $C_0,...,C_{s({\cal G},{\cal H})}$ have the same sign in the case {\bf (J1)},
\\
or
\\
{\bf (ii)}
conditions {\bf(P2)}, {\bf(S)} and {\bf(H1)} hold and
$C_0,...,C_{s({\cal G},{\cal H})-1}$ have the same sign in the case {\bf (J2)}.
\label{cor:sign}
\end{cor}

\begin{rmk}
{\rm
Corresponding to Remark~\ref{Rk:Sym},
Corollary~\ref{cor:sign}{\bf(ii)} still holds if
we exchange the positions of ${\cal G}$ and ${\cal H}$ in condition {\bf (P2)},
and replace the condition
{\bf(H1)} $\zeta(E_{s({\cal G})}({\cal G}))>\zeta(E_{s({\cal H})}({\cal H}))$ and $j_*\le \tilde{j}_*$
with the condition that
$\zeta(E_{s({\cal H})}({\cal H}))>\zeta(E_{s({\cal G})}({\cal G}))$ and $\tilde{j}_*\le j_*$
in the case {\bf(J2)}.
}
\label{Rk:corPS}
\end{rmk}

{\bf Proof of Corollary~\ref{cor:sign}.}
Comparing the conditions given in this corollary with the conditions
given in Theorems \ref{th:finite} and \ref{th:finiteJ2},
we see that it suffices to verify condition {\bf(Q)},
given just before Theorem~\ref{th:finite}, in cases {\bf(i)} and {\bf(ii)}.
In what follows,
we only prove {\bf(Q)} in case {\bf(ii)}
because the proof in case {\bf(i)} is similar and simpler.

Note that
condition {\bf(NP)} implies that
$\mathfrak{S}({\cal G})\cap\mathfrak{S}({\cal H})=\emptyset$.
It follows that
either {\bf(iia)} $\xi_k\in\mathfrak{S}({\cal G})\backslash\mathfrak{S}({\cal H})$,
or {\bf(iib)} $\xi_k\in\mathfrak{S}({\cal H})\backslash\mathfrak{S}({\cal G})$,
for each $k=1,...,s({\cal G},{\cal H})$.
In the case {\bf(iia)},
$\xi_k\in\mathfrak{S}({\cal G})$.
By \eqref{Dkf}, $D_k({\cal G})=V_{n(k)}({\cal G})$
for an integer $n(k)$ such that $\xi_k=\zeta(E_{n(k)}({\cal G}))$.
We see from \eqref{unionseq} that
$$
\zeta^-(V_{n(k)-1}({\cal G}))=\zeta(E_{n(k)-1}({\cal G}))
\le \xi_{k-1}
<\xi_k=\zeta(E_{n(k)}({\cal G}))=\zeta^+(V_{n(k)-1}({\cal G})).
$$
By \eqref{Dkf},
we have $D_{k-1}({\cal G})=V_{n(k)-1}({\cal G})$.
Thus $D_{k-1}({\cal G})$ and $D_k({\cal G})$ are two adjacent vertices and
$\zeta(\underline{D_{k-1}({\cal G})D_k({\cal G})})=\xi_k$.
Furthermore,
we see from definition~\eqref{defscomp} of $\xi_k$-component that
$\sigma({\cal G},\xi_k)=\underline{D_{k-1}({\cal G})D_k({\cal G})}$.
Noting that condition {\bf(T)} implies
${\cal N}_{\cal G}\cap\Delta({\cal G})=\Delta^V({\cal G})$,
by \eqref{defkxip} we get
\begin{eqnarray}
{\cal K}_{\xi_k}({\cal G})(\theta)=
a_{D_{k-1}({\cal G})}\theta^{j(D_{k-1}({\cal G}))}
+a_{D_k({\cal G})}\theta^{j(D_k({\cal G}))}.
\label{Z-G}
\end{eqnarray}
On the other hand,
$\xi_k\not\in\mathfrak{S}({\cal H})$ in {\bf(iia)}.
By \eqref{Dkf},
$\xi_k\in(\zeta^-(D_k({\cal H})),\zeta^+(D_k({\cal H})))$.
Since $\xi_{k-1}<\xi_k$,
we see from \eqref{unionseq} that
$\xi_{k-1}\in[\zeta^-(D_k({\cal H})),\zeta^+(D_k({\cal H})))$,
implying that $D_{k-1}({\cal H})=D_k({\cal H})$.
By definition~\eqref{defscomp} of $\xi_k$-component,
$\sigma({\cal H},\xi_k)=D_{k-1}({\cal H})=D_k({\cal H})$.
By \eqref{defkxip}, we get
\begin{eqnarray}
{\cal K}_{\xi_k}({\cal H})(\theta)=
a_{D_{k-1}({\cal H})}\theta^{j(D_{k-1}({\cal H}))}
=a_{D_k({\cal H})}\theta^{j(D_k({\cal H}))}.
\label{Z-H}
\end{eqnarray}
Then, it follows from \eqref{Z-G} and \eqref{Z-H} that
\begin{eqnarray}
[{\cal K}_{\xi_k}({\cal G}),{\cal K}_{\xi_k}({\cal H})]
=C_{k-1}\theta^{j(D_{k-1}({\cal G}))+j(D_{k-1}({\cal H}))-1}
+C_k\theta^{j(D_k({\cal G}))+j(D_k({\cal H}))-1},
\label{Gkcor51}
\end{eqnarray}
where $C_{k-1}$ and $C_k$ are defined in \eqref{defCk}.
Note that $C_{k-1}$ and $C_k$ have the same sign, as assumed in this corollary.
In what follows, we claim that
if $k\le s({\cal G},{\cal H})-1$
then both degrees
\begin{align}
j(D_{k-1}({\cal G}))+j(D_{k-1}({\cal H}))-1
~~~\mbox{and}~~~
j(D_k({\cal G}))+j(D_k({\cal H}))-1
\label{JDJDJDJD}
\end{align}
of the two terms in \eqref{Gkcor51} are even.
In fact,
$D_{k-1}({\cal G})$ and $D_k({\cal G})$ are both vertices of
${\cal N}_{\cal G}$.
Then
$j(D_{k-1}({\cal G})), j(D_{k}({\cal G}))\in\{j_0,...,j_{s({\cal G})}\}$.
Condition {\bf (H1)} implies that
$\zeta(E_{s({\cal G})}({\cal G}))>\zeta(E_{s({\cal H})}({\cal H}))$.
By the definition of $\xi_{s({\cal G},{\cal H})}$ given in \eqref{unionseq},
$\xi_{s({\cal G},{\cal H})}=\zeta(E_{s({\cal G})}({\cal G}))$,
implying that
$\xi_k\in[\zeta(E_{\iota(k)}({\cal G})),\zeta(E_{\iota(k)+1}({\cal G})))$
for an integer $\iota(k)$ less than $s({\cal G})$ since $k\le s({\cal G},{\cal H})-1$.
By \eqref{Dkf}, we have
$j(D_{k}({\cal G}))=j(V_{\iota(k)}({\cal G}))=j_{\iota(k)}\ne j_{s({\cal G})}$.
Similarly,
$j(D_{k-1}({\cal G}))\ne j_{s({\cal G})}$
since $k-1\le s({\cal G},{\cal H})-1$.
Note that
$\mathfrak{S}({\cal G})\cap\mathfrak{S}({\cal H})=\emptyset$.
By \eqref{def-Lambda},
$\Lambda_0({\cal G})=\vec{\Delta}^V({\cal G})$.
Then we see from {\bf(P2)} that $j_0,...,j_{s({\cal G})-1}$ are all odd, implying that
$j(D_{k-1}({\cal G}))$ and $j(D_{k}({\cal G}))$ are both odd.
Similarly,
$j(D_{k-1}({\cal H})),j(D_k({\cal H}))\in \{\tilde{j}_0,...,\tilde{j}_{s({\cal H})}\}$.
Again, from the fact
$\mathfrak{S}({\cal G})\cap\mathfrak{S}({\cal H})=\emptyset$, we also have
$\Lambda_0({\cal H})=\vec{\Delta}^V({\cal H})$.
It follows from {\bf(P2)} that
$\tilde{j}_0,...,\tilde{j}_{s({\cal H})}$ are all even, implying that
$j(D_{k-1}({\cal H}))$ and $j(D_k({\cal H}))$ are both even.
Therefore,
$j(D_{k-1}({\cal G}))+j(D_{k-1}({\cal H}))-1$ and
$j(D_k({\cal G}))+j(D_k({\cal H}))-1$
are both even and the claimed \eqref{JDJDJDJD} is proved.
Thus, from \eqref{Gkcor51} we see that
the Lie-bracket
$[{\cal K}_{\xi_k}({\cal G}),{\cal K}_{\xi_k}({\cal H})]$
has no nonzero real roots,
i.e., {\bf (Q)} is verified,
if $k\le s({\cal G},{\cal H})-1$.

In particular, for $k=s({\cal G},{\cal H})$,
we see from \eqref{Dkf} that $D_k({\cal G})=V_{s({\cal G})}({\cal G})$
since $\xi_k=\zeta(E_{s({\cal G})}({\cal G}))$.
Noticing that $\zeta(E_{s({\cal G})}({\cal G}))>\zeta(E_{s({\cal H})}({\cal H}))$
in {\bf(H1)},
we have
$$
\zeta^-(V_{s({\cal H})}({\cal H}))
=\zeta(E_{s({\cal H})}({\cal H}))
<\xi_k<0=\zeta^+(V_{s({\cal H})}({\cal H})).
$$
By \eqref{Dkf}, $D_k({\cal H})=V_{s({\cal H})}({\cal H})$.
In the case {\bf(J2)},
we have $j_{s({\cal G})}=\tilde{j}_{s({\cal H})}=0$,
implying that $C_k=0$ and $[{\cal K}_{\xi_k}({\cal G}),{\cal K}_{\xi_k}({\cal H})]$
is a monomial by \eqref{Gkcor51},
which clearly has no nonzero real roots,
i.e., {\bf (Q)} is verified,
if $k=s({\cal G},{\cal H})$.

In the above,
we proved {\bf (Q)} in the case {\bf(iia)}.
We can similarly prove {\bf (Q)} in the case {\bf(iib)}.
Consequently,
the results are obtained by Corollary \ref{cor:nonpara},
and the proof is completed.
\qquad$\Box$


Applying Corollary~\ref{cor:sign}{\bf(i)} to the general degenerate system
\begin{equation}
\big\{X_m(x,y)+\cdots+X_{m+n}(x,y)\big\}\frac{\partial }{\partial x}+
\big\{Y_m(x,y)+\cdots+Y_{m+n}(x,y)\big\}\frac{\partial }{\partial y}
\label{cor0ii}
\end{equation}
where $m\ge 2n^2+1$,
$n\ge 2$,
$X_k$ and $Y_k$ are homogeneous polynomials of degree $k$ for all $k=m,...,m+n$, we
obtain the following.

\begin{cor}
Let $\alpha_k:=2(n-k)(n-k-1)+1$, $\beta_k:=2(n-k-1)^2$,
$\gamma_k:=m+k-\alpha_k$ and $\delta_k:=m+k-\beta_k$,
where $n\ge 2$ and $k=0,...,n-1$.
Then in the exceptional direction $\theta=0$ of the singular point $O$ of system~\eqref{cor0ii},
we have either
{\bf(ia)}
$\eta_O=2n-1$, ${\cal S}_O^e=0$, ${\cal S}_O^h=2n-2$ and ${\cal S}_O^p=0$ if
$$
X_{m+k}(x,y)=(-1)^{k+1}x^{\delta_k}y^{\beta_k},~~~
Y_{m+k}(x,y)=(-1)^{k+1}x^{\delta_k-1}y^{\beta_k+1}+(-1)^kx^{\gamma_k}y^{\alpha_k}
$$
for all $k=0,...,n-1$, $X_{m+n}(x,y)=y^{m+n}$ and $Y_{m+n}(x,y)=0$,
or
{\bf(ib)}
$\eta_O=2n-2$, ${\cal S}_O^e=0$, ${\cal S}_O^h=2n-3$ and ${\cal S}_O^p=0$ if
$$
X_{m+k}(x,y)=(-1)^kx^{\delta_k-1}y^{\beta_k+1}+(-1)^{k+1}x^{\gamma_k}y^{\alpha_k},~~~
Y_{m+k}(x,y)=(-1)^{k+1}x^{\delta_k}y^{\beta_k}
$$
for all $k=0,...,n-1$,
$X_{m+n}(x,y)=0$ and $Y_{m+n}(x,y)=y^{m+n}$,
or
{\bf(ic)}
$\eta_O=+\infty$, ${\cal S}_O^e=2n-2$ and ${\cal S}_O^h={\cal S}_O^p=0$ if
$$
X_{m+k}(x,y)=(-1)^kx^{\delta_k}y^{\beta_k},~~~
Y_{m+k}(x,y)=(-1)^kx^{\delta_k-1}y^{\beta_k+1}+(-1)^kx^{\gamma_k}y^{\alpha_k}
$$
for all $k=0,...,n-1$,
$X_{m+n}(x,y)=y^{m+n}$ and $Y_{m+n}(x,y)=0$,
or {\bf(id)}
$\eta_O=+\infty$, ${\cal S}_O^e=2n-3$ and ${\cal S}_O^h={\cal S}_O^p=0$ if
$$
X_{m+k}(x,y)=(-1)^kx^{\delta_k-1}y^{\beta_k+1}+(-1)^kx^{\gamma_k}y^{\alpha_k},~~~
Y_{m+k}(x,y)=(-1)^{k+1}x^{\delta_k}y^{\beta_k}
$$
for all $k=0,...,n-1$,
$X_{m+n}(x,y)=0$ and $Y_{m+n}(x,y)=y^{m+n}$,
or {\bf(iia)}
$\eta_O=2n-1$, ${\cal S}_O^e=0$, ${\cal S}_O^h=2n-2$ and ${\cal S}_O^p=2$ if
$X_m(x,y)=x^{\delta_0}y^{\beta_0}$,
$Y_m(x,y)=x^{\delta_0-1}y^{\beta_0+1}$,
$$
X_{m+k}(x,y)=(-1)^kx^{\delta_k}y^{\beta_k},~~~
Y_{m+k}(x,y)=(-1)^kx^{\delta_k-1}y^{\beta_k+1}+(-1)^kx^{\gamma_{k-1}+1}y^{\alpha_{k-1}}
$$
for all $k=1,...,n-1$,
$X_{m+n}(x,y)=(-1)^{n+1}x^{m+n}$ and $Y_{m+n}(x,y)=y^{m+n}$,
or {\bf(iib)}
$\eta_O=2n-2$, ${\cal S}_O^e=0$, ${\cal S}_O^h=2n-3$ and ${\cal S}_O^p=2$ if
$X_m(x,y)=x^{\gamma_1-1}y^{\alpha_1}$,
$Y_m(x,y)=x^{\gamma_1-2}y^{\alpha_1+1}$,
{\small
$$
X_{m+k}(x,y)\!=\!(-1)^kx^{\gamma_{k+1}-1}y^{\alpha_{k+1}},~
Y_{m+k}(x,y)\!=\!(-1)^kx^{\gamma_{k+1}-2}y^{\alpha_{k+1}+1}
\!+\!(-1)^kx^{\delta_{k-1}+1}y^{\beta_{k-1}}
$$
}for all $k=1,...,n-2$,
$X_{m+n-1}(x,y)=0$,
$Y_{m+n-1}(x,y)=(-1)^{n-1}x^{m+n-3}y^2$,
$X_{m+n}(x,y)=y^{m+n}$ and $Y_{m+n}(x,y)=(-1)^nx^{m+n}$,
or {\bf(iic)}
$\eta_O=+\infty$, ${\cal S}_O^e=2n-2$, ${\cal S}_O^h=0$ and ${\cal S}_O^p=0$ if
$X_m(x,y)=x^{\delta_0}y^{\beta_0}$, $Y_m(x,y)=x^{\delta_0-1}y^{\beta_0+1}$,
$$
X_{m+k}(x,y)=(-1)^kx^{\delta_k}y^{\beta_k},~~~
Y_{m+k}(x,y)=(-1)^kx^{\delta_k-1}y^{\beta_k+1}+(-1)^{k+1}x^{\gamma_{k-1}+1}y^{\alpha_{k-1}}
$$
for all $k=1,...,n-1$, $X_{m+n}(x,y)=(-1)^nx^{m+n}$ and $Y_{m+n}(x,y)=y^{m+n}$,
or {\bf(iid)}
$\eta_O=+\infty$, ${\cal S}_O^e=2n-1$, ${\cal S}_O^h=0$ and ${\cal S}_O^p=0$ if
$X_m(x,y)=x^{\gamma_0}y^{\alpha_0}$,
$Y_m(x,y)=x^{\gamma_0-1}y^{\alpha_0+1}$,
$$
X_{m+k}(x,y)=(-1)^kx^{\gamma_k}y^{\alpha_k},~~~
Y_{m+k}(x,y)=(-1)^kx^{\gamma_k-1}y^{\alpha_k+1}+(-1)^kx^{\delta_{k-1}+1}y^{\beta_{k-1}}
$$
for all $k=1,...,n-1$, $X_{m+n}(x,y)=y^{m+n}$ and $Y_{m+n}(x,y)=(-1)^nx^{m+n}$.
\label{cor-n}
\end{cor}

\begin{rmk}
{\rm
Corollary~\ref{cor-n} contains two classes of results: {\bf(ia)}-{\bf(id)} and {\bf(iia)}-{\bf(iid)}.
In one of cases {\bf(ia)}-{\bf(id)},
$\theta=0$ is an isolated exceptional direction, i.e.,
$G_0(0)=0$ but $G_0(\theta)\not\equiv 0$.
In one of cases {\bf(iia)}-{\bf(iid)},
$\theta=0$ is a non-isolated exceptional direction, i.e.,
$G_0(\theta)\equiv 0$.
}
\label{Rk:cori-ii}
\end{rmk}

{\bf Proof of Corollary~\ref{cor-n}.}
We only check the result {\bf(ia)}, the remaining results can be similarly discussed.
It suffices to verify conditions in Corollary~\ref{cor:sign}{\bf(i)}.
As defined just below \eqref{equ:polar system},
we compute that
$G_k(\theta)=(-1)^k\theta^{\alpha_k}+O(\theta^{\alpha_k+1})$ and
$H_k(\theta)=(-1)^{k+1}\theta^{\beta_k}+O(\theta^{\beta_k+1})$ for each $k=0,...,n-1$,
$G_n(\theta)=O(\theta^{m+n+1})$ and $H_n(\theta)=O(\theta^{m+n})$.
Then,
$$
\vec{\Delta}^V({\cal G})=(V_0({\cal G}),...,V_{n-1}({\cal G})),~~~
\vec{\Delta}^V({\cal H})=(V_0({\cal H}),...,V_{n-1}({\cal H})),
$$
where $V_k({\cal G}):(k,\alpha_k)$ and $V_k({\cal H}):(k,\beta_k)$ for all $k=0,...,n-1$.
Clearly, condition {\bf(J1)} holds because $\alpha_{n-1}+\beta_{n-1}=1>0$.
Moreover,
for each $k=1,...,n-1$,
the $k$-th edges of ${\cal N}_{\cal G}$ and ${\cal N}_{\cal H}$
have slopes $-4(n-k)$ and $-4(n-k)+2$ respectively.
Thus,
\begin{eqnarray}
\zeta(E_1({\cal G}))<\zeta(E_1({\cal H}))
<\cdots<
\zeta(E_{n-1}({\cal G}))<\zeta(E_{n-1}({\cal H})),
\label{zEGH}
\end{eqnarray}
and therefore,
$$
\sharp\vec{\mathfrak{S}}({\cal G})\cup\vec{\mathfrak{S}}({\cal H})=2n-2
~~~\mbox{and}~~~
\sharp\vec{\mathfrak{S}}({\cal G})\cap\vec{\mathfrak{S}}({\cal H})=0.
$$
Since the intersection sequence is an empty set,
we see from definition \eqref{def-Lambda} that
$\Lambda_0({\cal G})=(V_0({\cal G}),...,V_{n-1}({\cal G}))$ and
$\Lambda_0({\cal H})=(V_0({\cal H}),...,V_{n-1}({\cal H}))$.
Then condition {\bf (P1)} holds because
the ordinate $\alpha_k$ of $V_k({\cal G})$ is odd and
the ordinate $\beta_k$ of $V_k({\cal H})$ is even for all $k=0,...,n-1$.
It is clear that condition {\bf(NP)} holds by \eqref{zEGH}.
Moreover, we see from abscissas of vertices of
${\cal N}_{\cal G}$ and ${\cal N}_{\cal H}$ that
${\cal N}_{\cal G}\cap\Delta({\cal G})=\Delta^V({\cal G})$ and
${\cal N}_{{\cal H}}\cap\Delta({\cal H})=\Delta^V({\cal H})$,
i.e., condition {\bf(T)} holds.

According to the slope order \eqref{zEGH} and the cardinality of the union sequence,
we see from definition \eqref{Dkf} that
$D_{2k}({\cal G})=V_k({\cal G})$ and
$D_{2k}({\cal H})=V_k({\cal H})$
for each $k=0,...,n-1$,
and $D_{2k+1}({\cal G})=V_{k+1}({\cal G})$ and
$D_{2k+1}({\cal H})=V_k({\cal H})$
for each $k=0,...,n-2$.
Further, by definition \eqref{defCk},
for each $k=0,...,n-1$,
$$
C_{2k}
=(j(D_{2k}({\cal G}))-j(D_{2k}({\cal H})))a_{D_{2k}({\cal G})}a_{D_{2k}({\cal H})}
=(\alpha_k-\beta_k)(-1)^{2k+1}<0,
$$
and for each $k=0,...,n-2$,
$$
C_{2k+1}
=(j(D_{2k+1}({\cal G}))-j(D_{2k+1}({\cal H})))a_{D_{2k+1}({\cal G})}a_{D_{2k+1}({\cal H})}
=(\alpha_{k+1}-\beta_k)(-1)^{2k+2}<0.
$$
Consequently,
conditions in Corollary~\ref{cor:sign}{\bf(i)} are verified
and therefore, Theorem~\ref{th:finite} holds.
Note that
$
C_0=(j_0-\tilde{j}_0)a_{i_0,j_0}({\cal G})a_{\tilde{i}_0,\tilde{j}_0}({\cal H})
=1-2n<0
$ and
$G_0(\theta)\not\equiv0$.
Then, system \eqref{cor0ii} has $N({\cal G})$ orbits connecting with $O$
in the direction $\theta=0$,
no e-tsectors and p-tsectors, and $\max\{0,N({\cal G})-1\}$ h-tsectors in this direction.
Since condition {\bf(NP)} holds,
by corollary~\ref{cor:nonpara} and the fact that
ordinates of vertices of ${\cal N}_{\cal G}$ are all odd,
$$
N({\cal G})
=\Xi({\cal A}({\cal G}),{\cal A}_j({\cal G}))
=\Xi\big((-1)^0,...,(-1)^{n-1},-(-1)^{n-1},...,-(-1)^0\big)
=2n-1.
$$
Thus, the proof of this corollary is completed.
\qquad$\Box$

Applying Corollary~\ref{cor:sign}{\bf(ii)} to the general degenerate system~\eqref{cor0ii},
we similarly obtain the following.

\begin{cor}
Let
$\alpha_k:=2(n-k)(n-k-1)+1$, $\beta_k:=2(n-k)^2$,
$\gamma_k:=m+k-\alpha_k$ and $\delta_k:=m+k-\beta_k$,
where $n\ge 2$ and $k=0,...,n-1$.
Then in the direction $\theta=0$ of the singular point $O$ of system~\eqref{cor0ii},
we have either
{\bf(ia)}
$\eta_O=2n-1$, ${\cal S}_O^e=0$, ${\cal S}_O^h=2n-2$ and ${\cal S}_O^p=0$ if
$$
X_{m+k}(x,y)=(-1)^{k+1}x^{\gamma_k-1}y^{\alpha_k+1}+(-1)^kx^{\delta_k}y^{\beta_k},~~~
Y_{m+k}(x,y)=(-1)^kx^{\gamma_k}y^{\alpha_k}
$$
for all $k=0,...,n-1$,
$X_{m+n}(x,y)=(-1)^nx^{m+n}+(-1)^{n+1}x^{m+n-1}y$ and
$Y_{m+n}(x,y)=(-1)^nx^{m+n}+y^{m+n}$,
or
{\bf(ib)}
$\eta_O=2n-2$, ${\cal S}_O^e=0$, ${\cal S}_O^h=2n-3$ and ${\cal S}_O^p=0$ if
$X_m(x,y)=-x^{\delta_0-1}y^{\beta_0+1}$, $Y_m(x,y)=x^{\delta_0}y^{\beta_0}$,
$$
X_{m+k}(x,y)=(-1)^kx^{\gamma_k}y^{\alpha_k},~
Y_{m+k}(x,y)=(-1)^kx^{\gamma_k-1}y^{\alpha_k+1}+(-1)^{k-1}x^{\delta_k}y^{\beta_k}
$$
for all $k=1,...,n-1$,
$X_{m+n}(x,y)=(-1)^nx^{m+n}+y^{m+n}$ and $Y_{m+n}(x,y)=(-1)^{n-1}$ $x^{m+n}+(-1)^nx^{m+n-1}y$,
or
{\bf(ic)}
$\eta_O=+\infty$, ${\cal S}_O^e=2n-2$, ${\cal S}_O^h=0$ and ${\cal S}_O^p=0$ if
$$
X_{m+k}(x,y)=(-1)^kx^{\gamma_k-1}y^{\alpha_k+1}+(-1)^kx^{\delta_k}y^{\beta_k},~~~
Y_{m+k}(x,y)=(-1)^{k+1}x^{\gamma_k}y^{\alpha_k}
$$
for all $k=0,...,n-1$,
$X_{m+n}(x,y)=(-1)^nx^{m+n}+(-1)^nx^{m+n-1}y$ and
$Y_{m+n}(x,y)=(-1)^{n+1}x^{m+n}+y^{m+n}$,
or {\bf(id)}
$\eta_O=+\infty$, ${\cal S}_O^e=2n-3$, ${\cal S}_O^h=0$ and ${\cal S}_O^p=0$ if
$X_m(x,y)=-x^{\delta_0-1}y^{\beta_0+1}$, $Y_m(x,y)=x^{\delta_0}y^{\beta_0}$,
$$
X_{m+k}(x,y)=(-1)^{k-1}x^{\gamma_k}y^{\alpha_k},~
Y_{m+k}(x,y)=(-1)^{k-1}x^{\gamma_k-1}y^{\alpha_k+1}+(-1)^{k-1}x^{\delta_k}y^{\beta_k}
$$
for all $k=1,...,n-1$,
$X_{m+n}(x,y)=(-1)^{n-1}x^{m+n}+y^{m+n}$ and $Y_{m+n}(x,y)=(-1)^{n-1}$ $x^{m+n}+(-1)^{n-1}x^{m+n-1}y$,
or {\bf(iia)}
$\eta_O=2n-1$, ${\cal S}_O^e=0$, ${\cal S}_O^h=2n-2$ and ${\cal S}_O^p=2$ if
$X_m(x,y)=x^{\delta_0}y^{\beta_0}$,
$Y_m(x,y)=x^{\delta_0-1}y^{\beta_0+1}$,
$$
X_{m+k}(x,y)=(-1)^{k-1}x^{\delta_k}y^{\beta_k}+(-1)^kx^{\gamma_k-1}y^{\alpha_k+1},~~~
Y_{m+k}(x,y)=(-1)^{k-1}x^{\gamma_k}y^{\alpha_k}
$$
for all $k=1,...,n-1$,
$X_{m+n}(x,y)=(-1)^{n-1}x^{m+n}+(-1)^nx^{m+n-1}y$ and
$Y_{m+n}(x,y)=(-1)^{n-1}x^{m+n}+y^{m+n}$,
or {\bf(iib)}
$\eta_O=2n-2$, ${\cal S}_O^e=0$, ${\cal S}_O^h=2n-3$ and ${\cal S}_O^p=2$ if
$X_m(x,y)=-x^{\gamma_1-1}y^{\alpha_1}$,
$Y_m(x,y)=-x^{\gamma_1-2}y^{\alpha_1+1}$,
{\small
$$
X_{m+k}(x,y)\!=\!(-1)^{k+1}x^{\gamma_{k+1}-1}y^{\alpha_{k+1}},~~~
Y_{m+k}(x,y)\!=\!(-1)^{k-1}x^{\delta_k}y^{\beta_k}
\!+\!(-1)^{k+1}x^{\gamma_{k+1}-2}y^{\alpha_{k+1}+1}
$$
}for all $k=1,...,n-2$,
$X_{m+n-1}(x,y)=(-1)^nx^{m+n-1}$,
$Y_{m+n-1}(x,y)=(-1)^nx^{m+n-3}y^2+(-1)^nx^{m+n-2}y$,
$X_{m+n}(x,y)=y^{m+n}$ and $Y_{m+n}(x,y)=(-1)^{n-1}x^{m+n}$,
or {\bf(iic)}
$\eta_O=+\infty$, ${\cal S}_O^e=2n-2$, ${\cal S}_O^h=0$ and ${\cal S}_O^p=0$ if
$X_m(x,y)=x^{\delta_0}y^{\beta_0}$, $Y_m(x,y):=x^{\delta_0-1}y^{\beta_0+1}$,
$$
X_{m+k}(x,y)=(-1)^{k-1}x^{\delta_k}y^{\beta_k}+(-1)^{k-1}x^{\gamma_k-1}y^{\alpha_k+1},~~~
Y_{m+k}(x,y)=(-1)^kx^{\gamma_k}y^{\alpha_k}
$$
for all $k=1,...,n-1$,
$X_{m+n}(x,y)=(-1)^{n-1}x^{m+n}+(-1)^{n-1}x^{m+n-1}y$ and
$Y_{m+n}(x,y)=(-1)^nx^{m+n}+y^{m+n}$,
or {\bf(iid)}
$\eta_O=+\infty$, ${\cal S}_O^e=2n-1$, ${\cal S}_O^h=0$ and ${\cal S}_O^p=0$ if
$X_m(x,y)=x^{\gamma_0}y^{\alpha_0}$,
$Y_m(x,y)=x^{\gamma_0-1}y^{\alpha_0+1}$,
$$
X_{m+k}(x,y)=(-1)^kx^{\gamma_k}y^{\alpha_k},~~~
Y_{m+k}(x,y)=(-1)^kx^{\gamma_k-1}y^{\alpha_k+1}+(-1)^kx^{\delta_k}y^{\beta_k}
$$
for all $k=1,...,n-1$,
$X_{m+n}(x,y)=(-1)^nx^{m+n}+y^{m+n}$ and
$Y_{m+n}(x,y)=(-1)^nx^{m+n}+(-1)^nx^{m+n-1}y$.
\label{cor-n2}
\end{cor}

Opposite to the above case {\bf (NP)},
the Newton polygon ${\cal N}_{\cal G}$ has an edge parallel to an edge of
${\cal N}_{\cal H}$,
called the case of {\it parallel edges}.
In this case, we see from the proof of Corollary~\ref{cor:nonpara} that
the polynomial ${\cal K}_\xi({\cal G})$ (or ${\cal K}_\xi({\cal H})$)
may have more than one positive or negative roots.
There will be more difficulties in giving a computation method for
$N({\cal G})$ and $N({\cal H})$
which is simpler than the process given in Remark~\ref{Rk:poly}.

\section{Desingularization of analytic functions}
\setcounter{equation}{0}
\setcounter{lm}{0}
\setcounter{thm}{0}
\setcounter{rmk}{0}
\setcounter{df}{0}
\setcounter{cor}{0}

In order to prove Theorems~\ref{th:finite} and \ref{th:finiteJ2},
we need to discuss tendencies of orbits in Z-sectors (\cite{SC}),
bounded by curves on which $\dot \rho=0$.
However we did not assume that $dH_0(0)/d\theta\ne 0$,
the classic Implicit Function Theorem is not applicable to
the equation ${\cal H}(\rho,\theta)=0$ near $O$.
Hence,
we apply the desingularization of analytic functions
to finding those real branches of the equation ${\cal H}(\rho,\theta)=0$,
and further determine the signs of the function ${\cal G}(\rho,\theta)$
on those real branches
and the signs of the function
$$
[{\cal G},{\cal H}]_\theta:={\cal G}'_\theta{\cal H}-{\cal G}{\cal H}'_\theta,
$$
the Lie-bracket of ${\cal G}$ and ${\cal H}$ in the variable $\theta$,
inside those Z-sectors.

\subsection{Estimate of the number of real branches}

For the analytic function ${\cal H}$ given in \eqref{expGH},
we want to find continuous solutions $\theta(\rho)$ of
the equation ${\cal H}(\rho,\theta)=0$ near $O:(0,0)$ such that $\theta(0)=0$.
For convenience,
each real continuous solution $\theta(\rho)$ defined on $(-\epsilon,0]$ or $[0,\epsilon)$
for small $\epsilon>0$ is called
a {\it real branch} of the equation ${\cal H}(\rho,\theta)=0$ passing through $O$.
Without loss of generality, we assume that
${\cal H}$ does not have a factor $\rho$
because the equation $\rho{\cal H}(\rho,\theta)=0$ and
the equation ${\cal H}(\rho,\theta)=0$ have the same real branches.
Thus,
\begin{eqnarray}
{\cal H}(0,\theta)=a_{0,\ell}({\cal H})\theta^\ell+O(\theta^{\ell+1}),
\label{H0t}
\end{eqnarray}
where $\ell\ge 1$ is an integer and $a_{0,\ell}({\cal H})\ne0$.

Since ${\cal H}$ is analytic at $O$ and satisfies \eqref{H0t},
by the Weierstrass Preparation Theorem
(see \cite[Theorem 6.1.3]{krantz} or the Appendix A.1),
there exists a unique real function $U(\rho,\theta)$, analytic at $O$,
such that $U(0,0)\ne 0$ and
\begin{eqnarray}
U(\rho,\theta){\cal H}(\rho,\theta)=W(\rho,\theta):=\theta^\ell+w_{\ell-1}(\rho)\theta^{\ell-1}+\cdots+w_0(\rho),
\label{Wei}
\end{eqnarray}
where $w_i$\,s are analytic at $\rho=0$ and vanish at $\rho = 0$.
By the Puiseux's Theorem (see \cite[Theorem 4.2.7]{krantz} or the Appendix A.1),
there are integers $d\ge1$ and $n\in\{0,...,\ell\}$ such that
either $W(u^d,\theta)=Q(u,\theta)$ if $n=0$
or
\begin{eqnarray}
W(u^d,\theta)=(\theta-R_1(u))(\theta-R_2(u))\cdots(\theta-R_n(u))Q(u,\theta)
\label{WRQ}
\end{eqnarray}
if $n\ge 1$, where
$Q(u,\theta)$ is a polynomial of degree $\ell-n$ in $\theta$
with coefficients being real analytic at $u=0$
and has no real roots when $u$ is small and nonzero,
each $R_i$ ($i=1,...,n$) is real analytic at $u=0$.
Decomposition \eqref{WRQ} implies that
on the half-plane $\rho\ge0$,
if $\theta(\rho)$ is a real branch but not identical with $0$,
by the analyticity of $R_i$\,s,
there exists an integer $b\ge1$ such that
\begin{eqnarray}
\theta(\rho)=c_b\rho^{b/d}+\sum_{i\ge b+1}c_i\rho^{i/d},~~~\rho\in[0,\epsilon),
\label{csq}
\end{eqnarray}
a convergent series in fractional powers of $\rho$,
where $c_i$\,s are real and $c_b\ne0$.

The Newton polygon can be used to construct those real branches (\cite[section~1.4]{CA00}).
For a given edge $E$ of the Newton polygon ${\cal N}_{\cal H}$,
as defined in \eqref{defkxip},
we have the {\it edge-polynomial}
\begin{eqnarray}
{\cal H}_E(\theta):=\sum_{(i,j)\in E\cap \Delta({\cal H})}a_{i,j}({\cal H})\theta^j.
\label{red}
\end{eqnarray}
Assume that the edge $E$ lines on the line $qu+pv=\varsigma$,
where $p$ and $q$ are a pair of coprime positive integers.
Further, we call the function
\begin{eqnarray}
{\cal D}_c {\cal H}(\rho,\theta):=\frac{{\cal H}(\rho^q,\rho^p(c+\theta))}{\rho^{\varsigma}}
\label{des}
\end{eqnarray}
the {\it desingularized function of ${\cal H}$ by $c$},
where $c$ is a nonzero real root of ${\cal H}_E$.
As stated in \cite[Section 1.4]{CA00},
the function ${\cal D}_c {\cal H}$ is analytic at $O$ and does not have a factor $\rho$.

For real branches of the equation ${\cal H}(\rho,\theta)=0$,
we have the following Lemma,
the first part of which is given by \cite[Lemma~1.3.1]{CA00}
and the the second one is indicated in \cite[pp.123-124]{Hans}.

\begin{lm}
Let ${\cal H}$ be given in \eqref{expGH} and satisfy \eqref{H0t}.
\\
{\bf(i)} Each nontrivial real branch of the equation ${\cal H}(\rho,\theta)=0$
on the half-plane $\rho\ge 0$ is of the form
\begin{eqnarray}
\theta(\rho)=c\rho^\iota+o(\rho^\iota),
\label{ycxio}
\end{eqnarray}
where $\iota:=-1/\zeta(E)$ for an edge $E$ of ${\cal N}_{\cal H}$
and $c$ is a nonzero real root of ${\cal H}_E$.
\\
{\bf(ii)} If there is an edge $E$ of ${\cal N}_{\cal H}$ such that
${\cal H}_E$ has a nonzero real root $c$ of odd multiplicity,
then the equation ${\cal H}(\rho,\theta)=0$ has a real branch of the form \eqref{ycxio}
with $\iota:=-1/\zeta(E)$.
Especially, if $c$ is a simple real root,
then the real branch of the form \eqref{ycxio} with $\iota:=-1/\zeta(E)$ is unique.
\label{lm-basic}
\end{lm}

By Lemma~\ref{lm-basic},
we can use the number of sign changes in the algebraic coefficient sequences,
defined in \eqref{defAA},
to estimate the number of real branches of the equation ${\cal H}(\rho,\theta)=0$.

\begin{prop}
Let ${\cal H}$ be given in \eqref{expGH} with $s({\cal H})+1$ vertices.
Then the equation ${\cal H}(\rho,\theta)=0$ has at least $\Xi({\cal A}({\cal H}))$ real branches in the interior of the first quadrant,
where ${\cal A}({\cal H})$ is the coefficient sequence of ${\cal H}$ as defined in \eqref{defAA}
and $\Xi$ denotes the number of sign changes.
\label{lm-solution by NP}
\end{prop}

{\bf Proof.}
Consider a sign change in the coefficient sequence ${\cal A}({\cal H})$.
Suppose that
$a_{\tilde{i}_{k-1},\tilde{j}_{k-1}}({\cal H})a_{\tilde{i}_k,\tilde{j}_k}({\cal H})<0$
for some $k\in\{1,...,s({\cal H})\}$.
We see from \eqref{red} that
$$
{\cal H}_{E_k({\cal H})}(\theta)=
a_{\tilde{i}_{k-1},\tilde{j}_{k-1}}({\cal H})\theta^{\tilde{j}_{k-1}}
+\cdots+
a_{\tilde{i}_k,\tilde{j}_k}({\cal H})\theta^{\tilde{j}_k}.
$$
Thus, coefficients of the highest and the lowest degree terms of ${\cal H}_{E_k({\cal H})}$
have different signs,
implying that ${\cal H}_{E_k({\cal H})}$ has a positive root $c$ of odd multiplicity.
By Lemma~\ref{lm-basic},
the equation ${\cal H}(\rho,\theta)=0$ has a real branch of the form \eqref{ycxio}
with $\iota=-1/\zeta(E_k({\cal H}))$.
Consequently, the equation ${\cal H}(\rho,\theta)=0$ has at leat
$\Xi({\cal A}({\cal H}))$ real branches
in the interior of the first quadrant,
and this proposition is proved.
\qquad$\Box$

In Proposition~\ref{lm-solution by NP},
we consider real branches of the equation ${\cal H}(\rho,\theta)=0$
only in the interior of the first quadrant.
In the remaining quadrants,
we can use the transformation $\rho\to-\rho$,
or the transformation $\rho\to-\rho,\theta\to-\theta$,
or the transformation $\theta\to-\theta$ to convert them to the first quadrant.
Thus, we can similarly estimate the number of real branches in the interior of the second,
the third and the fourth quadrants.
Moreover,
if
$a_{\tilde{i}_{s({\cal H})},\tilde{j}_{s({\cal H})}}({\cal H})
(-1)^{\tilde{j}_{s({\cal H})}}
a_{\tilde{i}_{s({\cal H})},\tilde{j}_{s({\cal H})}}({\cal H})<0$,
then $\tilde{j}_{s({\cal H})}$ is odd,
implying that $\theta(\rho)\equiv 0$ is a trivial real branch.
Thus, the equation ${\cal H}(\rho,\theta)=0$ has at least
$\Xi({\cal A}({\cal H}),{\cal A}_j({\cal H}))$ real branches on the half-plane $\rho\ge 0$.

In order to show the use of Proposition~\ref{lm-solution by NP},
let us consider the equation
$$
{\cal H}(\rho,\theta)
=\theta^8-\rho\theta^5+\rho^2\theta^4+2\rho^2\theta^3
-\rho^3\theta^2+2\rho^4\theta+\rho^5
=0.
$$
Its Newton polygon ${\cal N}_{\cal H}$ has vertices
$V_0({\cal H}):(0,8),~V_1({\cal H}):(1,5),~V_2({\cal H}):(2,3)$ and $V_3({\cal H}):(5,0)$,
which are determined by the terms $\theta^8$, $-\rho\theta^5$, $2\rho^2\theta^3$ and $\rho^5$ respectively.
Thus, by definition \eqref{defAA}, we have
${\cal A}({\cal H})=(1,-1,2,1)$ and ${\cal A}_j({\cal H})=(1,-2,1,1)$,
and therefore,
$\Xi({\cal A}({\cal H}))=2$ and $\Xi({\cal A}_j({\cal H}))=2$.
By Proposition~\ref{lm-solution by NP},
it has at least 2 real branches
in the interior of the first quadrant and the fourth quadrant respectively.

\subsection{The sign of an analytic function on a curve}

Proposition~\ref{lm-solution by NP} helps us give some continuous curves,
which will be employed to
partition the neighborhood of a degenerate singular point into several regions.
The following proposition show how we can use ${\cal N}_{\cal G}$ to determine
the sign of ${\cal G}$ on a given curve.
For convenience,
if ${\cal G}$ has a definite sign on a curve $\Upsilon$,
then let
\begin{eqnarray*}
{\rm sgn}({\cal G}|_\Upsilon):=\left\{
\begin{array}{lll}
~~1 &&  \mbox{as~} {\cal G}(\rho,\theta)>0 \mbox{~for all~}(\rho,\theta)\in\Upsilon,\\
~~0 &&  \mbox{as~} {\cal G}(\rho,\theta)=0 \mbox{~for all~}(\rho,\theta)\in\Upsilon,\\
-1  &&  \mbox{as~} {\cal G}(\rho,\theta)<0 \mbox{~for all~}(\rho,\theta)\in\Upsilon.
\end{array}
\right.
\end{eqnarray*}
Assume that
$\Upsilon: \theta=\gamma(\rho),~\rho \in (0,\epsilon)$, is a continuous curve with
\begin{eqnarray}
\gamma(\rho):=\sum_{i=0}^{N} C_i\rho^{\nu_i},~~~~C_i\ne 0,~~0<\nu_0<\nu_1<\cdots,
\label{curvGama}
\end{eqnarray}
where $\epsilon$ is small,
$N\in\mathbb{Z}_+\cup\{+\infty\}$ and $\nu_i\in(0,+\infty)$.
The following proposition
states that ${\cal G}$ has a definite sign on the curve $\Upsilon$
and helps us determine the sign.

\begin{prop}
Let ${\cal G}$ be given in \eqref{expGH} with $s({\cal G})+1$ vertices and
$\Upsilon$ be a curve given by \eqref{curvGama}.
Then ${\cal G}$ has a definite sign on the curve $\Upsilon$ for small $\epsilon$.
Additionally,
\\
{\bf (i)}
in the case that $-1/\nu_0$ is not a slope of any edge of ${\cal N}_{\cal G}$,
i.e., $\zeta^-(V_k({\cal G}))<-1/\nu_0<\zeta^+(V_k({\cal G}))$ for some
$k\in\{0,1,...,s({\cal G})\}$,
where $\zeta^{\pm}(V_k({\cal G}))$ are the left-sided slope and right-sided slope
of the vertex $V_k({\cal G})$ as defined in \eqref{SP+-},
$$
{\rm sgn}({\cal G}|_\Upsilon)= {\rm sgn}(a_{i_k,j_k}({\cal G})C_0^{j_k});
$$
{\bf (ii)}
in the case that $-1/\nu_0$ is the slope of an edge $E$ of ${\cal N}_{\cal G}$,
$$
{\rm sgn}({\cal G}|_\Upsilon)=
\left\{
\begin{aligned}
&{\rm sgn}({\cal G}_E(C_0))                       && \mbox{as~} {\cal G}_E(C_0)\ne 0,\\
&{\rm sgn}({\cal D}_{C_0} {\cal G}|_{\Upsilon_1}) && \mbox{as~} {\cal G}_E(C_0)=0,
\end{aligned}
\right.
$$
where ${\cal G}_E$ is the $E$-polynomial of ${\cal G}$,
${\cal D}_{C_0} {\cal G}$ is the desingularized function of ${\cal G}$ by $C_0$,
$\Upsilon_1$ is the curve determined by $\theta=\gamma_1(\rho):=\gamma(\rho^q)/\rho^{p}-C_0$, and
$p,q$ are a pair of coprime positive integers such that $\nu_0=p/q$.
\label{lm-lowest term}
\end{prop}

{\bf Proof.}
We see from expression \eqref{csq} of a nontrivial real branch of the equation
${\cal G}(\rho,\theta)=0$ and
expression \eqref{curvGama} of $\gamma(\rho)$ that
those real branches and $\gamma(\rho)$ are all monotonic for small $\rho>0$.
Then, the curve $\Upsilon$ either intersects with
none of the finitely many real branches or
coincides with one of them for small $\rho>0$.
It follows directly that $f$ has a definite sign on $\Upsilon$.

For small $\rho>0$,
the sign of ${\cal G}(\rho,\gamma(\rho))$ is determined
by its lowest degree term.
In case {\bf(i)},
since $\zeta^-(V_k({\cal G}))<-1/\nu_0<\zeta^+(V_k({\cal G}))$,
we see from the definition of Newton polygon that
$i+\nu_0 j>i_k+\nu_0 j_k$ for all $(i,j)\in\Delta({\cal G})\backslash V_k({\cal G})$.
Then
\begin{eqnarray*}
{\cal G}(\rho,\gamma(\rho))
=a_{i_k,j_k}({\cal G})C_0^{j_k}\rho^{i_k+\nu_0 j_k}+o(\rho^{i_k+\nu_0 j_k}),
\end{eqnarray*}
and result {\bf(i)} of this proposition follows directly.

In case {\bf(ii)},
assume that the edge $E$ lies on the line $u+\nu_0 v=b$ or equivalently,
$qu+pv=qb$ since $\nu_0=p/q$.
By the definition of Newton polygon,
$i+\nu_0 j=b$ for all $(i,j)\in\Delta({\cal G})\cap E$ and
$i+\nu_0 j>b$ for all $(i,j)\in\Delta({\cal G})\backslash E$.
In the circumstance ${\cal G}_E(C_0)\ne0$,
we see from expression \eqref{red} of the $E$-polynomial that
\begin{eqnarray}
{\cal G}(\rho,\gamma(\rho))&=&
\sum_{(i,j)\in \{\Delta({\cal G})\cap E\}\cup\{\Delta({\cal G})\setminus E\}}
a_{i,j}({\cal G})\rho^i\{C_0\rho^{\nu_0}(1+o(1))\}^j
\nonumber
\\
&=&\sum_{(i,j)\in \Delta({\cal G})\cap E}a_{i,j}({\cal G})C_0^j\rho^b+o(\rho^b)
\nonumber
\\
&=&{\cal G}_E(C_0) \rho^b+o(\rho^b),
\label{eld}
\end{eqnarray}
implying that
${\rm sgn}({\cal G}|_\Upsilon)={\rm sgn}({\cal G}_E(C_0))$.
In the opposite circumstance, i.e., ${\cal G}_E(C_0)=0$,
we need to consider higher order terms in \eqref{eld}
since the lowest degree term ${\cal G}_E(C_0) \rho^b$ vanishes.
Correspondingly, consider the higher order terms of $\gamma(\rho)$
and let $\psi(\rho):=\gamma(\rho)/\rho^{\nu_0}-C_0$.
Similar to the computation of \eqref{eld},
\begin{eqnarray}
{\cal G}(\rho,\rho^{\nu_0}(C_0+\psi(\rho)))=\rho^b{\tilde {\cal G}}(\rho,\psi(\rho)),
\label{fff}
\end{eqnarray}
where
$$
{\tilde {\cal G}}(\rho,\theta)
:={\cal G}_E(C_0+\theta)
+\sum_{(i,j)\in\Delta({\cal G})\setminus E}a_{i,j}({\cal G})\rho^{i+\nu_0 j-b}(C_0+\theta)^j.
$$
Moreover, by definition \eqref{des} and the fact that $\nu_0=p/q=-1/\zeta(E)$,
$$
\tilde {\cal G}(\rho^q,\theta)
={{\cal G}(\rho^q,\rho^p(C_0+\theta))}/{\rho^{qb}}
={\cal D}_{C_0}{\cal G}(\rho,\theta).
$$
Then,
by the definition of $\gamma_1(\rho)$ and the fact that $\rho^q\in(0,\epsilon)$ for all $\rho\in(0,\epsilon)$,
\begin{align*}
{\rm sgn}({\cal G}|_\Upsilon)
&={\rm sgn}({\cal G}(\rho,\gamma(\rho)))
={\rm sgn}({\cal G}(\rho,\rho^{\nu_0}(C_0+\psi(\rho))))
={\rm sgn}({\tilde {\cal G}}(\rho,\psi(\rho)))
\\
&={\rm sgn}({\tilde {\cal G}}(\rho^q,\psi(\rho^q)))
={\rm sgn}({\cal D}_{C_0}{\cal G}(\rho,\psi(\rho^q)))
={\rm sgn}({\cal D}_{C_0} {\cal G}(\rho,\gamma_1(\rho)))
\\
&={\rm sgn}({\cal D}_{C_0} {\cal G}|_{\Upsilon_1}).
\end{align*}
Thus, the proof of this proposition is completed.
\qquad$\Box$

Clearly, it is easy to determine the sign of ${\cal G}$ on the curve
$\theta=\gamma(\rho)\equiv 0$, $\rho\in (0,\epsilon)$,
because it depends on the lowest degree term of ${\cal G}(\rho,0)$.
In case {\bf(i)} the sign of ${\cal G}$ on the curve $\Upsilon$ is determined.
In case {\bf(ii)},
we let ${\cal G}^{(0)},~\Upsilon_0,~\gamma_0(\rho),~E^{(0)},~({\cal G}_E)^{(0)},~p_0$ and $q_0$ denote ${\cal G},~\Upsilon,~\gamma(\rho),~E,~{\cal G}_E,~p$ and $q$ respectively.
If $({\cal G}_E)^{(0)}(C_0)\ne 0$,
the sign of ${\cal G}$ on the curve $\Upsilon$ is determined.
If $({\cal G}_E)^{(0)}(C_0)=0$ in case {\bf (ii)}, Proposition~\ref{lm-lowest term}
does not answer the sign of ${\cal G}^{(0)}$ on the curve $\Upsilon_0$
but converts it to the sign of a new analytic function ${\cal D}_{C_0}{\cal G}^{(0)}$, denoted by ${\cal G}^{(1)}$,
on a new continuous curve
$$
\Upsilon_1:y=\gamma_1(\rho):=\gamma_0(\rho^{q_0})/\rho^{p_0}-C_0,
$$
the same form as in \eqref{curvGama}.
For $k\ge 1$ we apply Proposition~\ref{lm-lowest term} to
determine the sign of ${\cal G}^{(k)}$ on $\Upsilon_{k}$.
If ${\cal G}^{(k)}$ and $\Upsilon_k$ satisfy conditions in case {\bf (ii)},
then we let $E^{(k)},~({\cal G}_E)^{(k)},~p_k$ and $q_k$ denote
$E,~{\cal G}_E,~p$ and $q$ correspondingly.
Moreover, if $({\cal G}_E)^{(k)}(C_k)=0$, then
$$
{\rm sgn}({\cal G}^{(k)}|_{\Upsilon_k})={\rm sgn}({\cal G}^{(k+1)}|_{\Upsilon_{k+1}}),
$$
where ${\cal G}^{(k+1)}$ is the desingularized function of ${\cal G}^{(k)}$ by $C_k$ and
$$
\Upsilon_{k+1}:~y=\gamma_{k+1}(\rho):=\gamma_k(\rho^{q_k})/\rho^{p_k}-C_k,
$$
the same form as in \eqref{curvGama}.
Thus, there are totally two possibilities:
\begin{description}
\item[{\bf (a)}]
there is an integer $k_0\ge 1$ such that ${\cal G}^{(k)}$ and $\Upsilon_{k}$ satisfy conditions
in case {\bf(ii)} and $({\cal G}_E)^{(k)}(C_{k})=0$ for all $k\le k_0-1$ but either
$\Upsilon_{k_0}$ is the curve $\theta= 0$,
or ${\cal G}^{(k_0)}$ and $\Upsilon_{k_0}$ satisfy conditions in case {\bf(i)},
or ${\cal G}^{(k_0)}$ and $\Upsilon_{k_0}$ satisfy conditions in case {\bf(ii)} and $({\cal G}_E)^{(k_0)}(C_{k_0})\ne 0$.

\item[{\bf (b)}]
${\cal G}^{(k)}$ and $\Upsilon_{k}$ satisfy conditions
in case {\bf(ii)} and $({\cal G}_E)^{(k)}(C_{k})=0$ for all $k\ge 0$.
\end{description}
When {\bf (a)} happens,
using Proposition~\ref{lm-lowest term} at most $k_0$ times
we determine the sign of ${\cal G}$ on $\Upsilon$.
When {\bf (b)} happens,
it follows from the inductive Newton-Puiseux algorithm (\cite[Section 1.4]{CA00}) that
$\theta=\gamma(\rho)$ is a real branch of the equation ${\cal G}(\rho,\theta)=0$ and therefore,
${\rm sgn}({\cal G}|_\Upsilon)=0$.
Thus,
the sign of ${\cal G}$ on $\Upsilon$ can always be determined
by Proposition~\ref{lm-lowest term}.
For instance, consider ${\cal G}(\rho,\theta):=\theta^2+\rho^2\theta-\rho^3$ and $\Upsilon:\theta=\gamma(\rho):=\rho^{3/2}-\rho^2/2$
for small $\rho>0$.
Clearly, $-1/\nu_0=-2/3$ is the slope of the only edge $E$ of ${\cal N}_{\cal G}$ and
$C_0=1$ is a real root of the edge-polynomial ${\cal G}_E(\theta)=\theta^2-1$.
Thus, by Proposition~\ref{lm-lowest term}{\bf(ii)}, we obtain
${\rm sgn}({\cal G}|_\Upsilon)={\rm sgn}({\cal D}_{C_0} {\cal G}|_{\Upsilon_1})$,
where ${\cal D}_{C_0} {\cal G}(\rho,\theta)=\rho+2\theta+\rho\theta+\theta^2$
and $\Upsilon_1$ is determined by
$\theta=\gamma_1(\rho):=\gamma(\rho^2)/\rho^3-1=-\rho/2$.
We repeat the above process to determine
the sign of ${\cal D}_{C_0} {\cal G}$ on $\Upsilon_1$.
Let ${\cal G}^{(1)}$ denote ${\cal D}_{C_0}{\cal G}$.
Similarly,
one can check that ${\cal G}^{(1)}$ and $\Upsilon_1$ satisfy conditions in case {\bf(ii)} and $({\cal G}_E)^{(1)}(C_1)=0$.
Thus, by Proposition~\ref{lm-lowest term}{\bf (ii)}, we obtain
${\rm sgn}({\cal G}^{(1)}|_{\Upsilon_1})={\rm sgn}({\cal G}^{(2)}|_{\Upsilon_2})$,
where
$
{\cal G}^{(2)}(\rho,\theta)
:={\cal D}_{C_1} {\cal G}^{(1)}(\rho,\theta)=-\rho/4+2\theta+\rho\theta^2
$
and $\Upsilon_2$ is the line $\theta= 0$.
Hence, for sufficiently small $\rho>0$,
we have
$
{\rm sgn}({\cal G}|_\Upsilon)={\rm sgn}({\cal G}^{(1)}|_{\Upsilon_1})
={\rm sgn}({\cal G}^{(2)}|_{\Upsilon_2})
={\rm sgn}(-\rho/4)<0.
$

\subsection{The sign of an analytic function on a region}

We consider the sign of a function on a curve in Proposition~\ref{lm-lowest term},
and the following proposition is devoted to the region
\begin{eqnarray*}
\Omega_\epsilon:=\{(\rho,\theta)\in\mathbb{R}^2:\, 0<\rho<\epsilon,~|\theta|<\epsilon\}
\end{eqnarray*}
for small $\epsilon>0$.
For convenience, let
\begin{eqnarray}
f:=[{\cal G},{\cal H}]_\theta={\cal G}'_\theta{\cal H}-{\cal G}{\cal H}'_\theta,
\label{fGH}
\end{eqnarray}
the Lie-bracket of ${\cal G}$ and ${\cal H}$ in the variable $\theta$.
We say that the real function $f$
is {\it semi-positive} (or {\it semi-negative})
if $f(\rho,\theta)\ge 0$ (or $\le 0$) in region $\Omega_\epsilon$ for an $\epsilon>0$.
Further, we call $f$ {\it semi-definite}
if $f$ is either semi-positive or semi-negative.

In order to use the theory of discriminants (\cite{YL}),
for a polynomial $\varphi(\theta):=a_0\theta^m+a_1\theta^{m-1}+\cdots+a_m$ ($a_0\ne 0, m\ge 1$),
matrix
$$
M(\varphi):=\left(
              \begin{array}{cccccccc}
                a_0 & a_1  &a_2     &\cdots &a_m     &       &      &       \\
                0   & ma_0 &(m-1)a_1&\cdots &a_{m-1} &       &      &       \\
                    & a_0  &a_1     &\cdots &a_{m-1} &a_m    &      &       \\
                    & 0    &ma_0    &\cdots &2a_{m-2}&a_{m-1}&      &       \\
                    &      &        &\cdots &\cdots  &\cdots &      &       \\
                    &      &        &\cdots &\cdots  &\cdots &      &       \\
                    &      &        &\cdots &a_0     &a_1    &\cdots&a_m    \\
                    &      &        &\cdots &0       &ma_0   &\cdots&a_{m-1}\\
              \end{array}
            \right),
$$
a $2m\times2m$ matrix with the coefficients of $\varphi$ on the odd rows and
the coefficients of the derivative $\varphi'_\theta$ on the even rows on scale,
is called the {\it discrimination matrix} of $\varphi$.
For each $k=1,...,m$,
let $M_k$ denote the $2k$-th order leading principal minor of $M(\varphi)$ and
call it the $k$-th {\it discriminant}.
Then let $\omega_1:={\rm sgn}(M_1)$;
for $2\le i \le m-1$
we define $\omega_i:=(-1)^{[(i-r+1)/2]}{\rm sgn}(M_r)$ if $M_i=0$ and neither
$\{j: M_j\ne0, 0< j < i\}=\emptyset$ nor $\{k: M_k\ne0, i< k \le m\}=\emptyset$,
where $r:=\max\{j: M_j\ne0, 0< j < i\}$ and
$[(i-r+1)/2]$ denotes the greatest integer being $\le(i-r+1)/2$,
otherwise we directly define $\omega_i={\rm sgn}(M_i)$;
finally, let $\omega_m={\rm sgn}(M_m)$.
Then we obtain the {\it revised sign list}
$(\omega_1,\omega_2,...,\omega_m)$, denoted by ${\rm rsl}(\varphi)$.
Further,
we use $\Phi({\rm rsl}(\varphi))$ and $\Xi({\rm rsl}(\varphi))$ to present
the number of non-vanishing members and the number of sign changes in the sequence ${\rm rsl}(\varphi)$ respectively
as in \cite{YL}.

In order to state the results conveniently, we need the following hypotheses:
\begin{description}
\item[(N1)]
  $f$ has a factor $\theta$ of odd multiplicity.

\item[(N2)]
  ${\cal N}_f$ has an edge $E$ such that $f_E$ has a nonzero real root of odd multiplicity.

\item[(Y)]
  Neither {\bf(N1)} nor {\bf(N2)} holds and ${\cal N}_f$ consists of either no edge
  {\rm(}i.e., a singleton{\rm)} or at least one edge but for each edge $E$
\begin{eqnarray}
\Phi({\rm rsl}(\widehat{f_E}))=2\,\Xi({\rm rsl}(\widehat{f_E})),
\label{Theta-Xi}
\end{eqnarray}
where $\widehat{f_E}(\theta):=f_E(\theta)/\theta^{j_E}$ and $j_E$ is the ordinate of the right end-point of $E$.

\item[(U)]
Neither {\bf(N1)} nor {\bf(N2)} holds and
  $$
  \daleth:=\{E\in{\cal E}(f): \eqref{Theta-Xi}~\text{is invalid}\}\ne\emptyset,
  $$
  where ${\cal E}(f):=\{E_1(f),...,E_{s(f)}(f)\}$, the edge set of ${\cal N}_f$.
\end{description}

\begin{prop}
Let $f$ be given by \eqref{fGH}. Then,
\\
{\bf(R1)}
$f$ is not semi-definite if either {\bf(N1)} or {\bf(N2)} holds.
\\
Otherwise, all coefficients corresponding to
vertices of ${\cal N}_f$ have the same sign, denoted by $\omega$. Moreover,\\
{\bf(R2)}
under condition {\bf(Y)}, $\omega f$ is semi-positive,\\
{\bf(R3)}
under condition {\bf(U)},
$\omega f$ is semi-positive if and only if
for each edge $E\in\daleth$ the function $\omega{\cal D}_c f$ is semi-positive
for all nonzero real roots $c$ of $f_E$.
\label{lm-region}
\end{prop}

{\bf Proof.}
Assume without loss of generality that $\rho$ is not a factor of $f$.
It follows from \eqref{Wei} and \eqref{WRQ} that $f$ is semi-definite
on the region $\Omega_\epsilon$ for small $\epsilon>0$
if and only if
the equation $f(\rho,\theta)=0$ does not have a real branch of odd multiplicity
on $\Omega_\epsilon$.
If {\bf(N1)} holds,
then the equation $f(\rho,\theta)=0$ has a trivial branch $\theta(\rho)\equiv0$
of odd multiplicity,
implying that $f$ is not semi-definite.
If {\bf(N2)} holds,
let $c$ be such a real root of $f_E$ and
consider curves $\Upsilon_{\pm}:\theta=(c\pm\delta)\rho^\iota,~\rho\in(0,\epsilon)$,
where $\delta>0$ is small and
$\iota$ is positive such that $-1/\iota$ is the slope of the edge $E$.
By Proposition~\ref{lm-lowest term}{\bf(ii)},
we have ${\rm sgn}(f|_{\Upsilon_{\pm}})={\rm sgn}(f_E(c\pm\delta))$.
Since $c$ is a root of $f_E$ of odd multiplicity,
${\rm sgn}(f_E(c+\delta)){\rm sgn}(f_E(c-\delta))=-1$.
Thus, ${\rm sgn}(f|_{\Upsilon_{+}})=-{\rm sgn}(f|_{\Upsilon_{-}})\ne0$, i.e.,
$f$ is not semi-definite,
and therefore, {\bf(R1)} holds.

In the following,
we consider that neither {\bf(N1)} nor {\bf(N2)} holds.
If the signs of the coefficients corresponding to vertices of ${\cal N}_f$ are not the same,
then there is an edge $E$ such that
the coefficients of its left end-point and the right one have different signs.
We see from \eqref{red} that
coefficients of the highest degree term of $f_E(\theta)$ and the lowest one have different signs.
Then,
$f_E$ has a nonzero real root of odd multiplicity,
i.e., {\bf(N2)} holds, a contradiction.
Therefore,
all coefficients corresponding to vertices of ${\cal N}_f$ have the same sign, denoted by $\omega$.

Under condition {\bf(Y)},
if ${\cal N}_f$ is a singleton $V_0:(0,j_0)$,
then $f(\rho,\theta)=\theta^{j_0}(a_{0,j_0}(f)+o(1))$
because $j\ge j_0$ for all $(i,j)\in\Delta(f)$.
Note that $j_0$ is even because {\bf(N1)} does not hold.
Then, $\omega f$ is semi-definite since $\omega={\rm sgn}(a_{0,j_0}(f))$.
Oppositely,
we consider the situation that
${\cal N}_f$ consists of at least one edge.
For each edge $E\in{\cal E}(f)$,
the number of distinct real roots of $\widehat{f_E}$ is equal to
$\Phi({\rm rsl}(\widehat{f_E}))-2\,\Xi({\rm rsl}(\widehat{f_E}))$
by Theorem~2.1 of \cite{YL}.
Therefore,
\eqref{Theta-Xi} implies that $f_E$ has no nonzero real roots.
Lemma~\ref{lm-basic}{\bf(i)} further implies that
the equation $f(\rho,\theta)=0$ has at most the trivial branch $\theta(\rho)\equiv0$.
Clearly, the trivial branch is of even multiplicity because {\bf(N1)} does not hold.
Consequently, $f$ is semi-definite.
Additionally,
consider a vertex $V_k:(i_k,j_k)$ of ${\cal N}_f$ and
a curve $\Upsilon:\theta=\rho^\nu,~\rho\in(0,\epsilon)$,
for a positive $\nu$ such that $-1/\nu\in(\zeta^-(V_k),\zeta^+(V_k))$.
By Proposition~\ref{lm-lowest term}{\bf(i)},
we have ${\rm sgn}(f|_\Upsilon)={\rm sgn}(a_{i_k,j_k}(f))=\omega$
and therefore, $\omega f$ is semi-positive.

Under condition {\bf(U)},
if $\omega f$ is semi-positive,
then we see from definition \eqref{des} of the desingularized function of $f$ by $c$ that
$\omega {\cal D}_c f$ is also semi-positive, which proves the necessity.
In order to prove the sufficiency,
for a reduction to absurdity, assume that
$f$ is not semi-definite.
Then,
the equation $f(\rho,\theta)=0$ has a real branch $\theta_0(\rho)$ of odd multiplicity.
Clearly, $\theta_0(\rho)\not\equiv0$ since {\bf (N1)} does not hold.
By Lemma~\ref{lm-basic}{\bf(i)},
$\theta_0(\rho)=c\rho^\iota+o(\rho^\iota)$,
where $\iota$ is rational such that $-1/\iota$ is the slope of an edge $E$ of ${\cal N}_f$
and $c$ is a nonzero real root of $f_E$.
Moreover,
we see from \eqref{des} that
the equation ${\cal D}_c f(\rho,\theta)=0$ has a real branch of the same multiplicity as $\theta_0(\rho)$.
Then,
${\cal D}_c f$ is not semi-definite,
a contradiction to the sufficient condition.
Therefore, we have proved that $f$ is semi-definite.
Furthermore,
by definition \eqref{des} of ${\cal D}_c f$,
we obtain that $\omega f$ is semi-positive since $\omega {\cal D}_c f$ is semi-positive.
Thus, the proof is completed.
\qquad$\Box$

\begin{figure}[!h]
    \centering
     \subcaptionbox{%
     }{\includegraphics[height=1in]{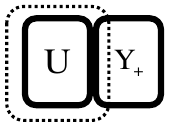}}~~~~~~~~~~~~
     \subcaptionbox{%
     }{\includegraphics[height=1in]{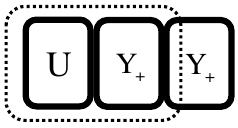}}
    \caption{Symbol sequences (a) $\Theta^{(1)}(f)$ and (b) $\Theta^{(2)}(f)$.}
    \label{fig:SS}
\end{figure}

Under condition {\bf(U)},
Proposition~\ref{lm-region} does not answer to the sign of $f$ on $\Omega_\epsilon$ yet
but converts it to the question of signs of desingularized functions ${\cal D}_c f$,
which will be discussed in the same procedure as $f$.
For convenience, let
$f^{(1)}_1, ..., f^{(1)}_{k}$ denote ${\cal D}_{c_1} f, ..., {\cal D}_{c_k} f$ respectively,
where $c_1,..., c_k$ are nonzero real roots
of edge-polynomials of $f$ for all edges of ${\cal N}_f$.
Let ${\cal D}^{(1)}(f):=(f^{(1)}_1, ..., f^{(1)}_{k})$, the sequence of the desingularized functions.
We match ${\cal D}^{(1)}(f)$ with the symbol sequence
$$
\Theta^{(1)}(f):=(\vartheta^{(1)}_1,\vartheta^{(1)}_2,...,\vartheta^{(1)}_k),
$$
where each $\vartheta^{(1)}_i\in  \{{\bf N}, {\bf Y_+}, {\bf Y_-}, {\bf U}\}$ is defined as
\begin{equation*}
\vartheta^{(1)}_i:=\left\{
\begin{array}{llllll}
{\bf N} &\mbox{if}~f_i^{(1)}~\mbox{satisfies {\bf(N1)} or {\bf(N2)} as}~f,
\\
{\bf Y_+} &\mbox{if}~f_i^{(1)}~\mbox{satisfies {\bf(Y)} as}~f~\mbox{and is semi-positive},
\\
{\bf Y_-} &\mbox{if}~f_i^{(1)}~\mbox{satisfies {\bf(Y)} as}~f~\mbox{and is  semi-negative},
\\
{\bf U} &\mbox{if}~f_i^{(1)}~\mbox{satisfies {\bf(U)} as}~f.
\end{array}
\right.
\end{equation*}
By Proposition~\ref{lm-region},
if the sequence $\Theta^{(1)}(f)$ contains a symbol ${\bf N}$ then
$f$ is not semi-definite;
if either $\Theta^{(1)}(f)=({\bf Y_+},...,{\bf Y_+})$
or $\Theta^{(1)}(f)=({\bf Y_-},...,{\bf Y_-})$ then
$f$ is semi-definite.
To the opposite,
the sequence $\Theta^{(1)}(f)$ contains at least a symbol ${\bf U}$ but no {\bf N}, that is,
$\Theta^{(1)}(f)$ contains either both ${\bf U}$ and ${\bf Y}_+$ only,
or both ${\bf U}$ and ${\bf Y}_-$ only,
or ${\bf U}$ only,
because
${\bf Y_+}$ and ${\bf Y_-}$ do not coexist in $\Theta^{(1)}(f)$ as proved in the Appendix~A.2.
For each $\vartheta^{(1)}_i$ being ${\bf U}$ in the sequence $\Theta^{(1)}(f)$,
we return to the beginning of this paragraph and
continue the procedure of desingularization for the corresponding $f^{(1)}_i$.
Let ${\cal D}^{(2)}_i(f):=(f^{(2)}_{i1},...,f^{(2)}_{i\tilde{k}})$,
where $f^{(2)}_{ir}$ simply denotes ${\cal D}_{\tilde{c}_r}f^{(1)}_{i}$,
and $\tilde{c}_1,..., \tilde{c}_{\tilde{k}}$ are nonzero real roots of edge-polynomials of $f_i^{(1)}$ for all edges of ${\cal N}_{f_i^{(1)}}$.
Similarly to the above, we match the sequence ${\cal D}_i^{(2)}(f)$
with a symbol sequence $\Theta_i^{(2)}(f_i^{(1)})$
defined in the same manner as we did for $\Theta^{(1)}(f)$.
Finally, we replace the symbol ${\bf U}$ at $\vartheta_i^{(1)}$ in $\Theta^{(1)}(f)$
with the subsequence $\Theta_i^{(2)}(f_i^{(1)})$ for all possible $i$
and create a new symbol sequence $\Theta^{(2)}(f)$,
called the {\it second symbol sequence}.
For intuitively speaking,
all unknown cards in $\Theta^{(1)}(f)$ are opened to show their suits
(as shown in Figure \ref{fig:SS}),
which become new cards in $\Theta^{(2)}(f)$.


If a new card is ${\bf N}$ in $\Theta^{(2)}(f)$,
we make sure that $f$ is not semi-definite.
If cards of $\Theta^{(2)}(f)$ are all ${\bf Y}_+$ (or all ${\bf Y}_-$),
then $f$ is semi-definite.
Note that by Appendix~A.2, symbols ${\bf Y}_+$ and ${\bf Y}_-$ do not coexist in
$\Theta^{(2)}(f)$ as in $\Theta^{(1)}(f)$.
Thus, as $\Theta^{(1)}(f)$, the sequence $\Theta^{(2)}(f)$ again has the same three subcases in which the sign of $f$ is unknown.
Therefore, we need to define and
discuss the $j$-th symbol sequence $\Theta^{(j)}(f)$ similarly for $j\ge 3$ if the symbol ${\bf U}$ still appears
but {\bf N} does not appear.
Therefore,
our discussion may be terminated in finitely many steps
with a symbol sequence $\Theta^{(j)}(f)$ which contains
either ${\bf Y}_+$ only or ${\bf Y}_-$ only or a symbol ${\bf N}$.

In order to show the process of identifying whether $f$ is semi-definite
on region $\Omega_\epsilon$ for an $\epsilon$,
we consider the polynomial
\begin{equation}
\begin{split}
f(\rho,\theta):=&\theta^6+(6\rho^2-4\rho)\theta^5+(11\rho^4-16\rho^3+6\rho^2)\theta^4
-(23\rho^6+10\rho^5
\\
&-12\rho^4+4\rho^3)\theta^3+(41\rho^7-12\rho^6+\rho^4)\theta^2-(12\rho^{10}+6\rho^9
\\
&+12\rho^8-10\rho^7+2\rho^6)\theta-17\rho^{11}+7\rho^{10}-6\rho^9+\rho^8
\end{split}
\label{fpro43}
\end{equation}
as an example.
Its Newton polygon ${\cal N}_f$ consists of three vertices $V_0:(0,6),V_1:(4,2)$ and $V_2:(8,0)$
and two edges $E_1$ and $E_2$, which link $V_0$ with $V_1$ and $V_1$ with $V_2$ respectively.
The $E_1$-polynomial and the $E_2$-polynomial are given by
$$
f_{E_1}(\theta)=\theta^2(\theta-1)^4~~~\mbox{and}~~~f_{E_2}(\theta)=(\theta-1)^2.
$$
Clearly, $f_{E_1}$ has one nonzero real root $c_1=1$ and
$f_{E_2}$ has one nonzero real root $c_2=1$.
Let $f^{(1)}_1:={\cal D}_{c_1} f$ and $f^{(1)}_2:={\cal D}_{c_2} f$.
Thus, by \eqref{des} we obtain
\begin{align*}
f^{(1)}_1(\rho,\theta)=&\theta^6+(6\rho+2)\theta^5+(11\rho^2+14\rho+1)\theta^4
+(-23\rho^3+34\rho^2+8\rho)\theta^3
\\
&-(28\rho^3-24\rho^2)\theta^2-(12\rho^5+6\rho^4-\rho^3)\theta-29\rho^5+\rho^4,
\\
f^{(1)}_2(\rho,\theta)=&\rho^4\theta^6+(12\rho^4-4\rho^3)\theta^5
+(56\rho^4-36\rho^3+6\rho^2)\theta^4+(101\rho^4-114\rho^3
\\
&+36\rho^2-4\rho)\theta^3+(72\rho^4-125\rho^3+60\rho^2-12\rho+1)\theta^2
\\
&-(\rho^4+38\rho^3-24\rho^2+2\rho)\theta-17\rho^4-12\rho^3+\rho^2.
\end{align*}
The polygon ${\cal N}_{f^{(1)}_1}$ has exactly one edge, denoted by $\tilde{E}$,
linking vertex $(0,4)$ with vertex $(4,0)$.
The $\tilde{E}$-polynomial of $f^{(1)}_1$ satisfies that
$$
(f^{(1)}_1)_{\tilde{E}}(\theta)=\theta^4+8\theta^3+24\theta^2+\theta+1.
$$
Let $g(\theta):=(f^{(1)}_1)_{\tilde{E}}(\theta)/\theta^{j_{\tilde{E}}}$,
where $j_{\tilde{E}}$ is the ordinate of the right end-point of $\tilde{E}$.
As we noted just before Proposition~\ref{lm-region}, for $k=1,2,3,4$,
we compute the $2k$-th discriminant of the polynomial $g$ and
obtain the sign sequence of discriminants
$$
({\rm sgn}(M_1),{\rm sgn}(M_2),{\rm sgn}(M_3),{\rm sgn}(M_4))=(1,0,-1,1),
$$
where $M_k$ is the $2k$-th order leading principal minor of the discrimination matrix of polynomial $g$.
Then the revised sign list of $g$ satisfies that
$$
{\rm rsl}(g)=(1,-1,-1,1).
$$
Clearly, $\Phi({\rm rsl}(g))=4$ and $\Xi({\rm rsl}(g))=2$,
and therefore,
$f^{(1)}_1$ satisfies the same condition as $f$ in {\bf(Y)}.
Moreover, $f^{(1)}_1$ is semi-positive
because coefficients corresponding to vertices of ${\cal N}_{f^{(1)}_1}$ are all positive.
On the other hand,
polygon ${\cal N}_{f^{(1)}_2}$ has exactly one edge linking vertex $(0,2)$ with vertex $(2,0)$,
and the edge-polynomial of $f^{(1)}_2$ is $(\theta-1)^2$.
Thus, $f^{(1)}_2$ satisfies the same condition as $f$ in {\bf(U)}.
Hence, we obtain the first symbol sequence
\begin{eqnarray}
\Theta^{(1)}(f)=({\bf Y}_+,{\bf U}).
\label{1st-seq}
\end{eqnarray}
With {\bf(U)}, we continue to desingularize function $f^{(1)}_2$.
Its edge-polynomial has exactly one nonzero real root $\tilde{c}_1=1$.
Let $f^{(2)}_1:={\cal D}_{\tilde{c}_1} f_2^{(1)}$.
Thus, by \eqref{des} we obtain
$$
f^{(2)}_1(\rho,\theta)=\theta^2(1+O(\theta))+\rho^2(1+O(\rho)).
$$
The polygon ${\cal N}_{f^{(2)}_1}$ has exactly one edge
linking vertex $(0,2)$ with vertex $(2,0)$,
and the edge-polynomial of $f^{(2)}_1$ is $\theta^2+1$.
Similarly to the above discussion of $f^{(1)}_1$,
one can check that $f^{(2)}_1$ satisfies the same condition as $f$ in {\bf(Y)} and
is semi-positive.
Hence, the unknown card {\bf U} appearing in the first symbol sequence \eqref{1st-seq}
is opened to show its suits: {\bf Y}$_+$.
Therefore, the second symbol sequence is
$
\Theta^{(2)}(f)=({\bf Y}_+,{\bf Y}_+),
$
which implies that the polynomial $f$ given in \eqref{fpro43} is semi-positive
on region $\Omega_\epsilon$ for an $\epsilon$.

The following proposition is devoted to the case that every $\Theta^{(j)}(f)$, $j=1,2,...$,
contains a symbol ${\bf U}$ but no ${\bf N}$,
that is, $\Theta^{(j)}(f)$ contains either both ${\bf U}$ and ${\bf Y}_+$ only
or both ${\bf U}$ and ${\bf Y}_-$ only or ${\bf U}$ only for all $j\ge1$.

\begin{prop}
Let $f$ be given by \eqref{fGH} and satisfy condition {\bf(U)},
given just before Proposition~\ref{lm-region}.
If the $j$-th symbol sequence $\Theta^{(j)}(f)$ contains
a symbol ${\bf U}$ but no ${\bf N}$ for all $j\ge1$,
then $\omega f$ is semi-positive,
where $\omega$ is defined in Proposition~\ref{lm-region}.
\label{lm-positive}
\end{prop}

{\bf Proof.}
As indicated at the beginning of the proof of Proposition~\ref{lm-region},
if $\omega f$ is not semi-positive,
then the equation $f(\rho,\theta)=0$ has a real branch $\theta^{(0)}(\rho)$
of odd multiplicity $\mu$.
We have $\theta^{(0)}(\rho)\not\equiv0$ since $f$ does not satisfy {\bf (N1)}.
By \eqref{csq}, $\theta^{(0)}(\rho)$ can be expanded as
$$
\theta^{(0)}(\rho)=\sum_{i=0}^{N} c_i\rho^{\nu_i^{(0)}},
$$
where $N\in\mathbb{Z}_+\cup\{+\infty\}$,
$c_i$\,s are non-vanished and $\nu_i^{(0)}$\,s are rational such that $0<\nu_0^{(0)}<\nu_1^{(0)}<\cdots$.
We denote $f$ by $f^{(0)}$ for convenience.
For $k=0,...,N$,
Lemma~\ref{lm-basic}{\bf(i)} implies that
${\cal N}_{f^{(k)}}$ has an edge $E^{(k)}$ with slope $-q_k/p_k$
for a pair of coprime positive integers $p_k$ and $q_k$ such that $\nu^{(k)}_k=p_k/q_k$,
and $c_k$ is a real root of $f^{(k)}_{E^{(k)}}$,
and we see from \eqref{des} that the equation $f^{(k+1)}(\rho,\theta)=0$ has the real branch
\begin{eqnarray}
\theta=\theta^{(k+1)}(\rho):=\frac{\theta^{(k)}(\rho^{q_{k}})}{\rho^{p_{k}}}-c_{k}
\left\{
\begin{aligned}
&=\sum_{i=k+1}^{N} c_i\rho^{\nu_i^{(k+1)}} &&\mbox{if}~k<N,\\
&\equiv 0 &&\mbox{if}~k=N<+\infty
\end{aligned}
\right.
\label{ykk}
\end{eqnarray}
of multiplicity $\mu$,
where $f^{(k+1)}:={\cal D}_{c_k}f^{(k)}$,
and $\nu^{(k+1)}_i:=q_k\nu^{(k)}_i -p_k$ for all $i=k+1,...,N$.

In the case $N<+\infty$,
we see from \eqref{ykk} that $\theta^{(N+1)}(\rho)\equiv0$.
Then, $f^{(N+1)}$ has a factor $\theta$ of odd multiplicity $\mu$,
implying that $f^{(N+1)}$ satisfies {\bf(N1)} as $f$,
and therefore,
the $(N+1)$-th symbol sequence contains a symbol {\bf N},
a contradiction to the assumption of this proposition.

In the case $N=+\infty$,
let $V_0^{(k)}:(i_0^{(k)},j_0^{(k)})$ be the left-most vertex of ${\cal N}_{f^{(k)}}$.
Then,
we have $i_0^{(k)}=0$ since $f^{(k)}$ does not have a factor $\rho$,
as indicated just below \eqref{des}.
On the other hand,
the sequence $\{j_0^{(k)}\}_{k=0}^{+\infty}$ of positive integers is decreasing.
In fact,
it is proved by \cite[Lemma~1.4.1]{CA00} that
$j_0^{(k+1)}$ equals the multiplicity of $c_k$ as a real root of
the polynomial $f^{(k)}_{E^{(k)}}$,
whose degree is not bigger than $j_0^{(k)}$.
It follows that $j_0^{(k)}\ge j_0^{(k+1)}>0$.
Hence, there is an integer $m\ge0$ such that
$j_*:=j_0^{(m)}=j_0^{(m+1)}=\cdots$.
We claim that
$\theta^{(m)}(\rho)$ is a real branch of the equation $f^{(m)}(\rho,\theta)=0$
of multiplicity $j_*$ (see \cite[Exercise~1.5]{CA00}).
Actually, Lemma~1.5.1 of \cite{CA00} indicates that
the Newton polygon ${\cal N}_{f^{(k)}}$ has exactly one edge and
its slope is $-1/p_k$, i.e., $q_k=1$, for all $k \ge m$.
Similar to the computation of \eqref{eld},
$$
f^{(k)}(\rho,\theta)=a^{(k)}_{0,j_*}(\theta-c_k\rho^{p_k})^{j_*}
+\sum_{i+jp_k> j_*p_k}a_{i,j}^{(k)}\rho^i\theta^j~~~\mbox{for all}~~~k\ge m.
$$
Let $g^{(k)}:=\partial^{j_*-1}f^{(k)}/\partial\theta^{j_*-1}$ for all $k\ge m$.
Then
$$
g^{(k)}(\rho,\theta)=j_*!a^{(k)}_{0,j_*}(\theta-c_k\rho^{p_k})+o(\theta)+o(\rho^{p_k}).
$$
Clearly, ${\cal N}_{g^{(k)}}$ has exactly one edge, whose slope is $-1/p_k$,
and the edge-polynomial has exactly one nonzero real root $c_k$.
Then the desingularized function
\begin{align*}
g_{c_k}^{(k)}(\rho,\theta)
&=\frac{g^{(k)}(\rho,\rho^{p_k}(c_k+\theta))}{\rho^{p_k}}
=\frac{\partial^{j_*-1}f^{(k)}(\rho,\rho^{p_k}(c_k+\theta))}{\rho^{p_kj_*}\partial\theta^{j_*-1}}
=\frac{\partial^{j_*-1}f^{(k+1)}(\rho,\theta)}{\partial \theta^{j_*-1}}
\\
&=g^{(k+1)}(\rho,\theta).
\end{align*}
It follows that $\theta^{(m)}(\rho)$ is a real branch of the equation $g^{(m)}(\rho,\theta)=0$
and therefore, our claim is proved.
This claim implies that $j_0^{(m)}$ equals $\mu$ and is odd.
Moreover, by \cite[Lemma~1.4.1]{CA00},
$f^{(m-1)}_{E^{(m-1)}}$ has a nonzero real root $c_{m-1}$ of odd multiplicity $j_0^{(m)}$,
implying that $f^{(m-1)}$ satisfies condition {\bf(N2)} as $f$.
Thus, the $(m-1)$-th symbol sequence contains a symbol {\bf N},
a contradiction to the assumption.

As a consequence of above two cases,
the equation $f(\rho,\theta)=0$ does not have a real branch of odd multiplicity.
Therefore, $\omega f$ is semi-definite and
the proof of this proposition is completed.
\qquad$\Box$

\section{Computation with Newton polygons}
\setcounter{equation}{0}
\setcounter{lm}{0}
\setcounter{thm}{0}
\setcounter{rmk}{0}
\setcounter{df}{0}
\setcounter{cor}{0}
\setcounter{pro}{0}

In last section
the Newton polygon and edge-polynomials of the Lie-bracket
$[{\cal G},{\cal H}]_\theta={\cal G}'_\theta{\cal H}-{\cal G}{\cal H}'_\theta$
are used to determine the sign of $[{\cal G},{\cal H}]_\theta$
on the region
$\Omega_\epsilon=\{(\rho,\theta)\in\mathbb{R}^2:0<\rho<\epsilon,|\theta|<\epsilon\}$
for small $\epsilon>0$,
but we need to know how
the Newton polygon and edge-polynomials of $[{\cal G},{\cal H}]_\theta$
can be computed with those of ${\cal G}$ and ${\cal H}$.
In order to do so,
we need knowledge on
{\it addition, multiplication} and {\it differentiation} of Newton polygons.
However, we only found the latter two from references.
In the book \cite{CA00}
both Newton polygons and edge-polynomials are computed for differentiation
but only edge-polynomials are obtained for multiplication.
Computation of Newton polygons for multiplication can be found from
\cite{Tei}, \cite[p.480, Lemma~3]{BK86} and \cite[pp.180-181]{JP00}.
However,
we found neither Newton polygons nor edge-polynomials for addition
from the above mentioned references \cite{BK86, CA00, JP00, Tei} or other literature.
In this section,
we first introduce those known results on differentiation and multiplication of Newton polygons.
Then, we investigate addition of Newton polygons.

For convenience,
we need some terminologies.
Similar to the $\xi$-componential polynomial defined in \eqref{defkxip},
a polynomial ${\cal G}$ restricted its exponents to a closed linear segment $\Gamma$,
whose end-points are both nonnegative lattice points on the $(u,v)$-plane,
is defined as
\begin{eqnarray}
{\cal G}_\Gamma(\theta)
:=\sum_{(i,j)\in\Gamma\cap\Delta({\cal G})} a_{i,j}({\cal G}) \theta^j,
\label{fonE}
\end{eqnarray}
called the {\it $\Gamma$-polynomial} of ${\cal G}$.
In particular,
when $\Gamma$ is a lattice point,
${\cal G}_\Gamma$ is a nonzero monomial if $\Gamma\in\Delta({\cal G})$;
otherwise, ${\cal G}_\Gamma=0$.
When $\Gamma$ is an edge of ${\cal N}_{\cal G}$,
then ${\cal G}_\Gamma$ is an edge-polynomial, defined by \eqref{red}.

The lower convex semi-hull of a lattice ${\cal Q}$,
a set consisting of points with integer coordinates, on the $(u,v)$-plane is defined as
\begin{eqnarray}
\wp({\cal Q}):={\rm conv}\{{\cal Q}+{\mathbb R}_+^2\}.
\label{wpQ}
\end{eqnarray}
We use ${\cal N}({\cal Q})$ to present the principal boundary of $\wp({\cal Q})$.
Endowed with the natural $u$-order, ${\cal N}({\cal Q})$ is called
the {\it ordered principal boundary} and denoted by $\vec{\cal N}({\cal Q})$.
Let ${\cal V}(\wp({\cal Q}))$ consist of all vertices of $\wp({\cal Q})$.
Then, $\wp({\cal Q})=\wp({\cal V}(\wp({\cal Q})))$.

\subsection{Differentiation}

In this subsection,
we show how to give the Newton polygon of the differentiation of
the nonzero analytic function ${\cal G}$.
Since only the $\theta$-derivative will be used in the following sections
and the case of the $\rho$-derivative can be discussed similarly,
we only deal with the $\theta$-derivative.
Computing the derivative in \eqref{expGH}, we get
\begin{alignat*}{2}
{\cal G}_\theta'(\rho,\theta)
=\sum_{(i,j)\in\Delta({\cal G})}ja_{i,j}({\cal G})\rho^i\theta^{j-1}
=\sum_{(i,j)\in\Delta({\cal G}),~j\ge 1}ja_{i,j}({\cal G})\rho^i\theta^{j-1}.
\end{alignat*}
Since $ja_{i,j}({\cal G})\ne 0$ for all $(i,j)\in\Delta({\cal G})$ with $j\ge 1$,
we obtain that $\Delta({\cal G}'_\theta)=\{(i,j)\in\Delta({\cal G}):j\ge 1\}-(0,1)$.
Hence,
we can write
$$
{\cal G}'_\theta(\rho,\theta)
=\sum_{(i,j)\in\Delta({\cal G}'_\theta)}a_{i,j}({\cal G}'_\theta)\rho^i\theta^j
$$
with $a_{i,j}({\cal G}'_\theta):=(j+1)a_{i,j+1}({\cal G})$.
The following proposition is summarized from \cite[pp.232-235]{CA00}.

\begin{pro}
Suppose that ${\cal G}$ is a nonzero analytic function of the form \eqref{expGH}.
Then the Newton polygon of the derivative ${\cal G}'_\theta$ satisfies that
$
{\cal N}_{{\cal G}'_\theta}={\cal N}(\Delta^*({\cal G}))-(0,1),
$
where $\Delta^*({\cal G}):=\{(i,j)\in\Delta({\cal G}):j\ge1\}$
and ${\cal N}(\Delta^*({\cal G}))$ is the principal boundary of the lower convex semi-hull of $\Delta^*({\cal G})$, defined just below \eqref{wpQ}.
The edge-polynomial of the derivative ${\cal G}'_\theta$ satisfies that
$({\cal G}'_\theta)_{E_k({\cal G}'_\theta)}=({\cal G}_{N_k})'_\theta~\mbox{for all}~k=1,...,s({\cal G}'_\theta),
$
where $N_k$ is the $k$-th edge of the ordered principal boundary
$\vec{\cal N}(\Delta^*({\cal G}))$
and $s({\cal G}'_\theta)$ is the number of edges of ${\cal N}_{{\cal G}'_\theta}$.
Moreover,
$
{\cal K}_\xi ({\cal G}'_\theta)=({\cal K}_\xi ({\cal G}))'_\theta~\mbox{for all}~\xi \in(-\infty,\zeta(E_{s({\cal G})}({\cal G}))]
$
if $s({\cal G})\ge1$,
where $s({\cal G})$ is the number of edges of ${\cal N}_{\cal G}$,
${\cal K}_\xi ({\cal G}'_\theta)$ and ${\cal K}_\xi ({\cal G})$ are
$\xi$-componential polynomials,
and $\zeta(E_{s({\cal G})}({\cal G}))$ is the slope of
the last edge $E_{s({\cal G})}({\cal G})$ of ${\cal N}_{\cal G}$.
\label{lm-Diff}
\end{pro}

\begin{rmk}
{\rm
{\bf (i)}
When $j_{s({\cal G})}$,
the ordinate of the right-most vertex of ${\cal N}_{\cal G}$,
is larger than $1$,
we have ${\cal N}_{{\cal G}'_\theta}={\cal N}_{\cal G}-(0,1)$,
$({\cal G}'_\theta)_{E_k({\cal G}'_\theta)}
=({\cal G}_{E_k({\cal G})})'_\theta$ for all $k=1,...,s({\cal G})$,
and ${\cal K}_\xi ({\cal G}'_\theta)=({\cal K}_\xi ({\cal G}))'_\theta$ for all $\xi<0$.
{\bf (ii)}
When $j_{s({\cal G})}$ is equal to $0$ and $s({\cal G})\ge 1$,
let $V_*({\cal G}):(i_*,j_*)$ be the second to the last valid index on $\vec{\cal N}_{\cal G}$.
If $\{(i,j)\in\Delta({\cal G}):0< j< j_*\}=\emptyset$,
i.e., condition {\bf(S)} given just before Theorem~\ref{th:finiteJ2} holds,
then
$
{\cal N}_{{\cal G}'_\theta}
=\{{\cal N}_{\cal G}\backslash\underline{V_*({\cal G})V_{s({\cal G})}({\cal G})}\}\cup V_*({\cal G})-(0,1),
$
where $\underline{V_*({\cal G})V_{s({\cal G})}({\cal G})}$ is
the closed linear segment linking $V_*({\cal G})$ with $V_{s({\cal G})}({\cal G})$.
Moreover,
$({\cal G}'_\theta)_{E_k({\cal G}'_\theta)}=({\cal G}_{E_k({\cal G})})'_\theta$ for all
$k=1,...,s({\cal G})-1+{\rm sgn}(j_{s({\cal G})-1}-j_*)$,
where $j_{s({\cal G})-1}$ is the ordinate of
the second to the last vertex of $\vec{\cal N}_{\cal G}$ and
${\rm sgn}$ is the signum function.
}\label{rmk:diff}
\end{rmk}

\subsection{Multiplication}

In order to obtain the Newton polygon of the multiplication of
nonzero analytic functions ${\cal G}$ and ${\cal H}$,
we introduce {\it joining} of convex polygons on the $(u,v)$-plane (\cite{GMN,Tei}).
We call a closed linear segment on the $(u,v)$-plane
a {\it basic segment} if the segment has a negative slope and
its end-points have nonnegative integer coordinates.
We also refer each point in $\mathbb{Z}^2_+$ to as a (trivial) basic segment,
which actually lies in the case that the left end-point coincides with the right one.
Let ${\cal S}$ be the set of all basic segments on the $(u,v)$-plane.
For a basic segment $S\in{\cal S}$,
we use ${\cal L}(S)$ and ${\cal R}(S)$ to denote
its {\it left end-point} and {\it right end-point} respectively.
Moreover,
we call the lengths of $S$ projected to the horizontal axis and the vertical one
respectively
the {\it length} $\ell(S)$ and the {\it height} $h(S)$ of $S$.
The {\it joining} $S_1 \uplus\cdots\uplus S_\varrho$
of basic segments $S_1,...,S_\varrho\in{\cal S}$
is the convex polygon
starting from the lattice point ${\cal L}(S_1)+\cdots+{\cal L}(S_\varrho)$
and joining all segments of positive length in the order of increasing slopes.

Let ${\cal P}$ denote the set of all Newton polygons, i.e.,
$
{\cal P}:=\{{\cal N}_{\cal G}:{\cal G}\in\mathbb{R}\{\rho,\theta\}\backslash\{0\}\},
$
where $\mathbb{R}\{\rho,\theta\}$ is the ring of real analytic functions at $O:(0,0)$.
Clearly,
$$
{\cal P}=\{S_1 \uplus\cdots\uplus S_\varrho:S_1,...,S_\varrho\in{\cal S},\varrho\ge1\},
$$
consisting of
all polygons obtained by joining finitely many basic segments in ${\cal S}$.
Newton polygons ${\cal N}_{\cal G}$ and ${\cal N}_{{\cal H}}$ have the following unique decompositions
\begin{eqnarray*}
{\cal N}_{\cal G}=S_0\uplus S_1\uplus\cdots \uplus S_{s({\cal G})}~~~\mbox{and}~~~
{\cal N}_{{\cal H}}=\tilde{S}_0\uplus \tilde{S}_1\uplus\cdots \uplus \tilde{S}_{s({\cal H})}
\end{eqnarray*}
respectively, where
each $S_k$, $k=1,...,s({\cal G})$,
is a basic segment with length $\ell(E_k({\cal G}))$ and
height $h(E_k({\cal G}))$ and end-points on the axes,
each $\tilde{S}_k$, $k=1,...,s({\cal H})$,
is a basic segment with length $\ell(E_k({\cal H}))$ and
height $h(E_k({\cal H}))$ and end-points on the axes, and
$S_0:(i_0,j_{s({\cal G})})$ and
$\tilde{S}_0:(\tilde{i}_0,\tilde{j}_{s({\cal H})})$
are lattice points.
Then the {\it joining} of Newton polygons ${\cal N}_{\cal G}$ and ${\cal N}_{\cal H}$
can be presented as
$$
{\cal N}_{\cal G}\uplus{\cal N}_{\cal H}
:=S_0\uplus S_1\uplus\cdots \uplus S_{s({\cal G})} \uplus
\tilde{S}_0\uplus \tilde{S}_1\uplus\cdots \uplus \tilde{S}_{s({\cal H})}.
$$

\begin{pro}
Suppose that ${\cal G}$ and ${\cal H}$ are both
nonzero analytic functions of the form \eqref{expGH}.
Then, ${\cal N}_{{\cal GH}}={\cal N}_{\cal G}\uplus{\cal N}_{\cal H}$, and
${\cal K}_\xi ({\cal GH})={\cal K}_\xi ({\cal G}){\cal K}_\xi ({\cal H})$ for all $\xi<0$,
where ${\cal K}_\xi$ is defined by \eqref{defkxip}.
\label{lm-MIL}
\end{pro}

The first conclusion of Proposition~\ref{lm-MIL} can be found in \cite[p.618]{Tei}
and also seen in \cite[p.480, Lemma~3]{BK86} and \cite[pp.180-181]{JP00}.
The second conclusion can be found in \cite[p.231]{CA00}.

\begin{rmk}
{\rm
{\bf(i)} The slope sequence of ${\cal GH}$ satisfies that
$\vec{\mathfrak{S}}({\cal GH})
=\vec{\mathfrak{S}}({\cal G})\cup\vec{\mathfrak{S}}({\cal H})$,
the sequence union defined by \eqref{unionseq},
and moreover, ${\cal N}_{\cal GH}$ has $s({\cal G},{\cal H})$ edges,
where $s({\cal G},{\cal H})$ is the cardinality of the union
$\mathfrak{S}({\cal G})\cup\mathfrak{S}({\cal H})$.
{\bf(ii)}
The left end-point $V_0({\cal GH})$ and
the right one $V_{s({\cal G},{\cal H})}({\cal GH})$ of ${\cal N}_{\cal GH}$
satisfy that
{\small
\begin{align*}
V_0({\cal GH})
&={\cal L}(S_0)+\cdots+{\cal L}(S_{s({\cal G})})+{\cal L}(\tilde{S}_0)
+\cdots+{\cal L}(\tilde{S}_{s({\cal H})})
=V_0({\cal G})+V_0({\cal H}),
\\
V_{s({\cal G},{\cal H})}({\cal GH})
&={\cal R}(S_0)+\cdots+{\cal R}(S_{s({\cal G})})+{\cal R}(\tilde{S}_0)
+\cdots+{\cal R}(\tilde{S}_{s({\cal H})})
=V_{s({\cal G})}({\cal G})+V_{s({\cal H})}({\cal H}).
\end{align*}
}Moreover,
if $E_k({\cal G})$ is not parallel to an edge of ${\cal N}_{{\cal H}}$,
then ${\cal N}_{{\cal GH}}$ has one edge $E({\cal GH})$ such that
$\zeta(E({\cal GH}))=\zeta(E_{k}({\cal G}))$ and
$h(E({\cal GH}))=h(E_{k}({\cal G}))$.
Oppositely,
if $E_{k}({\cal G})$ is parallel to
an edge $E_{\tilde{k}}({\cal H})$ of ${\cal N}_{\cal G}$,
then ${\cal N}_{\cal GH}$ has one edge $E({\cal GH})$ such that
$\zeta(E({\cal GH}))=\zeta(E_{k}({\cal G}))=\zeta(E_{\tilde{k}}({\cal H}))$
and
$h(E({\cal GH}))=h(E_{k}({\cal G}))+h(E_{\tilde{k}}({\cal H})).
$
}\label{rmk:mul}
\end{rmk}

\subsection{Addition}

The Newton polygons of ${\cal G}$ and ${\cal H}$ are determined by
the lower convex semi-hulls of $\Delta({\cal G})$ and $\Delta({\cal H})$ respectively.
Generally, the Newton polygon of the addition ${\cal G}+{\cal H}$ is determined by
the lower convex semi-hull of $\Delta({\cal G})\cup\Delta({\cal H})$.
However, in the case that
some common valid indices of ${\cal N}_{\cal G}$ and ${\cal N}_{\cal H}$ are vanished,
the Newton polygon of the addition ${\cal G}+{\cal H}$ cannot be determined by
the lower convex semi-hull of $\Delta({\cal G})\cup\Delta({\cal H})$.
The corresponding results are hardly found in the literature.

The so called {\it vanishing problem} is that
a common valid index of ${\cal G}$ and ${\cal H}$ is
no longer a valid index of ${\cal G}+{\cal H}$.
For example,
if ${\cal G}(\rho,\theta)$ contains the term $\rho^2\theta$ and
${\cal H}(\rho,\theta)$ contains the term $-\rho^2\theta$,
then point $(2,1)$ is a common valid index of them, however,
it is not a valid index of ${\cal G}+{\cal H}$.
The common valid index $(i,j)\in\Delta({\cal G})\cap\Delta({\cal H})$
is vanished under addition, i.e.,
$(i,j)\notin\Delta({\cal G}+{\cal H})$,
if $a_{i,j}({\cal G})+a_{i,j}({\cal H})=0$.
Otherwise, it is non-vanished if
\begin{eqnarray}
a_{i,j}({\cal G})+a_{i,j}({\cal H})\ne0.
\label{Anon}
\end{eqnarray}
Inequality \eqref{Anon} is called the {\it non-vanishing condition for addition}.

Whether common vertices of ${\cal G}$ and ${\cal H}$ are vanished or not,
it affects the Newton polygon of ${\cal G}+{\cal H}$.
Let
$$
\Lambda^V({\cal G}\backslash{\cal H}):=\{(i,j)\in\Delta^V({\cal G})\cap\Delta({\cal H}):
a_{i,j}({\cal G})+a_{i,j}({\cal H})=0\}.
$$
Similarly, we can define $\Lambda^V({\cal H}\backslash {\cal G})$.
They are the vanishing vertex sets of ${\cal G}$ and ${\cal H}$ respectively.
When common vertices vanish, i.e.,
$\Lambda^V({\cal G}\backslash{\cal H})
\cap\Lambda^V({\cal H}\backslash {\cal G})\ne\emptyset$,
consider
\begin{eqnarray*}
&{\cal G}_1(\rho,\theta):=
{\cal G}(\rho,\theta)-\sum_{(i,j)\in\Lambda^V({\cal G}\backslash{\cal H})\cap\Lambda^V({\cal H}\backslash {\cal G})}a_{i,j}({\cal G})\rho^i\theta^j,
\\
&{\cal H}_1(\rho,\theta):={\cal H}(\rho,\theta)-\sum_{(i,j)\in\Lambda^V({\cal G}\backslash{\cal H})\cap\Lambda^V({\cal H}\backslash {\cal G})}a_{i,j}({\cal H})\rho^i\theta^j.
\end{eqnarray*}
By the vanishing condition given just before \eqref{Anon},
we have ${\cal G}+{\cal H}={\cal G}_1+{\cal H}_1$,
which implies that the vanishing situation can be converted into the non-vanishing one.
In what follows, we only consider the situation that common vertices are non-vanished.
The following lemma states the relation of valid indices.

\begin{lm}
Suppose that ${\cal G}$ and ${\cal H}$ are
both nonzero analytic functions of the form \eqref{expGH}.
Then, $\Delta({\cal G}+{\cal H})\subset\Delta({\cal G})\cup\Delta({\cal H})$.
On the contrary, $\{\Delta({\cal G})\cup\Delta({\cal H})\}\backslash\{\Delta({\cal G})\cap\Delta({\cal H})\}\subset\Delta({\cal G}+{\cal H})$,
and for $(i,j)\in\Delta({\cal G})\cap\Delta({\cal H})$, $(i,j)\in\Delta({\cal G}+{\cal H})$ if and only if $(i,j)$ satisfies the non-vanishing condition \eqref{Anon} for addition.
\label{lm-Dff}
\end{lm}

{\bf Proof.}
Since ${\cal G}(\rho,\theta)=\sum_{(i,j)\in\Delta({\cal G})}a_{i,j}({\cal G})\rho^i\theta^j$,
${\cal H}(\rho,\theta)=\sum_{(i,j)\in\Delta({\cal H})} a_{i,j}({\cal H})\rho^i\theta^j$,
we have
\begin{equation}
\begin{aligned}
({\cal G}+{\cal H})(\rho,\theta)
=&\sum_{(i,j)\in\Delta({\cal G})\backslash\Delta({\cal H})}a_{i,j}({\cal G})\rho^i\theta^j
+\sum_{(i,j)\in\Delta({\cal H})\backslash\Delta({\cal G})}a_{i,j}({\cal H})\rho^i\theta^j
\\
&+\sum_{(i,j)\in\Delta({\cal G})\cap\Delta({\cal H})}\{a_{i,j}({\cal G})
+a_{i,j}({\cal H})\}\rho^i\theta^j.
\end{aligned}
\label{f+tf}
\end{equation}
It follows that
$$
\Delta({\cal G}+{\cal H})\subset
\{\Delta({\cal G})\backslash\Delta({\cal H})\} \cup
\{\Delta({\cal H})\backslash\Delta({\cal G})\} \cup
\{\Delta({\cal G})\cap\Delta({\cal H})\}
=\Delta({\cal G})\cup \Delta({\cal H}).
$$
On the contrary,
note that $(i,j)\in\{\Delta({\cal G})\cup\Delta({\cal H})\}\backslash
\{\Delta({\cal G})\cap\Delta({\cal H})\}
=\{\Delta({\cal G})\backslash\Delta({\cal H})\}\cup
\{\Delta({\cal H})\backslash\Delta({\cal G})\}$.
If $(i,j)\in\Delta({\cal G})\backslash\Delta({\cal H})$,
then \eqref{f+tf} implies that $a_{i,j}({\cal G}+{\cal H})=a_{i,j}({\cal G})\ne 0$,
and therefore, $(i,j)\in\Delta({\cal G}+{\cal H})$.
The case $(i,j)\in\Delta({\cal H})\backslash\Delta({\cal G})$ can be similarly discussed.
Hence,
$$
\{\Delta({\cal G})\cup\Delta({\cal H})\}\backslash
\{\Delta({\cal G})\cap\Delta({\cal H})\}
\subset\Delta({\cal G}+{\cal H}).
$$
Additionally, for any $(i,j)\in\Delta({\cal G})\cap\Delta({\cal H})$,
it follows from \eqref{f+tf} that
$a_{i,j}({\cal G}+{\cal H})=a_{i,j}({\cal G})+a_{i,j}({\cal H})$.
Then, $(i,j)\in\Delta({\cal G}+{\cal H})$ if and only if
$(i,j)$ satisfies the non-vanishing condition~\eqref{Anon}.
The proof of this lemma is completed.
\qquad$\Box$

\begin{pro}
Suppose that ${\cal G}$ and ${\cal H}$ are both nonzero analytic functions of the form \eqref{expGH} and
their common vertices satisfy the non-vanishing condition \eqref{Anon} for addition,
i.e., $\Lambda^V({\cal G}\backslash{\cal H})\cap
\Lambda^V({\cal H}\backslash {\cal G})=\emptyset$.
Then, ${\cal G}+{\cal H}$ has the Newton polygon
\begin{eqnarray}
{\cal N}_{{\cal G}+{\cal H}}={\cal N}(\Delta^V({\cal G})\cup\Delta^V({\cal H})),
\label{ADDN}
\end{eqnarray}
the principal boundary of
the lower convex semi-hull of $\Delta^V({\cal G})\cup\Delta^V({\cal H})$,
defined just below \eqref{wpQ}.
Moreover, the edge-polynomial of ${\cal G}+{\cal H}$ satisfies that
\begin{eqnarray}
({\cal G}+{\cal H})_{E_k({\cal G}+{\cal H})}={\cal G}_{N_k}+{\cal H}_{N_k}~\mbox{for all}~k=1,...,s({\cal G}+{\cal H}),
\label{ADDfE}
\end{eqnarray}
where $N_k$ is the $k$-th edge of the ordered principal boundary
$\vec{\cal N}(\Delta^V({\cal G})\cup\Delta^V({\cal H}))$,
and $s({\cal G}+{\cal H})$ is the number of edges of ${\cal N}_{{\cal G}+{\cal H}}$.
\label{ADD}
\end{pro}

{\bf Proof.}
For equality \eqref{ADDN},
we only need to show that
\begin{eqnarray}
\wp(\Delta({\cal G}+{\cal H}))=\wp(\Delta^V({\cal G})\cup\Delta^V({\cal H})).
\label{wpwp1}
\end{eqnarray}
By the inclusion $\Delta({\cal G}+{\cal H})\subset\Delta({\cal G})\cup\Delta({\cal H})$ given in Lemma~\ref{lm-Dff}
and the property of the union of convex hulls (\cite[Theorem 3.2]{Soltan}),
i.e., ${\rm conv}\{{\cal A}\cup{\cal B}\}
={\rm conv}\{{\rm conv}\{{\cal A}\}\cup{\rm conv}\{{\cal B}\}\}$
for any two sets ${\cal A},{\cal B}\subset\mathbb{R}^2$,
we obtain that
\begin{eqnarray*}
\wp(\Delta({\cal G}+{\cal H}))
&=&{\rm conv}\{\Delta({\cal G}+{\cal H})+{\mathbb R}_+^2\}
\\
&\subset&{\rm conv}\{\Delta({\cal G})\cup\Delta({\cal H})+{\mathbb R}_+^2\}
\\
&=&{\rm conv}\{{\rm conv}\{\Delta({\cal G})+{\mathbb R}_+^2\}
\cup{\rm conv}\{\Delta({\cal H})+{\mathbb R}_+^2\}\}
\\
&=&{\rm conv}\{{\rm conv}\{\Delta^V({\cal G})+{\mathbb R}_+^2\}
\cup{\rm conv}\{\Delta^V({\cal H})+{\mathbb R}_+^2\}\}
\\
&=&{\rm conv}\{\Delta^V({\cal G})\cup\Delta^V({\cal H})+{\mathbb R}_+^2\}
\\
&=&\wp(\Delta^V({\cal G})\cup\Delta^V({\cal H})),
\end{eqnarray*}
which proves the ``$\subset$''-part of \eqref{wpwp1}.
In what follows, we prove the ``$\supset$''-part of \eqref{wpwp1}.
It suffices to show the following inclusion of the  vertex set
\begin{eqnarray}
{\cal V}_*:={\cal V}(\wp(\Delta^V({\cal G})\cup\Delta^V({\cal H})))
\subset\Delta({\cal G}+{\cal H})
\label{VPVV}
\end{eqnarray}
because $\wp(\Delta^V({\cal G})\cup\Delta^V({\cal H}))=\wp({\cal V}_*)$.
If \eqref{VPVV} does not hold,
then there exists $P\in{\cal V}_*\subset\Delta^V({\cal G})\cup\Delta^V({\cal H})$ such that $P\notin\Delta({\cal G}+{\cal H})$.
Assume without loss of generality that $P\in\Delta^V({\cal G})$.
Then, $P\in\Delta^V({\cal G})\cap\Delta({\cal H})$ by Lemma \ref{lm-Dff}.
Since those common vertices of ${\cal G}$ and ${\cal H}$ are non-vanished,
$P\in\Delta({\cal H})\backslash\Delta^V({\cal H})$.
Therefore,
\begin{eqnarray}
P\in\partial\wp(\Delta({\cal H}))\backslash\Delta^V({\cal H}).
\label{PPDV}
\end{eqnarray}
Otherwise, noticing that $P\in\Delta({\cal H})\backslash\Delta^V({\cal H})\subset
\wp(\Delta({\cal H}))\backslash\Delta^V({\cal H})$,
we obtain
$$
P\in{\rm int}\wp(\Delta({\cal H}))
={\rm int}\wp(\Delta^V({\cal H}))
\subset{\rm int}\wp(\Delta^V({\cal G})\cup\Delta^V({\cal H})),
$$
where ${\rm int}$ represents the interior of a set,
a contradiction to the fact that $P\in{\cal V}_*$,
the vertex set of $\wp(\Delta^V({\cal G})\cup\Delta^V({\cal H}))$.
Thus, \eqref{PPDV} holds.
Note that $P$ is a vertex of $\wp(\Delta^V({\cal G})\cup\Delta^V({\cal H}))$.
Then, $P\notin\underline{Q_1Q_2}$ for any $Q_1,Q_2\in\wp(\Delta^V({\cal G})\cup\Delta^V({\cal H}))\backslash\{P\}$,
where $\underline{Q_1Q_2}$ is the linear closed segment linking $Q_1$ with $Q_2$.
We see from \eqref{PPDV} that either
$P\in\underline{V_{k-1}({\cal H})V_k({\cal H})}$ for a $k\in\{1,2,...,s({\cal H})\}$,
or $P\in\underline{Q^*V_0({\cal H})}$ for a $Q^*\in L_{{\cal H}}^v\backslash\{P\}$,
or $P\in\underline{Q_*V_{s({\cal H})}({\cal H})}$ for a $Q_*\in L_{{\cal H}}^h\backslash\{P\}$,
where $L_{{\cal H}}^v$ and $L_{{\cal H}}^h$ are
the vertical boundary and the horizontal one of $\wp(\Delta({\cal H}))$ respectively.
Since $\wp(\Delta({\cal H}))=\wp(\Delta^V({\cal H}))\subset\wp(\Delta^V({\cal G})\cup\Delta^V({\cal H}))$,
we have $V_{k-1}({\cal H}),V_k({\cal H}),Q_*,Q^*\in
\wp(\Delta^V({\cal G})\cup\Delta^V({\cal H}))$.
Whatever the case, it contradicts to the fact that
$P$ is a vertex of $\wp(\Delta^V({\cal G})\cup\Delta^V({\cal H}))$.
This proves \eqref{VPVV}, implying the ``$\supset$''-part of \eqref{wpwp1}.
Thus equality \eqref{wpwp1} holds and equality \eqref{ADDN} follows directly.

Finally, we consider equality \eqref{ADDfE}.
It follows form \eqref{ADDN} that $E_k({\cal G}+{\cal H})=N_k$
for all $k=1,...,s({\cal G}+{\cal H})$.
By definition \eqref{fonE},
\begin{align*}
({\cal G}+{\cal H})_{E_k({\cal G}+{\cal H})}(\theta)
&=\sum_{(i,j)\in N_k}a_{i,j}({\cal G}+{\cal H})\theta^j
\\
&=\sum_{(i,j)\in N_k}\{a_{i,j}({\cal G})\theta^j+a_{i,j}({\cal H})\theta^j\}
\\
&={\cal G}_{N_k}(\theta)+{\cal H}_{N_k}(\theta),
\end{align*}
that is, equality \eqref{ADDfE} holds.
This completes the proof of Proposition~\ref{ADD}.
\qquad$\Box$

\begin{rmk}
{\rm
The set $\Delta^V({\cal G})\cup\Delta^V({\cal H})$ in \eqref{ADDN} can be reduced to the set
$$
\Psi({\cal G},{\cal H}):=\{\Delta^V({\cal G})\backslash
\Lambda^V({\cal G}\backslash{\cal H})\}
\cup\{\Delta^V({\cal H})\backslash\Lambda^V({\cal H}\backslash {\cal G})\}.
$$
Actually, it follows from \eqref{wpwp1} that
$
\wp(\Psi({\cal G},{\cal H}))\subset
\wp(\Delta^V({\cal G})\cup\Delta^V({\cal H}))=\wp(\Delta({\cal G}+{\cal H})).
$
On the other hand,
we have ${\cal V}_*\subset\Psi({\cal G},{\cal H})$ because of formula \eqref{VPVV},
and therefore, we see from \eqref{wpwp1} that
$
\wp(\Delta({\cal G}+{\cal H}))=\wp({\cal V}_*)\subset\wp(\Psi({\cal G},{\cal H})).
$
Consequently, $\wp(\Delta({\cal G}+{\cal H}))=\wp(\Psi({\cal G},{\cal H}))$,
implying that ${\cal N}_{{\cal G}+{\cal H}}={\cal N}(\Psi({\cal G},{\cal H}))$.
}
\label{rmk:add}
\end{rmk}


\begin{figure}[!h]
    \centering
    \subcaptionbox{%
     }{\includegraphics[height=1.6in]{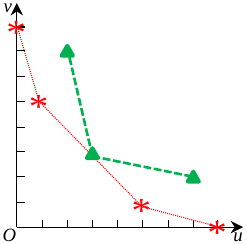}}~~~~~~~~~~~~
    \subcaptionbox{ %
     }{\includegraphics[height=1.6in]{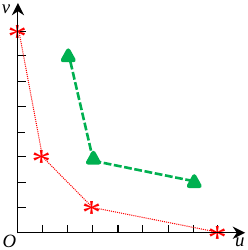}}
    \caption{(a) ${\cal N}_{\cal G}\le {\cal N}_{\cal H}$,
    (b) ${\cal N}_{\cal G}< {\cal N}_{\cal H}$, with
    dots for ${\cal N}_{\cal G}$ and dashes for ${\cal N}_{\cal H}$.
    }
    \label{fletf}
\end{figure}


Conclusions in Proposition~\ref{ADD} can be simplified in some special cases.
Consider prolongations
${\cal N}^\infty_{\cal G}:={\cal N}_{\cal G}\cup L^h_{\cal G}$ and
${\cal N}^\infty_{\cal H}:={\cal N}_{\cal H}\cup L^h_{\cal H}$,
where $L^h_{\cal G}$ and $L^h_{\cal H}$ are the horizontal boundaries of
$\wp(\Delta({\cal G}))$ and $\wp(\Delta({\cal H}))$ respectively.
Suppose that their expressions are
$v={\cal N}^\infty_{\cal G}(u)$ for $u\in[i_0,+\infty)$ and
$v={\cal N}^\infty_{\cal H}(u)$ for $u\in[\tilde{i}_0,+\infty)$ respectively,
where $i_0$ and $\tilde{i}_0$ are abscissas of left-most vertices of
${\cal N}_{\cal G}$ and ${\cal N}_{\cal H}$ respectively.
We say that ${\cal N}_{{\cal H}}$ {\it lies above} ${\cal N}_{\cal G}$,
simply denoted by ${\cal N}_{\cal G}\le {\cal N}_{\cal H}$,
if $i_0\le \tilde{i}_0 $ and
\begin{eqnarray}
{\cal N}^\infty_{{\cal G}}(u)\le {\cal N}^\infty_{{\cal H}}(u)~
\mbox{for all}~u\in [\tilde{i}_0,+\infty),
\label{NfleNtf}
\end{eqnarray}
as shown in Figure~\ref{fletf}(a).
Further, if either $\tilde{i}_0>i_{s({\cal G})}$ or
the inequality in \eqref{NfleNtf} is strictly $<$
for all $u\in[\tilde{i}_0,i_{s({\cal G})}]$,
as shown in Figure~\ref{fletf}(b),
then we call ${\cal N}_{\cal H}$ {\it strictly lies above} ${\cal N}_{\cal G}$,
denoted by ${\cal N}_{\cal G}<{\cal N}_{\cal H}$.

\begin{cor}
Suppose that ${\cal G}$ and ${\cal H}$ are both nonzero analytic functions
of the form \eqref{expGH}
and their common vertices satisfy the non-vanishing condition \eqref{Anon} for addition,
i.e., $\Lambda^V({\cal G}\backslash{\cal H})\cap
\Lambda^V({\cal H}\backslash {\cal G})=\emptyset$.
If ${\cal N}_{\cal G}\le {\cal N}_{\cal H}$ then
${\cal N}_{{\cal G}+{\cal H}}={\cal N}_{\cal G}$.
Further, $({\cal G}+{\cal H})_{E_k({\cal G}+{\cal H})}={\cal G}_{E_k({\cal G})}+{\cal H}_{E_k({\cal G})}$ for all $k=1,...,s({\cal G})$,
where $s({\cal G})$ is the number of edges of ${\cal N}_{\cal G}$.
Especially, if ${\cal N}_{\cal G}< {\cal N}_{\cal H}$ then
$({\cal G}+{\cal H})_{E_k({\cal G}+{\cal H})}={\cal G}_{E_k({\cal G})}$
for all $k=1,...,s({\cal G})$.
\label{cor-ADD}
\end{cor}

{\bf Proof.}
If ${\cal N}_{\cal G}\le{\cal N}_{\cal H}$ then $\wp(\Delta({\cal H}))\subset\wp(\Delta({\cal G}))$.
Similarly to the computation given just below \eqref{wpwp1}, we have
\begin{eqnarray*}
\wp(\Delta^V({\cal G})\cup\Delta^V({\cal H}))
&=&{\rm conv}\{{\rm conv}\{\Delta^V({\cal G})+{\mathbb R}_+^2\}
\cup{\rm conv}\{\Delta^V({\cal H})+{\mathbb R}_+^2\}\}
\\
&=&{\rm conv}\{\wp(\Delta({\cal G}))\cup\wp(\Delta({\cal H}))\}.
\\
&=&{\rm conv}\{\wp(\Delta({\cal G}))\}
\\
&=&\wp(\Delta({\cal G})).
\end{eqnarray*}
Then, we see from \eqref{ADDN} that
$
{\cal N}_{{\cal G}+{\cal H}}={\cal N}(\Delta^V({\cal G})\cup\Delta^V({\cal H}))={\cal N}(\wp({\cal H}))={\cal N}_{\cal G}.
$
It follows that $\vec{\cal N}(\Delta^V({\cal G})\cup\Delta^V({\cal H}))$ has exactly $s({\cal G})$ edges and
its $k$-th edge $N_k$ is $E_k({\cal G})$.
By \eqref{ADDfE},
$
({\cal G}+{\cal H})_{E_k({\cal G}+{\cal H})}
={\cal G}_{E_k({\cal G})}+{\cal H}_{E_k({\cal G})}~\mbox{for all}~k=1,...,s({\cal G}).
$
Especially, if ${\cal N}_{\cal G}<{\cal N}_{\cal H}$
then $E_k({\cal G})\cap\Delta({\cal H})\subset
{\cal N}_{\cal G}\cap\wp(\Delta({\cal H}))=\emptyset$.
By definition \eqref{fonE}, ${\cal H}_{E_k({\cal G})}=0$ in the above equality.
Thus, this corollary is proved.
\qquad$\Box$

\section{Proofs for main theorems}
\setcounter{equation}{0}
\setcounter{lm}{0}
\setcounter{thm}{0}
\setcounter{rmk}{0}
\setcounter{df}{0}
\setcounter{cor}{0}

In order to prove Theorems \ref{th:finite} and \ref{th:finiteJ2},
we first show that if
either conditions {\bf(P1)} and {\bf(Q)} hold in the one-above case {\bf (J1)},
or conditions {\bf(P2)}, {\bf(Q)}, {\bf(S)}
and one of conditions {\bf (H1)}, {\bf (H2)}, {\bf (H3)} and {\bf(H4)}
of Theorem~\ref{th:finiteJ2} hold in the two-below case {\bf (J2)},
then the function
$$
f:=[{\cal G},{\cal H}]_\theta={\cal G}'_\theta{\cal H}-{\cal G}{\cal H}'_\theta,
$$
the Lie-bracket of ${\cal G}$ and ${\cal H}$ in the variable $\theta$,
is semi-definite on the region
$$
\Omega_\epsilon:=\{(\rho,\theta)\in\mathbb{R}^2:0<\rho<\epsilon,|\theta|<\epsilon\}
$$
for an $\epsilon>0$.
For this purpose,
as indicated in section~3.3,
we need to know its Newton polygon and edge-polynomials,
which will be discussed in the following subsection.

\subsection{Newton polygon of the Lie-bracket $[{\cal G},{\cal H}]_\theta$}

\subsubsection{The one-above case}

\begin{lm}
In the one-above case {\bf (J1)},
suppose that
common vertices of ${\cal G}'_\theta{\cal H}$ and $-{\cal G}{\cal H}'_\theta$ satisfy
the non-vanishing condition~\eqref{Anon} for addition.
Then the Newton polygon ${\cal N}_f$ of the
{\rm Lie}-bracket
$f:=[{\cal G},{\cal H}]_\theta$ of ${\cal G}$ and ${\cal H}$ in the variable $\theta$
is equal to ${\cal N}_{{\cal G}{\cal H}}-(0,1)$,
which has $s({\cal G},{\cal H})$ edges,
and its $k$-th edge-polynomial $f_{E_k(f)}$ is equal to
$[{\cal K}_{\xi_k}({\cal G}),{\cal K}_{\xi_k}({\cal H})]$
for all $k=1,...,s({\cal G},{\cal H})$,
where
${\cal K}_{\xi_k}({\cal G})$ and ${\cal K}_{\xi_k}({\cal H})$
are $\xi_k$-componential polynomials,
$\xi_k$ and $s({\cal G},{\cal H})$ are defined in \eqref{unionseq}.
\label{J1Ng}
\end{lm}

\noindent{\bf Proof.}
Note that $j_{s({\cal G})}+\tilde{j}_{s({\cal H})}>0$ in {\bf(J1)}.
We only give the proof in the case $j_{s({\cal G})},\tilde{j}_{s({\cal H})}>0$,
and the case $j_{s({\cal G})}>0$ and $\tilde{j}_{s({\cal H})}=0$,
because the remaining case $j_{s({\cal G})}=0$ and $\tilde{j}_{s({\cal H})}>0$
is similar to the second case.
In the first case,
Remark~\ref{rmk:diff}{\bf(i)} shows that
${\cal N}_{\theta{\cal G}'_\theta}={\cal N}_{\cal G}$ and ${\cal N}_{\theta{\cal H}'_\theta}={\cal N}_{\cal H}$.
By Proposition~\ref{lm-MIL},
${\cal N}_{\theta{\cal G}'_\theta{\cal H}}
={\cal N}_{{\cal G}\theta{\cal H}'_\theta}
={\cal N}_{{\cal G}{\cal H}}$ and therefore
${\cal N}_{{\cal G}'_\theta{\cal H}}
={\cal N}_{-{\cal G}{\cal H}'_\theta}={\cal N}_{{\cal G}{\cal H}}-(0,1)$.
Under the assumption that common vertices of ${\cal G}'_\theta{\cal H}$ and $-{\cal G}{\cal H}'_\theta$ satisfy
the non-vanishing condition~\eqref{Anon} for addition,
Corollary~\ref{cor-ADD} shows that
$$
{\cal N}_f={\cal N}_{{\cal G}'_\theta{\cal H}}={\cal N}_{-{\cal G}{\cal H}'_\theta}={\cal N}_{{\cal G}{\cal H}}-(0,1).
$$
By Remark~\ref{rmk:mul}{\bf(i)},
if $s({\cal G},{\cal H})>0$,
then ${\cal N}_f$ has $s({\cal G},{\cal H})$ edges
and $\xi_k$ is the slope of the $k$-th edge $E_k({\cal GH})$ of ${\cal N}_{{\cal GH}}$.
Then, we see from the above equalities and definition \eqref{defscomp} of $\xi$-component that
$E_k(f)=E_k({\cal G}'_\theta{\cal H})=E_k({\cal G}{\cal H}'_\theta)$,
$E_k({\cal G}'_\theta{\cal H})=\sigma({\cal G}'_\theta{\cal H},\xi_k)$ and
$E_k({\cal G}{\cal H}'_\theta)=\sigma({\cal G}{\cal H}'_\theta,\xi_k)$
for all $k=1,...,s({\cal G},{\cal H})$.
By Propositions~\ref{lm-Diff}, \ref{lm-MIL} and \ref{ADD}, Remark~\ref{rmk:diff},
definition \eqref{defkxip} of $\xi$-componential polynomial,
and definition \eqref{fonE} of $\Gamma$-polynomial, we compute that
\begin{eqnarray}
f_{E_k(f)}&=&({\cal G}'_\theta{\cal H})_{E_k(f)}-({\cal G}{\cal H}'_\theta)_{E_k(f)}
\nonumber
\\
&=&({\cal G}'_\theta{\cal H})_{E_k({\cal G}'_\theta{\cal H})}-({\cal G}{\cal H}'_\theta)_{E_k({\cal G}{\cal H}'_\theta)}
\nonumber
\\
&=&{\cal K}_{\xi_k}({\cal G}'_\theta{\cal H})-{\cal K}_{\xi_k}({\cal G}{\cal H}'_\theta)
\nonumber
\\
&=&{\cal K}_{\xi_k}({\cal G}'_\theta) {\cal K}_{\xi_k}({\cal H})
-{\cal K}_{\xi_k} ({\cal G}) {\cal K}_{\xi_k}({\cal H}'_\theta)
\nonumber
\\
&=&({\cal K}_{\xi_k} ({\cal G}))'_\theta {\cal K}_{\xi_k}({\cal H})
-{\cal K}_{\xi_k}({\cal G}) ({\cal K}_{\xi_k}({\cal H}))'_\theta
\nonumber
\\
&=&[{\cal K}_{\xi_k}({\cal G}),{\cal K}_{\xi_k}({\cal H})].
\label{gEkgp1}
\end{eqnarray}
Thus, this lemma holds in the case
$j_{s({\cal G})},\tilde{j}_{s({\cal H})}>0$.

In the case $j_{s({\cal G})}>0$ and $\tilde{j}_{s({\cal H})}=0$,
since $\Delta(\theta{\cal H}'_\theta)=\Delta({\cal H}'_\theta)+(0,1)\subset\Delta({\cal H})$,
we have $\wp(\Delta(\theta{\cal H}'_\theta))\subset\wp(\Delta({\cal H}))$.
Further, the following inclusion
$$
\wp(\Delta({\cal G}\theta{\cal H}'_\theta))=\wp(\Delta({\cal G}))+\wp(\Delta(\theta{\cal H}'_\theta))
\subset\wp(\Delta({\cal G}))+\wp(\Delta({\cal H}))=\wp(\Delta({\cal G}{\cal H}))
$$
holds,
where we used the equality
$\wp(\Delta({\cal G}{\cal H}))=\wp(\Delta({\cal G}))+\wp(\Delta({\cal H}))$
given by \cite[Proposition~2.2]{Lip88}.
Then, by definition given just before Corollary~\ref{cor-ADD},
${\cal N}_{{\cal G}{\cal H}}\le{\cal N}_{{\cal G}\theta{\cal H}'_\theta}$.
Similar to the above case,
${\cal N}_{\theta{\cal G}'_\theta{\cal H}}={\cal N}_{{\cal G}{\cal H}}$
since $j_{s({\cal G})}>0$.
It follows that
${\cal N}_{\theta{\cal G}'_\theta{\cal H}}={\cal N}_{{\cal G}{\cal H}}\le{\cal N}_{{\cal G}\theta{\cal H}'_\theta}
={\cal N}_{-{\cal G}\theta{\cal H}'_\theta}$ and therefore
${\cal N}_{{\cal G}'_\theta{\cal H}}\le{\cal N}_{-{\cal G}{\cal H}'_\theta}$.
Under the assumption that common vertices of
${\cal G}'_\theta{\cal H}$ and $-{\cal G}{\cal H}'_\theta$ satisfy
the non-vanishing condition \eqref{Anon} for addition,
Corollary~\ref{cor-ADD} yields the first part of this lemma, i.e.,
\begin{eqnarray}
{\cal N}_f={\cal N}_{{\cal G}'_\theta{\cal H}}={\cal N}_{{\cal G}{\cal H}}-(0,1).
\label{lmJ1NG}
\end{eqnarray}
By Remark~\ref{rmk:mul}{\bf(i)}, ${\cal N}_f$ has $s({\cal G},{\cal H})$ edges.
In order to obtain the edge-polynomials as $s({\cal G},{\cal H})>0$,
we need to know the more specific relative position of
${\cal N}_{-{\cal G}{\cal H}'_\theta}$ and ${\cal N}_{{\cal G}'_\theta{\cal H}}$.
For this purpose, we claim that
\begin{eqnarray}
\Delta^V({\cal G}\theta{\cal H}'_\theta)\cap\Delta^V({\cal G}{\cal H})=\{V_0({\cal G}{\cal H}),...,V_{k_*-1}({\cal G}{\cal H})\}
\label{J1Vcom}
\end{eqnarray}
and that if $k_*<s({\cal G},{\cal H})$ then
\begin{eqnarray}
\Delta({\cal G}\theta{\cal H}'_\theta)\cap E_k({\cal G}{\cal H})=\emptyset~~~
\mbox{for all}~k=k_*+1,...,s({\cal G},{\cal H}),
\label{DyfdtfEk}
\end{eqnarray}
where $k_*$ is either zero as $s({\cal H})=0$,
or the integer such that $\xi_{k_*}=\zeta(E_{s({\cal H})}({\cal H}))$ as $s({\cal H})>0$.

In order to prove the above claim,
we consider the function $\psi:=\theta{\cal H}'_\theta+{\cal H}$.
Clearly, ${\cal N}_{{\cal H}}\le{\cal N}_{\theta{\cal H}'_\theta}$ and their common vertices satisfy
the non-vanishing condition~\eqref{Anon} for addition.
By Corollary~\ref{cor-ADD}, ${\cal N}_{\psi}={\cal N}_{{\cal H}}$ and
\begin{eqnarray}
a_{V_k(\psi)}(\psi)=(j(V_k({\cal H}))+1)a_{V_k({\cal H})}({\cal H})~\mbox{for all}~k=0,...,s({\cal H}),
\label{ahfatf}
\end{eqnarray}
where $a_{V_k(\psi)}(\psi)$ denotes the coefficient of $\psi$
corresponding to the vertex $V_k(\psi)$ of ${\cal N}_\psi$
and $j(V_k({\cal H}))$ denotes the ordinate of the vertex $V_k({\cal H})$.
We prove the above claim in the situation $k_*=0$,
the situation $k_*\in(0,s({\cal G},{\cal H}))$ and
the situation $k_*=s({\cal G},{\cal H})$ separately.

In the situation $k_*=0$,
we have $s({\cal H})=0$, i.e.,
${\cal N}_{\cal H}$ is a singleton $V_{s({\cal H})}({\cal H})$.
Then,
$\xi_k\in(\zeta^-(V_{s({\cal H})}({\cal H})),\zeta^+(V_{s({\cal H})}({\cal H})))=(-\infty,0)$
by \eqref{SP+-} for all $k=k_*+1,...,s({\cal G},{\cal H})$.
We see from the equality ${\cal N}_\psi={\cal N}_{{\cal H}}$
and definition \eqref{defscomp} of $\xi$-component that
\begin{eqnarray}
\sigma(\psi,\xi_k)=V_{s({\cal H})}(\psi)=V_{s({\cal H})}({\cal H})
=\sigma({\cal H},\xi_k)~\mbox{for all}~k=k_*+1,...,s({\cal G},{\cal H}).
\label{sVVs}
\end{eqnarray}
By \eqref{ahfatf}, \eqref{sVVs}, Proposition~\ref{lm-MIL} and
definition~\eqref{defkxip} of $\xi$-componential polynomial,
$$
{\cal K}_{\xi_k}({\cal G}\psi)
={\cal K}_{\xi_k}({\cal G}){\cal K}_{\xi_k}(\psi)
={\cal K}_{\xi_k}({\cal G})a_{V_{s({\cal H})}(\psi)}(\psi)
={\cal K}_{\xi_k}({\cal G})a_{V_{s({\cal H})}({\cal H})}({\cal H})
={\cal K}_{\xi_k}({\cal G}{\cal H}),
$$
where we used the fact that $j(V_{s({\cal H})}(\psi))=j(V_{s({\cal H})}({\cal H}))=\tilde{j}_{s({\cal H})}=0$.
On the other hand,
by the equality ${\cal N}_{\psi}={\cal N}_{\cal H}$ and Proposition~\ref{lm-MIL},
we have ${\cal N}_{{\cal G}\psi}={\cal N}_{{\cal G}{\cal H}}$.
Noticing that $\xi_k$ is the slope of edge $E_k({\cal GH})$,
we can rewrite ${\cal K}_{\xi_k}({\cal G}\psi)$ as
$$
{\cal K}_{\xi_k}({\cal G}\psi)
=({\cal G}\theta{\cal H}'_\theta+{\cal G}{\cal H})_{E_k({\cal G}{\cal H})}
=({\cal G}\theta{\cal H}'_\theta)_{E_k({\cal G}{\cal H})}+{\cal K}_{\xi_k}({\cal G}{\cal H}).
$$
It follows from the above two equalities about ${\cal K}_{\xi_k}({\cal G}\psi)$ that
\begin{eqnarray}
({\cal G}\theta{\cal H}'_\theta)_{E_k({\cal G}{\cal H})}=0
~~~\mbox{for all}~k=k_*+1,...,s({\cal G},{\cal H})
\label{GHEGH0}
\end{eqnarray}
and therefore \eqref{J1Vcom} and \eqref{DyfdtfEk} hold in the case $k_*=0$.

In the situation $k_*\in(0,s({\cal G},{\cal H}))$,
equality~\eqref{sVVs} still holds because
$$
\zeta^-(V_{s({\cal H})}({\cal H}))
=\zeta(E_{s({\cal H})}({\cal H}))
=\xi_{k_*}
<\xi_k<0
=\zeta^+(V_{s({\cal H})}({\cal H}))
$$
for all $k=k_*+1,...,s({\cal G},{\cal H})$.
The same computation of ${\cal K}_{\xi_k}({\cal G}\psi)$ as above shows that
\eqref{GHEGH0} still holds,
which implies \eqref{DyfdtfEk}.
In order to prove \eqref{J1Vcom},
for each $k=0,...,k_*-1$,
we choose $\nu_k\in(\zeta^-(V_k({\cal GH})),\zeta^+(V_k({\cal GH})))$.
Clearly, $\nu_k\notin\mathfrak{S}({\cal G})\cup\mathfrak{S}({\cal H})$
because of Remark~\ref{rmk:mul}{\bf(i)}.
Note that
\begin{align*}
\zeta^+(V_{s({\cal H})-1}({\cal H}))
&=\zeta(E_{s({\cal H})}({\cal H}))
=\xi_{k_*}
=\zeta(E_{k_*}({\cal GH}))
=\zeta^+(V_{k_*-1}({\cal GH}))
\\
&\ge\zeta^+(V_k({\cal GH}))
>\nu_k.
\end{align*}
Then, there is an $n(k)\in\{0,...,s({\cal H})-1\}$ such that
$\nu_k\in(\zeta^-(V_{n(k)}({\cal H})),\zeta^+(V_{n(k)}({\cal H})))$,
We see from \eqref{defscomp} that
$\sigma(\psi,\nu_k)=V_{n(k)}(\psi)=V_{n(k)}({\cal H})=\sigma({\cal H},\nu_k)$.
Similar to the computation given just below \eqref{sVVs},
\begin{align*}
{\cal K}_{\nu_k}({\cal G}\psi)
&={\cal K}_{\nu_k}({\cal G}){\cal K}_{\nu_k}(\psi)
\\
&={\cal K}_{\nu_k}({\cal G})a_{V_{n(k)}(\psi)}(\psi)\theta^{j(V_{n(k)}(\psi))}
\\
&={\cal K}_{\nu_k}({\cal G})(j(V_{n(k)}({\cal H}))+1)
a_{V_{n(k)}({\cal H})}({\cal H})\theta^{j(V_{n(k)}({\cal H}))}
\\
&=(j(V_{n(k)}({\cal H}))+1){\cal K}_{\nu_k}({\cal GH}).
\end{align*}
On the other hand, similar to the computation given just before \eqref{GHEGH0},
$$
{\cal K}_{\nu_k}({\cal G}\psi)
=({\cal G}\theta{\cal H}'_\theta)_{\sigma({\cal GH},\nu_k)}
+{\cal K}_{\nu_k}({\cal GH})
=a_{V_k({\cal GH})}({\cal G}\theta{\cal H}'_\theta)\theta^{j(V_k({\cal GH}))}
+{\cal K}_{\nu_k}({\cal GH}).
$$
Note that $j(V_{n(k)}({\cal H}))\ne0$ since $n(k)\le s({\cal H})-1$.
Then, the above two equalities about ${\cal K}_{\nu_k}({\cal G}\psi)$ indicate that
$a_{V_k({\cal GH})}({\cal G}\theta{\cal H}'_\theta)\ne0$,
i.e., $V_k({\cal GH})\in\Delta({\cal G}\theta{\cal H}'_\theta)$.
Since ${\cal N}_{\cal GH}\le {\cal N}_{{\cal G}\theta{\cal H}'_\theta}$,
we have
$\{V_0({\cal GH}),...,V_{k_*-1}({\cal GH})\}
\subset\Delta^V({\cal N}_{{\cal G}\theta{\cal H}'_\theta})$.
Then, \eqref{J1Vcom} follows form \eqref{DyfdtfEk}.

In the situation $k_*=s({\cal G},{\cal H})$,
similar to the above,
$\{V_0({\cal GH}),...,V_{k_*-1}({\cal GH})\}\subset
\Delta^V({\cal N}_{{\cal G}\theta{\cal H}'_\theta})$.
It suffices to show that $V_{k_*}({\cal GH})\notin\Delta({\cal G}\theta{\cal H}'_\theta)$.
Noticing that
$\xi_{k_*}=\zeta(E_{s({\cal H})}({\cal H}))=\zeta^-(V_{s({\cal H})}({\cal H}))$ and
$\xi_{k_*}=\zeta(E_{s({\cal G},{\cal H})}({\cal GH}))=\zeta^-(V_{s({\cal G},{\cal H})}({\cal GH}))$,
by \eqref{SP+-} and \eqref{defscomp},
we can choose $\nu_{k_*}\in(\xi_{k_*},0)$ such that
$$
\sigma(\psi,\nu_{k_*})=V_{s({\cal H})}(\psi)
=V_{s({\cal H})}({\cal H})=\sigma({\cal H},\nu_{k_*}),
~
\sigma({\cal G}\psi,\nu_{k_*})=\sigma({\cal G}{\cal H},\nu_{k_*})=V_{k_*}({\cal G}{\cal H}).
$$
Similar to the computation given in the above situation,
$$
{\cal K}_{\nu_{k_*}}({\cal G}\psi)
={\cal K}_{\nu_{k_*}}({\cal G}{\cal H}),
~~~
{\cal K}_{\nu_{k_*}}({\cal G}\psi)
=a_{V_{k_*}({\cal G}{\cal H})}({\cal G}\theta{\cal H}'_\theta)\theta^{j(V_{k_*}({\cal G}{\cal H}))}+{\cal K}_{\nu_{k_*}}({\cal G}{\cal H}),
$$
where we used the fact that $j(V_{s({\cal H})}({\cal H}))=0$ in the first equality.
It follows that $a_{V_{k_*}({\cal G}{\cal H})}({\cal G}\theta{\cal H}'_\theta)=0$,
i.e., $V_{k_*}({\cal G}{\cal H})\notin\Delta({\cal G}\theta{\cal H}'_\theta)$.
Consequently, \eqref{J1Vcom} and \eqref{DyfdtfEk} hold when $k_*=s({\cal G},{\cal H})$.
Thus, the claim given just below \eqref{lmJ1NG} is proved.

Having \eqref{J1Vcom} and \eqref{DyfdtfEk},
we further discuss those edge-polynomials of $[{\cal G},{\cal H}]_\theta$ in the second case,
i.e., $j_{s({\cal G})}>0$ and $\tilde{j}_{s({\cal H})}=0$.
We only consider the circumstance $s({\cal H})>0$,
and the circumstance $s({\cal H})=0$ can be similarly discussed.
When $s({\cal H})>0$,
for each $k=1,...,k_*-1$,
we have $\xi_k<\xi_{k_*}=\zeta(E_{s({\cal H})}({\cal H}))$,
and $E_k(f)=E_k({\cal G}'_\theta{\cal H})=E_k({\cal G}{\cal H}'_\theta)$ because of (\ref{lmJ1NG}) and (\ref{J1Vcom}).
Then, we can compute similarly to (\ref{gEkgp1}) that
$f_{E_k(f)}=[{\cal K}_{\xi_k}({\cal G}),{\cal K}_{\xi_k}({\cal H})]$.

For $k=k_*$,
we claim that
\begin{eqnarray}
\Delta({\cal G}{\cal H}'_\theta)\cap E_{k_*}(f)=
\Delta({\cal G}{\cal H}'_\theta)\cap \sigma({\cal G}{\cal H}'_\theta,\xi_{k_*}).
\label{equ:GHEGHs}
\end{eqnarray}
Actually,
the inequality
${\cal N}_f={\cal N}_{{\cal G}'_\theta{\cal H}}\le{\cal N}_{{\cal G}{\cal H}'_\theta}$
given just before \eqref{lmJ1NG} ensues that
\begin{eqnarray}
\wp(\Delta({\cal G}{\cal H}'_\theta))\subset\wp(\Delta(f)).
\label{equ:GHwpf}
\end{eqnarray}
Moreover, we see from \eqref{lmJ1NG} and \eqref{J1Vcom} that
$V_{k_*-2}(f):=V_{k_*-2}({\cal GH})-(0,1)$ and
$V_{k_*-1}(f):=V_{k_*-1}({\cal GH})-(0,1)$
are the ($k_*-2$)-th vertex and the ($k_*-1$)-th vertex of
${\cal N}_{{\cal G}{\cal H}'_\theta}$ and ${\cal N}_f$.
Then, \eqref{equ:GHwpf} implies that either
$$
V_{k_*}({\cal G}{\cal H}'_\theta)\in E_{k_*}(f)=E_{k_*}({\cal GH})-(0,1)
~~~\mbox{or}~~~
V_{k_*}({\cal G}{\cal H}'_\theta)\in {\rm int}\wp(\Delta(f)),
$$
where ${\rm int}$ is the interior of a set.
In the first situation,
since $\zeta(\underline{V_{k_*-1}(f)V_{k_*}({\cal G}{\cal H}'_\theta)})
=\zeta(E_{k_*}({\cal GH}))=\xi_{k_*}$,
we see from \eqref{defscomp} that
$\sigma({\cal G}{\cal H}'_\theta,\xi_{k_*})
=\underline{V_{k_*-1}(f)V_{k_*}({\cal G}{\cal H}'_\theta)}
\subset E_{k_*}(f)$.
By the definition of Newton polygon,
$\Delta({\cal G}{\cal H}'_\theta)\cap
\{E_{k_*}(f)\setminus \sigma({\cal G}{\cal H}'_\theta,\xi_{k_*})\}
=\emptyset$,
implying that \eqref{equ:GHEGHs} holds in the first situation.
In the second situation,
we have $\zeta(\underline{V_{k_*-2}(f)V_{k_*-1}(f)})
<\xi_{k_*}<\zeta(\underline{V_{k_*-1}(f)V_{k_*}({\cal G}{\cal H}'_\theta)})$.
By \eqref{defscomp},
$\sigma({\cal G}{\cal H}'_\theta,\xi_{k_*})=V_{k_*-1}(f)$.
Moreover, $\Delta({\cal G}{\cal H}'_\theta)\cap E_{k_*}(f)=V_{k_*-1}(f)$
since the slope of $E_{k_*}(f)$ is $\xi_{k_*}$ and its left end-point is $V_{k_*-1}(f)$.
Thus, \eqref{equ:GHEGHs} holds in the second situation
and therefore the claimed \eqref{equ:GHEGHs} is proved.
Then,
by definitions \eqref{defkxip} and \eqref{fonE},
Propositions \ref{lm-Diff} and \ref{lm-MIL},
we compute that
\begin{eqnarray*}
({\cal G}{\cal H}'_\theta)_{E_{k_*}(f)}
&=&\sum_{(i,j)\in \Delta({\cal G}{\cal H}'_\theta)\cap E_{k_*}(f)}a_{i,j}({\cal G}{\cal H}'_\theta)\theta^j
\\
&=&\sum_{(i,j)\in \Delta({\cal G}{\cal H}'_\theta)\cap \sigma({\cal G}{\cal H}'_\theta,\xi_{k_*})}a_{i,j}({\cal G}{\cal H}'_\theta)\theta^j
\\
&=&{\cal K}_{\xi_{k_*}}({\cal G}{\cal H}'_\theta)
\\
&=&{\cal K}_{\xi_{k_*}}({\cal G})({\cal K}_{\xi_{k_*}}({\cal H}))'_\theta.
\end{eqnarray*}
Further, similarly to (\ref{gEkgp1}),
we obtain that $f_{E_{k_*}(f)}=[{\cal K}_{\xi_{k_*}}{\cal G},{\cal K}_{\xi_{k_*}}{\cal H}]$.

Finally,
for each $k=k_*+1,...,s({\cal G},{\cal H})$,
\eqref{lmJ1NG} and \eqref{DyfdtfEk} imply that $({\cal G}{\cal H}'_\theta)_{E_k(f)}=0$.
Similar to (\ref{gEkgp1}),
$f_{E_k(f)}=({\cal K}_{\xi_k}({\cal G}))'_\theta{\cal K}_{\xi_k}({\cal H})$.
Since $\xi_k>\xi_{k_*}=\zeta(E_{s({\cal H})}({\cal H}))$,
\eqref{defscomp} shows that
$\sigma({\cal H},\xi_k)=V_{s({\cal H})}({\cal H})$.
Moreover, equality $j(V_{s({\cal H})}({\cal H}))=\tilde{j}_{s({\cal H})}=0$
implies by \eqref{defkxip} that
${\cal K}_{\xi_k}({\cal H})=a_{V_{s({\cal H})}({\cal H})}({\cal H})$
and therefore, $({\cal K}_{\xi_k}({\cal H}))'_\theta=0$.
Then,
$$
f_{E_k(f)}
=({\cal K}_{\xi_k} ({\cal G}))'_\theta{\cal K}_{\xi_k}({\cal H})
-{\cal K}_{\xi_k} ({\cal G})({\cal K}_{\xi_k}({\cal H}))'_\theta
=[{\cal K}_{\xi_k}({\cal G}),{\cal K}_{\xi_k}({\cal H})],
$$
which proves the conclusion on edge-polynomials.
Consequently, the proof of this lemma is completed.
\qquad$\Box$

\subsubsection{The two-below case}

\begin{lm}
In the two-below case {\bf (J2)},
suppose that condition {\bf(S)}, given just before Theorem \ref{th:finiteJ2}, holds,
and common vertices of ${\cal G}'_\theta{\cal H}$ and $-{\cal G}{\cal H}'_\theta$
satisfy the non-vanishing condition \eqref{Anon} for addition.
Let $f:=[{\cal G},{\cal H}]_\theta$ as defined in Lemma~\ref{J1Ng}.
Then
\\
{\bf(i)}
in the case $\zeta(E_{s({\cal G})}({\cal G}))>\zeta(E_{s({\cal H})}({\cal H}))$ and $j_*\le \tilde{j}_*$,
\begin{eqnarray}
{\cal N}_f=\{{\cal N}_{{\cal G}{\cal H}}
\backslash\underline{Q_*V_{s({\cal G},{\cal H})}({\cal G}{\cal H})}\}\cup Q_*-(0,1),
\label{J2iNg}
\end{eqnarray}
where $Q_*:=(i_*+\tilde{i}_{s({\cal H})},j_*)$,
$s({\cal G},{\cal H})$ is given in \eqref{unionseq},
$V_{s({\cal G},{\cal H})}({\cal G}{\cal H})$
is the $s({\cal G},{\cal H})$-th vertex of ${\cal N}_{\cal GH}$,
and $\underline{Q_*V_{s({\cal G},{\cal H})}({\cal G}{\cal H})}$ is
the closed linear segment linking
$Q_*$ with $V_{s({\cal G},{\cal H})}({\cal G}{\cal H})$.
Moreover, $s(f)$, the number of edges of ${\cal N}_f$,
is equal to $s({\cal G},{\cal H})-1+{\rm sgn}(j_{s({\cal G})-1}-j_*)$,
and the edge-polynomial of $f_{E_k(f)}$ satisfies that
\begin{eqnarray}
f_{E_k(f)}=[{\cal K}_{\xi_k}({\cal G}),{\cal K}_{\xi_k}({\cal H})]~\mbox{for all}~
k=1,...,s({\cal G},{\cal H})-1+{\rm sgn}(j_{s({\cal G})-1}-j_*),
\label{J2igEkg}
\end{eqnarray}
where $[\cdot,\cdot]$ is the {\rm Lie}-bracket,
${\cal K}_{\xi_k}({\cal G})$ and ${\cal K}_{\xi_k}({\cal H})$ are
$\xi_k$-componential polynomials,
$\xi_k$ is defined in \eqref{unionseq},
and ${\rm sgn}$ is the signum function;
\\
{\bf(ii)}
in the case $\zeta(E_{s({\cal G})}({\cal G}))>\zeta(E_{s({\cal H})}({\cal H}))$ and $j_*>\tilde{j}_*$,
\begin{eqnarray}
{\cal N}_f=\{{\cal N}_{{\cal G}{\cal H}}\backslash
\underline{Q_*V_{s({\cal G},{\cal H})}({\cal G}{\cal H})}\}
\cup \underline{Q_*\tilde{Q}_*}-(0,1),
\label{J2iiNg}
\end{eqnarray}
where $\tilde{Q}_*:=(i_{s({\cal G})}+\tilde{i}_*,\tilde{j}_*)$.
Moreover, $s(f)=s({\cal G},{\cal H})+{\rm sgn}(j_{s({\cal G})-1}-j_*)$,
and edge-polynomials of $f$ are given by \eqref{J2igEkg} and
\begin{eqnarray}
f_{E_{s(f)}(f)}
=j_*a_{i_*,j_*}({\cal G})
a_{\tilde{i}_{s({\cal H})},\tilde{j}_{s({\cal H})}}({\cal H})\theta^{j_*-1}
-\tilde{j}_*a_{i_{s({\cal G})},j_{s({\cal G})}}({\cal G})
a_{\tilde{i}_*,\tilde{j}_*}({\cal H})\theta^{\tilde{j}_*-1},
\label{gElast}
\end{eqnarray}
{\bf(iii)}
in the case $\zeta(E_{s({\cal G})}({\cal G}))=\zeta(E_{s({\cal H})}({\cal H}))$ and
$j_*\le \tilde{j}_*$,
$f$ has the same Newton polygon ${\cal N}_f$ as in \eqref{J2iNg}.
Moreover, $s(f)=s({\cal G},{\cal H})$ and
$f_{E_k(f)}=[{\cal K}_{\xi_k}({\cal G}),{\cal K}_{\xi_k}({\cal H})]$
for all $k=1,...,s({\cal G},{\cal H})$.
\label{precondiJ2}
\end{lm}

\begin{rmk}
{\rm
In the opposite cases, i.e.,
{\bf(i$'$)}
$\zeta(E_{s({\cal H})}({\cal H}))>\zeta(E_{s({\cal G})}({\cal G}))$ and $\tilde{j}_*\le j_*$,
{\bf(ii$'$)}
$\zeta(E_{s({\cal H})}({\cal H}))>\zeta(E_{s({\cal G})}({\cal G}))$ and $\tilde{j}_*>j_*$,
and
{\bf(iii$'$)}
$\zeta(E_{s({\cal H})}({\cal H}))=\zeta(E_{s({\cal G})}({\cal G}))$ and $\tilde{j}_*< j_*$,
we can obtain the Newton polygon and edge-polynomials (with opposite signs)
by directly exchanging ${\cal G}$ with ${\cal H}$ in Lemma~\ref{precondiJ2}
because $[{\cal G},{\cal H}]_\theta=-[{\cal H},{\cal G}]_\theta$,
where we exchange $j_*$ with $\tilde{j}_*$,
$i_{s({\cal G})}$ with $\tilde{i}_{s({\cal H})}$,
and $j_{s({\cal G})}$ with $\tilde{j}_{s({\cal H})}$ correspondingly.
}
\label{Rk:NJ2S}
\end{rmk}

{\bf Proof.}
In order to obtain ${\cal N}_f$,
we need to know the relative position of Newton polygons of $\theta{\cal G}'_\theta{\cal H}$ and $-{\cal G}\theta{\cal H}'_\theta$.
We claim that
if $\zeta(E_{s({\cal G})}({\cal G}))\ge\zeta(E_{s({\cal H})}({\cal H}))$, then
\begin{equation}
{\cal N}^\infty_{\theta{\cal G}'_\theta{\cal H}}(u)\left\{
\begin{array}{lllll}
={\cal N}^\infty_{{\cal G}{\cal H}}(u),  &u\in[i_0+\tilde{i}_0,i_*+\tilde{i}_{s({\cal H})}),
\\
=j_*,                                    &u\in[i_*+\tilde{i}_{s({\cal H})},+\infty),
\end{array}
\right.
\label{NydftfNftf}
\end{equation}
\begin{equation}
{\cal N}^\infty_{{\cal G}\theta{\cal H}'_\theta}(u)\left\{
\begin{array}{lllll}
\ge{\cal N}^\infty_{{\cal G}{\cal H}}(u),  &u\in[i_0+\tilde{i}_0,i_{s({\cal G})}+\tilde{i}_*),
\\
=\tilde{j}_*,                              &u\in[i_{s({\cal G})}+\tilde{i}_*,+\infty),
\end{array}
\right.
\label{NfydtfNftf}
\end{equation}
where ${\cal N}^\infty_{\theta{\cal G}'_\theta{\cal H}}(u)$,
${\cal N}^\infty_{{\cal G}\theta{\cal H}'_\theta}(u)$ and
${\cal N}^\infty_{{\cal G}{\cal H}}(u)$ are defined just before \eqref{NfleNtf}.
We first consider (\ref{NydftfNftf}).
By \eqref{unionseq},
we have
$\xi_{s({\cal G},{\cal H})}=\zeta(E_{s({\cal G})}({\cal G}))$
since $\zeta(E_{s({\cal G})}({\cal G}))\ge\zeta(E_{s({\cal H})}({\cal H}))$.
Similar to (\ref{J1Vcom}),
$\Delta^V(\theta{\cal G}'_\theta{\cal H})\cap\Delta^V({\cal GH})
=\{V_0({\cal GH}),...,V_{s({\cal GH})-1}({\cal GH})\}$,
and therefore,
\begin{eqnarray}
\Delta^V(\theta{\cal G}'_\theta{\cal H})
\supset
\{V_0({\cal G}{\cal H}),...,V_{s({\cal G},{\cal H})-1}({\cal G}{\cal H})\}.
\label{PVydftf}
\end{eqnarray}
Noticing that $j_{s({\cal G})}=0$ and condition {\bf(S)} holds,
we see from Remark~\ref{rmk:diff}{\bf(ii)} that
\begin{eqnarray}
{\cal N}_{\theta{\cal G}'_\theta}=\{{\cal N}_{\cal G}
\backslash\underline{V_*({\cal G})V_{s({\cal G})}({\cal G})}\}\cup V_*({\cal G}),
\label{Nydf}
\end{eqnarray}
where $V_*({\cal G}):(i_*,j_*)$ is the second to the last valid index on $\vec{\cal N}_{{\cal G}}$.
By Proposition~\ref{lm-MIL} and Remark~\ref{rmk:mul}{\bf(iii)},
the left end-point and right end-point of ${\cal N}_{\theta{\cal G}'_\theta{\cal H}}$ are
\begin{equation}
\begin{split}
V_0({\cal GH})&=V_0({\cal G})+V_0({\cal H})=(i_0+\tilde{i}_0,j_0+\tilde{j}_0),
\\
Q_*&=V_*({\cal G})+V_{s({\cal H})}({\cal H})=(i_*+\tilde{i}_{s({\cal H})},j_*),
\end{split}
\label{lrV0Q}
\end{equation}
respectively.
In order to prove \eqref{NydftfNftf},
we consider three circumstances:
{\bf(C1)}
$\zeta(E_{s({\cal G})}({\cal G}))>\zeta(E_{s({\cal H})}({\cal H}))$ and
$j_*=j_{s({\cal G})-1}$,
{\bf(C2)}
$\zeta(E_{s({\cal G})}({\cal G}))>\zeta(E_{s({\cal H})}({\cal H}))$ and
$j_*<j_{s({\cal G})-1}$,
and
{\bf(C3)} $\zeta(E_{s({\cal G})}({\cal G}))=\zeta(E_{s({\cal H})}({\cal H}))$.
In {\bf(C1)},
we have
$\mathfrak{S}(\theta{\cal G}'_\theta)
=\mathfrak{S}({\cal G})\backslash\{\zeta(E_{s({\cal G})}({\cal G}))\}$.
Since $\zeta(E_{s({\cal G})}({\cal G}))\notin
\mathfrak{S}({\cal G})\cap\mathfrak{S}({\cal H})$,
$$
\sharp\mathfrak{S}(\theta{\cal G}'_\theta)\cup\mathfrak{S}({\cal H})
=\sharp\mathfrak{S}({\cal G})\cup\mathfrak{S}({\cal H})-1=s({\cal G},{\cal H})-1,
$$
implying by Remark~\ref{rmk:diff}{\bf(i)} that
$\sharp\Delta^V(\theta{\cal G}'_\theta{\cal H})=s({\cal G},{\cal H})$.
We see from \eqref{PVydftf} that
$$
\vec{\Delta}^V(\theta{\cal G}'_\theta{\cal H})
=(V_0({\cal G}{\cal H}),...,V_{s({\cal G},{\cal H})-1}({\cal G}{\cal H})),
$$
and therefore,
$
Q_*=V_{s({\cal G},{\cal H})-1}({\cal G}{\cal H})\in
E_{s({\cal G},{\cal H})}({\cal G}{\cal H})
$.
Then, \eqref{NydftfNftf} holds in {\bf(C1)}.
Either in {\bf(C2)} or {\bf(C3)},
we have
$\mathfrak{S}(\theta{\cal G}'_\theta)\cup\mathfrak{S}({\cal H})
=\mathfrak{S}({\cal G})\cup\mathfrak{S}({\cal H})$.
By Remark~\ref{rmk:diff}{\bf(i)},
$\sharp\Delta^V(\theta{\cal G}'_\theta{\cal H})=s({\cal G},{\cal H})+1$.
Then, we see from \eqref{PVydftf} and \eqref{lrV0Q} that
$$
\vec{\Delta}^V(\theta{\cal G}'_\theta{\cal H})
=(V_0({\cal G}{\cal H}),...,V_{s({\cal G},{\cal H})-1}({\cal G}{\cal H}),Q_*).
$$
It follows that
$
\zeta(\underline{V_{s({\cal G},{\cal H})-1}({\cal G}{\cal H})Q_*})
=\xi_{s({\cal G},{\cal H})}
=\zeta(\underline{V_{s({\cal G},{\cal H})-1}({\cal G}{\cal H})
V_{s({\cal G},{\cal H})}({\cal G}{\cal H})}).
$
Note that $j(Q_*)=j_*>0$ and
by Remark~\ref{rmk:mul}{\bf(iii)},
$j(V_{s({\cal G},{\cal H})}({\cal G}{\cal H}))
=j_{s({\cal G})}+\tilde{j}_{s({\cal H})}=0$.
Then $Q_*\in E_{s({\cal G},{\cal H})}({\cal G}{\cal H})$,
which proves \eqref{NydftfNftf}.

For (\ref{NfydtfNftf}),
since $\tilde{j}_{s({\cal H})}=0$ and condition {\bf(S)} holds,
similar to (\ref{Nydf}),
the right end-point of ${\cal N}_{\theta{\cal H}'_\theta}$ is
$V_*({\cal H}):(\tilde{i}_*,\tilde{j}_*)$.
By Remark~\ref{rmk:mul}{\bf(iii)},
the left end-point and right end-point of
${\cal N}_{{\cal G}\theta{\cal H}'_\theta}$ are
\begin{equation}
\begin{split}
V_0({\cal GH})
&=V_0({\cal G})+V_0(\theta{\cal H}'_\theta)
=(i_0+\tilde{i}_0,j_0+\tilde{j}_0),
\\
\tilde{Q}_*
&=V_{s({\cal G})}({\cal G})+V_*({\cal H})
=(i_{s({\cal G})}+\tilde{i}_*,\tilde{j}_*)
\end{split}
\label{lrV0tQ}
\end{equation}
respectively.
Since
$\Delta(\theta{\cal H}'_\theta)
=\Delta({\cal H}'_\theta)+(0,1)\subset\Delta({\cal H})$,
we have
$\wp(\Delta(\theta{\cal H}'_\theta))\subset\wp(\Delta({\cal H}))$.
By the equality
$\wp(\Delta({\cal G}{\cal H}))=\wp(\Delta({\cal G}))+\wp(\Delta({\cal H}))$,
given by \cite[Proposition~2.2]{Lip88},
we have
$$
\wp(\Delta({\cal G}\theta{\cal H}'_\theta))
=\wp(\Delta({\cal G}))+\wp(\Delta(\theta{\cal H}'_\theta))
\subset\wp(\Delta({\cal G}))+\wp(\Delta({\cal H}))
=\wp(\Delta({\cal G}{\cal H})).
$$
Then ${\cal N}_{{\cal G}{\cal H}}\le{\cal N}_{{\cal G}\theta{\cal H}'_\theta}$,
which implies by \eqref{lrV0tQ} that \eqref{NfydtfNftf} holds.

We see from the proof of \eqref{NydftfNftf} that
the right end-point $Q_*$ of ${\cal N}_{\theta{\cal G}'_\theta{\cal H}}$
lies on the edge $E_{s({\cal G},{\cal H})}({\cal G}{\cal H})$
when $\zeta(E_{s({\cal G})}({\cal G}))\ge\zeta(E_{s({\cal H})}({\cal H}))$.
Then \eqref{NydftfNftf} ensures that
\begin{eqnarray}
{\cal N}_{{\cal G}'_\theta{\cal H}}
=\{{\cal N}_{{\cal G}{\cal H}}\backslash
\underline{Q_*V_{s({\cal G},{\cal H})}({\cal G}{\cal H})}\}\cup Q_*-(0,1).
\label{J1Ndftf}
\end{eqnarray}
On the other hand,
similar to equalities \eqref{J1Vcom} and \eqref{DyfdtfEk} in the proof of Lemma~\ref{J1Ng},
\begin{equation}
\begin{split}
\Delta^V({\cal G}\theta{\cal H}'_\theta)\cap\Delta^V({\cal G}{\cal H})&=\{V_0({\cal G}{\cal H}),...,V_{k_*-1}({\cal G}{\cal H})\},
\\
\Delta({\cal G}\theta{\cal H}'_\theta)\cap E_k({\cal G}{\cal H})&=\emptyset
~~~\mbox{for all}~~~k=k_*+1,...,s({\cal G},{\cal H}),
\end{split}
\label{dtfVcom}
\end{equation}
where $k_*$ is the integer such that $\xi_{k_*}=\zeta(E_{s({\cal H})}({\cal H}))$.

{\bf(i)}
In the case $\zeta(E_{s({\cal G})}({\cal G}))>\zeta(E_{s({\cal H})}({\cal H}))$ and
$j_*\le \tilde{j}_*$,
we see from (\ref{NydftfNftf}) and (\ref{NfydtfNftf}) that
${\cal N}_{{\cal G}'_\theta{\cal H}}\le{\cal N}_{-{\cal G}{\cal H}'_\theta}$.
Noticing that common vertices of
${\cal G}'_\theta{\cal H}$ and $-{\cal G}{\cal H}'_\theta$ satisfy
the non-vanishing condition \eqref{Anon} for addition,
we conclude by Corollary~\ref{cor-ADD} that
${\cal N}_f={\cal N}_{{\cal G}'_\theta{\cal H}}$.
By Remark \ref{rmk:mul}{\bf(i)},
we see from \eqref{J1Ndftf} that
$s(f)=s({\cal G},{\cal H})-1+{\rm sgn}(j_{s({\cal G})-1}-j_*)$.
Similarly to the analysis of edge-polynomials
in the case $\tilde{j}_{s({\cal H})}=0$
in the proof of Lemma~\ref{J1Ng},
we obtain equality \eqref{J2igEkg} by \eqref{dtfVcom}.

{\bf(ii)} In the case $\zeta(E_{s({\cal G})}({\cal G}))>\zeta(E_{s({\cal H})}({\cal H}))$ and $j_*>\tilde{j}_*$,
we claim that
\begin{eqnarray}\label{wpyg}
\wp(\Delta(\theta f))=\wp({\cal V}_*\cup Q_*\cup \tilde{Q}_*),
\end{eqnarray}
where ${\cal V}_*:=\{V_0({\cal G}{\cal H}),...,V_{s({\cal G},{\cal H})-1}({\cal G}{\cal H})\}$.
In fact,
it follows from (\ref{J1Ndftf}) that
\begin{eqnarray}\label{Vydftf}
\Delta^V(\theta{\cal G}'_\theta{\cal H})={\cal V}_*\cup Q_*.
\end{eqnarray}
Since common vertices of ${\cal G}'_\theta{\cal H}$ and $-{\cal G}{\cal H}'_\theta$ satisfy
the non-vanishing condition \eqref{Anon} for addition,
those of $\theta{\cal G}'_\theta{\cal H}$ and $-{\cal G}\theta{\cal H}'_\theta$ have the same property.
Then, using \eqref{Vydftf}, Proposition~\ref{ADD} and the fact that
$\tilde{Q}_*$ is the right end-point of ${\cal N}_{-{\cal G}\theta{\cal H}'_\theta}$,
we have
$$
\wp({\cal V}_*\cup Q_*\cup \tilde{Q}_*)
\subset\wp(\Delta^V(\theta{\cal G}'_\theta{\cal H})\cup\Delta^V(-{\cal G}\theta{\cal H}'_\theta))
=\wp(\Delta(\theta f)),
$$
implying the ``$\supset$''-part of \eqref{wpyg}.

Next,
we prove the ``$\subset$''-part of \eqref{wpyg}.
Condition {\bf(S)} and Remark~\ref{rmk:diff}{\bf(ii)} ensure that
$\mathfrak{S}(\theta{\cal H}'_\theta)\subset\mathfrak{S}({\cal H})$.
Note that
$\zeta(E_{s({\cal G})}({\cal G}))$ is the last element in
$\vec{\mathfrak{S}}({\cal G})\cup\vec{\mathfrak{S}}(\theta{\cal H}'_\theta)$ and
$\zeta(E_{s({\cal G})}({\cal G}))\notin
\mathfrak{S}({\cal G})\cap\mathfrak{S}(\theta{\cal H}'_\theta)$
since $\zeta(E_{s({\cal G})}({\cal G}))>\zeta(E_{s({\cal H})}({\cal H}))$.
By Remark~\ref{rmk:mul},
the last edge $E$ of ${\cal N}_{{\cal G}\theta{\cal H}'_\theta}$ satisfies that
$\zeta(E)=\zeta(E_{s({\cal G})}({\cal G}))$ and $h(E)=h(E_{s({\cal G})}({\cal G}))$,
and the right end-point of ${\cal N}_{{\cal G}\theta{\cal H}'_\theta}$ is
$\tilde{Q}_*:(i_{s({\cal G})}+\tilde{i}_*,\tilde{j}_*)$.
Then, the second to the last vertex on $\vec{\cal N}_{{\cal G}\theta{\cal H}'_\theta}$ is $(i_{s({\cal G})-1}+\tilde{i}_*,j_{s({\cal G})-1}+\tilde{j}_*)$.
Since $j_{s({\cal G})-1}+\tilde{j}_*>j_*$,
we see from (\ref{NydftfNftf}), (\ref{NfydtfNftf})
and (\ref{Vydftf}) that
\begin{eqnarray}
\Delta^V({\cal G}\theta{\cal H}'_\theta)\backslash \tilde{Q}_*
\subset\wp(\Delta^V(\theta{\cal G}'_\theta{\cal H}))
=\wp({\cal V}_*\cup Q_*).
\end{eqnarray}
By the property of the union of convex hulls (\cite[Theorem 3.2]{Soltan}),
i.e., ${\rm conv}\{{\cal A}\cup{\cal B}\}={\rm conv}\{{\rm conv}\{{\cal A}\}\cup{\rm conv}\{{\cal B}\}\}$
for any two sets ${\cal A},{\cal B}\subset\mathbb{R}^2$,
we obtain that
\begin{eqnarray*}
\wp(\Delta^V({\cal G}\theta{\cal H}'_\theta))
&=&\wp(\{\Delta^V({\cal G}\theta{\cal H}'_\theta)\backslash \tilde{Q}_*\}\cup \tilde{Q}_*)
\\
&=&{\rm conv}\{\wp(\Delta^V({\cal G}\theta{\cal H}'_\theta)\backslash \tilde{Q}_*)\cup\wp(\tilde{Q}_*)\}
\\
&\subset&{\rm conv}\{\wp({\cal V}_*\cup Q_*)\cup\wp(\tilde{Q}_*)\}
\\
&=&\wp({\cal V}_*\cup Q_*\cup \tilde{Q}_*).
\end{eqnarray*}
On the other hand,
(\ref{Vydftf}) shows that
$\wp(\Delta^V(\theta{\cal G}'_\theta{\cal H}))=\wp({\cal V}_*\cup Q_*)
\subset\wp({\cal V}_*\cup Q_*\cup \tilde{Q}_*)$.
Further, using Proposition~\ref{ADD} and
the property of the union of convex hulls (\cite[Theorem 3.2]{Soltan}),
we obtain that
\begin{eqnarray*}
\wp(\Delta(\theta f))
&=&\wp(\Delta^V(\theta{\cal G}'_\theta{\cal H})\cup\Delta^V(-{\cal G}\theta{\cal H}'_\theta))
\\
&=&{\rm conv}\{\wp(\Delta^V(\theta{\cal G}'_\theta{\cal H}))\cup\wp(\Delta^V({\cal G}\theta{\cal H}'_\theta))\}
\\
&\subset&{\rm conv}\{\wp({\cal V}_*\cup Q_*\cup \tilde{Q}_*)\}
\\
&=&\wp({\cal V}_*\cup Q_*\cup \tilde{Q}_*).
\end{eqnarray*}
This proves the ``$\subset$''-part of (\ref{wpyg}) and completes the proof of \eqref{wpyg}.

Having \eqref{wpyg},
we further show that
\begin{eqnarray}
\vec\Delta^V(\theta f)=(V_0({\cal G}{\cal H}),...,
V_{s({\cal G},{\cal H})-1}({\cal G}{\cal H}),Q_*,\tilde{Q}_*).
\label{Vtg}
\end{eqnarray}
Since $\zeta(E_{s({\cal G})}({\cal G}))>\zeta(E_{s({\cal H})}({\cal H}))$,
we have $k_*<s({\cal G},{\cal H})$.
Then, equalities in \eqref{dtfVcom} yield that
$\Delta({\cal G}\theta{\cal H}'_\theta)\cap
E_{s({\cal G},{\cal H})}({\cal G}{\cal H})=\emptyset$.
It follows that
$\tilde{Q}_*$ lies above the line containing the edge
$E_{s({\cal G},{\cal H})}({\cal G}{\cal H})$
since $\tilde{Q}_*\in\wp(\Delta({\cal G}\theta{\cal H}'_\theta))
\subset\wp(\Delta({\cal G}{\cal H}))$.
Moreover,
$\xi_{s({\cal G},{\cal H})}
=\zeta(E_{s({\cal G},{\cal H})}({\cal G}{\cal H}))
<\zeta(\underline{Q_*\tilde{Q}_*})<0$ because
$Q_*\in E_{s({\cal G},{\cal H})}({\cal G}{\cal H})$ and $j(Q_*)=j_*>\tilde{j}_*=j(\tilde{Q}_*)$.
When $j_*=j_{s({\cal G})-1}$,
we have $Q_*=V_{s({\cal G},{\cal H})-1}({\cal G}{\cal H})$,
and
$$
\zeta(\underline{V_{s({\cal G},{\cal H})-2}({\cal G}{\cal H})
V_{s({\cal G},{\cal H})-1}({\cal G}{\cal H})})
<\xi_{s({\cal G},{\cal H})}
<\zeta(\underline{V_{s({\cal G},{\cal H})-1}({\cal G}{\cal H})\tilde{Q}_*}).
$$
Then, \eqref{wpyg} yields \eqref{Vtg}.
When $j_*<j_{s({\cal G})-1}$,
we have $Q_*\ne V_{s({\cal G},{\cal H})-1}({\cal G}{\cal H})$,
and $\zeta(\underline{V_{s({\cal G},{\cal H})-1}({\cal G}{\cal H})Q_*})
=\xi_{s({\cal G},{\cal H})}<\zeta(\underline{Q_*\tilde{Q}_*})$.
Then, \eqref{wpyg} also yields \eqref{Vtg}.
This proves \eqref{Vtg},
implying that the Newton polygon of $f$ is given by \eqref{J2iiNg}.

Newton polygon ${\cal N}_f$ in case {\bf(ii)} has one more edge
than that in case {\bf(i)}.
It is similar to show that
the first $s(f)-1$ edge-polynomials of $f$ are given by \eqref{J2igEkg}.
Note that $\underline{Q_*\tilde{Q}_*}$ is the last edge of ${\cal N}_{\theta f}$.
We see from \eqref{Nydf} that
$\mathfrak{S}(\theta{\cal G}'_\theta)\subset\mathfrak{S}({\cal G})$
and therefore,
$\mathfrak{S}(\theta{\cal G}'_\theta{\cal H})\subset
\mathfrak{S}({\cal G}{\cal H})$ because of Remark~\ref{rmk:mul}{\bf(i)}.
Then,
as the right end-point of ${\cal N}_{\theta{\cal G}'_\theta{\cal H}}$,
$Q_*$ satisfies that $\zeta^+(Q_*)=0$ and
$\zeta^-(Q_*)\le \xi_{s({\cal G},{\cal H})}$
since $\xi_{s({\cal G},{\cal H})}$ is the smallest element in
$\mathfrak{S}({\cal G}{\cal H})$.
Similarly,
being the right end-point of ${\cal N}_{{\cal G}\theta{\cal H}'_\theta}$,
$\tilde{Q}_*$ satisfies that $\zeta^+(\tilde{Q}_*)=0$
and $\zeta^-(\tilde{Q}_*)\le \xi_{s({\cal G},{\cal H})}$.
Since it is proved that
$\xi_{s({\cal G},{\cal H})}<\zeta(\underline{Q_*\tilde{Q}_*})<0$ in the above paragraph,
$\zeta(\underline{Q_*\tilde{Q}_*})\in(\zeta^-(Q_*),\zeta^+(Q_*))$ and
$\zeta(\underline{Q_*\tilde{Q}_*})\in(\zeta^-(\tilde{Q}_*),\zeta^+(\tilde{Q}_*))$.
Then,
\begin{eqnarray}
\underline{Q_*\tilde{Q}_*}\cap\Delta(\theta{\cal G}'_\theta{\cal H})=Q_*~~~\mbox{and}~~~
\underline{Q_*\tilde{Q}_*}\cap\Delta({\cal G}\theta{\cal H}'_\theta)=\tilde{Q}_*.
\label{QQQQ}
\end{eqnarray}
Note that
$V_*({\cal G})$, $V_*({\cal H})$, $V_{s({\cal G})}({\cal G})$,
$V_{s({\cal H})}({\cal H})$, $Q_*$ and $\tilde{Q}_*$
are the $\xi$-components of ${\cal N}_{\theta{\cal G}'_\theta}$,
${\cal N}_{\theta{\cal H}'_\theta}$, ${\cal N}_{\cal G}$,
${\cal N}_{{\cal H}}$, ${\cal N}_{\theta{\cal G}'_\theta{\cal H}}$ and
${\cal N}_{{\cal G}\theta{\cal H}'_\theta}$ respectively
for small $\xi<0$.
Moreover, $Q_*=V_*({\cal G})+V_{s({\cal H})}({\cal H})$ and
$\tilde{Q}_*=V_{s({\cal G})}({\cal G})+V_*({\cal H})$.
By Proposition~\ref{lm-MIL},
\begin{eqnarray}
a_{Q_*}(\theta{\cal G}'_\theta{\cal H})
=j_*a_{i_*,j_*}({\cal G})a_{\tilde{i}_{s({\cal H})},\tilde{j}_{s({\cal H})}}({\cal H}),~
a_{\tilde{Q}_*}({\cal G}\theta{\cal H}'_\theta)
=\tilde{j}_*a_{i_{s({\cal G})},j_{s({\cal G})}}({\cal G})
a_{\tilde{i}_*,\tilde{j}_*}({\cal H}).
\label{aQfaQtf}
\end{eqnarray}
Then, we see from \eqref{fonE}, \eqref{QQQQ} and \eqref{aQfaQtf} that
the edge-polynomial corresponding to the last edge
$E_{s(f)}(f)=\underline{Q_*\tilde{Q}_*}-(0,1)$ is given by \eqref{gElast}.

{\bf(iii)}
In the case $\zeta(E_{s({\cal G})}({\cal G}))=\zeta(E_{s({\cal H})}({\cal H}))$ and
$j_*\le \tilde{j}_*$,
similar to \eqref{NydftfNftf},
\eqref{NfydtfNftf} can be rewritten as
\begin{equation}\label{NfydtfNftf2}
{\cal N}^\infty_{{\cal G}\theta{\cal H}'_\theta}(u)
=\left\{
\begin{array}{lllll}
{\cal N}^\infty_{{\cal G}{\cal H}}(u),  &u\in[i_0+\tilde{i}_0,i_{s({\cal G})}+\tilde{i}_*),
\\
\tilde{j}_*,                            &u\in[i_{s({\cal G})}+\tilde{i}_*,+\infty).
\end{array}
\right.
\end{equation}
Then we see from
\eqref{NydftfNftf} and \eqref{NfydtfNftf2} that
${\cal N}_{\theta{\cal G}'_\theta{\cal H}}\le{\cal N}_{{\cal G}\theta{\cal H}'_\theta}$.
Since common vertices of ${\cal G}'_\theta{\cal H}$  and $-{\cal G}{\cal H}'_\theta$ satisfy the non-vanishing condition \eqref{Anon} for addition,
common vertices of $\theta{\cal G}'_\theta{\cal H}$ and $-{\cal G}\theta{\cal H}'_\theta$ have the same property.
By Corollary \ref{cor-ADD},
we have ${\cal N}_{\theta f}={\cal N}_{\theta{\cal G}'_\theta{\cal H}}$ and therefore,
${\cal N}_f={\cal N}_{{\cal G}'_\theta{\cal H}}$.
Further, the Newton polygon of $f$ follows from (\ref{J1Ndftf}).
By \eqref{Nydf}, we have either
$\mathfrak{S}({\cal G}'_\theta)=\mathfrak{S}({\cal G})$ or
$\mathfrak{S}({\cal G}'_\theta)
=\mathfrak{S}({\cal G})\backslash\{\zeta(E_{s({\cal G})}({\cal G}))\}$.
Then, the equality $\zeta(E_{s({\cal G})}({\cal G}))=\zeta(E_{s({\cal H})}({\cal H}))$ implies that $\sharp\mathfrak{S}({\cal G}'_\theta)\cup\mathfrak{S}({\cal H})=\sharp\mathfrak{S}({\cal G})\cup\mathfrak{S}({\cal H})=s({\cal G},{\cal H})$.
Hence, $s(f)=s({\cal G},{\cal H})$ by Remark~\ref{rmk:mul}{\bf(i)}.
Similarly to the analysis of edge-polynomials of $f$ in the case $\tilde{j}_{s({\cal H})}=0$
in the proof of Lemma~\ref{J1Ng},
we have
$f_{E_k(f)}=[{\cal K}_{\xi_k}({\cal G}),{\cal K}_{\xi_k}({\cal H})]$
for all $k=1,...,s({\cal G},{\cal H})$.
Thus, the proof of this lemma is completed.
\qquad$\Box$

\subsection{Semi-definiteness of the Lie-bracket $[{\cal G},{\cal H}]_\theta$}

By those results on Newton polygon of the function $[{\cal G},{\cal H}]_\theta$
given in the above subsection,
the following lemma further shows its semi-definiteness.

\begin{lm}
Function $[{\cal G},{\cal H}]_\theta$,
the Lie-bracket of ${\cal G}$ and ${\cal H}$ in the variable $\theta$,
is semi-definite on $\Omega_\epsilon:=\{(\rho,\theta)\in\mathbb{R}^2:0<\rho<\epsilon,|\theta|<\epsilon\}$
for an $\epsilon>0$ if either
{\bf(L1):} conditions {\bf(P1)} and {\bf (Q)} hold in the one-above case {\bf (J1)},
or {\bf(L2):} conditions {\bf(P2)}, {\bf(Q)} and {\bf(S)},
and one of conditions {\bf(H1)}, {\bf(H2)}, {\bf(H3)} and {\bf(H4)}
of Theorem~\ref{th:finiteJ2} hold in the two-below case {\bf (J2)}.
Moreover, it is semi-positive {\rm(}or semi-negative{\rm)}
if $C_0>0$ {\rm(}or $<0${\rm)}, where $C_0$ is defined just before Theorem~\ref{th:finite}.
\label{defsignth1}
\end{lm}

\begin{rmk}
{\rm
Corresponding to Remark~\ref{Rk:NJ2S},
if we exchange ${\cal G}$ with ${\cal H}$
in conditions {\bf(P2)} and
{\bf(H1)}, {\bf(H2)}, {\bf(H3)} and {\bf(H4)} of Theorem~\ref{th:finiteJ2},
then Lemma~\ref{defsignth1} still holds,
where we exchange $j_*$ with $\tilde{j}_*$,
$i_{s({\cal G})}$ with $\tilde{i}_{s({\cal H})}$,
and $j_{s({\cal G})}$ with $\tilde{j}_{s({\cal H})}$ correspondingly.
\label{Rk:NJ2S5}
}
\end{rmk}

\noindent{\bf Proof.}
We first consider this lemma when conditions in {\bf(L1)} are satisfied.
Condition {\bf(J1)} implies that $j_{s({\cal G})}+\tilde{j}_{s({\cal H})}>0$,
and we assume without loss of generality that $j_{s({\cal G})}>0$.
In order to use Lemma~\ref{J1Ng},
we show that
common vertices of ${\cal G}'_\theta{\cal H}$ and $-{\cal G}{\cal H}'_\theta$ satisfy
the non-vanishing condition \eqref{Anon} for addition if condition {\bf(P1)} holds.
The proof is divided into the case $\tilde{j}_{s({\cal H})}>0$
and the case $\tilde{j}_{s({\cal H})}=0$.
In the first case,
Remark~\ref{rmk:diff}{\bf(i)} shows that
${\cal N}_{\theta{\cal G}'_\theta}={\cal N}_{\cal G}$
and
${\cal N}_{-\theta{\cal H}'_\theta}={\cal N}_{\cal H}$.
By Proposition~\ref{lm-MIL},
${\cal N}_{\theta{\cal G}'_\theta{\cal H}}
={\cal N}_{-{\cal G}\theta{\cal H}'_\theta}
={\cal N}_{{\cal GH}}$.
Moreover, Remark~\ref{rmk:mul}{\bf(i)} shows that
$$
\Delta^V(\theta{\cal G}'_\theta{\cal H})\cap\Delta^V(-{\cal G}\theta{\cal H}'_\theta)
=\Delta^V({\cal GH})
=\{V_0({\cal G}{\cal H}),....,V_{s({\cal G},{\cal H})}({\cal G}{\cal H})\}.
$$
For each $k=0,...,s({\cal G},{\cal H})$,
choosing $\nu_k\in(\zeta^-(V_k({\cal G}{\cal H})),\zeta^+(V_k({\cal G}{\cal H})))$,
we have $\nu_k\notin\mathfrak{S}({\cal G})\cup\mathfrak{S}({\cal H})$
by Remark~\ref{rmk:mul}{\bf(i)}.
Then,
there are
$n(k)\in\{0,...,s({\cal G})\}$ and $\tilde{n}(k)\in\{0,...,s({\cal H})\}$ such that
$V_{n(k)}({\cal G})$ and $V_{\tilde{n}(k)}({\cal H})$ are the $\nu_k$-components of
${\cal N}_{\cal G}$ and ${\cal N}_{{\cal H}}$ respectively.
By definition~\eqref{defkxip}, Proposition \ref{lm-MIL} and Remark~\ref{rmk:diff}{\bf(i)},
\begin{eqnarray*}
a_{V_k({\cal G}{\cal H})}(\theta{\cal G}'_\theta{\cal H})\theta^{j(V_k({\cal G}{\cal H}))}
&=&{\cal K}_{\nu_k}(\theta{\cal G}'_\theta{\cal H})
\\
&=&{\cal K}_{\nu_k}(\theta)({\cal K}_{\nu_k}({\cal G}))'_\theta{\cal K}_{\nu_k}({\cal H})
\\
&=&\theta j(V_{n(k)}({\cal G}))a_{V_{n(k)}({\cal G})}({\cal G})\theta^{j(V_{n(k)}({\cal G}))-1}{\cal K}_{\nu_k}({\cal H})
\\
&=&j(V_{n(k)}({\cal G})){\cal K}_{\nu_k}({\cal G}{\cal H}),
\end{eqnarray*}
and similarly,
$a_{V_k({\cal G}{\cal H})}(-{\cal G}\theta{\cal H}'_\theta)
\theta^{j(V_k({\cal G}{\cal H}))}
=-j(V_{\tilde{n}(k)}({\cal H})){\cal K}_{\nu_k}({\cal G}{\cal H})$,
where $a_Q({\cal G})$ represents the coefficient of function ${\cal G}$
corresponding to the lattice point $Q$ and $j(Q)$ represents the ordinate of $Q$.
We claim that
\begin{eqnarray}
j(V_{n(k)}({\cal G}))\ne j(V_{\tilde{n}(k)}({\cal H}))
~~~\mbox{for all}~~~k=0,...,s({\cal G},{\cal H}).
\label{jDkfnejDktf}
\end{eqnarray}
In fact,
since $\nu_k\in(\zeta^-(V_k({\cal G}{\cal H})),\zeta^+(V_k({\cal G}{\cal H})))$,
there exists $\varrho(k)\in\{0,...,r({\cal G},{\cal H})\}$ such that $\nu_k\in(\xi_{s_{\varrho(k)}},\xi_{s_{\varrho(k)+1}})$,
where $\xi_{s_{\varrho(k)}}$ and $\xi_{s_{\varrho(k)+1}}$ are elements
in the intersection
$\mathfrak{S}({\cal G})\cap\mathfrak{S}({\cal H})$
defined in \eqref{interseq},
and we complementarily define
$\xi_{s_0}:=-\infty$ and $\xi_{s_{r({\cal G},{\cal H})+1}}:=0$.
Then, we obtain that
\begin{eqnarray*}
&&\xi_{s_{\varrho(k)}}\le\zeta^-(V_{n(k)}({\cal G}))<
\zeta^+(V_{n(k)}({\cal G}))\le\xi_{s_{\varrho(k)+1}},
\\
&&\xi_{s_{\varrho(k)}}\le\zeta^-(V_{\tilde{n}(k)}({\cal H}))<
\zeta^+(V_{\tilde{n}(k)}({\cal H}))\le\xi_{s_{\varrho(k)+1}}.
\end{eqnarray*}
By definition (\ref{def-Lambda}),
$V_{n(k)}({\cal G})\in\Lambda_{\varrho(k)}({\cal G})$ and
$V_{\tilde{n}(k)}({\cal H})\in\Lambda_{\varrho(k)}({\cal H})$.
Condition {\bf(P1)} further implies that
$j(V_{n(k)}({\cal G}))$ and $j(V_{\tilde{n}(k)}({\cal H}))$
have different parities, and therefore, the claimed \eqref{jDkfnejDktf} holds.
It follows that
{\small
$$
\big(a_{V_k({\cal G}{\cal H})}(\theta{\cal G}'_\theta{\cal H})
\!+\!a_{V_k({\cal G}{\cal H})}(-{\cal G}\theta{\cal H}'_\theta)\big)\theta^{j(V_k({\cal GH}))}
=(j(V_{n(k)}({\cal G}))\!-\!j(V_{\tilde{n}(k)}({\cal H}))){\cal K}_{\nu_k}({\cal G}{\cal H})
\ne0.
$$
}Then,
common vertices of $\theta{\cal G}'_\theta{\cal H}$ and $-{\cal G}\theta{\cal H}'_\theta$ satisfy the non-vanishing condition \eqref{Anon} for addition
in the case $\tilde{j}_{s({\cal H})}>0$.

In the case $\tilde{j}_{s({\cal H})}=0$,
equality \eqref{J1Vcom} holds,
that is,
$\Delta^V({\cal G}\theta{\cal H}'_\theta)\cap\Delta^V({\cal G}{\cal H})
=\{V_0({\cal G}{\cal H}),...,V_{k_*-1}({\cal G}{\cal H})\}$,
where $k_*$ is either zero as $s({\cal H})=0$,
or the integer such that $\xi_{k_*}=\zeta(E_{s({\cal H})}({\cal H}))$ as $s({\cal H})>0$.
Since ${\cal N}_{\theta{\cal G}'_\theta{\cal H}}={\cal N}_{{\cal G}{\cal H}}$,
we obtain that
$$
\Delta^V(\theta{\cal G}'_\theta{\cal H})\cap\Delta^V(-{\cal G}\theta{\cal H}'_\theta)=\{V_0({\cal G}{\cal H}),...,V_{k_*-1}({\cal G}{\cal H})\}.
$$
If $s({\cal H})=0$ then $k_*=0$ and therefore,
$\Delta^V(\theta{\cal G}'_\theta{\cal H})\cap\Delta^V(-{\cal G}\theta{\cal H}'_\theta)=\emptyset$.
If $s({\cal H})>0$,
it is similar to the case $\tilde{j}_{s({\cal H})}>0$ that
common vertices of ${\cal N}_{\theta{\cal G}'_\theta{\cal H}}$ and
${\cal N}_{-{\cal G}\theta{\cal H}'_\theta}$
satisfy the non-vanishing condition \eqref{Anon} for addition because of {\bf(P1)}.
Consequently,
common vertices of ${\cal G}'_\theta{\cal H}$ and $-{\cal G}{\cal H}'_\theta$
satisfy the non-vanishing condition \eqref{Anon} for addition.

Let $f$ be the Lie-bracket $[{\cal G},{\cal H}]_\theta$.
Lemma~\ref{J1Ng} further shows that ${\cal N}_f={\cal N}_{{\cal G}{\cal H}}-(0,1)$,
having $s({\cal G},{\cal H})$ edges,
and $f_{E_k(f)}=[{\cal K}_{\xi_k}({\cal G}),{\cal K}_{\xi_k}({\cal H})]$
for all $k=1,...,s({\cal G},{\cal H})$.
By Remark~\ref{rmk:mul}{\bf(iii)},
the right end-point $V_{s({\cal G},{\cal H})}(f)$ is equal to
$V_{s({\cal G})}({\cal G})+V_{s({\cal H})}({\cal H})-(0,1)$.
We see from definition (\ref{def-Lambda}) that
$V_{s({\cal G})}({\cal G})\in\Lambda_{r({\cal G},{\cal H})}({\cal G})$ and
$V_{s({\cal H})}({\cal H})\in\Lambda_{r({\cal G},{\cal H})}({\cal H})$.
Then,
condition {\bf(P1)} implies that the ordinate of $V_{s({\cal G},{\cal H})}(f)$ is even,
and therefore, $f$ does not have a factor $\theta$ of odd multiplicity.
If $s({\cal G},{\cal H})=0$,
then $f$ is clearly semi-definite on $\Omega_\epsilon$ for an $\epsilon>0$.
If $s({\cal G},{\cal H})>0$,
then condition {\bf(Q)} ensures that
the edge-polynomial $f_{E_k(f)}$
has no nonzero real roots for all $k=1,...,s({\cal G},{\cal H})$.
Hence, $f$ is semi-definite on $\Omega_\epsilon$ for an $\epsilon>0$
by Proposition~\ref{lm-region}.
This proves the semi-definiteness of the Lie-bracket $[{\cal G},{\cal H}]_\theta$
when conditions in {\bf(L1)} are satisfied.

Next, we prove the semi-definiteness when conditions in {\bf(L2)} are satisfied.
Under conditions {\bf(S)} and {\bf(H1)} of Theorem~\ref{th:finiteJ2},
we have $k_*<s({\cal G},{\cal H})$,
where $k_*$ is the integer such that $\xi_{k_*}=\zeta(E_{s({\cal H})}({\cal H}))$.
Then \eqref{J1Ndftf} and \eqref{dtfVcom} in the proof of Lemma~\ref{precondiJ2} show that
\begin{eqnarray}
\Delta^V(\theta{\cal G}'_\theta{\cal H})\cap\Delta^V(-{\cal G}\theta{\cal H}'_\theta)
=\{V_0({\cal G}{\cal H}),...,V_{k_*-1}({\cal G}{\cal H})\}.
\label{L2GHV}
\end{eqnarray}
Similar to the proof in {\bf(L1)} of this lemma,
condition {\bf(P2)} indicates that
common vertices of ${\cal G}'_\theta{\cal H}$ and $-{\cal G}{\cal H}'_\theta$
satisfy the non-vanishing condition \eqref{Anon} for addition.
Then Lemma~\ref{precondiJ2}{\bf(i)} ensures that
the ordinate of the right end-point of ${\cal N}_f$ is $j_*-1$.

We claim that $j_*$ is odd.
Actually,
we have $\xi_k=\zeta(E_{s({\cal G})}({\cal G}))$
and $\zeta^-(V_{s({\cal H})}({\cal H}))
=\zeta(E_{s({\cal H})}({\cal H}))
<\xi_k<0=\zeta^+(V_{s({\cal H})}({\cal H}))$
for $k=s({\cal G},{\cal H})$.
By definition \eqref{defkxip},
$$
{\cal K}_{\xi_k}({\cal G})(\theta)
=a_{V_{s({\cal G})-1}({\cal G})}\theta^{j_{s({\cal G})-1}}
+\cdots
+a_{i_*,j_*}({\cal G})\theta^{j_*}+a_{V_{s({\cal G})}({\cal G})}
$$
and $K_{\xi_k}({\cal H})(\theta)=a_{V_{s({\cal H})}({\cal H})}$.
Then,
{\small
$$
[{\cal K}_{\xi_k}({\cal G}),{\cal K}_{\xi_k}({\cal H})]
=j_{s({\cal G})-1}a_{V_{s({\cal G})-1}({\cal G})}a_{V_{s({\cal H})}({\cal H})}
\theta^{j_{s({\cal G})-1}-1}
+\cdots
+j_*a_{i_*,j_*}({\cal G})a_{V_{s({\cal H})}({\cal H})}\theta^{j_*-1}.
$$
}Condition {\bf(P2)} indicates that $j_{s({\cal G})-1}$ is odd.
Condition {\bf(Q)} further implies that $j_*$ is odd.
Then, $f$ does not have a factor $\theta$ of odd multiplicity.
Moreover, by Lemma~\ref{precondiJ2},
all edge-polynomials are given by
$f_{E_k(f)}=[{\cal K}_{\xi_k}({\cal G}),{\cal K}_{\xi_k}({\cal H})]$,
$k=1,...,s({\cal G},{\cal H})-1+{\rm sgn}(j_{s({\cal G})-1}-j_*)$.
It follows from condition {\bf(Q)} and Proposition~\ref{lm-region} that
$f$ is semi-definite on $\Omega_\epsilon$ for an $\epsilon>0$.

Under conditions {\bf(S)} and {\bf(H2)} of Theorem~\ref{th:finiteJ2},
equality \eqref{L2GHV} still hods,
and condition {\bf(P2)} similarly implies that
common vertices of ${\cal G}'_\theta{\cal H}$ and $-{\cal G}{\cal H}'_\theta$ satisfy
the non-vanishing condition \eqref{Anon} for addition.
By Lemma~\ref{precondiJ2}{\bf(ii)},
all edge-polynomials of $f$ are given by
$[{\cal K}_{\xi_k}({\cal G}),{\cal K}_{\xi_k}({\cal H})]$, $k=1,...,s(f)-1$,
and $f_{E_{s(f)}(f)}$ given in \eqref{gElast}.
Condition {\bf(Q)} ensures that
$[{\cal K}_{\xi_k}({\cal G}),{\cal K}_{\xi_k}({\cal H})]$
has no nonzero real roots for all $k=1,...,s(f)-1$.
Similar to the above paragraph,
conditions {\bf(P2)} and {\bf(Q)} imply that $j_*$ is odd.
Since $\tilde{j}_*$ is also odd,
those two monomials of $f_{E_{s(f)}(f)}$ are both of even degrees,
and therefore, $f_{E_{s(f)}(f)}$ has no nonzero real roots because
$a_{i_*,j_*}({\cal G})a_{i_{s({\cal G})},j_{s({\cal G})}}({\cal G})
a_{\tilde{i}_*,\tilde{j}_*}({\cal H})
a_{\tilde{i}_{s({\cal H})},\tilde{j}_{s({\cal H})}}({\cal H})<0$.
Then,
each edge-polynomial of $f$ has no nonzero real roots.
Moreover,
by Lemma~\ref{precondiJ2}{\bf(ii)},
the ordinate of the right end-point of ${\cal N}_f$ is $\tilde{j}_*-1$,
implying that $f$ does not have a factor $\theta$ of odd multiplicity.
Consequently,
$f$ is semi-definite on $\Omega_\epsilon$ for an $\epsilon>0$
by Proposition~\ref{lm-region}.

Under conditions {\bf(S)} and {\bf(H3)} of Theorem~\ref{th:finiteJ2},
since $\zeta(E_{s({\cal G})}({\cal G}))=\zeta(E_{s({\cal H})}({\cal H}))$,
we see from \eqref{J1Ndftf} that
$\Delta^V(\theta{\cal G}'_\theta{\cal H})=
\{V_0({\cal G}{\cal H}),...,V_{s({\cal G},{\cal H})-1}({\cal G}{\cal H})\}\cup Q_*$.
Similarly, we have
$\Delta^V(\theta{\cal G}{\cal H}'_\theta)=
\{V_0({\cal G}{\cal H}),...,V_{s({\cal G},{\cal H})-1}({\cal G}{\cal H})\}\cup\tilde{Q}_*$.
Then,
$$
\Delta^V(\theta{\cal G}'_\theta{\cal H})\cap\Delta^V(\theta{\cal G}{\cal H}'_\theta)
=\{V_0({\cal G}{\cal H}),...,V_{s({\cal G},{\cal H})-1}({\cal G}{\cal H})\}
$$
because $j(Q_*)=j_*<\tilde{j}_*=j(\tilde{Q}_*)$.
Similar to the proof in {\bf(L1)} of this lemma,
condition {\bf(P2)} ensures that
common vertices of ${\cal G}'_\theta{\cal H}$ and $-{\cal G}{\cal H}'_\theta$ satisfy
the non-vanishing condition \eqref{Anon} for addition.
By Lemma~\ref{precondiJ2}{\bf(iii)},
all edge-polynomials are given by
$f_{E_k(f)}=[{\cal K}_{\xi_k}({\cal G}),{\cal K}_{\xi_k}({\cal H})]$,
$k=1,...,s({\cal G},{\cal H})$,
and the ordinate of the right end-point of ${\cal N}_f$ is $j_*-1$.
Similar to the above analysis, $j_*$ is odd,
and therefore, $f$ does not have a factor $\theta$ of odd multiplicity.
Then, condition {\bf(Q)} and Proposition~\ref{lm-region} similarly imply that
$f$ is semi-definite on $\Omega_\epsilon$ for some $\epsilon>0$.

Under conditions {\bf(S)} and {\bf(H4)} of Theorem~\ref{th:finiteJ2},
we have $Q_*=\tilde{Q}_*$ because $j_*=\tilde{j}_*$.
Thus,
$$
\Delta^V(\theta{\cal G}'_\theta{\cal H})\cap\Delta^V(\theta{\cal G}{\cal H}'_\theta)
=\{V_0({\cal G}{\cal H}),...,V_{s({\cal G},{\cal H})-1}({\cal G}{\cal H}),Q_*\}.
$$
Similarly,
common vertices $V_0({\cal G}{\cal H}),...,V_{s({\cal G},{\cal H})-1}({\cal G}{\cal H})$
satisfy the non-vanishing condition \eqref{Anon} for addition because of {\bf(P2)}.
For the common vertex $Q_*$,
noticing that \eqref{aQfaQtf} still holds, we have
$$
a_{Q_*}(\theta{\cal G}'_\theta{\cal H})+a_{Q_*}(-{\cal G}\theta{\cal H}'_\theta)
=j_*\big\{a_{i_*,j_*}({\cal G})a_{\tilde{i}_{s({\cal H})},\tilde{j}_{s({\cal H})}}({\cal H})
-a_{i_{s({\cal G})},j_{s({\cal G})}}({\cal G})a_{\tilde{i}_*,\tilde{j}_*}({\cal H})\big\}
\ne0,
$$
implying that $Q_*$ satisfies the non-vanishing condition \eqref{Anon} for addition.
Consequently, common vertices of
$\theta{\cal G}'_\theta{\cal H}$ and $-{\cal G}\theta{\cal H}'_\theta$
satisfy the non-vanishing condition \eqref{Anon} for addition,
so do common vertices of ${\cal G}'_\theta{\cal H}$ and $-{\cal G}{\cal H}'_\theta$.
We can further show similarly to the above paragraph that $f$ is semi-definite on $\Omega_\epsilon$ for an $\epsilon>0$.

Finally, we show that $C_0f$ is semi-positive.
In fact, since $V_0({\cal G}):(i_0,j_0)$ and $V_0({\cal H}):(\tilde{i}_0,\tilde{j}_0)$
are the left end-points of ${\cal N}_{\cal G}$ and ${\cal N}_{\cal H}$ respectively,
we have
\begin{eqnarray*}
{\cal G}(\rho,\theta)
=a_{i_0,j_0}({\cal G})\rho^{i_0}\theta^{j_0}+O(\rho^{i_0+1})+O(\theta^{j_0+1}),
\\
{\cal H}(\rho,\theta)
=a_{\tilde{i}_0,\tilde{j}_0}({\cal H})\rho^{\tilde{i}_0}\theta^{\tilde{j}_0}
+O(\rho^{\tilde{i}_0+1})+O(\theta^{\tilde{j}_0+1}).
\end{eqnarray*}
Simple computation shows that
$$
[{\cal G},{\cal H}]_\theta=C_0\rho^{i_0+\tilde{i}_0}\theta^{j_0+\tilde{j}_0-1}
+O(\rho^{i_0+\tilde{i}_0+1})+O(\theta^{j_0+\tilde{j}_0}),
$$
where $C_0$ are defined just before Theorem~\ref{th:finite}.
It follows that $C_0$ is the coefficient of $f$
corresponding to the left end-point of ${\cal N}_f$.
Proposition~\eqref{lm-region} further implies that $C_0f$ is semi-positive.
Thus, the proof of this lemma is completed.
\qquad$\Box$

\subsection{Proofs of main theorems}

Having the above lemmas given in subsections 5.1 and 5.2,
we are ready to prove Theorems \ref{th:finite} and \ref{th:finiteJ2}.

{\bf Proof of Theorem~\ref{th:finite}.}
In order to prove this theorem,
we consider the region
$$
{\cal U}_{\epsilon,\delta}:=\{(\rho,\theta)\in {\mathbb R}^2:0<\rho<\epsilon,|\theta|<\delta\},
$$
a neighborhood of the exceptional direction $\theta=0$ near $O:(0,0)$,
for small $\epsilon$ and $\delta$.
At each point in the region,
the motion is composition of expansion (or contraction) in the $\rho$-direction
and rotation in the $\theta$-direction.
First,
we divide the region into subregions in each of which
the motion in the $\rho$-direction is definitely expansion or contraction.
More concretely,
we need to divide ${\cal U}_{\epsilon,\delta}$ into Z-sectors (\cite[pp.71-83]{SC})
by real branches of the equation ${\cal H}(\rho,\theta)=0$
passing though the origin $O$.

\begin{lm}
If functions ${\cal G}$ and ${\cal H}$ satisfy conditions {\bf(P1)} and {\bf(Q)}
in the one-above case {\bf(J1)},
then equations ${\cal G}(\rho,\theta)=0$ and ${\cal H}(\rho,\theta)=0$
have $N({\cal G})$ and $N({\cal H})$ real branches
on the half-plane $\rho\ge0$ respectively,
where $N({\cal G})$ and $N({\cal H})$ are defined just before Theorem~\ref{th:finite}.
\label{Assertion1}
\end{lm}

{\bf Proof.}
We only give the proof for the equation ${\cal G}(\rho,\theta)=0$
since the proof for the equation ${\cal H}(\rho,\theta)=0$ is similar.
By Lemma~\ref{lm-basic}{\bf(i)},
each nontrivial real branch of the equation ${\cal G}(\rho,\theta)=0$
on the half-plane $\rho\ge 0$
is of the form $\theta=c\rho^\iota+o(\rho^\iota)$,
where $\iota$ is rational such that
$-1/\iota$ is the slope of an edge $E$ of ${\cal N}_{\cal G}$
and $c$ is a nonzero real root of the edge-polynomial ${\cal G}_E$.
For an edge $E$ of ${\cal N}_{\cal G}$,
let $\xi$ be its slope.
Then $\xi\in\mathfrak{S}({\cal G})$,
the slope set of ${\cal N}_{\cal G}$.
As indicated in the beginning of the proof of Lemma~\ref{defsignth1},
{\bf(P1)} implies that
common vertices of Newton polygons of
${\cal G}'_\theta{\cal H}$ and $-{\cal G}{\cal H}'_\theta$
satisfy the non-vanishing condition \eqref{Anon} for addition.
Then Lemma~\ref{J1Ng} holds and ensures that
\begin{eqnarray}
{\cal K}_\xi([{\cal G},{\cal H}]_\theta)
=[{\cal K}_\xi({\cal G}),{\cal K}_\xi({\cal H})]
=({\cal K}_\xi({\cal G}))'_\theta{\cal K}_\xi({\cal H})
-{\cal K}_\xi({\cal G})({\cal K}_\xi({\cal H}))'_\theta
\label{KGKHK4}
\end{eqnarray}
since $\xi\in\mathfrak{S}({\cal G})\subset\mathfrak{S}({\cal G})\cup\mathfrak{S}({\cal H})$.
Note that the polynomial ${\cal K}_\xi([{\cal G},{\cal H}]_\theta)$
has no nonzero real roots by {\bf(Q)}.
Then, all those nonzero real roots of ${\cal K}_\xi({\cal G})$ are of multiplicity 1.
We see from \eqref{defkxip} and \eqref{fonE} that
${\cal K}_\xi({\cal G})={\cal G}_E$ since $\xi$ is the slope of the edge $E$.
Hence, all those nonzero real roots of ${\cal G}_E$ are of multiplicity 1,
and each of them determines one nontrivial real branch of the equation
${\cal G}(\rho,\theta)=0$ by Lemma~\ref{lm-basic}{\bf(ii)}.
We see from \eqref{defChi} that
the equation ${\cal G}(\rho,\theta)=0$ has $\chi({\cal G})$ nontrivial real branches.
Additionally,
it has a trivial real branch $\theta(\rho)\equiv 0$ if and only if $j_{s({\cal G})}>0$
since $j_{s({\cal G})}$ is the ordinate of the right-most vertex of ${\cal N}_{\cal G}$.
By definitions of $N({\cal G})$ and $\chi({\cal G})$,
given just before Theorem~\ref{th:finite},
the equation ${\cal G}(\rho,\theta)=0$ has $N({\cal G})$ real branches
on the half-plane $\rho\ge0$.
Thus, the lemma is proved.
\qquad$\Box$

\begin{figure}[h]
  \centering
  \includegraphics[width=2.5in]{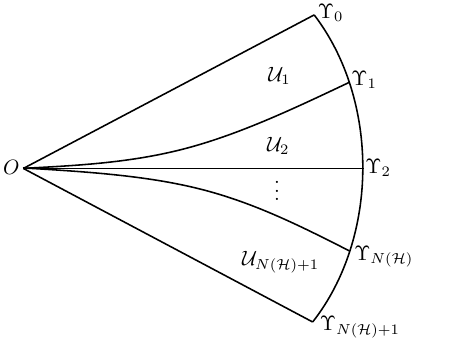}\\
  \caption{Partitions of the region ${\cal U}_{\epsilon,\delta}$.}
  \label{fig:PUed}
\end{figure}

By Lemma~\ref{Assertion1},
we assume that
$\theta_1(\rho),...,\theta_{N({\cal H})}(\rho)$
are all real branches of the equation ${\cal H}(\rho,\theta)=0$
on the half-plane $\rho\ge 0$,
which are ranked as
\begin{eqnarray}
\theta_1(\rho)>\cdots>\theta_{N({\cal H})}(\rho),~~~0<\rho<\epsilon.
\label{thetatheta}
\end{eqnarray}
Assume without loss of generality that
for each $k=1,...,N({\cal H})$
we have $|\theta_k(\rho)|<\delta$ for all $\rho\in[0,\epsilon)$ since $\theta_k(0)=0$.
Define $\theta_0(\rho):=\delta$ and $\theta_{N({\cal H})+1}(\rho):=-\delta$ complementarily.
Then, those $N({\cal H})$ zeros $\theta=\theta_{k}(\rho)$, $k=1,...,N({\cal H})$,
are $n$ branches of the curve ${\cal H}(\rho,\theta)=0$,
denoted by $\Upsilon_k$ respectively,
between the radial segments
$$
\Upsilon_0: \theta=\theta_0(\rho),~0<\rho<\epsilon,~~~\mbox{and}~~~
\Upsilon_{N({\cal H})+1}: \theta=\theta_{N({\cal H})+1}(\rho),~0<\rho<\epsilon.
$$
Those branches divide the region ${\cal U}_{\epsilon,\delta}$ into $N({\cal H})+1$ subregions.
Let ${\cal U}_k$ be the interior of the subregion between $\Upsilon_{k-1}$ and $\Upsilon_{k}$ for each $k=1,...,N({\cal H})+1$,
as shown in Figure~\ref{fig:PUed}.
Note that for $k=2,...,N({\cal H})$ the sector ${\cal U}_k$
contains the exceptional direction $\theta=0$ only because
$\theta_{k-1}(0)=\theta_k(0)=0$, i.e.,
its edges $\Upsilon_{k-1}$ and $\Upsilon_k$ connect with $O$ in the direction $\theta=0$.
However,
as explained in \cite{T-Z1} for Z-sectors and generalized normal sectors,
the sector ${\cal U}_1$ (resp. ${\cal U}_{N({\cal H})+1}$)
may contain another in addition to the exceptional direction $\theta=0$
because its edges $\Upsilon_0$ (resp. $\Upsilon_{N({\cal H})+1}$) and $\Upsilon_1$ (resp. $\Upsilon_{N({\cal H})}$)
connect with $O$ in directions $\theta=\delta$ (resp. $-\delta$) and $\theta=0$ respectively.

With the above division in the region ${\cal U}_{\epsilon,\delta}$,
we are in the position to investigate
the distribution of orbits of system~\eqref{equ:initial}
in the exceptional direction $\theta=0$.
We discuss in subcases {\bf (ia)}, {\bf (ib)}, {\bf (iia)} and {\bf (iib)}
of the theorem.

{\bf Step 1:} {\it Give distribution of orbits
in the exceptional direction $\theta=0$ in subcase~{\bf(ia)},
i.e., $G_0(\theta)\not\equiv 0$ and $C_0<0$.}

In this case,
the exceptional direction $\theta=0$ is isolated.
If an orbit of system~\eqref{equ:initial} connects with $O$ in ${\cal U}_k$
($k=1,...,N({\cal H})+1$),
by Theorem~3.1 of \cite[p.60]{ZZF},
which says that no orbits can approach or leave $O$
in a direction other than an exceptional direction,
it connects with $O$ in the direction $\theta=0$.
Then we only need to determine the class of each Z-sector and
further the number of orbits connecting with $O$ in this Z-sector.
So we consider the sign of $\dot\theta$ on edges of each Z-sector
for the motion in the $\theta$-direction.
It suffices to know
the multiplicities of real branches of the equations
${\cal G}(\rho,\theta)=0$ and ${\cal H}(\rho,\theta)=0$
and the relative positions of those branches.

\noindent{\bf Claim~1.1.}
{\it
Suppose that
the Lie-bracket $[{\cal G},{\cal H}]_\theta$ is semi-negative
on ${\cal U}_{\epsilon,\delta}$.
Then,
\\
{\bf(F1) Decrease:}
for each $k=1,...,N({\cal H})+1$ and each $\rho\in(0,\epsilon)$,
the ratio ${\cal G}(\rho,\theta)/{\cal H}(\rho,\theta)$ is decreasing in
$\theta\in(\theta_k(\rho),\theta_{k-1}(\rho))$;
\\
{\bf(F2) Odd Multiplicity:}
the multiplicity of each real branch of the equations ${\cal G}(\rho,\theta)=0$ and
${\cal H}(\rho,\theta)=0$ on the half-plane $\rho\ge0$ is odd; and
\\
{\bf(F3) Alternation:}
there is one real branch of the equation ${\cal G}(\rho,\theta)=0$
$($and ${\cal H}(\rho,\theta)=0$$)$ between each two adjacent real branches of the equation
${\cal H}(\rho,\theta)=0$ $($and ${\cal G}(\rho,\theta)=0$$)$.
}

In fact,
Lemma~\ref{defsignth1} ensures that
$[{\cal G},{\cal H}]_\theta\le 0$ on ${\cal U}_{\epsilon,\delta}$
since {\bf(P1)} and {\bf(Q)} hold in the one-above case {\bf(J1)} and $C_0<0$ in {\bf(ia)}.
For each $k=1,...,N({\cal H})+1$ we have
\begin{eqnarray}
\frac{\partial}{\partial \theta}\left(\frac{{\cal G}}{{\cal H}}\right)
=\frac{{\cal G}'_\theta{\cal H}-{\cal G}{\cal H}'_\theta}{{\cal H}^2}
=\frac{[{\cal G},{\cal H}]_\theta}{{\cal H}^2}\le0
~~~\mbox{for all}~~~(\rho,\theta) \in {\cal U}_k.
\label{PGH}
\end{eqnarray}
Note that ${\cal U}_k$ is bounded by
$\Upsilon_{k-1}:\theta=\theta_{k-1}(\rho)$ and $\Upsilon_k:\theta=\theta_k(\rho)$.
Then {\bf(F1)} follows directly from \eqref{PGH}.
For {\bf(F2)},
if it is not true, then we assume without loss of generality that in (\ref{thetatheta})
the real branch $\theta_k(\rho)$ of the equation ${\cal H}(\rho,\theta)=0$
is of even multiplicity $2\nu$ ($>0$).
It follows from \eqref{Wei}, \eqref{WRQ} and \eqref{csq} that
$$
{\cal H}(\rho,\theta)=(\theta-\theta_k(\rho))^{2\nu}\widetilde{\cal H}(\rho,\theta),
$$
where $\widetilde{\cal H}\in\mathbb{R}\{\rho^{1/d},\theta\}$ for a positive integer $d$ and
$\widetilde{\cal H}(\rho,\theta_k(\rho))\ne0$.
We compute that
$$
[{\cal G},{\cal H}]_\theta=-2\nu(\theta-\theta_k(\rho))^{2\nu-1}{\cal G}\widetilde{\cal H}
+(\theta-\theta_k(\rho))^{2\nu}[{\cal G},\widetilde{\cal H}]_\theta,
$$
implying that $\theta_k(\rho)$ is a real branch of $[{\cal G},{\cal H}]_\theta$ of multiplicity at least $2\nu-1$.
On the other hand,
$\theta_k(\rho)$ is not a real branch of ${\cal G}$ because ${\cal H}(\rho,\theta_k(\rho))\equiv0$
and
\begin{eqnarray}
\rho^{2m}({\cal G}^2(\rho,\theta)+{\cal H}^2(\rho,\theta))={\cal X}^2(x,y)+{\cal Y}^2(x,y)\ne0
\label{equ:GHXY}
\end{eqnarray}
in a vicinity of the isolated singular point $O$ of system~\eqref{equ:initial}.
Thus,
$\theta_k(\rho)$ is a real branch of $[{\cal G},{\cal H}]_\theta$ of odd multiplicity $2\nu-1$,
which contradicts to the semi-negativeness of $[{\cal G},{\cal H}]_\theta$
and proves {\bf(F2)}.
For {\bf(F3)},
if it is not true,
we assume without loss of generality that
between two adjacent real branches $\theta_k(\rho)$ and $\theta_{k+1}(\rho)$
of the equation ${\cal H}(\rho,\theta)=0$, seen in (\ref{thetatheta}),
there are either no or at least two real branches of the equation ${\cal G}(\rho,\theta)=0$.
In the first situation,
both ${\cal G}$ and ${\cal H}$ have definite signs in the region ${\cal U}_{k+1}$,
implying that for each $\rho\in(0,\epsilon)$,
$$
\lim_{\theta\to\theta_k(\rho)-0}\frac{{\cal G}(\rho,\theta)}{{\cal H}(\rho,\theta)}=+\infty(\mbox{or}-\infty)
~~~\mbox{and}~~~
\lim_{\theta\to\theta_{k+1}(\rho)+0}\frac{{\cal G}(\rho,\theta)}{{\cal H}(\rho,\theta)}=+\infty(\mbox{or}-\infty)
$$
since ${\cal G}$ and ${\cal H}$ do not have a common root as indicated in \eqref{equ:GHXY},
a contradiction to {\bf{(F1)}}.
In the second situation,
${\cal G}(\rho,\theta)/{\cal H}(\rho,\theta)=0$
at those real branches of the equation ${\cal G}(\rho,\theta)=0$ in ${\cal U}_{k+1}$,
which also contradicts to {\bf(F1)} and proves {\bf (F3)}.
Consequently, {\bf Claim~1.1} is proved.

In the following, we determine the class of each Z-sector.

\noindent{\bf Claim~1.2.}
{\it
For each $k=2,...,N({\cal H})$
the Z-sector ${\cal U}_k$ is of Class~II
and contains a unique orbit connecting with $O$.
The Z-sector ${\cal U}_1$ $($and ${\cal U}_{N({\cal H})+1}$$)$
is of either Class~II,
containing a unique orbit in connect with $O$,
or Class~III,
containing no orbits in connect with $O$.
}

Actually,
as indicated in \cite[pp.81-83]{SC},
there are totally three classes of Z-sector.
For each $k=1,...,N({\cal H})+1$,
lying between $\Upsilon_{k-1}$ and $\Upsilon_k$,
the Z-sector ${\cal U}_k$ is of Class I if
\begin{eqnarray}
\dot\rho|_{{\cal U}_k}\dot\theta|_{\Upsilon_{k-1}}>0
~~~\mbox{and}~~~
\dot\rho|_{{\cal U}_k}\dot\theta|_{\Upsilon_k}<0,
\label{con-CI}
\end{eqnarray}
is of Class II if
\begin{eqnarray}
\dot\rho|_{{\cal U}_k}\dot\theta|_{\Upsilon_{k-1}}<0
~~~\mbox{and}~~~
\dot\rho|_{{\cal U}_k}\dot\theta|_{\Upsilon_k}>0,
\label{con-CII}
\end{eqnarray}
and is of Class III if
\begin{eqnarray}
\dot\theta|_{\Upsilon_{k-1}}\dot\theta|_{\Upsilon_k}>0.
\label{con-CIII}
\end{eqnarray}
By results~{\bf(b)}, {\bf(c)} and {\bf(d)} of \cite[pp.81-83]{SC},
there are infinitely many orbits connecting with $O$
if ${\cal U}_k$ is of Class I;
either infinitely many orbits or exact one orbit connecting with $O$
if ${\cal U}_k$ is of Class II;
and
either infinitely many orbits or no orbits connecting with $O$
if ${\cal U}_k$ is of Class III.
On the other hand,
in each ${\cal U}_k$
we have ${\cal H}(\rho,\theta)\ne0$.
Further,
\eqref{PGH} ensures that
system~\eqref{equ:polar system} has at most one orbit connecting with $O$ in ${\cal U}_k$.
In fact,
assume that $\tilde{\theta}_1(\rho)$ and $\tilde{\theta}_2(\rho)$ are such two orbits
and $\tilde{\theta}_1(\rho)<\tilde{\theta}_2(\rho)$ for all $\rho\in(0,\epsilon)$.
By the Mean Value Theorem and \eqref{PGH},
there is $\xi(\rho)\in(\tilde{\theta}_1(\rho),\tilde{\theta}_2(\rho))$ such that
{\small
$$
\frac{d(\tilde{\theta}_2(\rho)-\tilde{\theta}_1(\rho))}{d\rho}
=\left.\frac{{\cal G}}{\rho{\cal H}}\right|_{(\rho,\tilde{\theta}_2(\rho))}-
\left.\frac{{\cal G}}{\rho{\cal H}}\right|_{(\rho,\tilde{\theta}_1(\rho))}
=\left.\frac{\partial}{\partial \theta}\left(\frac{{\cal G}}{\rho{\cal H}}\right)\right|_{(\rho,\xi(\rho))}(\tilde{\theta}_2(\rho)-\tilde{\theta}_1(\rho))
\le 0.
$$
}Associated with the initial data
$\tilde{\theta}_2(0)=\tilde{\theta}_1(0)=0$,
it gives that
$\tilde{\theta}_2(\rho)\le\tilde{\theta}_1(\rho)$ for all $\rho\in(0,\epsilon)$,
a contradiction to our assumption.
Then
${\cal U}_k$ can only be of either Class II or Class III for each $k=1,...,N({\cal H})$.
Moreover,
each of ${\cal U}_1$ and ${\cal U}_{N({\cal H})+1}$
contains a unique orbit connecting with $O$ in the case of Class II
or no orbits connecting with $O$ in the case of Class III.
On the other hand,
for each $k=2,...,N({\cal H})$,
fact {\bf(F3)} implies that
the equation ${\cal G}(\rho,\theta)=0$ has exactly one real branch
(odd multiplicity by {\bf(F2)}) in ${\cal U}_k$,
which implies that $\dot\theta|_{\Upsilon_{k-1}}\dot\theta|_{\Upsilon_k}<0$.
It follows from \eqref{con-CIII} that ${\cal U}_k$ is of Class~II rather than Class~III
and it contains a unique orbit connecting with $O$.
Thus {\bf Claim~1.2} is proved.

Having above two claims,
we are in the position to complete the proof in the subcase {\bf(ia)}.
By Lemma \ref{Assertion1}, we only need the following.

\noindent{\bf Claim~1.3.}
{\it
The number of orbits of system \eqref{equ:initial}
connecting with $O$ in the direction $\theta=0$ is equal to
the number of real branches of the equation ${\cal G}(\rho,\theta)=0$
on the half-plane $\rho\ge 0$.
}

If {\bf Claim~1.3} is true,
then by Lemma \ref{Assertion1}
system \eqref{equ:initial} has exactly $N({\cal G})$ orbits connecting with $O$
in the direction $\theta=0$,
making $\max\{0,N({\cal G})-1\}$ h-tsectors in this direction,
and therefore, conclusions of {\bf(ia)} hold.
In order to prove {\bf Claim~1.3},
we consider the following four circumstances:

\begin{center}
\begin{tabular}{ll}
{\bf(C1)} $a_{i_0,j_0}({\cal G})a_{\tilde{i}_0,\tilde{j}_0}({\cal H})<0$, $j_0$ is even,
& {\bf(C2)} $a_{i_0,j_0}({\cal G})a_{\tilde{i}_0,\tilde{j}_0}({\cal H})<0$, $j_0$ is odd, \\[5pt]
{\bf(C3)} $a_{i_0,j_0}({\cal G})a_{\tilde{i}_0,\tilde{j}_0}({\cal H})>0$, $j_0$ is even,
& {\bf(C4)} $a_{i_0,j_0}({\cal G})a_{\tilde{i}_0,\tilde{j}_0}({\cal H})>0$, $j_0$ is odd,
\end{tabular}
\end{center}
where $a_{i_0,j_0}({\cal G})$ (resp. $a_{\tilde{i}_0,\tilde{j}_0}({\cal H})$)
is the coefficient of ${\cal G}$ (resp. ${\cal H}$)
corresponding to the left-most vertex $(i_0,j_0)$
(resp. $(\tilde{i}_0,\tilde{j}_0)$) of ${\cal N}_{\cal G}$ (resp. ${\cal N}_{\cal H}$).
We only give the proof in {\bf(C1)} because the proofs in others are similar.
Assume without loss of generality that
\begin{align}
a_{i_0,j_0}({\cal G})<0~~~\mbox{and}~~~a_{\tilde{i}_0,\tilde{j}_0}({\cal H})>0
\label{assume-c1}
\end{align}
in {\bf(C1)}; otherwise we can similarly discussed.
In order to determine the number of orbits connecting with $O$,
we need to investigate the sign of $\dot \theta$ at each real branches $\Upsilon_k$,
which is given by
\begin{equation}
\dot\theta|_{\Upsilon_k}
\left\{
\begin{array}{lllll}
<0  &\mbox{for even}~k\in\{0,...,N({\cal H})\},
\\
>0  &\mbox{for odd}~k\in\{0,...,N({\cal H})\}
\end{array}
\right.
~~~\mbox{and}~~~
\dot\theta|_{\Upsilon_{N({\cal H})+1}}<0.
\label{iat0n0n}
\end{equation}
In fact,
since $(i_0,j_0)$ is the left-most vertex of ${\cal N}_{\cal G}$,
the left edge of $(i_0,j_0)$ is a vertical ray.
So we have
$\Delta({\cal G})=\Delta_0({\cal G})\cup \Delta_1({\cal G})$, where
$\Delta_0({\cal G}):=\{(i,j)\in\Delta({\cal G}):i=i_0,j\ge j_0\}$ and
$\Delta_1({\cal G}):=\{(i,j)\in\Delta({\cal G}):i\ge i_0+1\}$.
Then we see from \eqref{expGH} and the definition of $\Delta({\cal G})$,
given just below \eqref{expGH},
that
\begin{align}
{\cal G}(\rho,\theta)
&=\sum_{(i,j)\in\Delta({\cal G})} a_{i,j}({\cal G})\rho^i\theta^j
=\sum_{(i,j)\in\Delta_0({\cal G})} a_{i,j}({\cal G})\rho^i\theta^j
+\sum_{(i,j)\in\Delta_1({\cal G})} a_{i,j}({\cal G})\rho^i\theta^j
\nonumber
\\
&=\rho^{i_0}\{a_{i_0,j_0}({\cal G})\theta^{j_0}+O(\theta^{j_0+1})+O(\rho)\}.
\label{Gi0j0}
\end{align}
It follows that $\dot\theta|_{\Upsilon_0}={\cal G}(\rho,\delta)<0$
because of \eqref{assume-c1}.
Next, we determine the sign of $\dot\theta|_{\Upsilon_1}$.
By \eqref{assume-c1},
we compute similarly to \eqref{Gi0j0} that
$$
{\cal H}(\rho,\theta)|_{\Upsilon_0}
={\cal H}(\rho,\delta)
=\rho^{\tilde{i}_0}\{a_{\tilde{i}_0,\tilde{j}_0}({\cal H})\delta^{\tilde{j}_0}
+O(\delta^{\tilde{j}_0+1})+O(\rho)\}
>0
$$
and therefore by continuity ${\cal H}(\rho,\theta)>0$ in ${\cal U}_1$,
the first Z-sector given just before {\bf Step 1}
bounded by $\Upsilon_0$ and $\Upsilon_1$.
Note that $\theta_1(\rho)$, the expression of $\Upsilon_1$,
is a zero of ${\cal H}$.
By fact {\bf(F1)} of {\bf Claim~1.1},
\begin{eqnarray}
\lim_{\theta\to\theta_1(\rho)+0}
\frac{{\cal G}(\rho,\theta)}{{\cal H}(\rho,\theta)}=+\infty
\label{limit-infty}
\end{eqnarray}
since ${\cal G}$ and ${\cal H}$ do not have a common root as indicated in \eqref{equ:GHXY}.
It follows that $\dot\theta|_{\Upsilon_1}={\cal G}(\rho,\theta_1(\rho))>0$.
For the same reason, we also have the same limit as (\ref{limit-infty}) as $\theta\to\theta_k(\rho)+0$,
for each $k=2,...,N({\cal H})$.
By fact {\bf(F2)} of {\bf Claim~1.1},
each real branch $\theta_k(\rho)$ of the equation ${\cal H}(\rho,\theta)=0$
is of odd multiplicity,
implying that
${\cal H}(\rho,\theta)$ in ${\cal U}_k$ has an opposite sign to that in ${\cal U}_{k+1}$.
It follows from the inequality ${\cal H}(\rho,\theta)>0$ in ${\cal U}_1$,
given just before \eqref{limit-infty}, that
\begin{equation}
{\cal H}(\rho,\theta)
\left\{
\begin{array}{lllll}
>0 &\mbox{in}~{\cal U}_k~\mbox{for odd}~k\in\{1,...,N({\cal H})+1\},
\\
<0 &\mbox{in}~{\cal U}_k~\mbox{for even}~k\in\{1,...,N({\cal H})+1\}.
\end{array}
\right.
\label{equ:Huk}
\end{equation}
Thus
the first inequality in \eqref{iat0n0n} is obtained similarly for all $k=2,...,N({\cal H})$.
For a complement,
we see from \eqref{Gi0j0} that
$$
\dot\theta|_{\Upsilon_{N({\cal H})+1}}
={\cal G}(\rho,-\delta)
=\rho^{i_0}\{a_{i_0,j_0}({\cal G})(-\delta)^{j_0}+O(\delta^{j_0+1})+O(\rho)\}
<0
$$
since \eqref{assume-c1} holds and $j_0$ is even.
Up to now, the claimed \eqref{iat0n0n} is proved.

We see from \eqref{iat0n0n} and \eqref{equ:Huk} that
$\dot\rho|_{{\cal U}_1}\dot\theta|_{\Upsilon_0}<0$
and
$\dot\rho|_{{\cal U}_1}\dot\theta|_{\Upsilon_1}>0$.
By \eqref{con-CII},
the Z-sector ${\cal U}_1$ is of Class~II rather than Class~III.
Then {\bf Claim~1.2} implies that
${\cal U}_k$ has one orbit connecting with $O$ for each $k=1,...,N({\cal H})$.
However,
the last Z-sector ${\cal U}_{N({\cal H})+1}$
contains either no orbits or exact one orbit connecting with $O$.
In fact,
\begin{description}
\item[(S1)]
when $N({\cal H})$ is even,
the inequality
$\dot\theta|_{\Upsilon_{N({\cal H})}}\dot\theta|_{\Upsilon_{N({\cal H})+1}}>0$
in \eqref{iat0n0n} implies by \eqref{con-CIII} that
${\cal U}_{N({\cal H})+1}$ is of Class~III and therefore
${\cal U}_{N({\cal H})+1}$ contains no orbits connecting with $O$
by {\bf Claim~1.2};

\item[(S2)]
when $N({\cal H})$ is odd,
we see from \eqref{iat0n0n} and \eqref{equ:Huk} that
$\dot\rho|_{{\cal U}_{N({\cal H})+1}}\dot\theta|_{\Upsilon_{N({\cal H})}}<0$
and
$\dot\rho|_{{\cal U}_{N({\cal H})+1}}\dot\theta|_{\Upsilon_{N({\cal H})+1}}>0$,
which implies that ${\cal U}_{N({\cal H})+1}$ contains one orbit
connecting with $O$ as we saw in ${\cal U}_1$ above.
\end{description}
Consequently,
system~\eqref{equ:initial} has totally $N({\cal H})$
(or $N({\cal H})+1$) orbits connecting with $O$ in the direction $\theta=0$
in the case {\bf(S1)} (or {\bf(S2)}).

Finally, we prove that the number of real branches of
the equation ${\cal G}(\rho,\theta)=0$ is equal to $N({\cal H})$ and $N({\cal H})+1$
in the case {\bf(S1)} and {\bf(S2)} respectively,
showing that the number or orbits connecting with $O$ is the number of real branches.
In fact, in the case {\bf(S1)} we have
$\dot\theta|_{\Upsilon_{k-1}}\dot\theta|_{\Upsilon_k}<0$ for all $k=1,...,N({\cal H})$
and
$\dot\theta|_{\Upsilon_{N({\cal H})}}\dot\theta|_{\Upsilon_{N({\cal H})+1}}>0$
in \eqref{iat0n0n},
which implies by {\bf (F2)} and {\bf (F3)} that
the equation ${\cal G}(\rho,\theta)=0$ has $N({\cal H})$ real branches;
in the case {\bf(S2)} we have
$\dot\theta|_{\Upsilon_{k-1}}\dot\theta|_{\Upsilon_k}<0$ for each $k=1,...,N({\cal H})+1$
in \eqref{iat0n0n},
which implies by {\bf (F2)} and {\bf (F3)} that
the equation ${\cal G}(\rho,\theta)=0$ has $N({\cal H})+1$ real branches.
Thus {\bf Claim~1.3} is proved and the proof in {\bf Step~1} is completed.

{\bf Step~2:} {\it Give distribution of orbits
in the exceptional direction $\theta=0$ in the subcase {\bf(ib)},
i.e., $G_0(\theta)\not\equiv 0$ and $C_0>0$.}

We use the same procedure as in {\bf Step~1} to discuss this subcase.
By Lemma~\ref{defsignth1},
$[{\cal G},{\cal H}]_\theta\ge 0$ on ${\cal U}_{\epsilon,\delta}$
since conditions {\bf(P1)} and {\bf(Q)} hold in the one-above case {\bf(J1)}
and $C_0>0$.
Correspondingly to {\bf Claims~1.1}-{\bf 1.3} in {\bf Step~1},
we have the following:


\noindent{\bf Claim~2.1.}
{\it
Facts {\bf(F2)} and {\bf(F3)} of {\bf Claim~1.1} still hold and moreover we have
\vspace{-15pt}
\begin{description}
 \item[(F1$'$) Increase:]
{\it
for each $k=1,...,N({\cal H})+1$ and each $\rho\in(0,\epsilon)$,
the ratio ${\cal G}(\rho,\theta)/{\cal H}$ $(\rho,\theta)$ is increasing
in $\theta\in(\theta_k(\rho),\theta_{k-1}(\rho))$.
}
\end{description}
}

\noindent{\bf Claim~2.2.}
{\it
For each $k=2,...,N({\cal H})$
the Z-sector ${\cal U}_k$ is of Class~I and
contains infinitely many orbits connecting with $O$.
The Z-sector ${\cal U}_1$ $ ($and ${\cal U}_{N({\cal H})+1}$$)$ is of
either Class~I, containing infinitely many orbits in connect with $O$,
or Class~III, containing no orbits in connect with $O$.
}

\noindent{\bf Claim~2.3.}
{\it
The number of orbits of system \eqref{equ:initial} connecting with $O$
in the direction $\theta=0$
is zero
if the equation ${\cal G}(\rho,\theta)=0$ has no real branches
on the half-plane $\rho\ge 0$.
Otherwise,
system \eqref{equ:initial} has infinitely many such orbits,
which make finitely many e-tsectors.
Moreover,
the number of  e-tsectors is equal to
the number of real branches of the equation ${\cal G}(\rho,\theta)=0$
on the half-plane $\rho\ge 0$ minus 1.
}

The proof of {\bf Claim~2.3} needs {\bf Claim~2.1} and {\bf Claim~2.2}.
If {\bf Claim~2.3} is true, then conclusions of {\bf(ib)} follow directly from Lemma~\ref{Assertion1}.

In what follows, we need to prove {\bf Claims~2.1}-{\bf 2.3}.
The semi-positiveness of the Lie-bracket $[{\cal G},{\cal H}]_\theta$ implies that
facts {\bf(F2)} and {\bf(F3)} of {\bf Claim~1.1} still hold. Moreover,
correspondingly to {\bf(F1)},
we have {\bf (F1$'$)} in {\bf Claim~2.1} because
the inequality `$\le$' in \eqref{PGH} becomes `$\ge$'
when $[{\cal G},{\cal H}]_\theta \ge 0$ on ${\cal U}_{\epsilon,\delta}$.

Next, we prove {\bf Claim~2.2}.
We first consider Z-sectors ${\cal U}_k$, $k=2,...,N({\cal H})$, each of which
is bounded by
$\Upsilon_{k-1}:\theta=\theta_{k-1}(\rho)$ and $\Upsilon_k:\theta=\theta_k(\rho)$,
the two adjacent real branches
of the equation ${\cal H}(\rho,\theta)=0$
as indicated just below \eqref{thetatheta}.
Then {\bf(F1$'$)} implies that
$$
\lim_{\theta\to\theta_{k-1}(\rho)-0}
\frac{{\cal G}(\rho,\theta)}{{\cal H}(\rho,\theta)}=+\infty
~~~\mbox{and}~~~
\lim_{\theta\to\theta_k(\rho)+0}
\frac{{\cal G}(\rho,\theta)}{{\cal H}(\rho,\theta)}=-\infty
$$
since ${\cal G}$ and ${\cal H}$ do not have a common root as indicated in \eqref{equ:GHXY}.
We see from \eqref{equ:polar system} and the above two limits that
$\dot\rho|_{{\cal U}_k}\dot\theta|_{\Upsilon_{k-1}}>0$
and
$\dot\rho|_{{\cal U}_k}\dot\theta|_{\Upsilon_k}<0$.
Then the Z-sector ${\cal U}_k$ is of Class~I by \eqref{con-CI}, which implies
as indicated just below \eqref{con-CIII}
that
infinitely many orbits connect with $O$ in ${\cal U}_k$
for each $k=2,...,N({\cal H})$.

Having known ${\cal U}_k$ for $k=2,...,N({\cal H})$,
we need to discuss ${\cal U}_k$ for $k=1$ and $N({\cal H})+1$ complementarily.
Region ${\cal U}_1$,
bounded by $\Upsilon_0$ and $\Upsilon_1$,
cannot be a Z-sector of Class II.
Otherwise, by the definition of Class~II given in \eqref{con-CII},
$\dot\rho|_{{\cal U}_1}\dot\theta|_{\Upsilon_0}<0$ and
$\dot\rho|_{{\cal U}_1}\dot\theta|_{\Upsilon_1}>0$.
It follows from \eqref{equ:polar system} that
${\cal G}/{\cal H}<0$ near $\Upsilon_0$ and
${\cal G}/{\cal H}>0$ near $\Upsilon_1$ in the Z-sector ${\cal U}_1$,
a contradiction to {\bf(F1$'$)}.
This proves that ${\cal U}_1$ is of either Class I or Class III.
Clearly,
if ${\cal U}_1$ is of Class~I,
then ${\cal U}_1$ contains infinitely many orbits connecting with $O$;
if ${\cal U}_1$ is of Class~III,
then ${\cal U}_1$ contains either infinitely many or no orbits connecting with $O$,
as indicated just below \eqref{con-CIII}.
In order to determine the number of orbits in Class III,
we divide ${\cal U}_1$ into two parts:
$$
{\cal U}_1^+:={\cal U}_1\cap\{(\rho,\theta)\in\mathbb{R}^2:\rho>0,\theta\ge0\},
$$
the part lying in the closure of the first quadrant,
and its complement ${\cal U}_1\backslash{\cal U}_1^+$,
the part lying in the fourth quadrant.
There are totally 4 situations:
{\bf (2a)} $\dot\rho|_{{\cal U}_1}>0$ and $\dot\theta|_{\Upsilon_0}<0$,
{\bf (2b)} $\dot\rho|_{{\cal U}_1}>0$ and $\dot\theta|_{\Upsilon_0}>0$,
{\bf (2c)} $\dot\rho|_{{\cal U}_1}<0$ and $\dot\theta|_{\Upsilon_0}<0$, and
{\bf (2d)} $\dot\rho|_{{\cal U}_1}<0$ and $\dot\theta|_{\Upsilon_0}>0$.
We only discuss in {\bf(2a)} since the others are similar.
By the inequality $\dot\rho|_{{\cal U}_1}>0$ and the continuity of ${\cal H}$,
$\dot\rho|_{\Upsilon_0}>0$.
Then
\begin{eqnarray}
\left.\frac{\rho{\cal H}(\rho,\theta)}{{\cal G}(\rho,\theta)}\right|_{\Upsilon_0}
=\left.\frac{\dot\rho}{\dot\theta}\right|_{\Upsilon_0}
<0.
\label{frac-right}
\end{eqnarray}
On the other hand,
we see from \eqref{Gi0j0} that
\begin{align*}
{\cal G}(\rho,\theta)
=\rho^{i_0}\{a_{i_0,j_0}({\cal G})\theta^{j_0}+O(\theta^{j_0+1})+O(\rho)\}
=a_{0,j_0}({\cal G})\theta^{j_0}+O(\theta^{j_0+1})+O(\rho),
\end{align*}
where we note $i_0=0$ because $G_0(\theta)\not\equiv 0$.
Similarly,
$
{\cal H}(\rho,\theta)=\rho^{\tilde{i}_0}\{a_{\tilde{i}_0,\tilde{j}_0}({\cal H})\theta^{\tilde{j}_0}+O(\theta^{\tilde{j}_0+1})+O(\rho)\}.
$
Choosing small $\delta>0$,
we see from (\ref{frac-right}) that
$$
\left.\frac{\rho{\cal H}(\rho,\theta)}{{\cal G}(\rho,\theta)}\right|_{\Upsilon_0}
=\left\{\frac{a_{\tilde{i}_0,\tilde{j}_0}({\cal H})\delta^{\tilde{j}_0}+O(\delta^{\tilde{j}_0+1})}
{a_{0,j_0}({\cal G})\delta^{j_0}+O(\delta^{j_0+1})}+O(\rho)\right\}\rho^{\tilde{i}_0+1}
>
-c \rho^{\tilde{i}_0+1}~~\mbox{for all}~\rho\in(0,\epsilon)
$$
for a constant $c>0$.
Then this result and \eqref{frac-right} imply that
\begin{eqnarray}
-c \rho
\le
-c \rho^{\tilde{i}_0+1}
<
\frac{\dot\rho}{\dot\theta}
<0
~~~
\mbox{for all}~(\rho,\theta)\in \Upsilon_0.
\label{equ:drrdt0}
\end{eqnarray}
Note that the Z-sector ${\cal U}_1$ is of Class III and has two edges $\Upsilon_0$ and $\Upsilon_1$.
Since $\dot\theta|_{\Upsilon_0}<0$ and $\dot\rho|_{{\cal U}_1}>0$ in situation {\bf (2a)},
we see from \eqref{con-CIII} that $\dot\theta|_{\Upsilon_1}<0$ and therefore
$\dot\rho/\dot\theta\le 0$ for all $(\rho,\theta)\in \Upsilon_1$.
This together with (\ref{equ:drrdt0}) implies that
\begin{eqnarray}
-c \rho
<
\frac{\dot\rho}{\dot\theta}
<0
~~~
\mbox{for all}~(\rho,\theta)\in {\cal U}_1
\label{equ:drrdt}
\end{eqnarray}
by the monotonicity given in {\bf(F1$'$)}.
The right hand side of the above shows that
$\dot \theta<0$ on ${\cal U}_1^+$ because $\dot \rho|_{{\cal U}_1}>0$
in situation {\bf (2a)}.
This implies that $\rho(t)$ and $\theta(t)$ cannot approach simultaneously to $0$
in ${\cal U}_1^+$ as $t\to+\infty$ or $t\to -\infty$.
Hence system~\eqref{equ:polar system} has no orbits connecting with $O$ in ${\cal U}_1^+$.
In the complement ${\cal U}_1\setminus {\cal U}_1^+$,
we assume that system~\eqref{equ:polar system} has an orbit connecting with $O$.
Since the orbit $(\rho(t),\theta(t))$ lies in the fourth quadrant and
$\dot\rho \dot\theta<0$ on ${\cal U}_1\setminus {\cal U}_1^+$ by the right hand side of \eqref{equ:drrdt},
we can choose points $(\rho_0,\theta_0)$ and $(\rho,\theta)$ on this orbit
such that $\theta_0<\theta<0$ and $0<\rho<\rho_0$.
Integrating the both sides of the inequality on the left hand side of \eqref{equ:drrdt},
we obtain
$$
-c(\theta-\theta_0)<\ln\rho-\ln \rho_0,
$$
the right hand side of which tends to $-\infty$ as $\rho\to 0^+$.
This contradiction shows that system~\eqref{equ:polar system} has no orbits connecting with $O$
in ${\cal U}_1\setminus {\cal U}_1^+$.
Summarizing the above in ${\cal U}_1^+$ and ${\cal U}_1\setminus {\cal U}_1^+$,
we conclude that
${\cal U}_1$ contains no orbits connecting with $O$
if it is of Class~III.
Similarly to ${\cal U}_{1}$, we can prove that
the Z-sector ${\cal U}_{N({\cal H})+1}$ has the same conclusion as ${\cal U}_1$.
Thus, {\bf Claim~2.2} is proved.


Finally,
we prove {\bf Claim~2.3}.
Actually, we also need to discuss
in the same four circumstances {\bf (C1)}-{\bf (C4)} as listed in the proof of {\bf Claim~1.3}.
We only give the proof in {\bf (C1)} because the proofs in others are similar.
Assume without loss of generality that
$a_{i_0,j_0}({\cal G})<0$ and $a_{\tilde{i}_0,\tilde{j}_0}({\cal H})>0$ in {\bf(C1)},
i.e., \eqref{assume-c1} holds.
Similarly to \eqref{iat0n0n} and \eqref{equ:Huk},
we can use {\bf(F1$'$)} and {\bf(F2)} to show that
\begin{equation}
\dot \theta|_{\Upsilon_0}<0,~
\dot\theta|_{\Upsilon_k}
\left\{
\begin{array}{lllll}
<0~\mbox{for odd}~k\in\{1,...,N({\cal H})\},
\\
>0~\mbox{for even}~k\in\{1,...,N({\cal H})\},
\end{array}
\right.
~\mbox{and}~
\dot\theta|_{\Upsilon_{N({\cal H})+1}}<0,
\label{equ:ttkt}
\end{equation}
and
\begin{equation}
\dot \rho|_{{\cal U}_k}
\left\{
\begin{array}{llll}
>0 &\mbox{for odd}~k\in\{1,...,N({\cal H})+1\},
\\
<0 &\mbox{for even}~k\in\{1,...,N({\cal H})+1\}.
\end{array}
\right.
\label{equ:rhosignib}
\end{equation}
We discuss in the two cases: either the equation ${\cal G}(\rho,\theta)=0$ has no real branches,
or the equation ${\cal G}(\rho,\theta)=0$ has at least one real branch.

In the first case,
{\bf(F3)} implies that the equation ${\cal H}(\rho,\theta)=0$
has at most 1 real branches, i.e., $N({\cal H})=0$ or $1$.
When $N({\cal H})=0$,
we have ${\cal U}_{\epsilon,\delta}={\cal U}_1$,
as indicated in Figure \ref{fig:PUed}.
It follows from \eqref{equ:ttkt} that
\begin{eqnarray}
\dot\theta|_{\Upsilon_0}\dot\theta|_{\Upsilon_1}>0,
\label{theta-NH0}
\end{eqnarray}
implying by \eqref{con-CIII} that
the Z-sector ${\cal U}_{\epsilon,\delta}$ is of Class~III.
By {\bf Claim~2.2},
no orbits connect with $O$ in ${\cal U}_{\epsilon,\delta}$.
On the other hand,
when $N({\cal H})=1$,
the region ${\cal U}_{\epsilon,\delta}$ is divided into
Z-sectors ${\cal U}_1$ and ${\cal U}_2$ (i.e., ${\cal U}_{N({\cal H})+1}$).
We see from \eqref{equ:ttkt} that
\begin{eqnarray}
\dot\theta|_{\Upsilon_0}\dot\theta|_{\Upsilon_1}>0
~~~\mbox{and}~~~
\dot\theta|_{\Upsilon_1}\dot\theta|_{\Upsilon_2}>0,
\label{theta-NH1}
\end{eqnarray}
implying by \eqref{con-CIII} that
${\cal U}_1$ and ${\cal U}_2$ are both of Class~III.
By {\bf Claim~2.2},
orbits connect with $O$ in neither ${\cal U}_1$ nor ${\cal U}_2$,
implying that
no orbits connect with $O$ in the direction $\theta=0$.
Consequently, we obtain the first part of {\bf Claim~2.3}.

Oppositely, in the second case, i.e.,
the equation ${\cal G}(\rho,\theta)=0$ has at least 1 real branches,
we have $N({\cal H})\ge 2$.
Otherwise, if $N({\cal H})=0$,
as above we have \eqref{theta-NH0},
which together with {\bf(F2)} and {\bf (F3)} implies that
the equation ${\cal G}(\rho,\theta)=0$ has no real branches
in ${\cal U}_{\epsilon,\delta}$, a contradiction;
if $N({\cal H})=1$,
as above we have \eqref{theta-NH1},
which together with {\bf(F2)} and {\bf (F3)} leads to the same contradiction.
Since we have $N({\cal H})\ge 2$,
by {\bf Claim~2.2}
the Z-sector ${\cal U}_2$ is of Class~I.
Moreover,
both edges of ${\cal U}_2$ are tangent to the positive $x$-axis at $O$.
It implies that
system~\eqref{equ:initial} has infinitely many orbits
connecting with $O$ in the direction $\theta=0$.

In what follows,
we further determine the number of e-tsectors
formed by those orbits connecting with $O$, which only exist in the second case.
For even $k\in\{2,...,N({\cal H})-1\}$,
starting from each point $P_k$ on $\Upsilon_k$, the orbit $\phi(t,P_k)$
satisfies that
\begin{align}
\phi(t,P_k)\to O~\mbox{in}{~\cal U}_k~\mbox{as}~t \to +\infty
~~~\mbox{and}~~~
\phi(t,P_k)\to O~\mbox{in}~{\cal U}_{k+1}~\mbox{as}~t \to -\infty
\label{phitpk}
\end{align}
by \eqref{equ:ttkt} and \eqref{equ:rhosignib}.
Thus, this orbit from $-\infty$ to $+\infty$ forms a homoclinic orbit.
Similarly, for odd $k\in\{2,...,N({\cal H})-1\}$,
starting from each point on $\Upsilon_k$,
the orbit also forms a homoclinic orbit in the opposite direction.
Since there are totally $N({\cal H})-2$ curves $\Upsilon_k$\,s, the system has
at least $N({\cal H})-2$ e-tsectors.

In order to determine the exact number of e-tsectors,
we need to know whether there are homoclinic orbits passing through $\Upsilon_1$ and $\Upsilon_{N({\cal H})}$,
for which a prerequisite condition is that there are orbits connecting with $O$
in ${\cal U}_1$ and ${\cal U}_{N({\cal H})+1}$.
Actually,
no orbits connect with $O$ in ${\cal U}_1$
by {\bf Claim~2.2}
because the inequality $\dot\theta|_{\Upsilon_0}\dot\theta|_{\Upsilon_1}>0$ in \eqref{equ:ttkt}
implies by \eqref{con-CIII} that
the Z-sector ${\cal U}_1$ is of Class~III.
In contrast,
${\cal U}_{N({\cal H})+1}$ contains either infinitely many or no orbits connecting with $O$.
In fact,
\begin{description}
\item[(S1$'$)]
when $N({\cal H})$ is even,
we see from \eqref{equ:ttkt} and \eqref{equ:rhosignib} that
$\dot\rho|_{{\cal U}_{N({\cal H})+1}}\dot\theta|_{\Upsilon_{N({\cal H})}}>0$
and
$\dot\rho|_{{\cal U}_{N({\cal H})+1}}\dot\theta|_{\Upsilon_{N({\cal H})+1}}<0$,
which implies by \eqref{con-CI} that the Z-sector ${\cal U}_{N({\cal H})+1}$ is of Class~I and therefore
contains infinitely many orbits connecting with $O$;

\item[(S2$'$)]
when $N({\cal H})$ is odd,
we have $\dot\theta|_{\Upsilon_{N({\cal H})}}\dot\theta|_{\Upsilon_{N({\cal H})+1}}>0$ in \eqref{equ:ttkt},
which implies that
no orbits connect with $O$ in ${\cal U}_{N({\cal H})+1}$
as we saw in ${\cal U}_1$ above.
\end{description}
Summarily, $\Upsilon_{1}$ contributes no e-tsectors,
but $\Upsilon_{N({\cal H})}$ either contributes exact one (or no) e-tsector in the case {\bf (S1$'$)} (or {\bf (S2$'$)})
for the same reasons as given in last paragraph.
Consequently, in the direction $\theta=0$ the system has
totally $N({\cal H})-1$ (or $N({\cal H})-2$) e-tsectors in the case {\bf (S1$'$)} (or {\bf (S2$'$)}).
Moreover, there are no other types of tsectors.

At last, we prove that the number of real branches
of the equation ${\cal G}(\rho,\theta)=0$ is equal to $N({\cal H})$ and $N({\cal H})-1$ in the cases
{\bf (S1$'$)} and {\bf (S2$'$)} respectively,
showing that the number of e-tsectors is the number of real branches minus 1.
In fact,
in the case {\bf (S1$'$)}
we have
$\dot\theta|_{\Upsilon_0}\dot\theta|_{\Upsilon_1}>0$
and
$\dot\theta|_{\Upsilon_k}\dot\theta|_{\Upsilon_{k+1}}<0$ for all $k=1,...,N({\cal H})$
by \eqref{equ:ttkt},
which implies by {\bf(F2)} and {\bf(F3)} that
the equation ${\cal G}(\rho,\theta)=0$ has $N({\cal H})$ real ranches;
in the case {\bf (S2$'$)}
we have
$\dot\theta|_{\Upsilon_0}\dot\theta|_{\Upsilon_1}>0$,
$\dot\theta|_{\Upsilon_k}\dot\theta|_{\Upsilon_{k+1}}<0$ for all $k=1,2,...,N({\cal H})-1$,
and $\dot\theta|_{\Upsilon_{N({\cal H})}}\dot\theta|_{\Upsilon_{N({\cal H})+1}}>0$
by \eqref{equ:ttkt},
which implies by {\bf(F2)} and {\bf(F3)} that
the equation ${\cal G}(\rho,\theta)=0$ has $N({\cal H})-1$ real branches.
Thus, {\bf Claim~2.3} is proved and the proof in {\bf Step~2} is completed.

{\bf Step 3:} {\it Give distribution of orbits in the exceptional direction $\theta=0$
in the subcase {\bf(iia)}, i.e.,
$G_0(\theta)\equiv 0$ and $C_0<0$.}

We still use the same procedure as in {\bf Step~1} to discuss the subcase {\bf(iia)}
but a difference from {\bf Step~1}, where the subcase {\bf (ia)} is discussed,
is that the exceptional direction $\theta=0$ is {\it not isolated} since $G_0(\theta)\equiv 0$.
In this situation, the two equations in (\ref{equ:polar system}) have a common factor $\rho$,
which determines a line of singular points.
Applying the time-rescaling $dt=\rho d\tau$,
we reduce system~\eqref{equ:polar system} to the following
\begin{equation}
\dot \rho={\cal H}(\rho,\theta)=H_0(\theta)+O(\rho),~~~
\dot \theta={\cal G}(\rho,\theta)/\rho=G_1(\theta)+O(\rho),
\label{equ:polar system2}
\end{equation}
where $\tau$ is replaced with $t$ for simplification and we still use
$\dot \rho$ and $\dot \theta$ to denote $d\rho/dt$ and $d\theta/dt$ respectively.
Clearly,
system~\eqref{equ:polar system2} has no lines of singular points
and the point $(\rho,\theta)=(0,0)$ is either a regular point if $G_1(0)\ne 0$
or an isolated singular point if $G_1(0)=0$ since
$$
\rho^{2m}H_0^2(\theta)=\rho^{2m}(G_0^2(\theta)+H_0^2(\theta))
=X_m^2(x,y)+Y_m^2(x,y)\not\equiv 0
$$
as assumed just below \eqref{equ:initial}.
So it is more convenient to investigate system~\eqref{equ:polar system2}
near $(\rho,\theta)=(0,0)$ in ${\cal U}_{\epsilon,\delta}$
because
system~\eqref{equ:initial} has an orbit connecting with $O$
in the direction $\theta=0$ if and only if
system~\eqref{equ:polar system2} has an orbit
connecting with the point $(\rho,\theta)=(0,0)$
in ${\cal U}_{\epsilon,\delta}$.

Correspondingly to {\bf Claims 1.1}-{\bf 1.3} in {\bf Step~1},
we have the following:

\noindent{\bf Claim~3.1.}
{\it
Facts {\bf (F1)}-{\bf (F3)} in {\bf Claim~1.1} still hold.
}

\noindent{\bf Claim~3.2.}
{\it
For each $k=1,...,N({\cal H})+1$
the Z-sector ${\cal U}_k$
contains a unique orbit connecting with $O$ in the direction $\theta=0$.
}

\noindent{\bf Claim~3.3.}
{\it
The number of orbits of system~\eqref{equ:initial} connecting with $O$
in the direction $\theta=0$ is equal to
the number of real branches of the equation ${\cal H}(\rho,\theta)=0$
on the half-plane $\rho\ge 0$ plus 1.
}

The proof of {\bf Claim~3.3} needs {\bf Claim~3.1} and {\bf Claim~3.2}.
If {\bf Claim~3.3} is true,
then by Lemma \ref{Assertion1} there are $N({\cal H})$ h-tsectors.
Moreover,
since $H_0(0)=0$ as assumed in \eqref{G0H00} and
$H_0(\theta)\not\equiv 0$ as indicated just below \eqref{equ:polar system2},
we have $H_0(\theta)\ne 0$ for all $\theta\in(-\delta,0)\cup(0,\delta)$.
By Theorem 3.3 of \cite[p.63]{ZZF},
which says that in each exceptional direction $\theta=\theta_*$ such that $H_0(\theta_*)\ne 0$
there is exactly one orbit connecting with $O$,
in each direction $\theta=\vartheta\in(-\delta,0)\cup(0,\delta)$
system~\eqref{equ:initial} has exactly one orbit connecting with $O$.
It follows from the definition of ${\cal U}_1$ and ${\cal U}_{N({\cal H})+1}$,
given just below \eqref{thetatheta},
that
both of them contain infinitely many orbits of system~\eqref{equ:initial} connecting with $O$.
Then ${\cal S}_O^p=2$ and therefore, conclusions of {\bf(iia)} are obtained.

In what follows, we need to prove {\bf Claims~3.1}-{\bf 3.3}.
It is the same as in the subcase {\bf(ia)} that
$[{\cal G},{\cal H}]_\theta\le 0$ on ${\cal U}_{\epsilon,\delta}$ in the subcase {\bf(iia)}.
Then {\bf Claim 3.1} holds for the same reason as {\bf Claim 1.1}.

\begin{figure}[!h]
    \centering
     \subcaptionbox{%
     }{\includegraphics[height=1.6in]{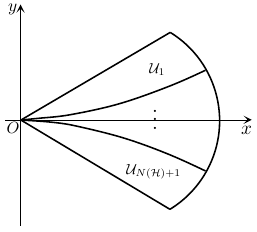}}~~~~~~~~~~~~~~~~
     \subcaptionbox{%
     }{\includegraphics[height=1.6in]{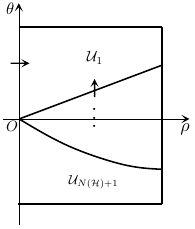}}
    \caption{Subregions on $(x,y)$-plane and $(\rho,\theta)$-plane.}
    \label{fig:Sxyrt}
\end{figure}

Next, we prove {\bf Claim 3.2}.
As in {\bf Claim 1.2},
this claim also needs to be discussed in two parts, i.e., for $k=2,...,N({\cal H})$
and for $k=1, N({\cal H})+1$ separately.
The first part can be proved by the same arguments given for {\bf Claim~1.2}.
However, the second part is quite different from that of {\bf Claim 1.2}.
The Z-sector ${\cal U}_1$ (or ${\cal U}_{N({\cal H})+1}$) in {\bf Claim 1.2} contains
the exceptional direction $\theta=0$ only,
but it in {\bf Claim 3.2} contains
infinitely many exceptional directions
because the exceptional directions are non-isolated.
For the second part,
the `at most' result for $k=2,...,N({\cal H})$
in the discussion on the first part
is also true for $k=1$ and $N({\cal H})+1$ for the same reason.
Further,
we claim that ${\cal U}_1$ (and ${\cal U}_{N({\cal H})+1}$)
contains exactly one orbit connecting with $O$.
Actually,
we consider system~\eqref{equ:polar system2} in the Cartesian coordinates and
the sector ${\cal U}_1$ on the $(x,y)$-plane
becomes a trapezoid on the $(\rho,\theta)$-plane,
as shown in Figure~\ref{fig:Sxyrt}.
Let $\Upsilon_0^*$ be the left-side vertical edge of ${\cal U}_1$, i.e.,
$$
\Upsilon_0^*:=\{(\rho,\theta)\in\mathbb{R}^2:\rho=0,0<\theta<\delta\}.
$$
Then $\dot \rho|_{\Upsilon_0^*}\ne 0$ for all $\theta\in(0,\delta)$ because
\begin{eqnarray}
\dot \rho|_{\Upsilon_0^*}=\dot \rho(0,\theta)={\cal H}(0,\theta)=H_0(\theta)\ne 0
~~~\mbox{for all}~~~\theta\in(0,\delta),
\label{dotrhoU*}
\end{eqnarray}
where we note that $H_0(\theta)\not\equiv 0$ when $G_0(\theta)\equiv 0$,
as indicated just below \eqref{equ:polar system2},
and that $H_0(0)=0$, as assumed in \eqref{G0H00}.
There are 2 situations:
{\bf(3a)} $\dot\rho|_{{\cal U}_1}>0$ and
{\bf(3b)} $\dot\rho|_{{\cal U}_1}<0$.
We only consider {\bf(3a)} since {\bf(3b)} can be discussed similarly.
The equality $\dot\rho|_{{\cal U}_1}>0$ in {\bf(3a)} implies by continuity that
$\dot \rho|_{\Upsilon_0^*}>0$.
Then all orbits of system \eqref{equ:polar system2}
starting from the left-side vertical edge $\Upsilon_0^*$ of ${\cal U}_1$ enter ${\cal U}_1$.
Additionally,
on the lower edge $\Upsilon_1$ of ${\cal U}_1$,
limit \eqref{limit-infty} still holds,
which implies that $\dot\theta|_{\Upsilon_1}=({\cal G}/\rho)|_{\Upsilon_1}>0$
since ${\cal H}|_{{\cal U}_1}=\dot\rho|_{{\cal U}_1}>0$ in {\bf(3a)}.
Then all orbits starting from $\Upsilon_1$ enter ${\cal U}_1$.
Note that
the orbit starting from each point $M\in \Upsilon_0^*$ intersects with
the compact part ${\cal B}_1$ of the boundary $\partial {\cal U}_1$
at a point $P_M$, where
$$
{\cal B}_1:=\partial {\cal U}_1\setminus\{\Upsilon_0^*\cup\Upsilon_1\cup O\}.
$$
By the uniqueness of orbits,
there is one point $P\in {\cal B}_1$ such that $\lim_{M\to O}P_M=P$.
Similarly,
there is one point $Q\in {\cal B}_1$ corresponding to $\Upsilon_1$.
The uniqueness of orbits ensures that
orbits starting from the segment on ${\cal B}_1$ between $P$ and $Q$
all negatively connect with $O$.
By the `at most' result mentioned just before our claim,
$P$ coincides with $Q$ and therefore exact one orbit connects with $O$ in ${\cal U}_1$.
Similarly, ${\cal U}_{N({\cal H})+1}$ has the same result.
Then, our claim given before \eqref{dotrhoU*} is proved,
which implies that
each of ${\cal U}_1$ and ${\cal U}_{N({\cal H})+1}$ contains
exactly one orbit of system~\eqref{equ:initial}
connecting with $O$ in the direction $\theta=0$.
Thus, {\bf Claim~3.2} is proved.

Finally, we prove {\bf Claim 3.3}.
By {\bf Claim~3.2},
there are $N({\cal H})+1$ orbits connecting with $O$ in the direction $\theta=0$.
Note that the equation ${\cal H}(\rho,\theta)=0$ has $N({\cal H})$ real branches
by Lemma~\ref{Assertion1}.
Then {\bf Claim~3.3} is proved and {\bf Step~3} is completed.

{\bf Step 4:} {\it Give distribution of orbits in the exceptional direction $\theta=0$
in the subcase {\bf(iib)}, i.e.,
$G_0(\theta)\equiv 0$ and $C_0<0$.}

We still use the same procedure as in {\bf Step~2}.
Correspondingly to {\bf Claims~2.1}-{\bf 2.3} in {\bf Step 2}, we have the following.

\noindent{\bf Claim~4.1.}
{\it
Facts {\bf (F1$'$)}, {\bf (F2)} and {\bf (F3)} hold.
}

\noindent{\bf Claim~4.2.}
{\it
For each $k=1,...,N({\cal H})+1$ the Z-sector ${\cal U}_k$
contains infinitely many orbits connecting with $O$.
}

\noindent{\bf Claim~4.3.}
{\it
System~\eqref{equ:initial} has one orbit connecting with $O$ in the direction $\theta=0$
if $N({\cal H})=0$, and there is one p-tsector and no other tsectors in this direction.
System~\eqref{equ:initial} has no orbits connecting with $O$ in the direction $\theta=0$
if $N({\cal H})=1$ and $\tilde{j}_{s({\cal H})}>0$,
and infinitely many orbits if
either $N({\cal H})=1$ and $\tilde{j}_{s({\cal H})}=0$ or $N({\cal H})> 1$.
Moreover, there are no h-tsectors and p-tsectors,
and the number of e-tsectors is equal to
the number of real branches of the equation ${\cal H}(\rho,\theta)=0$
on the half-plane $\rho\ge 0$ plus 1 if $N({\cal H})\ge 1$.
}

The proof of {\bf Claim~4.3} needs {\bf Claim~4.1} and {\bf Claim~4.2}.
If {\bf Claim~4.3} is true,
then conclusions in {\bf (iib)} are obvious by Lemma~\ref{Assertion1}.

It is the same as in the subcase {\bf (ib)} that
$[{\cal G},{\cal H}]_\theta\ge 0$ on ${\cal U}_{\epsilon,\delta}$ in the subcase {\bf (iib)}.
Then {\bf Claim 4.1} holds for the same reason as {\bf Claim 2.1}.

Next we prove {\bf Claim~4.2}.
This claim also needs to be discussed in two parts, i.e., for $k=2,...,N({\cal H})$
and for $k=1, N({\cal H})+1$ separately, as in {\bf Claim 2.2}.
The first part can be obtained in the same way as the first part of {\bf Claim~2.2}.
However,
the second part is quite different from that of {\bf Claim 2.2}.
The Z-sectors ${\cal U}_1$ (and ${\cal U}_{N({\cal H})+1}$) in {\bf Claim 2.2}
contains exactly one exceptional direction $\theta=0$,
but it in {\bf Claim 3.2} contains infinitely many ones
because the exceptional directions are non-isolated.
For the second part,
using the same argument in the paragraph given just below {\bf Claim~3.3},
we can obtain that
both ${\cal U}_1$ and ${\cal U}_{N({\cal H})+1}$ contain infinitely many orbits of system~\eqref{equ:initial} connecting with $O$.
Then {\bf Claim~4.2} is proved.

Finally, we prove {\bf Claim~4.3}.
If $N({\cal H})> 1$ then the two edges $\Upsilon_1$ and $\Upsilon_2$ of ${\cal U}_2$
are both tangent to the positive $x$-axis at $O$.
By {\bf Claim~4.2},
system~\eqref{equ:initial} has infinity many orbits connecting with $O$
in the direction $\theta=0$.
In order to determine the number of e-tsectors formed by those orbits connecting with $O$,
we also need to discuss in the same four circumstance {\bf(C1)}-{\bf(C4)}
as listed in the proof of {\bf Claim~1.3}.
We still only give the proof in {\bf (C1)} and
assume without loss of generality that
$a_{i_0,j_0}({\cal G})<0$ and $a_{\tilde{i}_0,\tilde{j}_0}({\cal H})>0$.
Then \eqref{equ:ttkt} and \eqref{equ:rhosignib} still hold.
Similarly to \eqref{phitpk} in the proof of {\bf Claim~2.3}, we see that
there are $N({\cal H})-2$ e-tsectors
each of which passes through the curve $\Upsilon_k$
for one $k\in\{2,...,N({\cal H})-1\}$.
Further,
we need to know whether there are homoclinic orbits passing through
$\Upsilon_1$ and $\Upsilon_{N({\cal H})}$.
By {\bf Claim~4.2},
we can choose an orbit $\Gamma$ connecting with $O$ in ${\cal U}_1$.
Moreover,
we can choose sufficiently small $\epsilon$ such that
the intersection point of the orbit $\Gamma$ and the boundary $\partial {\cal U}_1$
lies on the arc $\rho=\epsilon$.
Then, $\Gamma$ divides the region ${\cal U}_1$ into two connected subregions,
one of which between $\Gamma$ and $\Upsilon_1$ is denoted by ${\cal U}'_1$.
By \eqref{equ:ttkt} and \eqref{equ:rhosignib} we have $\dot \theta|_{\Upsilon_1}<0$ and $\dot\rho|_{{\cal U}'_1}>0$ respectively,
which implies that the orbit $\phi(t,P_1)$, passing through a point $P_1\in \Upsilon_1$,
satisfies that
\begin{align}
\phi(t,P_1)\to O~\mbox{in}~{\cal U}'_1~~~\mbox{as}~ t\to -\infty.
\label{phitp1}
\end{align}
On the other hand,
by \eqref{equ:ttkt} and \eqref{equ:rhosignib} we also have
inequalities $\dot \theta|_{\Upsilon_1}<0$ and $\dot\rho|_{{\cal U}_2}<0$,
which ensures that
\begin{align}
\phi(t,P_1)\to O~\mbox{in}~{\cal U}_2~~~\mbox{as}~ t\to +\infty
\label{phitp2}
\end{align}
since the Z-sector ${\cal U}_2$ is of Class I.
Then, by (\ref{phitp1}) and (\ref{phitp2})
there is one elliptic sector passing through $\Upsilon_1$.
Similarly,
there is also one elliptic sector passing through $\Upsilon_{N({\cal H})}$.
Consequently,
system~\eqref{equ:initial} has totally $N({\cal H})$ e-tsectors.
Note that the equation ${\cal H}(\rho,\theta)=0$ has $N({\cal H})$ real branches
by Lemma~\ref{Assertion1}.
Then {\bf Claim~4.3} holds in the situation $N({\cal H})>1$.
Moreover, there are no other types of tsectors.


In what follows,
we consider the remaining situations $N({\cal H})=0$ and $N({\cal H})=1$ separately.
In the situation $N({\cal H})=0$,
we only need to consider {\bf (C4)}, as listed in the proof of {\bf Claim~1.3}.
In fact,
the equation ${\cal H}(\rho,\theta)=0$ has no real branches by Lemma~\ref{Assertion1},
which implies that
\begin{eqnarray}
\tilde{j}_{s({\cal H})}=0~~~\mbox{and}~~~j_{s({\cal G})}>0
\label{jsH}
\end{eqnarray}
in case {\bf(J1)}. Furthermore, we claim that
\begin{description}
\item{\bf Claim H0:}
Ordinates $\tilde{j}_0,...,\tilde{j}_{s({\cal H})}$
of vertices of ${\cal N}_{\cal H}$ are all even
and
coefficients
$a_{\tilde{i}_0,\tilde{j}_0}({\cal H}),...,
a_{\tilde{i}_{s({\cal H})},\tilde{j}_{s({\cal H})}}({\cal H})$ corresponding to those vertices
have the same sign.

\item{\bf Claim G0:}
Ordinates $j_0,...,j_{s({\cal G})}$
of vertices of ${\cal N}_{\cal G}$ are all odd
and
coefficients
$a_{i_0,j_0}({\cal G}),...,a_{i_{s({\cal G})},j_{s({\cal G})}}({\cal G})$ corresponding to those vertices
have the same sign.
\end{description}
For an indirect proof to {\bf Claim H0},
assume that either two coefficients have different signs or
one of ordinates, denoted by $\tilde{j}_{k}$, is odd.
If the first one in `assume' is true,
then the number of sign changes in the vertex coefficient sequence
$
{\cal A}({\cal H})=
(a_{\tilde{i}_0,\tilde{j}_0}({\cal H}),...,
a_{\tilde{i}_{s({\cal H})},\tilde{j}_{s({\cal H})}}({\cal H})),
$
denoted by $\Xi({\cal A}({\cal H}))$ just before Corollary~\ref{cor:nonpara},
is greater than 0,
implying by Proposition~\ref{lm-solution by NP} that there is a real branch,
a contradiction.
If the first one is denied but the second one is true in `assume',
then $(-1)^{\tilde{j}_k}<0$ and $(-1)^{\tilde{j}_{s({\cal H})}}>0$ by \eqref{jsH}
and therefore the $j$-algebraic vertex coefficient sequence
$
{\cal A}_j({\cal H})=
((-1)^{\tilde{j}_{s({\cal H})}}
a_{\tilde{i}_{s({\cal H})},\tilde{j}_{s({\cal H})}}({\cal H}),...,$
$(-1)^{\tilde{j}_0}a_{\tilde{i}_0,\tilde{j}_0}({\cal H}))
$
satisfies that
$\Xi({\cal A}_j({\cal H}))\ge 1$,
implying by the remark given just below the proof of Proposition~\ref{lm-solution by NP} that there is a real branch,
the same contradiction.
Thus, {\bf Claim H0} is proved.
By {\bf Claim H0} and {\bf(P1)},
$j_0,...,j_{s({\cal G})}$ are all odd.
Moreover,
{\bf(F3)} of {\bf Claim~4.1} ensures that
$\theta=0$ is the only real branch of the equation ${\cal G}(\rho,\theta)=0$
because $j_{s({\cal G})}>0$ and
the equation ${\cal H}(\rho,\theta)=0$ has no real branches,
as indicated just before \eqref{jsH}.
Similar to {\bf Claim~H0},
{\bf Claim G0} follows from Proposition~\ref{lm-solution by NP}.
By {\bf Claim G0},
$j_0$ is odd, which implies that
neither {\bf(C1)} nor {\bf(C3)} is possible.
In {\bf(C2)},
similarly to \eqref{Gi0j0},
by {\bf Claim H0} and {\bf Claim G0} we have
\begin{align*}
\frac{{\cal G}(\rho,\delta)}{{\cal H}(\rho,\delta)}
&=\frac{\rho^{i_0}\{a_{i_0,j_0}({\cal G})\delta^{j_0}+O(\delta^{j_0+1})+O(\rho)\}}
{\rho^{\tilde{i}_0}\{a_{\tilde{i}_0,\tilde{j}_0}({\cal H})\delta^{\tilde{j}_0}
+O(\delta^{\tilde{j}_0+1})+O(\rho)\}}
<0,
\\
\frac{{\cal G}(\rho,-\delta)}{{\cal H}(\rho,-\delta)}
&=\frac{\rho^{i_0}\{a_{i_0,j_0}({\cal G})(-\delta)^{j_0}+O(\delta^{j_0+1})+O(\rho)\}}
{\rho^{\tilde{i}_0}\{a_{\tilde{i}_0,\tilde{j}_0}({\cal H})(-\delta)^{\tilde{j}_0}
+O(\delta^{\tilde{j}_0+1})+O(\rho)\}}
>0,
\end{align*}
a contradiction to the fact {\bf(F1$'$)} of {\bf Claim~4.1}.
Consequently, only {\bf(C4)} is possible.
In what follows,
we assume without loss of generality that
\begin{eqnarray}
a_{i_0,j_0}({\cal G})<0~~\mbox{ and }~~a_{\tilde{i}_0,\tilde{j}_0}({\cal H})<0
\label{assume-c4}
\end{eqnarray}
in {\bf(C4)};
otherwise, we can discussed similarly.
Note that $\theta=0$ is an orbit
since $\theta=0$ is a real branch of the equation ${\cal G}(\rho,\theta)=0$.
Then system~\eqref{equ:polar system2} has at least one orbit
connecting with $O$ in the region ${\cal U}_{\epsilon,\delta}$,
defined in the beginning of the proof of Theorem~\ref{th:finite}.
Next, we determine the exact number of orbits
connecting with $O$ in the region ${\cal U}_{\epsilon,\delta}$, sbut
the criteria given in \cite[p.220, Theorem~5]{SC} for Z-sectors of Class I cannot be used
because our exceptional direction $\theta=0$ satisfying \eqref{G0H00} is irregular.
Moreover,
the criterion given in \cite[Lemma~5]{T-Z1} for general sectors
is also invalid because
$$
\frac{\partial }{\partial \theta} \frac{{\cal G}(\rho,\theta)}{{\cal H}(\rho,\theta)}
=\frac{[{\cal G},{\cal H}]_\theta}{\rho{\cal H}^2}
\ge 0~~~\mbox{for all}~(\rho,\theta)\in{\cal U}_{\epsilon,\delta},
$$
where we note that $[{\cal G},{\cal H}]_\theta \ge 0$ as shown just below {\bf Claim~4.3}.
So we turn to claim that
\begin{align}
\frac{\partial }{\partial \theta} \frac{{\cal G}(\rho,\theta)}{{\cal H}(\rho,\theta)}
=\frac{[{\cal G},{\cal H}]_\theta}{{\cal H}^2}\le \rho^{1/n}
~~~\mbox{for all}~(\rho,\theta)\in{\cal U}_{\epsilon,\delta},
\label{GHH2U}
\end{align}
where $n$ ($>0$) is a large integer.
If it is true, let
$\vartheta_1(\rho)$ and $\vartheta_2(\rho)$ be two orbits
of system~\eqref{equ:polar system2} connecting with $O$.
Assume without loss of generality that
$0<\vartheta_1(\rho)<\vartheta_2(\rho)<\delta$ for all $\rho\in(0,\epsilon)$.
By the Mean Value Theorem and \eqref{GHH2U},
\begin{align*}
\frac{d(\vartheta_2(\rho)-\vartheta_1(\rho))}{d\rho}
=\frac{{\cal G}(\rho,\vartheta_2(\rho))}{\rho{\cal H}(\rho,\vartheta_2(\rho))}
-\frac{{\cal G}(\rho,\vartheta_1(\rho))}{\rho{\cal H}(\rho,\vartheta_1(\rho))}
\le \frac{1}{\rho^{(n-1)/n}}(\vartheta_2(\rho)-\vartheta_1(\rho)).
\end{align*}
Integrating both sides from $\rho$ to $\rho_0$ with $0<\rho<\rho_0$,
we obtain that
\begin{align}
\ln (\vartheta_2(\rho_0)-\vartheta_1(\rho_0))
-\ln (\vartheta_2(\rho)-\vartheta_1(\rho))
\le (\rho_0^{1/n}-\rho^{1/n})/n,
\label{wGHwp}
\end{align}
the left hand side of which tends to $+\infty$ as $\rho\to 0^+$.
This contradiction shows that
$\theta=0$ is the only orbit of system~\eqref{equ:polar system2}
connecting with $O$ in ${\cal U}_{\epsilon,\delta}$, i.e.,
system~\eqref{equ:initial} has one orbit connecting with $O$ in the direction $\theta=0$.
Moreover,
using the same argument in the paragraph just below {\bf Claim~3.3},
in each direction $\theta=\theta_*\in(-\delta,0)\cup(0,\delta)$
system~\eqref{equ:initial} has one orbit connecting with $O$.
Then there is one p-tsector and no other tsectors in the direction $\theta=0$.

Now we are going to prove (\ref{GHH2U}), which is equivalent to the inequality
\begin{align}
{\cal H}^2\ge \rho^{-1/n}[{\cal G},{\cal H}]_\theta
\label{H>[]}
\end{align}
on ${\cal U}_{\epsilon,\delta}$.
We simply let $f$
denote $[{\cal G},{\cal H}]_\theta$
and claim that
\begin{align}
{\cal N}_{{\cal H}^2}<{\cal N}_f,
\label{NH2Ng}
\end{align}
i.e.,
${\cal N}_f$ strictly lies above ${\cal N}_{{\cal H}^2}$,
defined just below \eqref{NfleNtf}.
In fact,
by Remark~\ref{rmk:mul}{\bf(ii)},
the left-most and the right-most vertices of ${\cal N}_{{\cal H}^2}$ are points
\begin{eqnarray}
(2\tilde{i}_0,2\tilde{j}_0)
~~~\mbox{and}~~~
(2\tilde{i}_{s({\cal H})},2\tilde{j}_{s({\cal H})})
\label{ijH}
\end{eqnarray}
respectively.
On the other hand, by Lemma~\ref{J1Ng} we see from \eqref{lmJ1NG} and \eqref{jsH} that
$$
{\cal N}_f={\cal N}_{{\cal G}'_\theta{\cal H}}.
$$
By Remarks~\ref{rmk:diff} and \ref{rmk:mul},
the left-most and the right-most vertices of ${\cal N}_f$ are points
\begin{align}
(i_0+\tilde{i}_0,j_0+\tilde{j_0}-1)
~~\mbox{ and }~~
(i_{s({\cal G})}+\tilde{i}_{s({\cal H})},j_{s({\cal G})}+\tilde{j}_{s({\cal H})}-1)
\label{iijj1ij1}
\end{align}
respectively.
By the definition of ``strictly lies above'' given just below \eqref{NfleNtf},
the claimed \eqref{NH2Ng} is equivalent to the statement that
either
${\cal N}_{{\cal H}^2}^{\infty}(u)\le {\cal N}_f^{\infty}(u)$
for all $u\in[i_0+\tilde{i}_0,+\infty)$ when $i_0+\tilde{i}_0>2\tilde{i}_{s({\cal H})}$ or
\begin{equation}
{\cal N}_{{\cal H}^2}^{\infty}(u)
\left\{
\begin{array}{lllll}
<{\cal N}_f^{\infty}(u)    &\mbox{for all} ~~~u\in[i_0+\tilde{i}_0,2\tilde{i}_{s({\cal H})}],
\\
\le {\cal N}_f^{\infty}(u) &\mbox{for all} ~~~u\in(2\tilde{i}_{s({\cal H})},+\infty)
\end{array}
\right.
\label{nnn}
\end{equation}
when $i_0+\tilde{i}_0\le 2\tilde{i}_{s({\cal H})}$.
Thus,
in the case $i_0+\tilde{i}_0> 2\tilde{i}_{s({\cal H})}$,
the ordinates of the second points in \eqref{ijH} and \eqref{iijj1ij1} imply that
$$
{\cal N}_{{\cal H}^2}^\infty(u)=2\tilde{j}_{s({\cal H})}=0
\le
j_{s({\cal G})}+\tilde{j}_{s({\cal H})}-1\le {\cal N}^\infty_f(u)
~~~\forall u\in[i_0+\tilde{i}_0,+\infty)
$$
because of \eqref{jsH}. This shows that the claimed \eqref{NH2Ng} is true.
In the opposite case $i_0+\tilde{i}_0\le 2\tilde{i}_{s({\cal H})}$,
we prove \eqref{nnn} indirectly.
First, we consider the first points in (\ref{ijH}) and (\ref{iijj1ij1}).
Because of (\ref{assume-c4}),
the inequality
$C_0=(j_0-\tilde{j}_0)a_{i_0,j_0}({\cal G})a_{\tilde{i}_0,\tilde{j}_0}({\cal H})>0$
given in {\bf (iib)} of Theorem~\ref{th:finite}
implies that $j_0>\tilde{j_0}$ and therefore,
\begin{align}
j_0+\tilde{j}_0-1\ge 2\tilde{j}_0.
\label{jj12j}
\end{align}
On the other hand,
$i_0\ge 1$
but
$\tilde{i}_0=0$
because $H_0(\theta)\not\equiv 0$ when $G_0(\theta)\equiv 0$
as indicated below \eqref{equ:polar system2}.
It follows that the left-most vertex of ${\cal N}_f$ strictly
lies above the polygon
${\cal N}^\infty_{{\cal H}^2}(u)$ for $u\in[2\tilde{i}_0,+\infty)$.
Thus, by the assumption of the indirect proof
we see that
the edge $E_{k'}({\cal H}^2)$
linking vertices $V_{k'-1}({\cal H}^2):(\hat{i}_{k'-1},\hat{j}_{k'-1})$ and
$V_{k'}({\cal H}^2):(\hat{i}_{k'},\hat{j}_{k'})$
and
the edge $E_{k''}(f)$
linking vertices $V_{k''-1}(f):(\breve{i}_{k''-1},\breve{j}_{k''-1})$ and
$V_{k''}(f):(\breve{i}_{k''},\breve{j}_{k''})$
intersect at a point, denoted by $(i^*,j^*)$, such that
\begin{align}
\hat{i}_{k'-1}< i^* \le \hat{i}_{k''},
~~~
\hat{j}_{k'-1}>\breve{j}_{k''}
~~~\mbox{and}~~~
\zeta(E_{k'}({\cal H}^2))>\zeta(E_{k''}(f)).
\label{iiijjEE}
\end{align}
Choose
$\nu'\in(\zeta(E_{k'}({\cal H}^2))-\varepsilon,\zeta(E_{k'}({\cal H}^2)))$ and
$\nu''\in(\zeta(E_{k''}(f)),\zeta(E_{k''}(f))+\varepsilon)$
for sufficiently small $\varepsilon>0$.
Clearly, $\nu'>\nu''$ because of the third inequality in \eqref{iiijjEE}.
By definition \eqref{defscomp},
the $\nu'$-component of ${\cal N}_{{\cal H}^2}$ is the vertex $V_{k'-1}({\cal H}^2)$ and
the $\nu''$-component of ${\cal N}_f$ is the vertex $V_{k''}(f)$.
By definition \eqref{defkxip},
the corresponding componential polynomials are the following
\begin{align*}
{\cal K}_{\nu'}({\cal H}^2)(\theta)
=a_{\hat{i}_{k'-1},\hat{j}_{k'-1}}({\cal H}^2)\theta^{\hat{j}_{k'-1}},
~~~
{\cal K}_{\nu''}({\cal G}'_\theta{\cal H})(\theta)
=a_{\breve{i}_{k''},\breve{j}_{k''}}({\cal G}'_\theta{\cal H})\theta^{\breve{j}_{k''}},
\end{align*}
where we note that ${\cal N}_f={\cal N}_{{\cal G}'_\theta{\cal H}}$,
as indicated just before \eqref{iijj1ij1}.
By Remark~\ref{rmk:mul},
${\mathfrak S}({\cal H}^2)={\mathfrak S}({\cal H})$ and
${\mathfrak S}(f)={\mathfrak S}({\cal G})\cup{\mathfrak S}({\cal H})$,
which implies that
$\nu',\nu''\notin {\mathfrak S}({\cal G})\cup{\mathfrak S}({\cal H})$.
Then the $\nu'$- and $\nu''$-components of ${\cal N}_{\cal H}$ are both vertices,
denoted by $V_{\varrho'}({\cal H}):(\tilde{i}_{\varrho'},\tilde{j}_{\varrho'})$ and
$V_{\varrho''}({\cal H}):(\tilde{i}_{\varrho''},\tilde{j}_{\varrho''})$ respectively,
and the $\nu''$-component of ${\cal N}_{\cal G}$ is also a vertex,
denoted by $V_{\sigma}({\cal G}):(i_\sigma,j_\sigma)$.
By definition \eqref{defkxip},
they also correspond to the componential polynomials
{\small
$$
{\cal K}_{\nu'}({\cal H})(\theta)
=a_{\tilde{i}_{\varrho'},\tilde{j}_{\varrho'}}({\cal H})\theta^{\tilde{j}_{\varrho'}},
~~
{\cal K}_{\nu''}({\cal H})(\theta)
=a_{\tilde{i}_{\varrho''},\tilde{j}_{\varrho''}}({\cal H})\theta^{\tilde{j}_{\varrho''}},
~~
{\cal K}_{\nu''}({\cal G})(\theta)
=a_{i_\sigma,j_\sigma}({\cal G})\theta^{j_\sigma}.
$$
}By Remark~\ref{rmk:diff} and Proposition~\ref{lm-MIL},
${\cal K}_{\nu'}({\cal H}^2)
={\cal K}_{\nu'}({\cal H}){\cal K}_{\nu'}({\cal H})$
and
${\cal K}_{\nu''}({\cal G}'_\theta{\cal H})
=({\cal K}_{\nu''}({\cal G}))'_\theta{\cal K}_{\nu''}({\cal H})$.
From the above expressions of polynomials, it follows that
\begin{align*}
a_{\hat{i}_{k'-1},\hat{j}_{k'-1}}({\cal H}^2)\theta^{\hat{j}_{k'-1}}
&=a_{\tilde{i}_{\varrho'},\tilde{j}_{\varrho'}}({\cal H})\theta^{\tilde{j}_{\varrho'}}
a_{\tilde{i}_{\varrho'},\tilde{j}_{\varrho'}}({\cal H})\theta^{\tilde{j}_{\varrho'}},
\\
a_{\breve{i}_{k''},\breve{j}_{k''}}({\cal G}'_\theta{\cal H})\theta^{\breve{j}_{k''}}
&=j_\sigma a_{i_\sigma,j_\sigma}({\cal G})\theta^{j_\sigma-1}
a_{\tilde{i}_{\varrho''},\tilde{j}_{\varrho''}}({\cal H})\theta^{\tilde{j}_{\varrho''}}.
\end{align*}
Thus, we see from the degrees of $\theta$ that
\begin{align}
\hat{j}_{k'-1}=\tilde{j}_{\varrho'}+\tilde{j}_{\varrho'}
~~~\mbox{and}~~~
\breve{j}_{k''}=j_\sigma-1+\tilde{j}_{\varrho''}.
\label{j6j6j6}
\end{align}
On the other hand,
${\cal K}_{\nu''}(f)
=[{\cal K}_{\nu''}({\cal G}),{\cal K}_{\nu''}({\cal H})]$
as indicated in Lemma~\ref{J1Ng}.
By the above expressions of polynomials, it follows that
{\small
$$
a_{\breve{i}_{k''},\breve{j}_{k''}}(f)\theta^{\breve{j}_{k''}}
=(j_\sigma-\tilde{j}_{\varrho''})
a_{i_\sigma,j_\sigma}({\cal G})
a_{\tilde{i}_{\varrho''},\tilde{j}_{\varrho''}}({\cal H})
\theta^{j_\sigma-1+\tilde{j}_{\varrho''}},
$$
}which implies from their coefficients that
\begin{eqnarray}
a_{\breve{i}_{k''},\breve{j}_{k''}}(f)
=(j_\sigma-\tilde{j}_{\varrho''})
a_{i_\sigma,j_\sigma}({\cal G})
a_{\tilde{i}_{\varrho''},\tilde{j}_{\varrho''}}({\cal H}).
\label{coef-a}
\end{eqnarray}
Note that $a_{\breve{i}_{k''},\breve{j}_{k''}}(f)$ is the coefficient of $f$
corresponding to the vertex $V_{k''}(f):(\breve{i}_{k''},\breve{j}_{k''})$.
According to Proposition~\ref{lm-lowest term}, the semi-positiveness of $f$ implies that
$
{\rm sgn}(a_{\breve{i}_{k''},\breve{j}_{k''}}(f))
={\rm sgn}(f|_\Gamma)
\ge 0,
$
where $\Gamma:\theta=\rho^{-1/\nu''},\rho\in(0,\epsilon)$, and
$\nu''$ is chosen just below \eqref{iiijjEE} such that
$\zeta^-(V_{k''}(f))=\zeta(E_{k''}(f))<\nu''<\zeta(E_{k''}(f))+\varepsilon=\zeta^+(V_{k''}(f))$
as defined in \eqref{SP+-}.
On the right hand side of \eqref{coef-a},
$
a_{i_\sigma,j_\sigma}({\cal G})
a_{\tilde{i}_{\varrho''},\tilde{j}_{\varrho''}}({\cal H})>0
$
by {\bf Claim H0} and {\bf Claim G0} because of (\ref{assume-c4}).
It follows from (\ref{coef-a}) that
$
j_\sigma>\tilde{j}_{\varrho''}.
$
Note that $\tilde{j}_{\varrho'}\le \tilde{j}_{\varrho''}$
because the inequality $\nu'>\nu''$ ensures that the vertex $V_{\varrho'}({\cal H}):(\tilde{i}_{\varrho'},\tilde{j}_{\varrho'})$ does not lie before
the vertex $V_{\varrho''}({\cal H}):(\tilde{i}_{\varrho''},\tilde{j}_{\varrho''})$
on $\vec{{\cal N}}_{{\cal H}}$.
Then we have $j_\sigma>\tilde{j}_{\varrho'}$ and therefore we see
from \eqref{j6j6j6} that $\breve{j}_{k''}\ge \hat{j}_{k'-1}$,
a contradiction to the second inequality in \eqref{iiijjEE}.
This shows that (\ref{nnn}) is true and the the claimed \eqref{NH2Ng} is proved.

Having \eqref{NH2Ng}, we return to prove (\ref{GHH2U}), i.e.,
the inequality ${\cal H}^2\ge \rho^{-1/n} f$ on ${\cal U}_{\epsilon,\delta}$.
In order to avoid the fractional power $-1/n$ in this inequality,
which prevents from reduction to Newton polygons,
we let
\begin{eqnarray}
\hat\rho:=\rho^{1/n}, \ \
\hat{\cal H}(\hat\rho,\theta):={\cal H}(\hat{\rho}^n,\theta),\ \
\hat f(\hat\rho,\theta):=f(\hat{\rho}^n,\theta)/\hat\rho.
\label{stretching-n}
\end{eqnarray}
Note that the valid index set $\Delta(\hat{\cal H})$
is obtained by stretching the set $\Delta({\cal H})$ horizontally by a multiplier $n$.
Similarly, $\Delta(\hat{f})$ is obtained by stretching $\Delta(f)$ horizontally by $n$ and then shifting 1 unit to the left.
From the definitions of the Newton polygons of $\hat{\cal H}$ and $\hat{f}$
we have
${\cal N}^\infty_{\hat{\cal H}^2}(u)={\cal N}^\infty_{{\cal H}^2}(nu)$
and
${\cal N}^\infty_{\hat{f}}(u)={\cal N}^\infty_{f}(nu+1)$
and, moreover,
the right-most vertex of ${\cal N}_{\hat{\cal H}^2}$
and the left-most vertex of ${\cal N}_{\hat{f}}$ are points
$(2n\tilde{i}_{s({\cal H})},2\tilde{j}_{s({\cal H})})$
and
$(ni_0+n\tilde{i}_0-1,j_0+\tilde{j}_0-1)$
respectively.
By \eqref{NH2Ng},
there is a large integer $n$ such that
either
${\cal N}_{{\cal H}^2}^{\infty}(u)=0\le {\cal N}_f^{\infty}(u+1/n)$
for all $u\in[i_0+\tilde{i}_0-1/n,+\infty)$
when $i_0+\tilde{i}_0-1/n>2\tilde{i}_{s({\cal H})}$ or
\begin{equation*}
{\cal N}_{{\cal H}^2}^{\infty}(u)
\left\{
\begin{array}{lllll}
<{\cal N}_f^{\infty}(u+1/n)
&\mbox{for all} ~~~u\in[i_0+\tilde{i}_0-1/n,2\tilde{i}_{s({\cal H})}],
\\
\le {\cal N}_f^{\infty}(u+1/n)
&\mbox{for all} ~~~u\in(2\tilde{i}_{s({\cal H})},+\infty)
\end{array}
\right.
\end{equation*}
when $i_0+\tilde{i}_0-1/n\le 2\tilde{i}_{s({\cal H})}$.
This implies that
\begin{align}
{\cal N}_{\hat{\cal H}^2}<{\cal N}_{\hat f},
\label{NH2Nf}
\end{align}
the same as (\ref{NH2Ng}) with $\hat{\cal H}$ and $\hat{f}$.
It follows by Corollary~\ref{cor-ADD} that
the Newton polygon of the function $\hat {\cal H}^2-\hat{f}$ and its edge-polynomials
are the same as those of $\hat {\cal H}^2$.
Although
edges of ${\cal N}_{\hat{\cal H}^2}$ are stretched from edges of ${\cal N}_{{\cal H}^2}$
as in (\ref{stretching-n}),
coefficients of their edge-polynomials with the corresponding valid indices are the same.
Thus $\hat {\cal H}^2$ and ${\cal H}^2$ have the same edge-polynomials
by definition \eqref{red}.
Note that
each edge-polynomial of ${\cal H}$ has no nonzero real roots;
otherwise,
the same argument given just below \eqref{KGKHK4} shows that
such a nonzero real root is of multiplicity 1 and
determines a real branch of the equation ${\cal H}(\rho,\theta)=0$ by Lemma~\ref{lm-basic},
a contradiction to the fact indicated just before \eqref{jsH} that
the equation ${\cal H}(\rho,\theta)=0$ has no real branches.
By Proposition~\ref{lm-MIL},
each edge-polynomial of ${\cal H}^2$ has no nonzero real roots.
It follows that each edge-polynomial of $\hat{\cal H}^2-\hat{f}$ has no nonzero real roots.
By Proposition~\ref{lm-region},
$\hat{\cal H}^2-\hat{f}\ge 0$
for all
$(\hat \rho,\theta)\in
\{(\hat \rho,\theta)\in\mathbb{R}^2:0<\hat\rho<\epsilon^{1/n},-\delta<\theta<\delta\}$
and therefore
(\ref{H>[]}) is true
for all $(\rho,\theta)\in{\cal U}_{\epsilon,\delta}$,
which implies (\ref{GHH2U}) as mentioned before (\ref{NH2Ng})
and completes the proof in the first situation $N({\cal H})=0$.

In the second situation $N({\cal H})=1$,
we only need to consider {\bf(C1)} and {\bf(C3)}.
In fact,
by Lemma~\ref{Assertion1'}
the equation ${\cal H}(\rho,\theta)=0$ has only one real branch,
i.e., $\Upsilon_1:\theta=\theta_1(\rho)$ in \eqref{thetatheta},
which is of odd multiplicity by {\bf(F2)} of {\bf Claim 4.1}.
Then
$$
{\cal H}(\rho,\delta){\cal H}(\rho,-\delta)
=(-1)^{\tilde{j}_0}\rho^{2\tilde{i}_0}
\{a^2_{\tilde{i}_0,\tilde{j}_0}({\cal H})\delta^{2\tilde{j}_0}
+O(\delta^{2\tilde{j}_0+1})+O(\rho)\}
<0,
$$
where the computation is similar to \eqref{Gi0j0}.
It follows that $\tilde{j}_0$ is odd.
By condition {\bf (P2)}, $j_0$ is even,
which implies that only {\bf(C1)} and {\bf (C3)} is possible.
We only consider {\bf(C3)} since {\bf(C1)} is similar and simpler.
In {\bf(C3)},
we assume without loss of generality that
$a_{i_0,j_0}({\cal G})<0$ and $a_{\tilde{i}_0,\tilde{j}_0}({\cal H})<0$,
i.e., (\ref{assume-c4}) is true for even $j_0$.

There are two cases:
{\bf(4a)} $\tilde{j}_{s({\cal H})}=0$ and {\bf(4b)} $\tilde{j}_{s({\cal H})}>0$.
In {\bf(4a)},
since $\tilde{j}_{s({\cal H})}=0$ in {\bf(J1)},
we have $j_{s({\cal G})}>0$ and therefore
$\theta=0$ is an orbit of system~\eqref{equ:polar system2}.
The only one real branch $\Upsilon_1$ divides the region ${\cal U}_{\epsilon,\delta}$
into two Z-sectors ${\cal U}_1$ and ${\cal U}_2$,
as indicated below \eqref{thetatheta}.
Similarly to \eqref{equ:ttkt} and \eqref{equ:rhosignib},
we have
\begin{align}
\dot\theta|_{\Upsilon_0}<0,~~~
\dot\theta|_{\Upsilon_1}>0,~~~
\dot\theta|_{\Upsilon_2}<0,~~~
\dot\rho|_{{\cal U}_1}<0,~~~
\dot\rho|_{{\cal U}_2}>0.
\label{3t2r}
\end{align}
Hence,
the region between $\Upsilon_1$ and the orbit $\theta=0$ is
a generalized normal sector of Class I (\cite{T-Z1}),
a sector containing no singular points and satisfying that
$\dot \rho<0$ (resp. $>0$) inside and
none of orbits starting from its sides leaves (resp. enters).
By Lemma~1 of \cite{T-Z1}, the sector contains infinitely many orbits connecting with $O$ and therefore
system~\eqref{equ:initial} has infinitely many orbits connecting with $O$
in the direction $\theta=0$.
Furthermore,
similar to \eqref{phitp1} in the situation $N({\cal H})> 1$,
\eqref{3t2r} implies that
the orbit $\phi(t,P_1)$ passing through a point $P_1\in \Upsilon_1$ satisfies that
\begin{align}
\phi(t,P_1)\to O~\mbox{in}~{\cal U}_1~\mbox{as}~t\to+\infty~~~\mbox{and}~~~
\phi(t,P_1)\to O~\mbox{in}~{\cal U}_2~\mbox{as}~t\to-\infty.
\end{align}
Then there is exactly one e-tsector, which passes through the unique real branch $\Upsilon_1$.
Moreover,
there are no other tsectors.

In {\bf(4b)}, we claim that
\begin{description}
\item[ Claim H1:]
$\tilde{j}_0,...,\tilde{j}_{s({\cal H})}$ are odd
and
$a_{\tilde{i}_0,\tilde{j}_0}({\cal H}),...,
a_{\tilde{i}_{s({\cal H})},\tilde{j}_{s({\cal H})}}({\cal H})<0$.

\item[ Claim G1:]
 $j_0,...,j_{s({\cal G})}$ are all even,
$a_{i_0,j_0}({\cal G}),...,$ $a_{i_{k_*-1},j_{k_*-1}}({\cal G})<0$
and
$a_{i_{k_*},j_{k_*}}({\cal G}),$ $...,a_{i_{s({\cal G})},j_{s({\cal G})}}({\cal G})>0$
for an integer $k_*$.
\end{description}
Actually,
{\bf Claim H1} can be proved similarly to {\bf Claim H0}.
By {\bf Claim H1} and {\bf (P1)},
$j_0,...,j_{s({\cal G})}$ are all even.
Furthermore,
since \eqref{3t2r} still holds,
the first three inequalities in \eqref{3t2r} imply that
the equation ${\cal G}(\rho,\theta)=0$ has two real branch
by {\bf(F2)} and {\bf(F3)} of {\bf Claim~4.1}.
By the remark given below Proposition~\ref{lm-solution by NP},
$\Xi({\cal A}({\cal G}),{\cal A}_j({\cal G}))\le 2$.
Since $\tilde{j}_{s({\cal H})}>0$ in {\bf (4b)},
the unique real branch $\Upsilon_1:\theta=\theta_1(\rho)$ is the line $\theta=0$.
We see from \eqref{equ:polar system2} and the second inequality in \eqref{3t2r} that
\begin{align}
0<\theta|_{\Upsilon_1}=
\{
a_{i_{s({\cal G})},j_{s({\cal G})}}({\cal G})
\rho^{i_{s({\cal G})}-1}\theta^{j_{s({\cal G})}}
+O(\theta^{j_{s({\cal G})}+1})
+O(\rho^{i_{s({\cal G})}}\theta^{j_{s({\cal G})}})
\}|_{\theta=0}.
\label{thUp1}
\end{align}
Then $j_{s({\cal G})}=0$ and
$a_{i_{s({\cal G})},j_{s({\cal G})}}({\cal G})>0$.
However,
$a_{i_0,j_0}({\cal G})<0$ as assumed in {\bf (C3)} in the first paragraph of the proof for the situation $N(\mathcal{H})=1$.
Thus $\Xi({\cal A}({\cal G}),{\cal A}_j({\cal G}))= 2$ and therefore in ${\cal A}({\cal G})$
the first $k_*$ vertex coefficients are negative but the rest are positive, which gives {\bf Claim G1}.

Knowing the two claims, we can discuss the number of orbits in the direction $\theta=0$.
In the case {\bf (4b)} of the situation $N(\mathcal{H})=1$
the equation $\mathcal{H}(\rho,\theta)=0$
has exactly one real branch $\Upsilon_1:\theta=\theta_1(\rho)\equiv 0$
as indicated before (\ref{thUp1}),
which
divides ${\cal U}_{\epsilon,\delta}$ into ${\cal U}_{\epsilon,\delta}={\cal U}_1\cup \Upsilon_1\cup {\cal U}_2$,
where
${\cal U}_1=\{(\rho,\theta)\in\mathbb{R}^2:0<\rho<\epsilon,0<\theta<\delta\}$
and
${\cal U}_2=\{(\rho,\theta)\in\mathbb{R}^2:0<\rho<\epsilon,-\delta<\theta<0\}$
are Z-sectors defined below \eqref{thetatheta}.
By {\bf Claim G1}, {\bf Claim~4.3} and Lemma~\ref{lm-basic},
the equation ${\cal G}(\rho,\theta)=0$ in ${\cal U}_1$  has exactly one real branch, which is of the form
\begin{eqnarray}
\theta=\theta_*(\rho):=c_*\rho^{\iota}+o(\rho^\iota)
\label{G-branch}
\end{eqnarray}
and denoted by $\Gamma_*$,
where $-1/\iota$ is the slope of the edge $E_{k_*}({\cal G})$ and
$c_*$ is a nonzero real root of the edge-polynomial ${\cal G}_{E_{k_*}({\cal G})}$.
This real branch (\ref{G-branch}) divides ${\cal U}_1$ into
${\cal U}_1={\cal U}'_1\cup \Gamma_*\cup{\cal U}''_1$, where
\begin{align*}
{\cal U}'_1:=\{(\rho,\theta)\in\mathbb{R}^2:0<\rho<\epsilon,\theta_*(\rho)<\theta<\delta\},~~~
{\cal U}''_1:={\cal U}_1\setminus\overline{{\cal U}'_1},
\end{align*}
and $\overline{{\cal U}'_1}$ is the closure of ${\cal U}'_1$.
In ${\cal U}''_1$,
we have $\dot\theta|_{{\cal U}''_1}>0$ and $\dot\rho|_{{\cal U}''_1}<0$ since \eqref{3t2r} still holds in the case {\bf (4b)},
which implies that no orbits connect with $O$ in ${\cal U}''_1$.
In ${\cal U}'_1$,
we claim that
\begin{align}
\frac{\partial }{\partial \theta} \frac{{\cal G}(\rho,\theta)}{{\cal H}(\rho,\theta)}
=\frac{[{\cal G},{\cal H}]_\theta}{{\cal H}^2}\le \rho^{1/n}
~~~\mbox{for all }~(\rho,\theta)\in {\cal U}'_1,
\label{GHH2U2}
\end{align}
which is the same as \eqref{GHH2U} but does not hold on the whole ${\cal U}_{\epsilon,\delta}$.
If this claim is true,
for the same reason
system~\eqref{equ:polar system2} has at most one orbit
connecting with $O$ in ${\cal U}'_1$,
but on the other hand
we see that ${\cal U}'_1$ contains either $0$ or infinitely many orbits connecting with $O$
as discussed for ${\cal U}_1$ in the proof of {\bf Claim 3.2}.
It follows that no orbits connect with $O$ in ${\cal U}'_1$.
By results on ${\cal U}'_1$ and ${\cal U}''_1$,
no orbits connect with $O$ in ${\cal U}_1$.
Similarly,
we can also prove that no orbits connect with $O$ in ${\cal U}_2$ and therefore in the whole ${\cal U}_{\epsilon,\delta}$,
implying that no orbits connect with $O$ in $\theta=0$.
Similar to the case {\bf (4a)},
there are exactly one e-tsector,
which passes through the only real branch $\Upsilon_1$ of the equation ${\cal H}(\rho,\theta)=0$,
but no other tsectors.

Finally,
we complete the proof of {\bf Claim~4.3} by proving the claimed \eqref{GHH2U2}.
In fact,
\eqref{GHH2U2} is equivalent to the same inequality as the one in (\ref{H>[]})
on the subregion ${\cal U}'_1$ of $\mathcal{U}_{\epsilon,\delta}$.
Let $f:=[{\cal G},{\cal H}]_\theta$.
By Lemma~\ref{J1Ng},
${\cal N}_f$ has $s({\cal G},{\cal H})$ vertices.
We claim that
\begin{align}
\{V_0(f),...,V_{k'_*-1}(f)\}\in \wp(\Delta({\cal H}^2)),~~~
\{V_{k'_*}(f),...,V_{s({\cal G},{\cal H})}(f)\}\notin \wp(\Delta({\cal H}^2))
\label{VfHVfH}
\end{align}
for an integer $k'_*$ satisfying
\begin{align}
\zeta(E_{k'_*}(f))\ge \zeta(E_{k_*}({\cal G})),
\label{EkfEkG}
\end{align}
where
$\wp(\Delta({\cal H}^2)$ is the lower convex semi-hull defined below \eqref{expGH}
and $k_*$ is given in {\bf Claim G1}.
Note that in the situation $N({\cal H})=0$ we only have the first inclusion in (\ref{VfHVfH}) for all vertices
$V_0(f),..., V_{s({\cal G},{\cal H})}(f)$ because
${\cal N}_f$ strictly lies above ${\cal N}_{{\cal H}^2}$ as indicated in \eqref{NH2Ng},
but
${\cal N}_f$ intersects with ${\cal N}_{{\cal H}^2}$ in the current situation $N({\cal H})=1$.
In order to prove the claimed \eqref{VfHVfH},
we notice in the current situation that
\eqref{ijH}, \eqref{iijj1ij1} and \eqref{jj12j} still hold and
$i_0\ge 1$ but $\tilde{i}_0=0$,
which implies that $V_0(f)$ lies in $\wp(\Delta({\cal H}^2))$.
On the other hand,
the second points in \eqref{ijH} and \eqref{iijj1ij1} imply that
$V_{s({\cal G},{\cal H})}(f)$ does not lie in $\wp(\Delta({\cal H}^2))$
since $\tilde{j}_{s({\cal H})}>0$ in {\bf (4b)} and
$j_{s({\cal G})}=0$, obtained below \eqref{thUp1}.
Thus there is an integer $k'_*$ such that
$V_{k'_*}\notin\wp(\Delta({\cal H}^2))$ and
the first inclusion in \eqref{VfHVfH} holds.
Further, we prove that $k'_*$ satisfies \eqref{EkfEkG}.
Otherwise,
\begin{align}
\zeta(E_{k'_*}(f))< \zeta(E_{k_*}({\cal G})).
\label{EkfEkG<}
\end{align}
However, as indicated before (\ref{EkfEkG<}),
$V_0(f)\in \wp(\Delta({\cal H}^2))$ but $V_{k'_*}\notin\wp(\Delta({\cal H}^2))$,
which implies that
the edge $E_{k'_*}(f)$ linking
vertices $V_{k'_*-1}(f):(\breve{i}_{k'_*-1},\breve{j}_{k'_*-1})$ and
$V_{k'_*}(f):(\breve{i}_{k'_*},\breve{j}_{k'_*})$
and the edge $E_{\ell}({\cal H})$ linking
vertices $V_{\ell-1}({\cal H}^2):(\hat{i}_{\ell-1},\hat{j}_{\ell-1})$
and $V_{\ell}({\cal H}^2):(\hat{i}_\ell,\hat{j}_\ell)$
intersect at a point, denoted by $(i',j')$, such that
\begin{align}
\hat{i}_{\ell-1}<i'\le\hat{i}_{\ell},~~~
\hat{j}_{\ell-1}>\breve{j}_{k'_*}~~~\mbox{and}~~~
\zeta(E_{\ell}({\cal H}^2))>\zeta(E_{k'_*}(f)).
\label{jjzz}
\end{align}
Choosing
$\mu'\in(\zeta(E_\ell({\cal H}^2))-\varepsilon,\zeta(E_\ell({\cal H}^2)))$ and
$\mu''\in(\zeta(E_{k'_*}(f)),\zeta(E_{k'_*}(f))+\varepsilon)$
for sufficiently small $\varepsilon>0$,
similarly to \eqref{j6j6j6} and \eqref{coef-a},
we obtain that
\begin{align}
&\hat{j}_{\ell-1}=\tilde{j}_{\varrho_1}+\tilde{j}_{\varrho_1},~~
\breve{j}_{k'_*}=\tilde{j}_{\sigma'}-1+\tilde{j}_{\varrho_2},~~
\label{j-a1}
\\
&a_{\breve{i}_{k'_*},\breve{j}_{k'_*}}(f)
=(j_{\sigma'}-\tilde{j}_{\varrho_2})
a_{i_{\sigma'},j_{\sigma'}}({\cal G})
a_{\tilde{i}_{\varrho_2},\tilde{j}_{\varrho_2}}({\cal H}),
\label{j-a2}
\end{align}
where points $(\tilde{i}_{\varrho_1},\tilde{j}_{\varrho_1})$ and $(\tilde{i}_{\varrho_2},\tilde{j}_{\varrho_2})$
are $\mu'$- and $\mu''$-components of ${\cal N}_{\cal H}$ and
point $(i_{\sigma'},j_{\sigma'})$ is the $\mu''$-component of ${\cal N}_{\cal G}$.
We see from \eqref{EkfEkG<} that $\mu''<\zeta(E_{k_*}({\cal G}))$ and therefore
$a_{i_{\sigma'},j_{\sigma'}}({\cal G})<0$ by {\bf Claim G1}.
Moreover,
$a_{\breve{i}_{k'_*},\breve{j}_{k'_*}}(f)>0$
for the same reason as indicated below \eqref{coef-a},
and $a_{\tilde{i}_{\varrho_2},\tilde{j}_{\varrho_2}}({\cal H})>0$ by {\bf Claim H1}.
It follows from (\ref{j-a2}) that $j_{\sigma'}>\tilde{j}_{\varrho_2}$.
On the other hand,
$\tilde{j}_{\varrho_1}\le \tilde{j}_{\varrho_2}$ since $\mu'>\mu''$
by the second inequality in \eqref{jjzz}.
Consequently, from (\ref{j-a1}) we get
$\hat{j}_{\ell-1}\le \breve{j}_{k'_*}$,
a contradiction to the second inequality in \eqref{jjzz}.
Thus we obtain \eqref{EkfEkG} by indirect proof.
At last,
we prove that those chosen $k'_*$ also satisfies
the second inclusion in \eqref{VfHVfH}.
Actually,
if it is not true,
then at least one of vertices listed in the second inclusion of \ref{VfHVfH} lies in $\wp(\Delta({\cal H}^2))$.
However, $V_{k'_*}(f)\notin \wp(\Delta({\cal H}^2))$ as shown just before (\ref{EkfEkG<}).
It follows that
the edge $E_{\tilde{k}'}({\cal H}^2)$
linking vertices $V_{\tilde{k}'-1}({\cal H}^2):
(\hat{i}_{\tilde{k}'-1},\hat{j}_{\tilde{k}'-1})$ and
$V_{\tilde{k}'}({\cal H}^2):(\hat{i}_{\tilde{k}'},\hat{j}_{\tilde{k}'})$ and
the edge $E_{\tilde{k}''}(f)$ linking vertices $V_{\tilde{k}''-1}(f):(\breve{i}_{\tilde{k}''-1},\breve{j}_{\tilde{k}''-1})$ and
$V_{\tilde{k}''}(f):(\breve{i}_{\tilde{k}''},\breve{j}_{\tilde{k}''})$,
where $\tilde{k}''\ge k'_*+1$,
intersect at a point, denoted by $(\tilde{i}^*,\tilde{j}^*)$, such that
\begin{align}
\breve{i}_{\tilde{k}''-1}< \tilde{i}^* \le \breve{i}_{\tilde{k}''},
~~~
\breve{j}_{\tilde{k}''-1}>\hat{j}_{\tilde{k}'}
~~~\mbox{and}~~~
\zeta(E_{\tilde{k}''}(f))>\zeta(E_{\tilde{k}'}({\cal H}^2)).
\label{ZJ}
\end{align}
Thus, the same arguments as given below (\ref{jjzz}) also lead to a contradiction to the second inequality in (\ref{ZJ})
and therefore the claimed \eqref{VfHVfH} is proved.

Having \eqref{VfHVfH}, we return to prove \eqref{GHH2U2}, i.e.,
the inequality ${\cal H}^2\ge \rho^{-1/n} f$ on ${\cal U}'_1$.
We see from \eqref{VfHVfH} that
\begin{equation}
\begin{split}
\wp(\Delta({\cal H}^2-f))
=\wp(\{V_0({\cal H}^2),...,V_{\ell}({\cal H}^2),
V_{k'_*}(f),...,V_{s({\cal G},{\cal H})}(f)\}),
\label{wpDH2f}
\end{split}
\end{equation}
where $\ell$ is defined just before \eqref{jjzz}.
Hence the vertex sequence of the Newton polygon of the function ${\cal H}^2-f$
satisfies that
$$
\vec{\Delta}^V({\cal H}^2-f)
=(V_0({\cal H}^2),...V_{\ell'}({\cal H}^2),V_{k''_*}(f),...,V_{s({\cal G},{\cal H})}(f))
$$
for integers $\ell'\in \{0,...,\ell\}$ and $k''_*\in\{k'_*,...,s({\cal G},{\cal H})\}$.
Then the point $V_{k''_*-1}(f)$ lies above the line linking points
$V_{\ell'}({\cal H}^2)$ and $V_{k''_*}(f)$.
If the the point $V_{k''_*-1}(f)$ exactly lies on the line,
then $k''_*>k'_*$
since $V_{k'_*-1}(f)\in \wp(\Delta({\cal H}^2))$ by \eqref{VfHVfH}.
It follows that
$$
\zeta(\underline{V_{\ell'}({\cal H}^2)V_{k''_*}(f)})
=\zeta(E_{k''_*}(f))
>\zeta(E_{k'_*}(f))
=-1/\iota.
$$
On the other hand,
if the point $V_{k''_*-1}(f)$ strictly lies above the line, then
$$
\zeta(\underline{V_{\ell'}({\cal H}^2)V_{k''_*}(f)})
>\zeta(E_{k''_*}(f))
\ge \zeta(E_{k'_*}(f))
=-1/\iota.
$$
Consequently,
each edge and its edge-polynomial of ${\cal H}^2-f$ are the same as those of ${\cal H}^2$
if the slope of the edge is not greater than $-1/\iota$.
As discussing just below (\ref{stretching-n}) in the situation $N({\cal H})=0$,
we similarly see that
each edge-polynomial of $\hat{\cal H}^2-\hat{f}$ is the same as that of $\hat{\cal H}^2$
if the slope of the edge is not greater than $-1/(n\iota)$.
Moreover,
by the same arguments as given at the end of the proof in the situation $N({\cal H})=0$,
each edge-polynomial of $\hat{\cal H}^2$ has no nonzero real roots.
It follows that each edge-polynomial $\hat{\cal H}^2-\hat{f}$ has no nonzero real roots
if the slope of the edge is not greater than $-1/(n\iota)$.
By Lemma~\ref{lm-basic},
the equation $\hat{\cal H}^2-\hat{f}=0$ has no real branches in the region
$$
\hat{\cal U}'_1:=\{(\hat{\rho},\theta)\in\mathbb{R}^2:
0<\hat{\rho}<\epsilon^{1/n},\theta_*(\hat{\rho}^n)<\theta<\delta\},
$$
where $\theta_*(\hat{\rho}^n)=c_*\hat\rho^{n\iota}+o(\hat\rho^{n\iota})$,
with the lowest order $n\iota$
according to the expression of $\theta_*(\rho)$ given before \eqref{GHH2U2}.
It follows that
$\hat{\cal H}^2-\hat{f}\ge 0$ on $\hat{\cal U}'_1$
and therefore
${\cal H}^2-\rho^{-1/n}[{\cal G},{\cal H}]_\theta\ge 0$ on ${\cal U}'_1$.
As indicated in the beginning of this paragraph,
the claimed \eqref{GHH2U2} is obtained.
This completes the proof of Theorem~\ref{th:finite}.
\qquad$\Box$

\begin{table}[h]
\begin{minipage}{0.5\linewidth}
\centering
\begin{tabular}{|l|l|l|}
\hline
$N({\cal G})$ & odd $N({\cal H})$   & even $N({\cal H})$  \\ \hline
{\bf(C1)}     & $N({\cal H})+1$     & $N({\cal H})$       \\ \hline
{\bf(C2)}     & $N({\cal H})$       & $N({\cal H})+1$     \\ \hline
{\bf(C3)}     & $N({\cal H})-1$     & $N({\cal H})$       \\ \hline
{\bf(C4)}     & $N({\cal H})$       & $N({\cal H})-1$     \\ \hline
\end{tabular}
\\
\vspace{5pt}
(a) in results {\bf(ia)} and {\bf(iia)}
\end{minipage}
\begin{minipage}{0.5\linewidth}
\centering
\begin{tabular}{|l|l|l|}
\hline
$N({\cal G})$ & odd $N({\cal H})$  & even $N({\cal H})$    \\ \hline
{\bf(C1)}     & $N({\cal H})-1$    & $N({\cal H})$         \\ \hline
{\bf(C2)}     & $N({\cal H})$      & $N({\cal H})-1$       \\ \hline
{\bf(C3)}     & $N({\cal H})+1$    & $N({\cal H})$         \\ \hline
{\bf(C4)}     & $N({\cal H})$      & $N({\cal H})+1$       \\ \hline
\end{tabular}
\\
\vspace{5pt}
(b) in results {\bf(ib)} and {\bf(iib)}
\end{minipage}
\caption{Relations between $N({\cal G})$ and $N({\cal H})$.}
\label{table-NGNH}
\end{table}


As remarked just after Theorem~\ref{th:finite},
results {\bf(ia)} and {\bf(ib)} are stated in terms of $N({\cal G})$
but {\bf(iia)} and {\bf(iib)} in terms of $N({\cal H})$.
Actually, both $N({\cal G})$ and $N({\cal H})$ can be used to state those results.
In fact,
only results on e-tsectors and h-tsectors are stated in $N({\cal G})$ or $N({\cal H})$.
By {\bf Claim~1.2} and {\bf Claim~3.2},
the two edges of each h-tsector lie in two adjacent Z-sectors separately and therefore
each h-tsector determines a unique real branch of the equation ${\cal H}(\rho,\theta)=0$.
As stated in the second last paragraph and the third last paragraph
in {\bf Step~2}
and the last paragraph in {\bf Step~4},
each e-tsector also determines a unique real branch of
the equation ${\cal H}(\rho,\theta)=0$.
By {\bf Claims} {\bf 1.1}, {\bf 2.1}, {\bf 3.1} and {\bf 4.1},
there is one real branch of the equation ${\cal G}(\rho,\theta)=0$
$($and ${\cal H}(\rho,\theta)=0$$)$ between each two adjacent real branches of the equation
${\cal H}(\rho,\theta)=0$ $($and ${\cal G}(\rho,\theta)=0$$)$.
Since
$N({\cal G})$ and $N({\cal H})$ are numbers of real branches of the equations
${\cal G}(\rho,\theta)=0$ and ${\cal H}(\rho,\theta)=0$ respectively by Lemma~\ref{Assertion1},
we see that
$$
N({\cal H})=N({\cal G})\pm 1
\  \ \mbox{  or  }\ \
N({\cal H})=N({\cal G}).
$$
More precisely,
in results {\bf(ia)} and {\bf(iia)} (resp. {\bf(ib)} and {\bf(iib)}), the relations between
$N({\cal G})$ and $N({\cal H})$ are shown in Table \ref{table-NGNH}(a) (resp. Table \ref{table-NGNH}(b)),
where {\bf(C1)}-{\bf(C4)} are given just at the beginning of the proof of {\bf Claim~1.3}.
Note that
the relation in circumstance {\bf(C1)} in Table~\ref{table-NGNH}(a)
(resp. Table~\ref{table-NGNH}(b)) is proved in the last paragraph
of the proof of {\bf Claim 1.3} (resp. {\bf Claim 2.3}).
The relations in other circumstances in Table~\ref{table-NGNH}(a)
(resp. Table~\ref{table-NGNH}(b)) can be proved similarly.
However,
the cases {\bf(ia)} and {\bf(ib)}, where the exceptional direction $\theta=0$ is isolated,
and the cases {\bf(iia)} and {\bf(iib)}, where $\theta=0$ is non-isolated,
are quite different.
This difference only occurs in the Z-sectors ${\cal U}_1$ and ${\cal U}_{N({\cal H})+1}$, the first one and the last one
as defined just below \eqref{thetatheta}.
In results {\bf(ia)} and {\bf(ib)} we prefer $N({\cal G})$ to $N({\cal H})$
because there are 8 situations if we state in terms of $N({\cal H})$.
In contrast,
in results {\bf(iia)} and {\bf(iib)} we prefer $N({\cal H})$ to $N({\cal G})$
because there are also 8 situations if we state in terms of $N({\cal G})$.


Next, we prove Theorem~\ref{th:finiteJ2}.
Corresponding to Lemma~\ref{Assertion1}, the following lemma is needed.

\begin{lm}
If conditions {\bf(P2)}, {\bf(Q)} and {\bf(S)} are satisfied
and one of conditions {\bf(H1)}, {\bf(H2)}, {\bf(H3)} and {\bf(H4)} of
Theorem~\ref{th:finiteJ2} holds,
then equations ${\cal G}(\rho,\theta)=0$ and ${\cal H}(\rho,\theta)=0$ have
$N({\cal G})$ and $N({\cal H})$ real branches on the half-plane $\rho\ge0$ respectively.
\label{Assertion1'}
\end{lm}

{\bf Proof.}
The proof is similar to that of Lemma~\ref{Assertion1}.
We only give the proof for the equation ${\cal G}(\rho,\theta)=0$
when conditions {\bf(P2)}, {\bf(Q)} and {\bf(S)} are satisfied and
condition {\bf(H1)} of Theorem~\ref{th:finiteJ2} holds,
and the proof is similar
when either {\bf(H2)} or {\bf(H3)} or {\bf(H4)} of Theorem~\ref{th:finiteJ2} holds.
Moreover,
it is similar to discuss the equation ${\cal H}(\rho,\theta)=0$.
As indicated just below \eqref{L2GHV},
common vertices of ${\cal G}'_\theta{\cal H}$ and $-{\cal G}{\cal H}'_\theta$
satisfy the non-vanishing condition \eqref{Anon} for addition because of {\bf(P2)}.
It follows that Lemma~\ref{precondiJ2}{\bf(i)} holds.
By {\bf(H1)}, Lemma~\ref{precondiJ2}{\bf(i)} contains exactly
the case $j_*<j_{s({\cal G})-1}$ and the case $j_*=j_{s({\cal G})-1}$.

In the first case,
by equality \eqref{J2igEkg} in Lemma~\ref{precondiJ2}{\bf(i)},
equality \eqref{KGKHK4}, given in the proof of Lemma~\ref{Assertion1}, still holds.
Then it is the same to show that
the equation ${\cal G}(\rho,\theta)=0$ has
$\chi({\cal G})$ nontrivial real branches on the half-plane $\rho\ge0$.
Additionally,
the equation ${\cal G}(\rho,\theta)=0$ does not have
a trivial real branch $\theta(\rho)\equiv0$ since $j_{s({\cal G})}=0$
in the two-below case {\bf(J2)}.
Thus,
by the relation of $N({\cal G})$ and $\chi({\cal G})$,
given just before Theorem~\ref{th:finite},
the equation ${\cal G}(\rho,\theta)=0$ has
$N({\cal G})$ real branches on the half-plane $\rho\ge0$,
i.e., Lemma \ref{Assertion1'} holds in the first case.

In the second case,
equality \eqref{J2igEkg} in Lemma~\ref{precondiJ2}{\bf(i)} implies that
\eqref{KGKHK4} still holds for all
$\xi\in\mathfrak{S}({\cal G})\backslash\{\zeta(E_{s({\cal G})}({\cal G}))\}$
since $\zeta(E_{s({\cal G})}({\cal G}))=\xi_{s({\cal G},{\cal H})}$
when $\zeta(E_{s({\cal G})}({\cal G}))>\zeta(E_{s({\cal H})}({\cal H}))$.
Thus,
we need to consider the situation $\xi<\zeta(E_{s({\cal G})}({\cal G}))$
and the situation $\xi=\zeta(E_{s({\cal G})}({\cal G}))$ separately.
In the first situation,
it is the same as the proof of Lemma~\ref{Assertion1} to show that
each nonzero real roots of $E$-polynomial ${\cal G}_E$ is of multiplicity 1 and
determines one nontrivial real branch of the equation ${\cal G}(\rho,\theta)=0$,
where $E$ is an edge of ${\cal N}_{\cal G}$ with slope $<\zeta(E_{s({\cal G})}({\cal G}))$.
In the second situation,
we see from \eqref{defkxip} and \eqref{fonE} that
\begin{align*}
{\cal K}_\xi({\cal G})(\theta)
&={\cal G}_{E_{s({\cal G})}({\cal G})}(\theta)
=\sum_{(i,j)\in E_{s({\cal G})}({\cal G})\cap\Delta({\cal G})}a_{i,j}({\cal G})\theta^j
\\
&=a_{i_{s({\cal G})-1},j_{s({\cal G})-1}}({\cal G})\theta^{j_{s({\cal G})-1}}
+a_{i_{s({\cal G})},j_{s({\cal G})}}({\cal G}).
\end{align*}
Actually,
the equality $j_*=j_{s({\cal G})-1}$ in the second case implies that
$E_{s({\cal G})}({\cal G})\cap\Delta({\cal G})
=\{V_{s({\cal G})-1}({\cal G}),V_{s({\cal G})}({\cal G})\}$
because $j_*$ is the ordinate of
the second to the last valid index on ${\cal N}_{\cal G}$
and $j_{s({\cal G})-1}$ is the ordinate of
the second to the last vertex of ${\cal N}_{\cal G}$.
Additionally, $j_{s({\cal G})-1}$ is odd by {\bf(P2)}.
Then,
the polynomial ${\cal G}_{E_{s({\cal G})}({\cal G})}$ has one nonzero real root,
a simple one,
which determines one real branch of the equation ${\cal G}(\rho,\theta)=0$
by Lemma~\ref{lm-basic}{\bf(ii)}.
As a consequence of the above two situations,
the equation ${\cal G}(\rho,\theta)=0$ also has $\chi({\cal G})$, defined in \eqref{defChi},
nontrivial real branches on the half-plane $\rho\ge0$.
Moreover,
the equation ${\cal G}(\rho,\theta)=0$ does not have
a trivial real branch $\theta(\rho)\equiv0$
since $j_{s({\cal G})}=0$ in the two-below case {\bf(J2)}.
Then, the equation ${\cal G}(\rho,\theta)=0$ has
$N({\cal G})$ real branches on the half-plane $\rho\ge0$,
i.e., Lemma \ref{Assertion1'} holds in the second case.
Thus, the proof of Lemma \ref{Assertion1'} is completed.
\qquad$\Box$


{\bf Proof of Theorem~\ref{th:finiteJ2}.}
Although this proof is aimed to case {\bf (J2)}, i.e., neither ${\cal N}_{\cal G}$ nor ${\cal N}_{\cal H}$ ends above the $u$-axis,
the procedure of this proof is the same as that of Theorem~\ref{th:finite}
but we use Lemma \ref{Assertion1'} instead of Lemma \ref{Assertion1}.
Unlike case {\bf (J1)}, where conditions {\bf (P1)} and {\bf (Q)} are enough,
we need
conditions {\bf(P2)}, {\bf(Q)}, {\bf(S)}
and one of conditions {\bf(H1)}, {\bf(H2)}, {\bf(H3)} and {\bf(H4)}, given in Theorem~\ref{th:finiteJ2} in case {\bf(J2)},
to ensure by Lemma~\ref{defsignth1} that
$[{\cal G},{\cal H}]_\theta$ is semi-definite on the region
${\cal U}_{\epsilon,\delta}:=\{(\rho,\theta)\in {\mathbb R}^2:0<\rho<\epsilon,|\theta|<\delta\}$.
Thus we repeat {\bf Steps 1}-{\bf 3} and {\bf Steps 4} with $N({\cal H})>1$ given in the proof of Theorem~\ref{th:finite}
to prove results of {\bf (ia)}, {\bf (ib)} and {\bf (iia)} and
the result of {\bf (iib)} with $N({\cal H})>1$ respectively.

In what follows,
we consider the situations $N({\cal H})=0$ and $N({\cal H})=1$ for results in {\bf (iib)},
appeared in {\bf Step~4} of the proof of Theorem~\ref{th:finite} but the discussion of which
is more complicated in case {\bf (J2)}.

In the situation $N({\cal H})=0$,
we also have the same circumstances {\bf (C1)}-{\bf (C4)}
as indicated in {\bf Step~4} of the proof of Theorem~\ref{th:finite},
but we claim that only {\bf (C4)} is possible, i.e.,
$a_{i_0,j_0}({\cal G})a_{\tilde{i}_0,\tilde{j}_0}({\cal H})>0$ and $j_0$ is odd.
In fact,
the equation ${\cal H}(\rho,\theta)=0$ has no real branches by Lemma~\ref{Assertion1'},
which implies that {\bf Claim H0} given in the proof of Theorem~\ref{th:finite} still holds.
On the other hand,
correspondingly to {\bf Claim G0}, we have
\begin{description}
\item[Claim G0$'$:]
For $\iota=0,...,s({\cal G})-1$,
$j_\iota$\,s are odd and
$a_{i_\iota,j_\iota}({\cal G})$\,s
have the same sign.
\end{description}
Actually,
$j_\iota$\,s are odd
by {\bf Claim H0} and {\bf(P2)}.
Furthermore,
noting before {\bf Claim G0$'$} that the equation ${\cal H}(\rho,\theta)=0$ has no real branches,
by {\bf(F3)} of {\bf Claim 4.1} in the proof of Theorem~\ref{th:finite}, which still holds now,
we see that the equation ${\cal G}(\rho,\theta)=0$ has at most one real branch.
It follows from Proposition~\ref{lm-solution by NP} and the remark below its proof that
$\Xi({\cal A}({\cal H}),{\cal A}_j({\cal H}))\le 1$.
Since $j_{s({\cal G})}=0$ in {\bf(J2)} and $j_{s({\cal G})-1}$ is odd,
either
$$
a_{i_{s({\cal G})-1},j_{s({\cal G})-1}}\!({\cal G})
a_{i_{s({\cal G})},j_{s({\cal G})}}\!({\cal G})\!<\!0
~or~
(-1)^{j_{s({\cal G})}}a_{i_{s({\cal G})},j_{s({\cal G})}}\!({\cal G})
(-1)^{j_{s({\cal G})-1}}a_{i_{s({\cal G})-1},j_{s({\cal G})-1}}\!({\cal G})\!<\!0.
$$
It follows that
$a_{i_0,j_0}({\cal G}),...,a_{i_{s({\cal G})-1},j_{s({\cal G})-1}}({\cal G})$
all have the same sign, which proves {\bf Claim G0$'$}.
By {\bf Claim H0} and {\bf Claim G0$'$},
which is similar to {\bf Claim G0}
but $j_{s({\cal G})}$ and $a_{i_{s({\cal G})},j_{s({\cal G})}}({\cal G})$ are not involved,
we also see that only {\bf (C4)} is possible because
the missed $j_{s({\cal G})}$ and $a_{i_{s({\cal G})},j_{s({\cal G})}}({\cal G})$
do not affect the proof.
This proves our claim on {\bf (C4)}.
Having {\bf (C4)}, without loss of generality,
we still assume that \eqref{assume-c4} holds.
Then
by {\bf Claim H0}, {\bf Claim G0$'$} and Lemma~\ref{precondiJ2}
we similarly have \eqref{GHH2U}, which implies
that
system~\eqref{equ:polar system2} has at most one orbit
connecting with $O$ in ${\cal U}_{\epsilon,\delta}$.
On the other hand,
noticing that $\tilde{j}_0$ is even by {\bf Claim H0}
and that $\tilde{i}_0=0$ since $H_0(\theta)\not\equiv 0$
when $G_0(\theta)\equiv 0$ as indicated below \eqref{equ:polar system2},
we see from \eqref{equ:polar system2}, \eqref{assume-c4} and {\bf Claim H0} that
$$
\dot \theta|_{\rho=0}={\cal H}(0,\theta)=H_0(\theta)
=a_{0,\tilde{j}_0}({\cal H})\theta^{\tilde{j}_0}+O(\theta^{\tilde{j}_0+1})<0
~~~\mbox{for all}~
\theta\in(-\delta,0)\cup(0,\delta).
$$
It follows that
orbits starting from the two linear segments
$$
\{(\rho,\theta)\in\mathbb{R}^2:\rho=0,0<\theta<\delta\}
~~~\mbox{and}~~~
\{(\rho,\theta)\in\mathbb{R}^2:\rho=0,-\delta<\theta<0\}
$$
on the boundary of ${\cal U}_{\epsilon,\delta}$
all leave ${\cal U}_{\epsilon,\delta}$.
Then,
similar to the discussion on ${\cal U}_1$ given just below \eqref{dotrhoU*},
system~\eqref{equ:polar system2} has either one or infinitely many orbits connecting with $O$
in ${\cal U}_{\epsilon,\delta}$.
Consequently,
system~\eqref{equ:polar system2} has exactly one orbit
connecting with $O$ in ${\cal U}_{\epsilon,\delta}$,
which implies that
system~\eqref{equ:initial} has exactly one orbit connecting with $O$ in the direction $\theta=0$.
Moreover,
similar to the discussion given just below \eqref{wGHwp},
there is one p-tsector and no other tsectors
in the situation $N({\cal H})=0$.
Thus,
those results in {\bf (iib)} with $N({\cal H})=0$ are proved.

In the other situation $N({\cal H})=1$,
similarly to the same situation in {\bf Step 4}
in the proof of Theorem~\ref{th:finite},
we can also see that $\tilde{j}_0$ is odd and $j_0$ is even and therefore
only circumstances {\bf (C1)} and {\bf(C3)} are possible.
In the following,
we only consider {\bf (C3)} because the discussion in {\bf (C1)} is similar and simpler.
In {\bf (C3)} we assume without loss of generality that
\begin{align}
a_{i_0,j_0}({\cal G})<0~~~\mbox{and}~~~a_{\tilde{i}_0,\tilde{j}_0}({\cal H})<0,
\label{assume-c3}
\end{align}
as in the end of the paragraph before \eqref{3t2r}.
Note that
in {\bf (iib)} of Theorem~\ref{th:finiteJ2} we only need to consider the two cases:
{\bf (4$_+$)} $a_{i_0,j_0}({\cal G})a_{i_{s({\cal G})},j_{s({\cal G})}}({\cal G})>0$
and
{\bf (4$_-$)} $a_{i_0,j_0}({\cal G})a_{i_{s({\cal G})},j_{s({\cal G})}}({\cal G})<0$
when $N({\cal H})=1$.

\begin{figure}[h]
    \centering
     \subcaptionbox{%
     }{\includegraphics[height=1.6in]{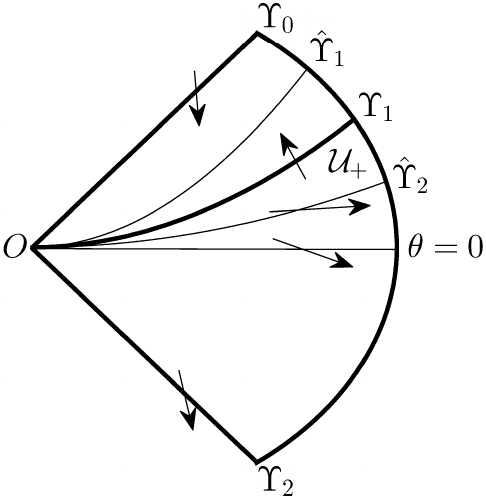}}~~~~~~~~~~~~~
     \subcaptionbox{%
     }{\includegraphics[height=1.6in]{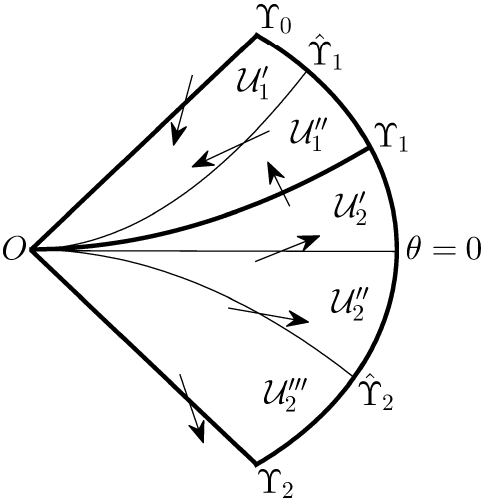}}
    \caption{Division of ${\cal U}_{\epsilon,\delta}$ in (a) case {\bf(4$_+$)} and (b) case {\bf(4$_-$)}.}
    \label{fig:th2nh1}
\end{figure}

In case {\bf (4$_+$)},
similar computation to \eqref{Gi0j0} shows that
\begin{align}
\dot\theta|_{\theta=0}={\cal G}(0,\rho)/\rho
=a_{i_{s({\cal G})},j_{s({\cal G})}}({\cal G})\rho^{i_{s({\cal G})}-1}
+O(\rho^{i_{s({\cal G})}})
<0
\label{dtheta0}
\end{align}
because of \eqref{assume-c3} and the inequality of {\bf (4$_+$)}.
On the other hand,
similar to \eqref{equ:ttkt} and \eqref{equ:rhosignib},
assumption \eqref{assume-c3} ensures that
\begin{align}
\dot\theta|_{\Upsilon_0}<0,~~~
\dot\theta|_{\Upsilon_1}>0,~~~
\dot\theta|_{\Upsilon_2}<0,~~~
\dot\rho|_{{\cal U}_1}<0~~~\mbox{and}~~~
\dot\rho|_{{\cal U}_2}>0,
\label{3theta2rho}
\end{align}
where ${\cal U}_1$ and ${\cal U}_2$ are Z-sectors
obtained from the division of ${\cal U}_{\epsilon,\delta}$ by the curve $\Upsilon_1:\theta=\theta_1(\rho)$ given in \eqref{thetatheta},
which
by Lemma~\ref{Assertion1'}
is the only real branch of the equation ${\cal H}(\rho,\theta)=0$.
By {\bf (F2)} and {\bf(F3)} of {\bf Claim 4.1},
which remains true in the present situation,
each of ${\cal U}_1$ and ${\cal U}_2$ contains exactly one real branch of
the equation ${\cal G}(\rho,\theta)=0$,
denoted by
\begin{align}
\hat\Upsilon_1:\theta=\hat\theta_1(\rho), ~~\rho\in(0,\epsilon)~~~\mbox{and}~~~
\hat\Upsilon_2:\theta=\hat\theta_2(\rho), ~~\rho\in(0,\epsilon)
\label{bhG}
\end{align}
respectively,
and moreover both of them are of odd multiplicity.
By \eqref{dtheta0} and the inequalities on $\Upsilon_0$ and $\Upsilon_2$ in \eqref{3theta2rho},
we see that the vector field rotates clockwise on curves $\Upsilon_0$, $\theta=0$ and $\Upsilon_2$.
It follows from the odd multiplicity that
the curves $\hat\Upsilon_1$ and $\hat\Upsilon_2$
both lie between either $\Upsilon_0$ and $\theta=0$ or $\theta=0$ and $\Upsilon_2$.
Since
$\hat\theta_2(\rho)<\theta_1(\rho)<\hat{\theta}_1(\rho)$
for all $\rho\in(0,\epsilon)$
by the division of ${\cal U}_{\epsilon,\delta}$,
the sector ${\cal U}_+$ (or ${\cal U}_-$) between $\Upsilon_1$ and $\hat\Upsilon_2$
(or between $\hat\Upsilon_1$ and $\Upsilon_1$)
lies in the first (or fourth) quadrant.
Further,
we see from \eqref{3theta2rho} that
orbits starting from $\Upsilon_1$ and $\hat\Upsilon_2$ (or $\hat\Upsilon_1$ and $\Upsilon_1$)
all leave (or enter) the region ${\cal U}_+$ (or ${\cal U}_-$),
as illustrated in Figure~\ref{fig:th2nh1}(a).
Note that curves $\Upsilon_1$, $\hat\Upsilon_1$ and $\hat\Upsilon_2$
all tangent to the $x$-axis at the origin $O$.
Then ${\cal U}_+$ (or ${\cal U}_-$)
is a generalized normal sector of Class I (\cite{T-Z1})
and therefore
infinitely many orbits of system~\eqref{equ:initial}
connect with $O$ in the direction $\theta=0$ by Lemma~1 of \cite{T-Z1}.
Moreover,
similar to the situation $N({\cal H})= 1$
in {\bf Step 4} of the proof of Theorem~\ref{th:finite},
there is exactly one e-tsector,
passing through the only real branch $\Upsilon_1$
of the equation ${\cal H}(\rho,\theta)=0$,
but no other tsectors.


In case {\bf (4$_-$)},
there are two circumstances:
{\bf (4$_-^+$)} $0<\theta_1(\rho)<\delta$ for all $\rho\in(0,\epsilon)$,
and
{\bf (4$_-^-$)} $-\delta<\theta_1(\rho)<0$ for all $\rho\in(0,\epsilon)$,
where $\theta=\theta_1(\rho)$ is the only real branch of the equation ${\cal H}(\rho,\theta)=0$,
as indicated just below \eqref{3theta2rho}.
We only discuss {\bf (4$_-^+$)} because the discussion of {\bf (4$_-^-$)} is similar.
In {\bf (4$_-^+$)},
we claim that
\begin{description}
\item[Claim H1$'$:]
$\tilde{j}_0,...,\tilde{j}_{q-1}$ are all odd and
$\tilde{j}_q,...,\tilde{j}_{s({\cal H})}$ are all even for an integer
$q\in\{1,...,s({\cal H})\}$,
$a_{\tilde{i}_0,\tilde{j}_0}({\cal H}),...,a_{\tilde{i}_{q-1},\tilde{j}_{q-1}}({\cal H})<0$
and
$a_{\tilde{i}_q,\tilde{j}_q}({\cal H}),...,
a_{\tilde{i}_{s({\cal H})},\tilde{j}_{s({\cal H})}}({\cal H})>0$.

\item[Claim G1$'$:]
$j_0,...,j_{q'-1}$ are all even and
$j_{q'},...,j_{s({\cal G})-1}$ are all odd for an integer $q'\in\{1,...,s({\cal G})-1\}$
such that $\zeta(E_{q'}({\cal G}))=\zeta(E_q({\cal H}))$,
$a_{i_0,j_0}({\cal G}),...,a_{i_{q'-1},j_{q'-1}}({\cal G})<0$
and
$a_{i_{q'},j_{q'}}({\cal G}),...,a_{i_{s({\cal G})},j_{s({\cal G})}}({\cal G})>0$.
\end{description}

We prove the two claims at the last of the proof.
By {\bf Claim H1$'$} and Lemma~\ref{lm-basic},
the only real branch
$\Upsilon_1: \theta=\theta_1(\rho)$ of the equation ${\cal H}(\rho,\theta)=0$
is of the form
\begin{align}
\theta=\theta_1(\rho)=c_1\rho^{-1/\xi}+o(\rho^{-1/\xi}),
\label{H1'theta1}
\end{align}
where $c_1$ is a positive root of the polynomial ${\cal K}_\xi({\cal H})$
and $\xi:=\zeta(E_q({\cal H}))$.
Similar to case {\bf(4$_+$)},
this branch divides ${\cal U}_{\epsilon,\delta}$ into
two Z-sectors ${\cal U}_1$ and ${\cal U}_2$.
On the other hand,
by {\bf Claim G1$'$} and Lemma~\ref{lm-basic},
the two real branches $\hat\Upsilon_1:\theta=\hat\theta_1(\rho)$ and $\hat\Upsilon_2:\theta=\hat\theta_2(\rho)$
of the equation ${\cal G}(\rho,\theta)=0$,
which are obtained in ${\cal U}_1$ and ${\cal U}_2$ as in \eqref{bhG},
are of the forms
$$
\theta=\hat\theta_1(\rho)=\hat{c}_1\rho^{-1/\xi}+o(\rho^{-1/\xi})
~~~\mbox{and}~~~
\theta=\hat\theta_2(\rho)=\hat{c}_2\rho^{-1/\iota}+o(\rho^{-1/\iota}),
$$
where $\hat{c}_1>0$ and $\hat{c}_2<0$ are real roots of
${\cal K}_\xi({\cal G})$ and ${\cal K}_\iota({\cal G})$ respectively,
$\xi$ is given below (\ref{H1'theta1}),
and $\iota:=\zeta(E_{s({\cal G})}({\cal G}))$.
Further,
the real branch $\hat\Upsilon_1$ divides
${\cal U}_1$ into two subregions
$$
{\cal U}'_1:=\{(\rho,\theta)\in\mathbb{R}^2:0<\rho<\epsilon,\hat\theta_1(\rho)<\theta<\delta\}
~~~\mbox{and}~~~
{\cal U}''_1:={\cal U}_1\backslash\overline{{\cal U}'_1},
$$
where $\overline{{\cal U}'_1}$ is the closure of ${\cal U}'_1$,
and the real branch $\hat\Upsilon_2$ and the line $\theta=0$
divide ${\cal U}_2$ into the three subregions
\begin{align*}
{\cal U}'_2&:=\{(\rho,\theta)\in\mathbb{R}^2:
0<\rho<\epsilon,0<\theta<\theta_1(\rho)\},
\\
{\cal U}''_2&:=\{(\rho,\theta)\in\mathbb{R}^2:
0<\rho<\epsilon,\hat\theta_2(\rho)<\theta<0\}
\end{align*}
and ${\cal U}'''_2:={\cal U}_2\backslash\overline{{\cal U}'_2\cup{\cal U}''_2}$,
as shown in Figure~\ref{fig:th2nh1}(b).
In what follows, we discuss the number of orbits
in the five sectors ${\cal U}'_1$, ${\cal U}''_1$, ${\cal U}'_2$, ${\cal U}''_2$ and ${\cal U}'''_2$
separately.

The region ${\cal U}''_1$ (resp. ${\cal U}''_2$) contains no orbits
connecting with the singular point $O$
because ${\cal U}''_1$ (resp. ${\cal U}''_2$) lies in the first (resp. fourth) quadrant and
$\dot\rho\dot\theta<0$ (resp. $>0$) in ${\cal U}''_1$ (resp. ${\cal U}''_2$)
by \eqref{3theta2rho}.

For the remaining regions ${\cal U}'_1$, ${\cal U}'_2$ and ${\cal U}'''_2$,
we consider
${\cal U}_\gamma:=\{(\rho,\theta)\in\mathbb{R}^2:0<\rho<\epsilon,
(c_1-\gamma)\rho^{-1/\xi}<\theta<(c_1+\gamma)\rho^{-1/\xi}\}$. Clearly,
${\cal U}'_1$, ${\cal U}'_2$ and ${\cal U}'''_2$ are all subsets of
${\cal U}_{\epsilon,\delta,\gamma}:={\cal U}_{\epsilon,\delta}\backslash {\cal U}_\gamma$
for sufficiently small $\gamma>0$.

\begin{lm}
For small $\gamma>0$ the following inequality holds:
\begin{align*}
\frac{\partial }{\partial \theta} \frac{{\cal G}(\rho,\theta)}{{\cal H}(\rho,\theta)}
=\frac{[{\cal G},{\cal H}]_\theta}{{\cal H}^2}\le \rho^{1/n}
~~~\mbox{for all}~(\rho,\theta)\in{\cal U}_{\epsilon,\delta,\gamma}.
\end{align*}
\label{Lie-th2.2}
\end{lm}

The above lemma will be proved after the proof of the theorem is completed.
This lemma clearly ensures that
neither ${\cal U}'_1$ nor ${\cal U}'''_2$ contains orbits of system~\eqref{equ:polar system2}
connecting with $O$
by the argument given just below \eqref{GHH2U}.
For ${\cal U}'_2$,
there are three cases:
{\bf(U$_+$)} $u({\cal G},\xi)>u({\cal H},\xi)-1/\xi$,
{\bf(U$_-$)} $u({\cal G},\xi)<u({\cal H},\xi)-1/\xi$ and
{\bf(U$_0$)} $u({\cal G},\xi)=u({\cal H},\xi)-1/\xi$,
where both $u({\cal G},\xi)$ and $u({\cal H},\xi)$ are defined just below \eqref{defscomp}.
Case {\bf(U$_0$)} is explained in Remark~\ref{Rk:th22u} because
the number of orbits in ${\cal U}'_2$ does not decided by the principal parts.

In case {\bf(U$_+$)},
consider the curve
$
\theta=\tilde{\theta}(\rho):=\tilde{c}\rho^{-1/\xi},~\rho\in(0,\epsilon),
$
for a constant $\tilde{c}\in(0,c_1-\gamma)$ and consider the region
$$
\widetilde{\cal U}'_2:=\{(\rho,\theta)\in\mathbb{R}^2:0<\rho<\epsilon,
\tilde{\theta}(\rho)<\theta<\theta_1(\rho)\}.
$$
Similarly to \eqref{eld}, we compute
\begin{eqnarray}
\left.\frac{d\theta}{d\rho}\right|_{\theta=\tilde{\theta}(\rho)}
=\left.\frac{{\cal G}(\rho,\theta)}
{\rho{\cal H}(\rho,\theta)}\right|_{\theta=\tilde{\theta}(\rho)}
<2M\rho^{u({\cal G},\xi)-u({\cal H},\xi)-1}
<-\frac{\tilde{c}}{\xi}\rho^{-1/\xi-1}
=\frac{d\tilde{\theta}(\rho)}{d\rho},
\label{Muu}
\end{eqnarray}
where $M:={\cal K}_\xi({\cal G})(\tilde{c})/{\cal K}_\xi({\cal H})(\tilde{c})$,
showing that
orbits starting from the curve $\theta=\tilde{\theta}(\rho)$
leave the region $\widetilde{\cal U}'_2$.
On the other hand,
the inequality $\dot\theta|_{\Upsilon_1}>0$ in \eqref{3theta2rho},
which still holds in the present circumstance,
implies that
orbits starting from the curve $\theta=\theta_1(\rho)$ also leave
the region $\widetilde{\cal U}'_2$.
Hence,
$\widetilde{\cal U}'_2$ is a generalized normal sector of Class I (\cite{T-Z1}),
implying by \cite[Lemma~1]{T-Z1} that
system~\eqref{equ:initial} has infinitely many orbits
connecting with $O$ in the direction $\theta=0$.
Moreover,
similar to the situation $N({\cal H})= 1$
in {\bf Step 4} of the proof of Theorem~\ref{th:finite},
there is one e-tsector,
passing through the only real branch $\theta=\theta_1(\rho)$
of the equation ${\cal H}(\rho,\theta)=0$,
but no other tsectors.

In case {\bf(U$_-$)},
the curve $\theta=\tilde\theta(\rho)$, considered in case {\bf(U$_+$)},
divides the region ${\cal U}'_2$
into two subregions $\widetilde{\cal U}'_2$ and ${\cal V}$,
where $\widetilde{\cal U}'_2$ is defined just before \eqref{Muu} and
$$
{\cal V}:=\{(\rho,\theta)\in\mathbb{R}^2:0<\rho<\epsilon,0<\theta<\tilde\theta(\rho)\}.
$$
For $\widetilde{\cal U}'_2$,
the lower edge of which is the curve $\theta=\tilde\theta(\rho)$,
by {\bf(F1$'$)} of {\bf Claim~4.1}
we see from \eqref{dtdrdttdr} that
\begin{align*}
\frac{d\theta}{d\rho}
=\frac{{\cal G}(\rho,\theta)}{\rho{\cal H}(\rho,\theta)}
>\frac{M}{2}\rho^{u({\cal G},\xi)-u({\cal H},\xi)-1}
~~~\forall (\rho,\theta)\in {\cal U}'_2\setminus\widetilde{\cal U}'_2.
\end{align*}
If system~\eqref{equ:polar system2} has an orbit
$\theta=\vartheta(\rho)$ connecting with $O$ in $\widetilde{\cal U}'_2$,
then the above inequality implies that
$$
\frac{d\vartheta(\rho)}{d\rho}>\frac{M}{2}\rho^{u({\cal G},\xi)-u({\cal H},\xi)-1}.
$$
However,
since $\tilde{\theta}(\rho),\theta_1(\rho)= O(\rho^{-1/\xi})$,
we have $\vartheta(\rho)= O(\rho^{-1/\xi})$ and therefore
$d\vartheta(\rho)/d\rho= O(\rho^{-1/\xi-1})$,
a contradiction to the above inequality
because $u({\cal G},\xi)-u({\cal H},\xi)-1<-1/\xi-1$ in case {\bf(U$_-$)}.
Then system~\eqref{equ:polar system2} has no orbits connecting with $O$
in $\widetilde{\cal U}'_2$.

For the other subregion ${\cal V}$,
similarly to (\ref{Muu}), we have
\begin{align}
\left.\frac{d\theta}{d\rho}\right|_{\theta=\tilde{\theta}(\rho)}
=\left.\frac{{\cal G}(\rho,\theta)}
{\rho{\cal H}(\rho,\theta)}\right|_{\theta=\tilde{\theta}(\rho)}
>\frac{M}{2}\rho^{u({\cal G},\xi)-u({\cal H},\xi)-1}
>-\frac{\tilde{c}}{\xi}\rho^{-1/\xi-1}
=\frac{d\tilde{\theta}(\rho)}{d\rho},
\label{dtdrdttdr}
\end{align}
which together with
the inequality $\dot\rho|_{{\cal U}_2}>0$ given in \eqref{3theta2rho}
ensures that orbits starting from the curve $\theta=\tilde{\theta}(\rho)$
leave the region ${\cal V}$.
Moreover,
since $a_{i_{s({\cal G})},j_{s({\cal G})}}({\cal G})>0$
in case {\bf(4$_-$)} because of \eqref{assume-c3},
similar to \eqref{dtheta0},
\begin{align*}
\dot\theta|_{\theta=0}={\cal G}(0,\rho)/\rho
=a_{i_{s({\cal G})},j_{s({\cal G})}}({\cal G})\rho^{i_{s({\cal G})}-1}
+O(\rho^{i_{s({\cal G})}})
>0,
\end{align*}
which implies that orbits starting from the curve $\theta=0$
enter the region ${\cal V}$.
Thus the region ${\cal V}$
is a generalized normal sector of Class III (\cite{T-Z1}),
implying by \cite[Lemma~3]{T-Z1}
that system~\eqref{equ:initial} has either no or infinitely many orbits
connecting with $O$ in ${\cal V}$.
Further, since $\tilde{c}\in(0,c_1-\gamma)$,
we have ${\cal V} \subset \widetilde{\cal U}_{\epsilon,\delta}$.
Then we can use Lemma~\ref{Lie-th2.2} to discuss
as in the argument given just below \eqref{GHH2U},
which ensures that ${\cal V}$ contains at most one orbit connecting with $O$.
Consequently,
no orbits connect with $O$ in ${\cal V}$.

Summarily,
system~\eqref{equ:initial} has no orbits connecting with $O$ in the direction $\theta=0$.
Moreover,
similar to the situation $N({\cal H})= 1$
in {\bf Step 4} of the proof of Theorem~\ref{th:finite},
there is one e-tsector,
passing through the only real branch $\theta=\theta_1(\rho)$
of the equation ${\cal H}(\rho,\theta)=0$,
but no other tsectors.
Thus the proof of Theorem~\ref{th:finiteJ2} is completed.


Finally, we prove {\bf Claim H1$'$} and {\bf Claim G1$'$}.
In fact,
the equation ${\cal H}(\rho,\theta)=0$ has no real branches in the fourth quadrant
since the only real branch $\theta=\theta_1(\rho)$ lies in the the first quadrant
in {\bf (4$_-^+$)}.
By the remark given just below Proposition~\ref{lm-solution by NP},
\begin{align}
\Xi({\cal A}_j({\cal H}))
=\Xi((-1)^{\tilde{j}_{s({\cal H})}}
a_{\tilde{i}_{s({\cal H})},\tilde{j}_{s({\cal H})}}({\cal H}),...,
(-1)^{\tilde{j}_0}a_{\tilde{i}_0,\tilde{j}_0}({\cal H}))
=0.
\label{XAJH}
\end{align}
Note that $\tilde{j}_0$ is odd and $\tilde{j}_{s({\cal H})}$ is even,
as indicated before \eqref{assume-c3} and in {\bf (J2)} respectively,
and that $a_{\tilde{i}_0,\tilde{j}_0}({\cal H})<0$ in \eqref{assume-c3}.
Then \eqref{XAJH} ensures that
$a_{\tilde{i}_{s({\cal H})},\tilde{j}_{s({\cal H})}}({\cal H})>0$.
If the parities of $\tilde{j}_k$\,s given in {\bf Claim H1$'$} are not true,
then there are integers $\tilde{\ell}, \tilde{\ell}'\in\{1,...,s({\cal H})-1\}$ such that
$\tilde{\ell}<\tilde{\ell}'$,
$\tilde{j}_{\tilde{\ell}}$ is even and $\tilde{j}_{\tilde{\ell}'}$ is odd.
We see from \eqref{XAJH} that
$a_{\tilde{i}_{\tilde{\ell}},\tilde{i}_{\tilde{\ell}}}({\cal H})>0$ and
$a_{\tilde{i}_{\tilde{\ell}'},\tilde{i}_{\tilde{\ell}'}}({\cal H})<0$.
Thus $\Xi({\cal A}({\cal H}))\ge 3$,
implying by Proposition~\ref{lm-solution by NP} that
the equation ${\cal H}(\rho,\theta)=0$ has at least 3 real branches,
a contradiction.
This shows that $\tilde{j}_0,...,\tilde{j}_{q-1}$ are odd and
$\tilde{j}_q,...,\tilde{j}_{s({\cal H})}$ are even for an integer $q$,
as stated in {\bf Claim H1$'$}.
Moreover,
since all elements in the sequence in \eqref{XAJH} are positive,
we obtain that
\begin{align*}
a_{\tilde{i}_0,\tilde{j}_0}({\cal H}),...,a_{\tilde{i}_{q-1},\tilde{j}_{q-1}}({\cal H})<0
~~~\mbox{and}~~~
a_{\tilde{i}_q,\tilde{j}_q}({\cal H}),...,
a_{\tilde{i}_{s({\cal H})},\tilde{j}_{s({\cal H})}}({\cal H})>0,
\end{align*}
which proves {\bf Claim H1$'$}.
By {\bf Claim H1$'$} and {\bf (P2)},
we obtain that
there is an integer $q'\in\{1,...,s({\cal G})-1\}$ such that
\begin{align}
\zeta(E_{q'}({\cal G}))=\zeta(E_q({\cal H})),
~~~
j_0,...,j_{q'-1}~\mbox{are even}
~~~\mbox{and}~~~
j_{q'},...,j_{s({\cal G})-1}~\mbox{are odd},
\label{CG1'q'}
\end{align}
as stated in {\bf Claim G1$'$}.
On the other hand,
we see from \eqref{assume-c3} that $a_{i_0,j_0}({\cal G})<0$
and therefore $a_{i_{s({\cal G})},j_{s({\cal G})}}({\cal G})>0$
since they have different signs in {\bf(4$_-^+$)}.
If there are integers $\ell,\ell'\in\{1,...,s({\cal G})-1\}$ such that $\ell<\ell'$,
$a_{i_{\ell},j_{\ell}}({\cal G})>0$ and $a_{i_{\ell'},j_{\ell'}}({\cal G})<0$,
then $\Xi({\cal A}({\cal G}))\ge 3$,
implying by Proposition~\ref{lm-solution by NP} that
the equation ${\cal G}(\rho,\theta)=0$ has at least 3 real branches,
a contradiction to the fact of exactly 2 real branches given in \eqref{bhG}.
Thus, there is an integer $q''$ such that
\begin{align*}
a_{i_0,j_0}({\cal G})<0,...,a_{i_{q''-1},j_{q''-1}}({\cal G})<0
~~~\mbox{and}~~~
a_{i_{q''},j_{q''}}({\cal G})>0,...,a_{i_{s({\cal G})},j_{s({\cal G})}}({\cal G})>0.
\end{align*}
Finally, we complete the proof of {\bf Claim G1$'$} by showing that $q''=q'$.
In fact,
if $q''<q'$ then the two inequalities
\begin{align*}
(-1)^{j_{s({\cal G})}}a_{i_{s({\cal G})},j_{s({\cal G})}}({\cal G})
(-1)^{j_{s({\cal G})}-1}a_{i_{s({\cal G})-1},j_{s({\cal G})-1}}({\cal G})
<0,
\\
(-1)^{j_{q'}}a_{i_{q'},j_{q'}}({\cal G})
(-1)^{j_{q''}}a_{i_{q''},j_{q''}}({\cal G})
<0
\end{align*}
imply that $\Xi({\cal A}_j({\cal G}))\ge 2$
since $q'\le s({\cal G})$ as indicated just before \eqref{CG1'q'}.
Note that the above existence of $q''$ indicates that $\Xi({\cal A}({\cal G}))=1$.
Then
$
\Xi({\cal A}({\cal G}),{\cal A}_j({\cal G}))
\ge
\Xi({\cal A}({\cal G}))+
\Xi({\cal A}_j({\cal G}))\ge 3,
$
the same contradiction as above.
In the other case $q''>q'$,
since
$a_{\tilde{i}_{q-1},\tilde{j}_{q-1}}({\cal H})
a_{\tilde{i}_q,\tilde{j}_q}({\cal H})<0$,
the edge-polynomial
$$
{\cal H}_{E_q({\cal H})}(\theta)=
a_{\tilde{i}_{q-1},\tilde{j}_{q-1}}({\cal H})\theta^{\tilde{j}_{q-1}}+\cdots+
a_{\tilde{i}_q,\tilde{j}_q}({\cal H})\theta^{\tilde{j}_q}
$$
has a simple positive real root $c_1$,
which implies by Lemma~\ref{lm-basic} that
the only real branch $\theta=\theta_1(\rho)$ of the equation ${\cal H}(\rho,\theta)=0$
is of the form
\begin{align}
\theta_1(\rho)=c_1\rho^{-1/\xi}+o(\rho^{-1/\xi})
\label{theta1H}
\end{align}
with $\xi=\zeta(E_{q}({\cal H}))$.
Similarly,
since $a_{i_{q''-1},j_{q''-1}}({\cal G})a_{i_{q''},j_{q''}}({\cal G})<0$
and
$$(-1)^{j_{s({\cal G})}}a_{i_{s({\cal G})},j_{s({\cal G})}}({\cal G})
(-1)^{j_{s({\cal G})-1}}a_{i_{s({\cal G})-1},j_{s({\cal G})-1}}({\cal G})<0,
$$
edge-polynomials
${\cal G}_{E_{q''}({\cal G})}$
and
${\cal G}_{E_{s({\cal G})}({\cal G})}$
have a simple positive real root $\hat{c}_1$ and
a simple negative real root $\hat{c}_2$ respectively.
By Lemma~\ref{lm-basic},
the two real branches $\theta=\theta_1(\rho)$ and $\theta=\theta_2(\rho)$,
which are obtained in ${\cal U}_1$ and ${\cal U}_2$ as in \eqref{bhG} respectively,
are of the forms
\begin{align}
\hat{\theta}_1(\rho)=\hat{c}_1\rho^{-1/\xi''}+o(\rho^{-1/\xi''})
~~~\mbox{and}~~~
\hat{\theta}_2(\rho)=\hat{c}_2\rho^{-1/\iota}+o(\rho^{-1/\iota})
\label{theta12G}
\end{align}
with $\xi'':=\zeta(E_{q''}({\cal H}))$ and $\iota:=\zeta(E_{s({\cal G})}({\cal G}))$ respectively.
Note that
$\zeta(E_{q}({\cal H}))=\zeta(E_{q'}({\cal G}))<\zeta(E_{q''}({\cal G}))<0$
in the case $q''>q'$.
Then we see from \eqref{theta1H} and the first equality of \eqref{theta12G} that
the curve $\theta=\hat{\theta}_1(\rho)$ does not lie in the region
${\cal U}_1=\{(\rho,\theta)\in\mathbb{R}^2:\theta_1(\rho)<\theta<\delta\}$,
a contradiction to \eqref{bhG}.
Consequently,
$q''=q'$ and therefore {\bf Claim G1$'$} is proved.
\qquad$\Box$


{\bf Proof of Lemma~\ref{Lie-th2.2}.}
Actually,
we only need to show that
${\cal H}^2-\rho^{-1/n}f\ge 0$ on ${\cal U}_{\epsilon,\delta,\gamma}$,
where $f:=[{\cal G},{\cal H}]_\theta$.
By {\bf Claim H1$'$}, {\bf Claim G1$'$} and Lemma~\ref{precondiJ2},
we can similarly prove that \eqref{NH2Ng} and \eqref{NH2Nf} still hold.
Then functions $\hat {\cal H}^2-\hat{f}$ and $\hat {\cal H}^2$
have the same Newton polygon and corresponding edge-polynomials,
where $\hat {\cal H}$ and $\hat{f}$ are defined in \eqref{stretching-n}.
Although
edges of ${\cal N}_{\hat{\cal H}^2}$ are stretched from edges of ${\cal N}_{{\cal H}^2}$
as in (\ref{stretching-n}),
coefficients of their edge-polynomials with the corresponding valid indices are the same.
Thus $\hat {\cal H}^2$ and ${\cal H}^2$ have the same edge-polynomials
by definition \eqref{red}.
The same argument given just below \eqref{KGKHK4} shows that
each nonzero real root
(if exists) of each edge-polynomial of ${\cal H}$ is of multiplicity 1,
which determines one real branch of the equation ${\cal H}(\rho,\theta)=0$
by Lemma~\ref{lm-solution by NP}.
It follows from \eqref{H1'theta1} that only the edge-polynomial of ${\cal H}$
corresponding the $q$-th edge with  slope $\xi$ has one nonzero real root $c_1$
but the others have no nonzero real roots.
By Proposition~\ref{lm-MIL},
only the edge-polynomial of ${\cal H}^2$
corresponding to the $q$-th edge with slope $\xi$
has only one nonzero real roots $c_1$
but the others have no nonzero real roots.
Then only the edge-polynomial of $\hat{\cal H}^2$ corresponding to
the $q$-th edge with slope $\xi/n$ has one nonzero real root $c_1$
but the others have no nonzero real roots.
By Lemma~\ref{lm-basic},
each real branch (if exists) of the equation
$\hat {\cal H}^2-\hat{f}=0$ is of the form
$$
\theta=c_1\hat\rho^{-n/\xi}+o(\hat\rho^{-n/\xi}),
$$
which implies that the function ${\cal H}^2-\rho^{-1/n}f$ has definite signs on regions
\begin{align*}
&{\cal V}_+:=
\{(\rho,\theta)\in\mathbb{R}^2:0<\rho<\epsilon,(c_1+\gamma)\rho^{-1/\xi}<\theta<\delta\},
\\
&{\cal V}_-:=
\{(\rho,\theta)\in\mathbb{R}^2:0<\rho<\epsilon,-\delta<\theta<(c_1-\gamma)\rho^{-1/\xi}\}.
\end{align*}
In order to determine the signs,
we investigate the signs of $g:={\cal H}^2-\rho^{-1/n}f$ on curves
$\Gamma_{\pm}:\theta=\pm\rho^{-1/\nu},\rho\in(0,\epsilon)$,
where $\nu$ is smaller than the slope of the first edge of ${\cal N}_g$.
According to Proposition~\ref{lm-lowest term},
$$
{\rm sgn}(g|_{\Gamma_+})={\rm sgn}(a_{\check{i}_0,\check{j}_0}(g))
~~\mbox{and}~~
{\rm sgn}(g|_{\Gamma_-})={\rm sgn}((-1)^{\check{j}_0}a_{\check{i}_0,\check{j}_0}(g)),
$$
where
$(\check{i}_0,\check{j}_0)$ is the left-most vertex of the Newton polygon ${\cal N}_g$.
Since ${\cal H}^2$ and $g$ have the same Newton polygon,
$(\check{i}_0,\check{j}_0)=(2\tilde{i}_0,2\tilde{j_0})$,
as indicated in \eqref{ijH}.
Moreover,
by Proposition \ref{lm-MIL},
$$
a_{\check{i}_0,\check{j}_0}(g)
=a_{\check{i}_0,\check{j}_0}({\cal H}^2)
=a^2_{\tilde{i}_0,\tilde{j}_0}({\cal H})>0,
$$
which implies that
$$
{\rm sgn}(g|_{\Gamma_+})={\rm sgn}(a^2_{\tilde{i}_0,\tilde{j}_0}({\cal H}))
> 0
~~\mbox{and}~~
{\rm sgn}(g|_{\Gamma_-})={\rm sgn}((-1)^{2\tilde{j}_0}a^2_{\tilde{i}_0,\tilde{j}_0}({\cal H}))
> 0.
$$
Since
$\Gamma_\pm\subset {\cal V}_\pm$,
we see from the above two inequalities that
${\cal H}^2-\rho^{-1/n}f=g\ge 0$ on
${\cal U}_{\epsilon,\delta,\gamma}={\cal V}_+\cup{\cal V}_-$
and therefore Lemma~\ref{Lie-th2.2} is proved.
\qquad$\Box$

\begin{rmk}
{\rm
We see from the above proofs of Theorems \ref{th:finite} and \ref{th:finiteJ2}
({\bf Claims 1.3}, {\bf 2.3}, {\bf 3.3} and {\bf 4.3})
that the numbers $\eta_O$, ${\cal S}_O^e$, ${\cal S}_O^h$ and ${\cal S}_O^p$
are determined by the numbers of real branches of
the equations ${\cal G}(\rho,\theta)=0$ and ${\cal H}(\rho,\theta)=0$
on the half-plane $\rho\ge 0$
in the case that
the Lie-bracket $[{\cal G},{\cal H}]_\theta$ is {\it semi-definite} on
$\Omega_\epsilon:=\{(\rho,\theta)\in\mathbb{R}^2:0<\rho<\epsilon,|\theta|<\epsilon\}$
for an $\epsilon>0$.
By Lemma~\ref{defsignth1},
the semi-definiteness is guaranteed by either conditions {\bf(P1)} and {\bf(Q)}
or
conditions {\bf(P2)}, {\bf(Q)}, {\bf(S)} and one of conditions {\bf(H1)}-{\bf(H4)}.
Further,
by Lemma~\ref{Assertion1} or \ref{Assertion1'},
the numbers of real branches are exactly equal to
$N({\cal G})$ and $N({\cal H})$ respectively.
However,
in the case that some of those conditions are violated,
it gets  more difficult to judge the semi-definiteness
although Propositions~\ref{lm-region} and \ref{lm-positive}
can still be used as done for the polynomial \eqref{fpro43} in section~4.3.
In fact,
remainder parts of ${\cal G}$ and ${\cal H}$ will be involved.
Additionally,
in this case Lemmas~\ref{Assertion1} and \ref{Assertion1'} may not be true
and therefore the numbers of real branches are no longer given by
$N({\cal G})$ and $N({\cal H})$ respectively
because those conditions are all defined by principal parts of ${\cal G}$ and ${\cal H}$.
\label{Rk:general}
}
\end{rmk}

\section{Applications}

\setcounter{equation}{0}
\setcounter{lm}{0}
\setcounter{thm}{0}
\setcounter{rmk}{0}
\setcounter{df}{0}
\setcounter{cor}{0}

\begin{eg}
{\rm
Krauskopf and Rousseau (\cite{K-R}) investigated in their Theorem 3
phase portraits of the degenerate system
\begin{align}
\dot x=-axy-x^3,~~~\dot y=-y^2+bx^4
\label{equ:KR2h1e}
\end{align}
on the half-plane $x\ge 0$.
In the case {\bf (i)} $1/2<a<1$ and $1/a^2<b<1/(2a-1)$
and the case {\bf (ii)} $1/2<a<1$ and $b>1/(2a-1)$,
blowing-up the degenerate singular point $O$,
they found that
system~\eqref{equ:KR2h1e} has infinitely many orbits connecting with $O$
in each direction of $\theta=0$ and $\pm\pi/2$ in case {\bf (i)}
and infinitely many orbits connecting with $O$ in each direction of $\theta=\pm\pi/2$ but
no orbits in the direction $\theta=0$ in case {\bf (ii)}.
The phase portraits in the two cases were given in the fifth and sixth sub-figures
in their Figure 2 respectively,
i.e., our Figures~\ref{fig:eg60KR}(i) and \ref{fig:eg60KR}(ii) respectively.


\begin{figure}[h!]
\centering
\subcaptionbox*{(i)}{\includegraphics[height=1.6in,width=1.3in]{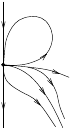}}
~~~~~~~
\subcaptionbox*{(ii)}{\includegraphics[height=1.6in,width=1.3in]{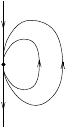}}
\caption{Phase portraits of \eqref{equ:KR2h1e} given in \cite{K-R}
         for cases {\bf (i)} and {\bf (ii)}.}
\label{fig:eg60KR}
\end{figure}


The phase portraits in Figures~\ref{fig:eg60KR}(i) and \ref{fig:eg60KR}(ii) show that
there are infinitely many orbits connecting with $O$ in the directions $\theta=\pm\pi/2$.
However,
system \eqref{equ:KR2h1e} is a special form of \eqref{equ:initial}
with $m=2$, $X_2(x,y)=-axy$ and $Y_2(x,y)=-y^2$,
implying by definitions given just below \eqref{equ:polar system} that
$G_0(\theta)=(a-1)\cos\theta\sin^2\theta$ and
$H_0(\theta)=-a\cos^2\theta\sin\theta-\sin^3\theta$.
Then system~\eqref{equ:KR2h1e} has exactly three exceptional directions
$\theta=\pm\pi/2$ and $\theta=0$ on the half-plane $x\ge 0$.
Moreover, $H_0(\pm\pi/2)\ne 0$ and $H_0(0)= 0$, implying that
the two ones $\theta=\pm\pi/2$ are regular exceptional directions but
the direction $\theta=0$ is irregular.
For $\theta=\pm\pi/2$, we have
$$
\frac{d}{d\theta}G_0(\pm\pi/2) H_0(\pm\pi/2)=a-1<0.
$$
By Theorem~6 of \cite[p.220]{SC} or Theorem~3.7 of \cite[p.70]{ZZF},
in each of the two directions there is a unique orbit connecting with $O$,
which lies on the $y$-axis.
This contradicts to the two phase portraits in Figure \ref{fig:eg60KR}.

In order to correct the two phase portraits, we discuss the direction $\theta=0$,
but none of results given in \cite{SC} and \cite{ZZF} can be applied to irregular ones.
Rewrite \eqref{equ:KR2h1e} in the polar coordinates as
$\dot \rho=\rho{\cal H}(\rho,\theta)$ and
$\dot \theta={\cal G}(\rho,\theta)$,
where
\begin{align*}
{\cal H}(\rho,\theta)
&=-a\theta-\rho+O(\theta^3)+O(\rho\theta^2)+O(\rho^2\theta),
\\
{\cal G}(\rho,\theta)
&=(a-1)\theta^2+\rho\theta+b\rho^2+O(\theta^4)+ O(\rho\theta^3)+O(\rho^2\theta^2),
\end{align*}
and compute the Newton polygons of ${\cal G}$ and ${\cal H}$.
Function ${\cal G}$ has two vertices $V_0({\cal G}):(0,2)$ and $V_1({\cal G}):(2,0)$,
which are linked by the Newton polygon ${\cal N}_{\cal G}$ as shown in Figure~\ref{fig:eg60NP}(b).
Similarly,
function ${\cal H}$ has two vertices $V_0({\cal H}):(0,1)$ and $V_1({\cal H}):(1,0)$,
which are linked by the Newton polygon ${\cal N}_{\cal H}$ as shown in Figure~\ref{fig:eg60NP}(c).
Noting that
ordinates of the right-most vertices of ${\cal N}_{\cal G}$ and ${\cal N}_{\cal H}$ are both zero,
we see that none of the two polygons ends above the $u$-axis, i.e.,
system \eqref{equ:KR2h1e} lies in the case {\bf (J2)}.
In order to determine the numbers of orbits and tsectors in the direction $\theta=0$,
we check conditions {\bf(P2)}, {\bf(Q)}, {\bf(S)} and {\bf(H4)} of our Theorem~\ref{th:finiteJ2}.

{\bf Assertion 1:} {\bf(P2)} is true.
Actually,
we have slope sequences
$\vec{\mathfrak{S}}({\cal G})=(-1)$ and $\vec{\mathfrak{S}}({\cal H})=(-1)$.
By \eqref{unionseq} and \eqref{interseq},
$$
\vec{\mathfrak{S}}({\cal G})\cup\vec{\mathfrak{S}}({\cal H})
=(-1)
~~~\mbox{and}~~~
\vec{\mathfrak{S}}({\cal G})\cap\vec{\mathfrak{S}}({\cal H})=(-1).
$$
Further \eqref{def-Lambda} shows that
vertex sequences
$\vec\Delta^V({\cal G})$ and $\vec\Delta^V({\cal H})$
are divided into subsequences
$\Lambda_0({\cal G})=(V_0({\cal G}))$,
$\Lambda_1({\cal G})=(V_1({\cal G}))$,
$\Lambda_0({\cal H})=(V_0({\cal H}))$ and
$\Lambda_1({\cal H})=(V_1({\cal H}))$.
Then {\bf(P2)} holds since
ordinates of lattice points in $\Lambda_0({\cal H})$
(and $\Lambda_1({\cal H})\setminus(V_{s({\cal H})}({\cal H}))$) are all odd (and odd)
and those in $\Lambda_0({\cal G})$ (and $\Lambda_1({\cal G})$) are all even (and even).

{\bf Assertion 2:} {\bf(Q)} is true.
Actually,
we only need to check that the polynomial
$[{\cal K}_\xi({\cal G}),{\cal K}_{\xi}({\cal H})]$
has no nonzero real roots for all
$\xi\in\vec{\mathfrak{S}}({\cal G})\cup\vec{\mathfrak{S}}({\cal H})=(-1)$.
Then {\bf(Q)} holds because the polynomial
$[{\cal K}_{-1}({\cal G}),{\cal K}_{-1}({\cal H})]=a(1-a)\theta^2+2(1-a)\theta+(ab-1)$
has no nonzero real roots in both cases {\bf(i)} and {\bf(ii)}.

{\bf Assertion 3:} {\bf(S)} is true.
Actually,
the second last valid index $(i_*,j_*)$ (and $(\tilde{i_*},\tilde{j}_*)$)
on ${\cal N}_{\cal G}$ (and ${\cal N}_{\cal H}$) is the point $(1,1)$ (and $(0,1)$).
Then {\bf(S)} holds because
we see from the above expansions of ${\cal G}$ and ${\cal H}$ that
$\{(i,j)\in\Delta({\cal G}):0<j<j_*\}
=\{(\tilde{i},\tilde{j})\in\Delta({\cal H}):0<\tilde{j}<\tilde{j}_*\}
=\emptyset.$

{\bf Assertion 4:} {\bf(H4)} is true.
Actually,
we have $\zeta(E_{s({\cal H})}({\cal H}))=-1=\zeta(E_{s({\cal G})}({\cal G}))$,
$\tilde{j}_*=j_*$ and
$
a_{\tilde{i}_*,\tilde{j}_*}({\cal H})
a_{i_{s({\cal G})},j_{s({\cal G})}}({\cal G})
=-ab
\ne
-1
=a_{\tilde{i}_{s({\cal H})},\tilde{j}_{s({\cal H})}}({\cal H})
a_{i_*,j_*}({\cal G}),
$
where $s({\cal H})=s({\cal G})=1$.

The above assertions show that
conditions of Theorem~\ref{th:finiteJ2} are verified fully by Remark~\ref{Rk:Sym}.
Thus Theorem~\ref{th:finiteJ2} indicates that
conclusions of subcases {\bf (ia)} and {\bf (ib)} in Theorem~\ref{th:finite} are true
because $G_0(\theta)\not\equiv 0$.
As defined just before Theorem~\ref{th:finite},
\begin{align*}
C_0&=(j_0-\tilde{j_0})a_{i_0,j_0}({\cal G})a_{\tilde{i}_0,\tilde{j}_0}({\cal H})=a(1-a)>0,
\\
N({\cal G})&=\chi({\cal G})=\sum_{\xi\in\mathfrak{S}({\cal G})}
\sharp\{\theta\in\mathbb{R}\backslash\{0\}:{\cal K}_\xi({\cal G})(\theta)=0\}=2,
\end{align*}
where $\mathfrak{S}({\cal G})=\{-1\}$ as indicated in {\bf Assertion~1} and
${\cal K}_{-1}({\cal G})(\theta)=(a-1)\theta^2+\theta+b$,
having two different nonzero real roots in both cases {\bf(i)} and {\bf(ii)}.
By Theorem~\ref{th:finite}{\bf (ib)},
system~\eqref{equ:KR2h1e} has infinitely many orbits connecting with $O$
in the direction $\theta=0$, which form exactly one e-tsector.
Hence,
as shown in Figure~\ref{fig:eg60NP}(a),
the phase portraits in the two cases {\bf(i)} and {\bf(ii)} are the same.
This is not consistent with those by Krauskopf and Rousseau.
Actually, although the phase portraits of blow-up systems in the two cases are different
as shown in Figure~3 of \cite{K-R}, they are the same after blowing-down.
}
\label{EG:KR}
\end{eg}

\begin{figure}[h]
  \centering
   \subcaptionbox{%
     }{\includegraphics[height=1.6in]{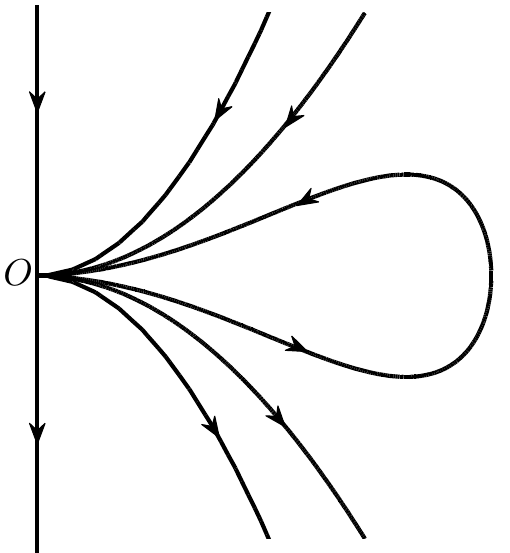}}~~~
   \subcaptionbox{%
     }{\includegraphics[height=1.6in]{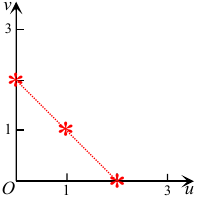}}~~~
   \subcaptionbox{%
     }{\includegraphics[height=1.6in]{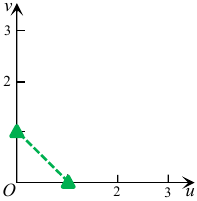}}
    \caption{(a) Phase portrait of \eqref{equ:KR2h1e} in both {\bf (i)} and {\bf (ii)}.
   (b) Newton polygon ${\cal N}_{\cal G}$. (c) Newton polygon ${\cal N}_{\cal H}$.}
  \label{fig:eg60NP}
\end{figure}

In addition to the 2-orbit case (cusp) given in the introduction,
we can use Corollary~\ref{cor-n2}{\bf(ia)} to construct a degenerate system
having exactly 3 orbits which connect with $O$ in the direction $\theta=0$
and make exactly 2 h-tsectors.

\begin{figure}[h!]
  \centering
  \includegraphics[height=1.6in]{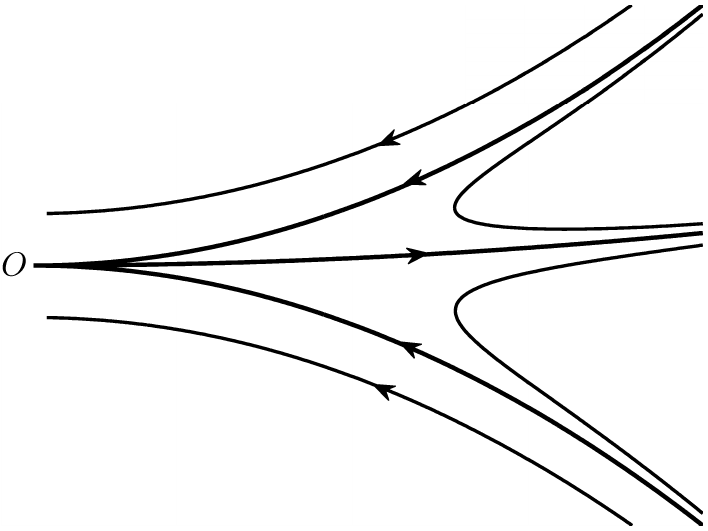}
  \caption{Orbits of system~\eqref{eg:3orb} in $\theta=0$.}
  \label{fig:eg61}
\end{figure}

\begin{eg}
{\rm
The degenerate system
\begin{eqnarray}
\dot x=-x^3y^6+xy^8+x^{11},~~~
\dot y=x^4y^5-x^9y+x^{11}
\label{eg:3orb}
\end{eqnarray}
has 3 orbits which connect with $O$ in the direction $\theta=0$ and
make 2 h-tsectors but no other tsectors,
as shown in Figure~\ref{fig:eg61}.
In fact,
this system is of the form \eqref{cor0ii}.
We choose $m=9$ and $n=2$,
then $\alpha_0=5$, $\alpha_1=1$,
$\beta_0=8$, $\beta_1=2$,
$\gamma_0=4$, $\gamma_1=9$,
$\delta_0=1$ and $\delta_1=8$, as defined in Corollary~\ref{cor-n2}.
By Corollary~\eqref{cor-n2}{\bf(ia)}, system~\eqref{eg:3orb} has $2n-1$ ($=3$) orbits
which connect with $O$ in the direction $\theta=0$ and
make 2 h-tsectors in this direction.
}
\label{EG:3orb}
\end{eg}


Note that Corollary~\ref{cor-n2} is a corollary of Corollary~\ref{cor:sign}{\bf(ii)},
which is aimed to the case {\bf(J2)} with condition {\bf(H1)} given in Theorem~\ref{th:finiteJ2}.
Actually, the above example obviously lies in this case and satisfies {\bf(H1)}.

The next example also lies in case {\bf(J2)} but does not satisfy condition {\bf(H1)}.
So Corollary~\ref{cor-n2} is not applicable.
Note that our Corollaries~\ref{cor:nonpara}-\ref{cor-n2} all require condition {\bf (NP)} (condition of no parallel edges),
but the example does not satisfy {\bf (NP)}.
Thus none of those corollaries can be applicable.
We will use Theorem~\ref{th:finiteJ2} directly.

\begin{figure}[h!]
  \centering
  \subcaptionbox{%
     }{\includegraphics[height=1.6in]{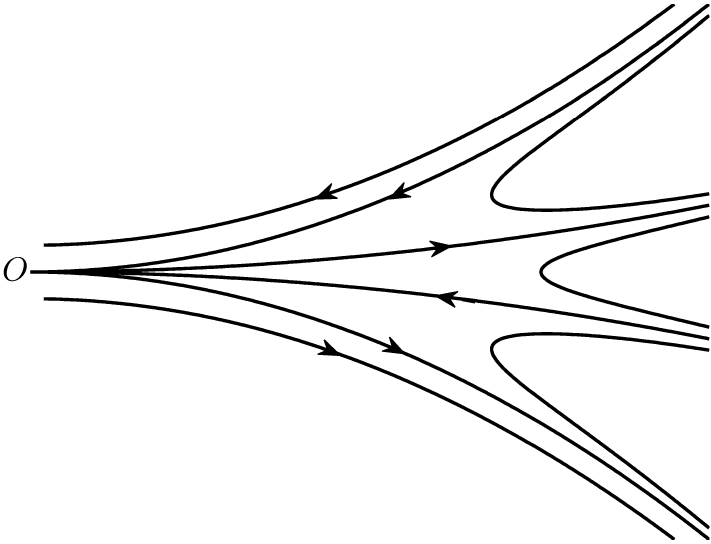}}~~~~~~
  \subcaptionbox{%
     }{\includegraphics[height=1.6in]{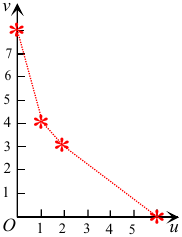}}~~~~~~
     \subcaptionbox{%
     }{\includegraphics[height=1.6in]{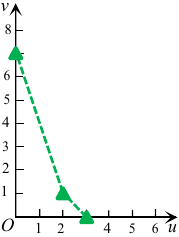}}
    \caption{(a) Orbits of system~\eqref{eg:4orb} in $\theta=0$,
    (b) ${\cal N}_{\cal G}$ and (c) ${\cal N}_{\cal H}$.}
  \label{fig:eg62}
\end{figure}

\begin{eg}
\label{EG:4orb}
{\rm
The degenerate system
\begin{eqnarray}
\left\{
\begin{array}{llll}
\dot x=-y^7+x^3y^5+x^8y-x^{10}/2-7x^{12}y,
\\
\dot y=-x^4y^4+x^6y^3+x^7y^2-x^9y/2+x^{13}
\end{array}
\right.
\label{eg:4orb}
\end{eqnarray}
has 4 orbits which connect with $O$ in the direction $\theta=0$ and
make 3 h-tsectors but no other tsectors,
as shown in Figure~\ref{fig:eg62}(a).
In fact,
similarly to \eqref{equ:polar system},
we rewrite \eqref{eg:4orb} in the polar coordinates as
$\dot \rho=\rho{\cal H}(\rho,\theta)$ and
$\dot \theta={\cal G}(\rho,\theta)$,
where
\begin{align*}
{\cal H}(\rho,\theta)
&=-\theta^7+\rho^2\theta-\rho^3/2
+O(\theta^9)+O(\rho^2\theta^3)+O(\rho^3\theta^2)+O(\rho^6\theta),
\\
{\cal G}(\rho,\theta)
&=\theta^8-\rho\theta^4+\rho^2\theta^3+\rho^6
+O(\theta^{10})+O(\rho\theta^6)+O(\rho^2\theta^5)+O(\rho^6\theta^4).
\end{align*}
The Newton polygon ${\cal N}_{\cal G}$ has $s({\cal G})$ ($=3$) edges,
linking vertices $V_0({\cal G}):(0,8)$,
$V_1({\cal G}):(1,4)$, $V_2({\cal G}):(2,3)$ and $V_3({\cal G}):(6,0)$ successively,
as shown in Figure~\ref{fig:eg62}(b).
The Newton polygon ${\cal N}_{\cal H}$ has $s({\cal H})$ ($=2$) edges,
linking vertices $V_0({\cal H}):(0,7)$, $V_1({\cal H}):(2,1)$ and
$V_2({\cal H}):(3,0)$ successively,
as shown in Figure~\ref{fig:eg62}(c).
Similar to Example~\ref{EG:KR},
system \eqref{eg:4orb} lies in the case {\bf (J2)}.
In order to determine the numbers of orbits and tsectors in the direction $\theta=0$,
we check conditions {\bf(P2)}, {\bf(Q)}, {\bf(S)} and {\bf(H2)} of our Theorem~\ref{th:finiteJ2}.

\begin{description}

\item[{\bf(P2)} is true.]
Actually, we have slope sequences
$\vec{\mathfrak{S}}({\cal G})=(-4,-1,-3/4)$ and $\vec{\mathfrak{S}}({\cal H})=(-3,-1)$.
By \eqref{unionseq} and \eqref{interseq},
$$
\vec{\mathfrak{S}}({\cal G})\cup\vec{\mathfrak{S}}({\cal H})
=(-4,-3,-1,-3/4)
~~~\mbox{and}~~~
\vec{\mathfrak{S}}({\cal G})\cap\vec{\mathfrak{S}}({\cal H})=(-1).
$$
Further \eqref{def-Lambda} shows that
the vertex sequences
$\vec\Delta^V({\cal G})$ and $\vec\Delta^V({\cal H})$
are divided into subsequences
$\Lambda_0({\cal G})=(V_0({\cal G}),V_1({\cal G}))$,
$\Lambda_1({\cal G})=(V_2({\cal G}),V_3({\cal G}))$,
$\Lambda_0({\cal H})=(V_0({\cal H}),V_1({\cal H}))$ and
$\Lambda_1({\cal H})=(V_2({\cal H}))$.
Then {\bf(P2)} holds since
ordinates of lattice points in $\Lambda_0({\cal G})$
(and $\Lambda_1({\cal G})\setminus(V_{s({\cal G})}({\cal G}))$) are all even (and odd)
and those in $\Lambda_0({\cal H})$ (and $\Lambda_1({\cal H})$) are all odd (and even).

\item[{\bf(Q)} is true.]
Actually,
we only need to check that
the polynomial $[{\cal K}_\xi({\cal G}),{\cal K}_{\xi}({\cal H})]$
has no nonzero real roots for all
$\xi\in\vec{\mathfrak{S}}({\cal G})\cup\vec{\mathfrak{S}}({\cal H})$.
We compute that
{\small
\begin{align*}
[{\cal K}_{-4}({\cal G}),{\cal K}_{-4}({\cal H})]&=-\theta^{10}(\theta^4+3),
&[{\cal K}_{-3}({\cal G}),{\cal K}_{-3}({\cal H})]&=-3\theta^4(\theta^6+1),
\\
[{\cal K}_{-1}({\cal G}),{\cal K}_{-1}({\cal H})]&=-\theta^2(3\theta^2-4\theta+3/2),
&[{\cal K}_{-\frac{3}{4}}({\cal G}),{\cal K}_{-\frac{3}{4}}({\cal H})]&=-3\theta^2/2,
\end{align*}
}each of which has no nonzero real roots, and therefore {\bf(Q)} holds.

\item[{\bf(S)} is true.]
Actually,
we see from the above expansions of ${\cal G}$ and ${\cal H}$ that
points $(i_*,j_*)=(2,3)$ and $(\tilde{i}_*,\tilde{j}_*)=(1,1)$
are the second last valid indices on
${\cal N}_{\cal G}$ and ${\cal N}_{\cal H}$ respectively.
Similar to Example~\ref{EG:KR}, {\bf(S)} holds.

\item[{\bf(H2)} is true.]
Actually, we have
$\zeta(E_{s({\cal G})}({\cal G}))=-3/4>-1=\zeta(E_{s({\cal H})}({\cal H}))$,
$j_*>\tilde{j}_*$, $\tilde{j}_*$ is odd, and
$
a_{i_*,j_*}({\cal G})
a_{\tilde{i}_{s({\cal H})},\tilde{j}_{s({\cal H})}}({\cal H})
a_{i_{s({\cal G})},j_{s({\cal G})}}({\cal G})
a_{\tilde{i}_*,\tilde{j}_*}({\cal H})
=-1/2<0.
$
\end{description}
Consequently,
conditions of Theorem~\ref{th:finiteJ2} are verified fully.
Thus Theorem~\ref{th:finiteJ2} indicates that
conclusions in Theorem~\ref{th:finite}{\bf (ia)} are true
because $G_0(\theta)\not\equiv 0$ and
$C_0=(j_0-\tilde{j}_0)a_{i_0,j_0}({\cal G})a_{\tilde{i}_0,\tilde{j_0}}({\cal H})=-1<0$.
Similarly to Example~\ref{EG:KR},
noticing that polynomials
${\cal K}_{-4}({\cal G})(\theta)=\theta^4(\theta^4-1)$,
${\cal K}_{-1}({\cal G})(\theta)=-\theta^3(\theta-1)$ and
${\cal K}_{-3/4}({\cal G})(\theta)=\theta^3+1$
have 2, 1 and 1 nonzero real roots respectively,
we see from definition of $N({\cal G})$ given just below \eqref{defChi} that $N({\cal G})=4$.
By Theorem~\ref{th:finite}{\bf (ia)},
system~\eqref{eg:4orb} has 4 orbits connecting with $O$ in the direction $\theta=0$
and forming 3 h-tsectors.
}
\end{eg}

The above example is devoted to {\bf(H2)} in case {\bf(J2)}.
The following example is concerning {\bf(H3)} in case {\bf(J2)}.

\begin{figure}[!h]
  \centering
  \subcaptionbox{%
     }{\includegraphics[height=1.6in]{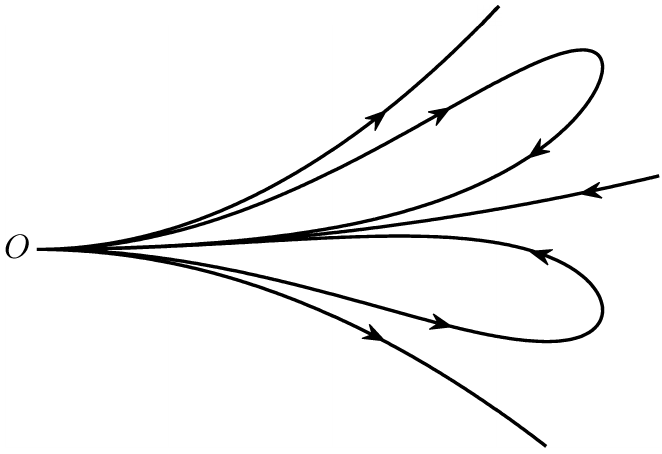}}~~~
       \subcaptionbox{%
     }{\includegraphics[height=1.6in]{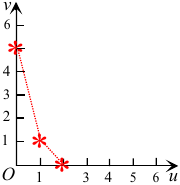}}~~~
     \subcaptionbox{%
     }{\includegraphics[height=1.6in]{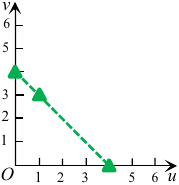}}
    \caption{(a) Orbits of system~\eqref{eg:2esector} in $\theta=0$,
    (b) ${\cal N}_{\cal G}$ and (c) ${\cal N}_{\cal H}$.}
  \label{fig:eg63}
\end{figure}

\begin{eg}
\label{EG:2esector}
{\rm
The degenerate system
\begin{eqnarray}
\dot x=x^2y^4+x^5y^2-x^4y^3-x^7y-x^{10},~~~
\dot y=2xy^5-x^6y+x^8-\frac{11}{2}x^9y+y^{11}
\label{eg:2esector}
\end{eqnarray}
has infinitely many orbits which connect with $O$ in the direction $\theta=0$ and
make 2 e-tsectors but no other tsectors,
as shown in Figure~\ref{fig:eg63}(a).
In fact, similarly to \eqref{equ:polar system},
we rewrite \eqref{eg:2esector} in the polar coordinates as
$\dot \rho=\rho{\cal H}(\rho,\theta)$ and
$\dot \theta={\cal G}(\rho,\theta)$,
where
\begin{align*}
{\cal H}(\rho,\theta)
=&\theta^4-\rho\theta^3-\rho^4+O(\theta^6)+O(\rho\theta^5)+O(\rho^4\theta^4),
\\
{\cal G}(\rho,\theta)
=&\theta^5-\rho\theta+\rho^2+O(\theta^7)+O(\rho\theta^3)+O(\rho^2\theta^2)+O(\rho^4\theta).
\end{align*}
The Newton polygon ${\cal N}_{\cal G}$ has $s({\cal G})$ ($=2$) edges,
linking vertices $V_0({\cal G}):(0,5)$,
$V_1({\cal G}):(1,1)$ and $V_2({\cal G}):(2,0)$ successively,
as shown in Figure~\ref{fig:eg63}(b).
The Newton polygon ${\cal N}_{\cal H}$ has $s({\cal H})$ ($=1$) edge
linking vertex $V_0({\cal G}):(0,4)$ with vertex $V_1({\cal G}):(4,0)$,
as shown in Figure~\ref{fig:eg63}(c).
Similar to Example~\ref{EG:4orb},
system \eqref{eg:2esector} lies in the case {\bf (J2)}.
In order to determine the numbers of orbits and tsectors in the direction $\theta=0$,
we check conditions {\bf(P2)}, {\bf(Q)}, {\bf(S)} and {\bf(H3)} of our Theorem~\ref{th:finiteJ2}.
We can verify {\bf(P2)} and {\bf(S)} similarly to Example~\ref{EG:4orb}.
Moreover,

\begin{description}

\item[{\bf(Q)} is true.]
Actually,
we have $\vec{\mathfrak{S}}({\cal G})\cup\vec{\mathfrak{S}}({\cal H})=(-4,-1)$
and all Lie-brackets are given by
$[{\cal K}_{-4}({\cal G}),{\cal K}_{-4}({\cal H})]=\theta^4(\theta^4+3)$
and
$[{\cal K}_{-1}({\cal G}),{\cal K}_{-1}({\cal H})]=3\theta^2(\theta-1)^2+1$,
both of which have no nonzero real roots.
Therefore, {\bf(Q)} holds.

\item[{\bf(H3)} is true.]
Actually,
$\zeta(E_{s({\cal G})}({\cal G}))=-1=\zeta(E_{s({\cal H})}({\cal H}))$ and
$j_*=1<3=\tilde{j}_*$.
\end{description}
Then conditions of Theorem~\ref{th:finiteJ2} are verified fully
and therefore Theorem~\ref{th:finite}{\bf(ib)} holds
because $G_0(\theta)\ne 0$ and
$C_0=(j_0-\tilde{j}_0)a_{i_0,j_0}({\cal G})a_{\tilde{i}_0,\tilde{j_0}}({\cal H})=1>0$.
Similarly to Example~\ref{EG:KR},
since polynomials
${\cal K}_{-4}({\cal G})(\theta)=\theta^5-\theta$ and
${\cal K}_{-1}({\cal G})(\theta)=-\theta+1$
have 2 and 1 nonzero real roots respectively,
we obtain $N({\cal G})=3$.
By Theorem~\ref{th:finite}{\bf(ib)},
there are 2 e-tsectors,
as shown in Figure~\ref{fig:eg63}(a)
}
\end{eg}

The following example is concerning {\bf(H4)} in case {\bf(J2)}
and shows 3 e-tsectors in an exceptional direction.

\begin{figure}[h!]
  \centering
  \subcaptionbox{%
     }{\includegraphics[height=1.6in]{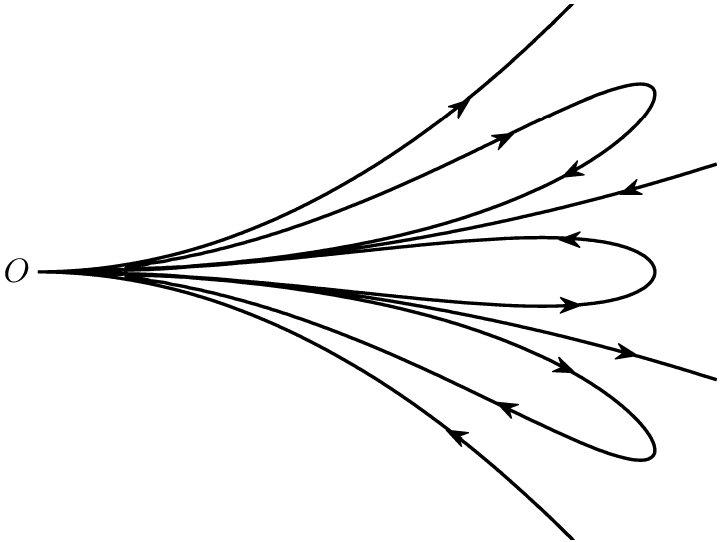}}~~~
  \subcaptionbox{%
     }{\includegraphics[height=1.6in]{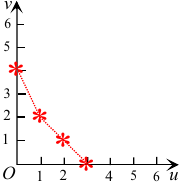}}~~~
     \subcaptionbox{%
     }{\includegraphics[height=1.6in]{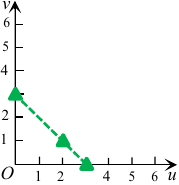}}
    \caption{(a) Orbits of system~\eqref{eg:3esector} in $\theta=0$,
    (b) ${\cal N}_{\cal G}$ and (c) ${\cal N}_{\cal H}$.}
  \label{fig:eg64}
\end{figure}

\begin{eg}
\label{EG:3esector}
{\rm
The degenerate system
\begin{eqnarray}
\dot x=x^2y^3+x^3y^3-6x^6y-x^8,~~~
\dot y=2xy^4-x^4y^2+x^6y+x^8+y^9
\label{eg:3esector}
\end{eqnarray}
has infinitely many orbits connecting with $O$ in the direction $\theta=0$
and making 3 e-tsectors and no other tsectors,
as shown in Figure~\ref{fig:eg64}(a).
In fact,
similarly to \eqref{equ:polar system},
we rewrite \eqref{eg:3esector} in the polar coordinates as
$\dot \rho=\rho{\cal H}(\rho,\theta)$ and
$\dot \theta={\cal G}(\rho,\theta)$,
where
\begin{align*}
{\cal H}(\rho,\theta)
=&\theta^3-6\rho^2\theta-\rho^3+O(\theta^7)+O(\rho^2\theta^2)+O(\rho^3\theta),
\\
{\cal G}(\rho,\theta)
=&\theta^4-\rho\theta^2+\rho^2\theta+\rho^3
+O(\theta^6)+O(\rho^2\theta^2)+O(\rho^3\theta).
\end{align*}
The Newton polygon ${\cal N}_{\cal G}$ has $s({\cal G})$ ($=2$) edges,
linking vertices $V_0({\cal G}):(0,4)$,
$V_1({\cal G}):(1,2)$ and $V_2({\cal G}):(3,0)$ successively,
as shown in Figure~\ref{fig:eg64}(b).
The Newton polygon ${\cal N}_{\cal H}$ has $s({\cal H})$ ($=1$) edge,
linking vertex $V_0({\cal G}):(0,3)$ with vertex $V_1({\cal G}):(3,0)$,
as shown in Figure~\ref{fig:eg64}(c).
Similar to Example~\ref{EG:4orb},
system \eqref{eg:3esector} lies in the case {\bf (J2)}.
In order to determine the numbers of orbits and tsectors in the direction $\theta=0$,
we check conditions {\bf(P2)}, {\bf(Q)}, {\bf(S)} and {\bf(H4)} of our Theorem~\ref{th:finiteJ2}.
We can verify {\bf(P2)} and {\bf(S)} similarly to Example~\ref{EG:4orb}.
Moreover,

\begin{description}

\item[{\bf(Q)} is true.]
Actually,
we have $\vec{\mathfrak{S}}({\cal G})\cup\vec{\mathfrak{S}}({\cal H})=(-2,-1)$
and all Lie-brackets are given by
$[{\cal K}_{-2}({\cal G}),{\cal K}_{-2}({\cal H})]=\theta^4(\theta^2+1)$
and
$[{\cal K}_{-1}({\cal G}),{\cal K}_{-1}({\cal H})]
=(\theta^2+\theta+1)(\theta^2-3\theta+5)$,
both of which have no nonzero real roots.
Then {\bf(Q)} holds.

\item[{\bf(H4)} is true.]
Actually, we have
$\zeta(E_{s({\cal G})}({\cal G}))=-1=\zeta(E_{s({\cal H})}({\cal H}))$ and
$j_*=1=\tilde{j}_*$, and
$
a_{i_*,j_*}({\cal G})a_{\tilde{i}_{s({\cal H})},\tilde{j}_{s({\cal H})}}({\cal H})
=-1\ne-6=
a_{i_{s({\cal G})},j_{s({\cal G})}}({\cal G})a_{\tilde{i}_*,\tilde{j}_*}({\cal H}).
$
\end{description}
Then conditions of Theorem~\ref{th:finiteJ2} are verified fully
and therefore, Theorem~\ref{th:finite}{\bf(ib)} holds
because $G_0(\theta)\ne 0$ and
$C_0=(j_0-\tilde{j}_0)a_{i_0,j_0}({\cal G})a_{\tilde{i}_0,\tilde{j_0}}({\cal H})=1>0$.
Similarly to Example~\ref{EG:KR},
since polynomials
${\cal K}_{-2}({\cal G})(\theta)=\theta^2(\theta^2-1)$ and
${\cal K}_{-1}({\cal G})(\theta)=-\theta^2+\theta+1$
have totally 4 nonzero real roots,
we obtain that $N({\cal G})=4$.
By Theorem~\ref{th:finite}{\bf(ib)},
there are 3 e-tsectors,
as shown in Figure~\ref{fig:eg64}(a).
}
\end{eg}

The above 5 examples all lie in case {\bf (J2)}. Actually, it is easier than in {\bf (J2)}
to give examples in case {\bf (J1)}. Similarly,
we prefer using Corollary~\ref{cor-n} immediately.
If none of conditions in Corollary~\ref{cor-n} is true,
we have to use Theorem~\ref{th:finite} directly
as we did in examples~\ref{EG:4orb}-\ref{EG:3esector}.

We finally give an example in case {\bf (J1)}, which
safisfises neither {\bf(NP)} (condition of no parallel edges)
nor {\bf(T)} (condition of no valid indices except for vertices).
Thus we cannot use Corollaries~\ref{cor:nonpara}-\ref{cor-n} but
have to use Theorem~\ref{th:finite} directly.
In comparison with example 6.2 and examples 6.3-6.4, which
does not satisfy {\bf(NP)} and
also satisfy neither {\bf(NP)} nor {\bf(T)} respectively,
the following example contains more valid indices on edges of the Newton polygon.
We need to use the Complete Discrimination System (\cite{YL}),
stated in the second paragraph in section~3.3,
to verify condition {\bf(Q)} and determine $N({\cal G})$ and $N({\cal H})$ as indicated in Remark~\ref{Rk:poly}.

\begin{figure}[!h]
    \centering
    \subcaptionbox{%
     }{\includegraphics[height=1.6in]{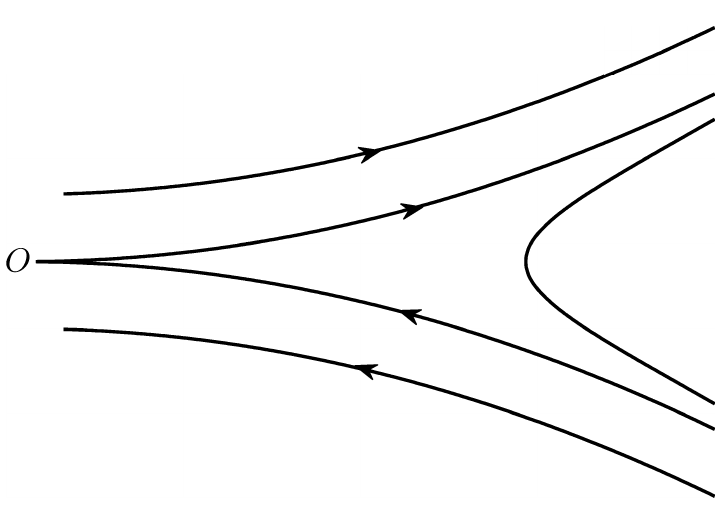}}~~~
     \subcaptionbox{%
     }{\includegraphics[height=1.6in]{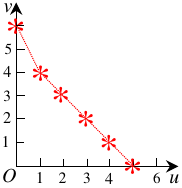}}~~~
     \subcaptionbox{%
     }{\includegraphics[height=1.6in]{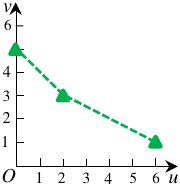}}
    \caption{(a) Orbits of system~\eqref{eg:2orb} in $\theta=0$,
    (b) ${\cal N}_{\cal G}$ and (c) ${\cal N}_{\cal H}$.}
    \label{fig:GHEg65}
\end{figure}

\begin{eg}
{\rm
The degenerate system
\begin{equation}
\left\{
\begin{array}{lllll}
\dot x=xy^5-x^2y^5+\frac{187}{1250}x^5y^3-\frac{13}{10}x^6y^3
-\frac{5}{8}x^8y^2-\frac{1489}{10000}x^{10}y+x^{11}y,
\\
\dot y=x^3y^4+\frac{3}{2}x^5y^3+\frac{13}{10}x^7y^2+\frac{5}{8}x^9y
+\frac{1489}{10000}x^{11}+x^{10}y^2+y^{13}
\end{array}
\right.
\label{eg:2orb}
\end{equation}
has 2 orbits which connect with $O$ in the direction $\theta=0$
and make 1 h-tsector but no other tsectors,
as shown in Figure~\ref{fig:GHEg65}(a).
In fact,
similarly to \eqref{equ:polar system},
we rewrite system~\eqref{eg:2orb} as
$\dot \rho=\rho{\cal H}(\rho,\theta)$ and $\dot \theta={\cal G}(\rho,\theta)$,
where
\begin{align*}
{\cal H}(\rho,\theta)
=&\theta^5+\frac{187}{1250}\rho^2\theta^3+\rho^6\theta
+O(\theta^7)+O(\rho^2\theta^4)+O(\rho^6\theta^3)+O(\rho^7\theta^{14}),
\\
{\cal G}(\rho,\theta)
=&-\theta^6
+\rho\theta^4
+\frac{3}{2}\rho^2\theta^3
+\frac{13}{10}\rho^3\theta^2
+\frac{5}{8}\rho^4\theta
+\frac{1489}{10000}\rho^5
+O(\theta^8)
\\
&
+O(\rho\theta^6)
+O(\rho^2\theta^4)
+O(\rho^3\theta^4)
+O(\rho^4\theta^3)
+O(\rho^5\theta^2)
+O(\rho^7\theta^{13}).
\end{align*}
The Newton polygon ${\cal N}_{\cal G}$ has two edges,
linking vertices $V_0({\cal G}):(0,6)$, $V_1({\cal G}):(1,4)$
and $V_2({\cal G}):(5,0)$ successively,
as shown in Figure~\ref{fig:GHEg65}(a).
The Newton polygon ${\cal N}_{\cal H}$ has two edges,
linking vertices $V_0({\cal H}):(0,5)$, $V_1({\cal H}):(2,3)$
and $V_2({\cal H}):(6,1)$ successively,
as shown in Figure~\ref{fig:GHEg65}(b).
Clearly,
${\cal N}_{\cal H}$ ends above the $u$-axis, i.e.,
system~\eqref{eg:2orb} lies in the case {\bf(J1)}.
In order to determine the numbers of orbits and tsectors in the direction $\theta=0$,
we check conditions {\bf(P1)} and {\bf(Q)} of Theorem~\ref{th:finite}.

\begin{description}

\item[{\bf(P1)} is true.]
Actually, we have
$\vec{\mathfrak{S}}({\cal G})=(-2,-1)$ and
$\vec{\mathfrak{S}}({\cal H})=(-1,-1/2)$.
By \eqref{unionseq} and \eqref{interseq},
$$
\vec{\mathfrak{S}}({\cal G})\cup\vec{\mathfrak{S}}({\cal H})=(-2,-1,-1/2)
~~~\mbox{and}~~~
\vec{\mathfrak{S}}({\cal G})\cap\vec{\mathfrak{S}}({\cal H})=(-1).
$$
Further \eqref{def-Lambda} shows that sequences
$\vec\Delta^V({\cal G})$ and $\vec\Delta^V({\cal H})$
are divided into subsequences
$\Lambda_0({\cal G})=(V_0({\cal G}),V_1({\cal G}))$,
$\Lambda_1({\cal G})=(V_2({\cal G}))$,
$\Lambda_0({\cal H})=(V_0({\cal H}))$ and
$\Lambda_1({\cal H})=(V_1({\cal H}),V_2({\cal H}))$.
Then {\bf(P1)} holds because
ordinates of lattice points in $\Lambda_0({\cal G})$ and $\Lambda_1({\cal G})$
are all even and those in $\Lambda_0({\cal H})$ and $\Lambda_1({\cal H})$ are all odd.

\item[{\bf(Q)} is true.]
Actually, we only need to check that the polynomial
$[{\cal K}_\xi({\cal G}),{\cal K}_{\xi}({\cal H})]$
has no nonzero real roots for all
$\xi\in\vec{\mathfrak{S}}({\cal G})\cup\vec{\mathfrak{S}}({\cal H})=(-2,-1,-1/2)$.
For $\xi=-2$, the polynomial
$[{\cal K}_{-2}({\cal G}),{\cal K}_{-2}({\cal H})]=-\theta^{10}-\theta^8$
has no nonzero real roots.
For $\xi=-1$, we have
$[{\cal K}_{-1}({\cal G}),{\cal K}_{-1}({\cal H})]=\theta^2P(\theta)$,
where
$$
P(\theta):=-\theta^6-3\theta^5-\frac{2344}{625}\theta^4-\frac{5}{2}\theta^3
-\frac{46949}{50000}\theta^2-\frac{187}{1000}\theta-\frac{835329}{12500000}.
$$
As defined in the second paragraph of section~3.3,
the revised sign list ${\rm rsl}(P)$ satisfies that
${\rm rsl}(P)=(1,-1,1,-1,-1,-1)$,
in which the number $\Phi({\rm rsl}(P))$ of non-vanishing members
and the number $\Xi({\rm rsl}(P))$ of sign changes are $6$ and $3$ respectively,
i.e., $\Phi({\rm rsl}(P))=2\Xi({\rm rsl}(P))$.
By Theorem~2.1 of \cite{YL},
$P$ has no real roots, i.e.,
$[{\cal K}_{-1}({\cal G}),{\cal K}_{-1}({\cal H})]$
has no nonzero real roots.
For $\xi=-1/2$,
the polynomial
$
[{\cal K}_{-\frac{1}{2}}({\cal G}),{\cal K}_{-\frac{1}{2}}({\cal H})]
=-\frac{1489}{10000}(\frac{561}{1250}\theta^2+1)
$
clearly has no nonzero real roots.
Thus, {\bf(Q)} holds.
\end{description}
Then conditions of Theorem~\ref{th:finite} are fully verified.
Since $G_0(\theta)\ne0$ and
$C_0=(j_0-\tilde{j}_0)a_{i_0,j_0}({\cal G})a_{\tilde{i}_0,\tilde{j_0}}({\cal H})=-1<0$,
by Theorem~\ref{th:finite}{\bf(ia)},
system~\eqref{eg:2orb} has $N({\cal G})$ orbits
which connect with $O$ in the direction $\theta=0$
and make $\max\{0,N({\cal G})-1\}$ h-tsectors.
Since ${\cal N}_{\cal G}$ ends on the $u$-axis,
as defined just before Theorem~\ref{th:finite},
$$
N({\cal G})=\sum_{\xi\in\mathfrak{S}({\cal G})}
\sharp\{\theta\in\mathbb{R}\backslash\{0\}:{\cal K}_\xi({\cal G})(\theta)=0\},
$$
where $\mathfrak{S}({\cal G})=\{-2,-1\}$.
For $\xi=-2$,
the polynomial ${\cal K}_{-2}({\cal G})(\theta)=\theta^4(1-\theta^2)$
has two nonzero real roots.
For $\xi=-1$,
the revised sign list of the polynomial
$$
{\cal K}_{-1}({\cal G})(\theta)
=\theta^4+\frac{3}{2}\theta^3+\frac{13}{10}\theta^2+\frac{5}{8}\theta+\frac{1489}{10000}
$$
is $(1,-1,-1,1)$,
in which numbers of non-vanishing members and sign changes are 4 and $2$ respectively.
Then ${\cal K}_{-1}({\cal G})$ has no real roots by Theorem~2.1 of \cite{YL}.
Thus, $N({\cal G})=2$, i.e.,
system~\eqref{eg:2orb} has 2 orbits connecting with $O$ in the direction $\theta=0$,
making 1 h-tsector.
}
\end{eg}


As indicated in Remark \ref{Rk:general},
if those conditions
{\bf(P1)} and {\bf(Q)} or
{\bf(P2)}, {\bf(Q)}, {\bf(S)} and one of conditions {\bf(H1)}-{\bf(H4)} are violated,
numbers of real branches of ${\cal G}(\rho,\theta)=0$ and ${\cal H}(\rho,\theta)=0$
are no longer determined by $N({\cal G})$ and $N({\cal H})$ respectively,
but remainder parts of ${\cal G}$ and ${\cal H}$ will be involved and things get complicated.
However,
by Proposition~\ref{lm-solution by NP} and the remark given just below its proof,
we can easily obtain lower bounds for numbers of the real branches
with the Newton polygons ${\cal N}_{\cal G}$ and ${\cal N}_{\cal H}$.
Actually, the equation ${\cal G}(\rho,\theta)=0$ (resp. ${\cal H}(\rho,\theta)=0$)
has at least $\Xi({\cal A}({\cal G}),{\cal A}_j({\cal G}))$
(resp. $\Xi({\cal A}({\cal H}),{\cal A}_j({\cal H}))$) real branches.

\appendix
\renewcommand{\thesection}{}
\renewcommand{\thesubsection}{A.\arabic{subsection}}
\section{Appendix}
{\small
\subsection{Useful theorems}

{\bf Definition} (\cite[p.152, Definition 6.1.2]{krantz}).
A real function $W(\boldsymbol{x},y)$,
analytic in a neighborhood of $({\bm 0},0)\in\mathbb{R}^n\times\mathbb{R}$,
is called a {\it Weierstrass polynomial} of degree $m$ if
$$
W(\boldsymbol{x},y)=y^m+w_{m-1}(\boldsymbol{x})y^{m-1}+\cdots+w_1(\boldsymbol{x})y+w_0(x),
$$
where $w_i$\,s are real analytic functions in a neighborhood of ${\bm 0}\in\mathbb{R}^n$ and
vanish at $\boldsymbol{x}={\bm 0}$.

Let $\Lambda(n)$ denote the set of multi-indices ${\footnotesize\bm\mu}=(\mu_1,...,\mu_n)$, where $\mu_i$\,s are all nonnegative integers.
For $\boldsymbol{x}=(x_1,...,x_n)\in \mathbb{R}^n$,
let $\boldsymbol{x}^{\footnotesize\bm\mu}:=x_1^{\mu_1}x_2^{\mu_2}\cdots x_n^{\mu_n}$.
As shown in \cite[pp.152-153, Theorem 6.1.3]{krantz}, we have the following.

\noindent
{\bf Theorem A1} {\bf (Weierstrass Preparation Theorem)}.
{\it
Let
$$
\Phi(\boldsymbol{x},y):=\sum_{\boldsymbol{\mu}\in\Lambda(n)}\sum_{j=0}^{+\infty}
\Phi_{\boldsymbol{\mu},j}\boldsymbol{x}^{\boldsymbol{\mu}} y^j,
$$
where $\boldsymbol{x}=(x_1,...,x_n)\in\mathbb{R}^n$ and $y\in\mathbb{R}$,
be real analytic in a neighborhood of $(\mathbf{0},0)\in\mathbb{R}^n\times\mathbb{R}$
and
suppose that there is a positive integer $m$ such that
$$
\Phi_{\boldsymbol{0},0}=\Phi_{\boldsymbol{0},1}=\cdots=\Phi_{\boldsymbol{0},m-1}=0,
~\mbox{and}~\Phi_{\boldsymbol{0},m}=1.
$$
Then
\begin{description}
  \item[(W1)]
for any real analytic $\Psi(x,y)$ in a neighborhood of $(\mathbf{0},0)\in\mathbb{R}^n\times\mathbb{R}$,
there are uniquely real analytic functions
\begin{eqnarray*}
Q(\boldsymbol{x},y)=\sum_{\boldsymbol{\mu}\in\Lambda(n)}\sum_{j=0}^{+\infty}
Q_{\boldsymbol{\mu},j}\boldsymbol{x}^{\boldsymbol{\mu}} y^j,
~~
R(\boldsymbol{x},y)=\sum_{\boldsymbol{\mu}\in\Lambda(n)}\sum_{j=0}^{+\infty}
R_{\boldsymbol{\mu},j}\boldsymbol{x}^{\boldsymbol{\mu}} y^j,
\end{eqnarray*}
where $R_{{\bm\mu},j}=0$ for all $j\ge m$ and all multi-indices ${\bm\mu}$, such that
$
\Psi=Q\Phi+R.
$

\item[(W2)]
there is a function $U(\boldsymbol{x},y)$,
real analytic and non-vanishing in a neighborhood $\Omega$ of $({\bm 0},0)\in\mathbb{R}^n\times\mathbb{R}$
such that
$
U\Phi=W,
$
a Weierstrass polynomial of degree $m$, in $\Omega$.
\end{description}
}

\noindent
{\bf Theorem A2} ({\bf Puiseux's Theorem}, \cite[p.100]{krantz}).
{\it
Let $W(x,y)$ be a polynomial in $y$ of the form
$$
W(x,y)=y^n+w_{n-1}(x)y^{n-1}+\cdots+w_0(x),
$$
where each $w_i$ is real analytic at $x=0$.
Then there is an integer $d>0$ such that
$$
W(u^d,y)=(y-R_1(u))(y-R_2(u))\cdots (y-R_m(u))Q(u,y),
$$
where $R_i$\,s are real analytic at $u=0$ and
$Q(u,y)$ is a polynomial of degree $n-m$ in $y$
whose coefficients are real functions of $u$ analytic at $u=0$
but has no real zeros for any non-vanished $u$ near $0$.
}

Note that integer $m\in\{0,...,n\}$, and if $m=0$ then $W(u^d,y)=Q(u,y)$.

\subsection{Symbol sequences}

{\bf Lemma}.
{\it
Symbols ${\bf Y_+}$ and ${\bf Y_-}$ do not coexist
in the symbol sequence $\Theta^{(j)}(f)$ for $j\ge1$.
}

{\bf Proof.}
Since $f$ satisfies condition {\bf (U)}, given just before Proposition~\ref{lm-region},
coefficients of vertices of ${\cal N}_f$ have the same sign by Proposition~\ref{lm-region}.
Assume without loss of generality that they are all positive.
Then, for each edge $E_k\in\daleth$,
the edge set defined just before Proposition~\ref{lm-region},
we have
$$
f_{E_k}(\theta)=a_{i_{k-1},j_{k-1}}(f)\theta^{j_{k-1}}+\cdots+a_{i_k,j_k}(f)\theta^{j_k}
\ge0~~~\mbox{for all}~\theta\in\mathbb{R}
$$
because of condition {\bf(U)}.
Moreover, for each nonzero real root $c_i$ of $f_{E_k}$,
similar to the computation in \eqref{eld},
$
f_i^{(1)}(0,\theta)=f_{E_k}(c_i+\theta),
$
which is semi-positive.
On the other hand,
$
f_i^{(1)}(0,\theta)=a_{\hat{i}_0,\hat{j}_0}\theta^{\hat{j}_0}+O(\theta^{\hat{j}_0+1}),
$
where $(\hat{i}_0,\hat{j}_0)$ is the left end-point of ${\cal N}_{f_i^{(1)}}$.
Hence,
the coefficient of the left end-point of ${\cal N}_{f_i^{(1)}}$ is positive.
By Proposition~\ref{lm-region},
$f_i^{(1)}$ cannot be semi-negative and therefore,
the $i$-th symbol $\vartheta_i^{(1)}$ of the sequence $\Theta^{(1)}(f)$
cannot be ${\bf Y}_-$,
i.e., $\Theta^{(1)}(f)$ does not contain the symbol ${\bf Y_-}$.

If $f_i^{(1)}$ satisfies the same condition as $f$ in {\bf (U)},
then coefficients corresponding to vertices of ${\cal N}_{f_i^{(1)}}$ are all positive
by Proposition~\ref{lm-region}.
Hence, the sequence $\Theta_i^{(2)}(f_i^{(1)})$ does not contain the symbol ${\bf Y}_-$
for the same reason as $\Theta^{(1)}(f)$.
Then  $\Theta^{(2)}(f)$ does not contain the symbol ${\bf Y}_-$ either,
because $\Theta^{(2)}(f)$ is obtained by replacing
the symbol ${\bf U}$ at $\vartheta_i^{(1)}$ in $\Theta^{(1)}(f)$
with the subsequence $\Theta_i^{(2)}(f_i^{(1)})$ for all possible $i$.
One can prove inductively that
the $j$-th symbol sequence $\Theta^{(j)}(f)$ does not contain the symbol ${\bf Y}_-$ for $j\ge3$,
since for those desingularized functions satisfying condition {\bf (U)} we just repeat the same procedure as for $f$.
Thus, the lemma is proved.
\qquad$\Box$
}


{\bf Acknowledgement}: The author Jun Zhang is supported by NSFC \# 12101087.
The author Xingwu Chen is supported by NSFC \# 11871355.
The author Weinian Zhang is supported by
National Key R\&D Program of China (2022YFA1005900)
and
NSFC \# 12171336 and \# 11831012.



\end{document}